\documentclass[reqno]{amsart}
\usepackage[margin = 1in]{geometry}
\usepackage{amsmath, amssymb, amsthm, fancyhdr, verbatim, graphicx}
\usepackage{enumerate}
\usepackage[all]{xy}
\usepackage[usenames,dvipsnames]{xcolor}
\usepackage{mathrsfs}
\usepackage{tikz-cd}
\usepackage{framed, hyperref}
\usepackage[titletoc]{appendix}
\usepackage{bbm}
\usepackage{lipsum}
\usepackage{adjustbox}

\numberwithin{equation}{section}

\newcommand{\LL}{\mathfrak{L}}
\newcommand{\LN}{\mathfrak{N}}

\newcommand{\rd}{\mathrm{d}}

\DeclareFontFamily{U}{matha}{\hyphenchar\font45}
\DeclareFontShape{U}{matha}{m}{n}{
      <5> <6> <7> <8> <9> <10> gen * matha
      <10.95> matha10 <12> <14.4> <17.28> <20.74> <24.88> matha12
      }{}
\DeclareSymbolFont{matha}{U}{matha}{m}{n}
\DeclareFontFamily{U}{mathx}{\hyphenchar\font45}
\DeclareFontShape{U}{mathx}{m}{n}{
      <5> <6> <7> <8> <9> <10>
      <10.95> <12> <14.4> <17.28> <20.74> <24.88>
      mathx10
      }{}
\DeclareSymbolFont{mathx}{U}{mathx}{m}{n}

\DeclareMathSymbol{\obot}         {2}{matha}{"6B}

\RequirePackage{xspace}

\newcommand{\BA}{\ensuremath{\mathbb{A}}\xspace}

\newcommand{\BC}{\ensuremath{\mathbb{C}}\xspace}

\newcommand{\BG}{\ensuremath{\mathbb{G}}\xspace}

\newcommand{\BQ}{\ensuremath{\mathbb{Q}}\xspace}

\newcommand{\BZ}{\ensuremath{\mathbb{Z}}\xspace}

\newcommand{\sA}{\ensuremath{\mathscr{A}}\xspace}

\newcommand{\sE}{\ensuremath{\mathscr{E}}\xspace}

\newcommand{\sG}{\ensuremath{\mathscr{G}}\xspace}
\newcommand{\sH}{\ensuremath{\mathscr{H}}\xspace}

\newcommand{\sL}{\ensuremath{\mathscr{L}}\xspace}
\newcommand{\sM}{\ensuremath{\mathscr{M}}\xspace}

\newcommand{\sS}{\ensuremath{\mathscr{S}}\xspace}
\newcommand{\sT}{\ensuremath{\mathscr{T}}\xspace}

\newcommand{\sX}{\ensuremath{\mathscr{X}}\xspace}
\newcommand{\sY}{\ensuremath{\mathscr{Y}}\xspace}
\newcommand{\sZ}{\ensuremath{\mathscr{Z}}\xspace}

\usepackage{color}

\newcommand{\F}{\mathbf{F}}

\newcommand{\CC}{\mathbf{C}}
\newcommand{\G}{\mathbf{G}}
\newcommand{\tr}[0]{\operatorname{tr}}
\newcommand{\wt}[1]{\widetilde{#1}}

\newcommand{\Q}{\mathbf{Q}}
\newcommand{\Z}{\mathbf{Z}}

\newcommand{\mf}[1]{\mathfrak{#1}}

\newcommand{\Gal}{\operatorname{Gal}}

\newcommand{\ul}[1]{\underline{#1}}
\newcommand{\ol}[1]{\overline{#1}}
\newcommand{\wh}[1]{\widehat{#1}}

\newcommand{\Cal}[1]{\mathcal{#1}}
\newcommand{\A}{\mathbf{A}}

\newcommand{\ft}{{}^{\tau}} 
\newcommand{\co}{\colon}
\newcommand{\mrm}[1]{\mathrm{#1}}

\newcommand{\bs}{\backslash}

\newcommand{\PP}{\mathbf{P}}
\newcommand{\scr}[1]{\mathscr{#1}}

\newcommand{\inj}{\hookrightarrow}
\newcommand{\surj}{\twoheadrightarrow}

\DeclareMathOperator{\GL}{GL}

\DeclareMathOperator{\Frob}{Frob}

\DeclareMathOperator{\coker}{coker}

\DeclareMathOperator{\N}{\mathbf{N}}
\DeclareMathOperator{\Tr}{Tr}

\DeclareMathOperator{\Hom}{Hom}
\DeclareMathOperator{\Ind}{Ind}
\DeclareMathOperator{\Int}{Int}

\DeclareMathOperator{\rank}{rank}

\DeclareMathOperator{\Aut}{Aut}

\DeclareMathOperator{\Nm}{Nm}
\DeclareMathOperator{\Spec}{Spec\,}

\DeclareMathOperator{\Lie}{Lie}
\DeclareMathOperator{\End}{End}

\DeclareMathOperator{\Res}{Res}

\DeclareMathOperator{\Bun}{Bun}
\DeclareMathOperator{\Ext}{Ext}
\DeclareMathOperator{\Pic}{Pic}

\DeclareMathOperator{\Id}{Id}

\DeclareMathOperator{\Sym}{Sym}

\DeclareMathOperator{\pt}{pt}

\DeclareMathOperator{\Sht}{Sht}

\DeclareMathOperator{\pr}{pr}

\DeclareMathOperator{\Hk}{Hk}

\DeclareMathOperator{\Eis}{Eis}
\DeclareMathOperator{\Ch}{Ch}

\DeclareMathOperator{\Std}{Std}
\DeclareMathOperator{\Herm}{Herm}
\DeclareMathOperator{\Lagr}{Lagr}

\DeclareMathOperator{\Den}{Den}
\DeclareMathOperator{\Map}{Map}

\DeclareMathOperator{\len}{len}

\DeclareMathOperator{\Vect}{V}
\DeclareMathOperator{\Prym}{Prym}
\DeclareMathOperator{\all}{all}
\DeclareMathOperator{\RHom}{RHom}

\DeclareMathOperator{\Perf}{Perf}

\DeclareMathOperator{\Tot}{Tot}
\DeclareMathOperator{\taut}{taut}

\DeclareMathOperator{\ns}{ns}

\newcommand\old{\textup{old}}

\newcommand{\sSect}{\mathscr{S}\hspace{-.05cm}ect}
\newcommand{\sMap}{\mathscr{M}\hspace{-.05cm}ap}

\newcommand{\sHk}{\mathscr{H}\hspace{-.05cm}k}

\newcommand{\bL}{\mathbf{L}}
\newcommand{\bT}{\mathbf{T}}

\DeclareMathOperator{\BM}{BM}
\newcommand{\mBM}{\mrm{H}^{\BM}}

\DeclareMathOperator{\Spr}{Spr}

\DeclareMathOperator{\Lang}{Lang}
\DeclareMathOperator{\AJ}{AJ}

\DeclareMathOperator{\Sect}{Sect}
\DeclareMathOperator{\Irr}{Irr}
\newcommand\ups{\upsilon}
\renewcommand\div{\textup{div}}

\newcommand{\incl}{\hookrightarrow}
\newcommand{\isom}{\stackrel{\sim}{\to}}

\newcommand{\Qlbar}{\overline{\Q}_\ell}

\renewcommand{\j}[1]{\langle{#1}\rangle}
\newcommand\un{\underline}
\newcommand{\bu}{\bullet}
\newcommand{\ov}{\overline}

\newcommand\sss{\subsubsection}
\newcommand\xr{\xrightarrow}
\newcommand\op{\oplus}
\newcommand\ot{\otimes}

\newcommand\one{\mathbf{1}}

\renewcommand\c{\circ}
\newcommand\vn{\varnothing}
\newcommand\upH{\textup{H}}
\newcommand{\cohog}[2]{\textup{H}^{#1}({#2})}     

\newcommand{\cHom}{\Cal{H}om}
\newcommand{\cExt}{\Cal{E}xt}
\newcommand{\cEnd}{\Cal{E}nd}
\newcommand{\cRHom}{\Cal{R}\Cal{H}om}
\newcommand{\cHerm}{\Cal{H}erm}
\newcommand{\cAut}{\Cal{A}ut}

\renewcommand\a\alpha
\renewcommand\b\beta
\newcommand\g\gamma
\renewcommand\d\delta
\newcommand\D\Delta
\newcommand{\e}{\epsilon}
\newcommand{\io}{\iota}
\renewcommand{\th}{\theta}
\newcommand{\Th}{\Theta}

\newcommand{\s}{\sigma}
\renewcommand{\t}{\tau}

\newcommand{\y}{\eta}
\newcommand{\z}{\zeta}
\newcommand{\ep}{\epsilon}
\newcommand{\vp}{\varpi}
\renewcommand{\l}{\lambda}
\renewcommand{\L}{\Lambda}
\newcommand{\om}{\omega}

\renewcommand{\r}{\rho}

\newcommand{\nai}{\mrm{naive}}

\newcommand\hs{\heartsuit}

\newcommand\sh{\sharp}

\newcommand\cA{\mathcal{A}}

\newcommand\cE{\mathcal{E}}
\newcommand\cF{\mathcal{F}}
\newcommand\cG{\mathcal{G}}
\newcommand\cH{\mathcal{H}}

\newcommand\cK{\mathcal{K}}
\newcommand\cL{\mathcal{L}}
\newcommand\cM{\mathcal{M}}

\newcommand\cO{\mathcal{O}}
\newcommand\cP{\mathcal{P}}
\newcommand\cQ{\mathcal{Q}}

\newcommand\cS{\mathcal{S}}
\newcommand\cT{\mathcal{T}}

\newcommand\cV{\mathcal{V}}

\newcommand\cX{\mathcal{X}}
\newcommand\cY{\mathcal{Y}}
\newcommand\cZ{\mathcal{Z}}

\newcommand\frD{\mathfrak{D}}

\newcommand\frL{\mathfrak{L}}
\newcommand\frM{\mathfrak{M}}
\newcommand\frN{\mathfrak{N}}

\newcommand\frZ{\mathfrak{Z}}

\newtheorem{thm}{Theorem}[section]
\newtheorem{lemma}[thm]{Lemma}
\newtheorem{prop}[thm]{Proposition}
\newtheorem{cor}[thm]{Corollary}

\newtheorem{conj}[thm]{Conjecture}
\newtheorem{defn-prop}[thm]{Definition-Proposition}

\theoremstyle{remark}
\newtheorem{remark}[thm]{Remark} 
\newtheorem{defn}[thm]{Definition}

\newtheorem{example}[thm]{Example}

\newtheorem{warning}[thm]{Warning}

\makeatletter
\def\th@remark{%
  \thm@headfont{\bfseries}%
  \normalfont 
  \thm@preskip \thm@preskip 
  \thm@postskip\thm@preskip
}
\def\imod#1{\allowbreak\mkern5mu({\operator@font mod}\,\,#1)}
\makeatother

\numberwithin{equation}{section}

\title[Higher theta series for unitary groups over function fields]{Higher theta series for unitary groups \\ over function fields \\ 
\vspace{.2cm}  S\'{e}ries th\^{e}ta sup\'{e}rieures pour les groupes unitaires \\ sur les corps de fonctions}

\author{Tony Feng}
\address{University of California, Berkeley, Department of Mathematics, 970 Evans Hall, Berkeley, CA 94720, USA}
\email{fengt@berkeley.edu}

\author{Zhiwei Yun}
\address{Massachusetts Institute of Technology, Department of Mathematics, 77 Massachusetts Avenue, Cambridge, MA 02139, USA}
\email{zyun@mit.edu}

\author{Wei Zhang}

\address{Massachusetts Institute of Technology, Department of Mathematics, 77 Massachusetts Avenue, Cambridge, MA 02139, USA}
\email{weizhang@mit.edu}

\begin{document}

\begin{abstract}
In previous work, we defined certain virtual fundamental classes for special cycles on the moduli stack of Hermitian shtukas, and related them to the higher derivatives of non-singular Fourier coefficients of Siegel-Eisenstein series. In the present article, we construct virtual fundamental classes in greater generality, including those expected to relate to the higher derivatives of \emph{singular} Fourier coefficients. We assemble these classes into ``higher'' theta series, which we conjecture to be modular. Two types of evidence are presented: structural properties affirming that the cycle classes behave as conjectured under certain natural operations such as intersection products, and verification of modularity in several special situations. One innovation underlying these results is a new approach to special cycles in terms of derived algebraic geometry. \\

\vspace{.2cm} 
Dans des travaux pr\'{e}c\'{e}dents, nous avons d\'{e}fini certaines classes fondamentales virtuelles pour des cycles sp\'eciaux sur les champs de chtoucas hermitiens, et les avons li\'{e}es aux d\'{e}riv\'{e}es sup\'{e}rieures des coefficients de Fourier non singuliers des s\'{e}ries de Siegel-Eisenstein. Dans cet article, nous construisons des classes fondamentales virtuelles dans des contextes plus g\'{e}n\'{e}raux, y compris celles qui sont cens\'{e}es \^{e}tre li\'{e}es aux d\'{e}riv\'{e}es sup\'{e}rieures des coefficients de Fourier singuliers. Nous assemblons ces classes en des s\'{e}ries th\^{e}ta ``sup\'{e}rieures'', et conjecturons que ces s\'{e}ries th\^{e}tas sont modulaires. Deux types d'indications sont pr\'{e}sent\'es en faveur de cette conjecture : des propri\'{e}t\'{e}s structurelles qui affirment que ces classes de cycles se comportent conform\'ement \`a cette conjecture sous certaines op\'{e}rations naturelles (par exemple des produits d'intersection), et la v\'{e}rification de la modularit\'{e} dans quelques situations sp\'{e}ciales. Ces r\'esultats s'appuient sur une nouvelle approche des cycles sp\'eciaux en termes de g\'{e}om\'{e}trie alg\'{e}brique d\'{e}riv\'{e}e.
\end{abstract}

\maketitle

\tableofcontents


\section{Introduction}

The earliest examples of theta functions were generating series for the number of representations of integers by quadratic forms. It has been known at least since the work of Jacobi that theta functions enjoy remarkable symmetry properties, which later became known as \emph{modularity}, that underlie many of their applications. An incarnation of theta functions in arithmetic algebraic geometry was discovered by Kudla, who named them \emph{arithmetic theta series}. This paper is about modularity in the context of arithmetic theta series. 

The earliest examples of arithmetic theta series were constructed by Kudla as generating series with coefficients being cycle classes in the Chow groups of Shimura varieties \cite{Kud04}. Kudla envisioned a conjectural \emph{arithmetic Siegel--Weil formula} \cite{Kudla1997}, which would further require extending the special cycles to good integral models of Shimura varieties. A significant difficulty is the task of defining the appropriate cycle classes in the arithmetic Chow group indexed by singular Fourier coefficients. For example, for unitary Shimura varieties Kudla and Rapoport constructed the cycle classes on their integral models indexed by non-singular Fourier coefficients in \cite{KRI,KRII}, while Li and the third author \cite{LZ1} proved an arithmetic Siegel--Weil formula for the non-singular Fourier coefficients (see also \cite{LZ2} for the orthogonal analog). However, the definition of the singular terms, and therefore also the full arithmetic theta series, remains open (except in some lower dimensional case, see \cite{KRY}).


In \cite{FYZ} we proposed a function field analogue of this story: we defined special cycles on the moduli stack of Hermitian shtukas, constructed certain virtual fundamental classes for the cycles indexed by \emph{non-singular} Fourier coefficients, and related them to the Taylor expansion of Fourier coefficients of corresponding Siegel-Eisenstein series. A novel feature of the function field version is that cycle classes can be defined for each non-negative integer $r$, and related to the $r^{\mrm{th}}$ derivative of the Fourier coefficients of Siegel-Eisenstein series, whereas only the cases $r=0$ and $r=1$ seem to be witnessed over number fields (at least for the time being). 

In this paper, we will construct virtual fundamental classes in general, going beyond the non-singular cases considered in \cite{FYZ}, and assemble them into full ``higher'' arithmetic theta series (so named because they are related to higher derivatives of Siegel-Eisenstein series). The form of the singular terms exhibits interesting complexities that will be discussed further in \S \ref{ssec: intro singular terms}. We formulate a conjecture about the modularity of such theta series, and then give evidence for this conjecture. 

\subsection{The modularity conjecture}\label{ssec: intro singular terms}

We now introduce notation so as to be able to describe our conjecture and the main results with more precision. Let $X$ be a smooth, proper and geometrically connected curve over $k=\F_q$ of characteristic $p\ne2$, and let $\nu \co X' \rightarrow X$ be a connected \'{e}tale double cover, with the non-trivial automorphism denoted $\sigma \in \Aut(X'/X)$. Let $F$ be the function field of $X$ and let $F'$ be the function field of $X'$. In \cite{FYZ} we defined the moduli stack $\Sht_{U(n)}^r$ parametrizing rank $n$ ``Hermitian shtukas'' with $r$ legs. We also defined certain special cycles $\Cal{Z}_{\cE}^r(a)$ indexed by $\cE$, a vector bundle of rank $m$ with $1\leq m\leq n$ on $X'$, and a Hermitian map $a \co \cE \rightarrow \sigma^* \cE^{\vee}$ where $\cE^{\vee}:= \cHom(\cE, \omega_{X'})$ is the Serre dual of $\cE$. The space of such $a$ was called $\cA_{\cE}^{\all}(k)$ in \cite{FYZ}, but is called $\cA_{\cE}(k)$ in this paper. (Everything in \cite{FYZ} works in a slightly more general setup allowing a similitude factor, but for simplicity we omit this from our introduction.) 


To define the higher theta series, we construct an appropriate virtual fundamental class $[\Cal{Z}_{\cE}^r(a)]\in\Ch_{r(n-m)}(\Cal{Z}_{\cE}^r(a))$ for {\em every} $a\in \cA_{\cE}(k)$. 

This was done in \cite{FYZ} when $a$ is {\em non-singular} (meaning that $a \co \cE \rightarrow \sigma^*  \cE^{\vee}$ is injective as a map of coherent sheaves) and either $\rank \cE = n$ or $\cE$ is a direct sum of line bundles, by taking derived intersections from the situation where $\rank \cE = 1$, following the ideas of \cite{KRII} in the number field case. 
However, even in the non-singular case, to handle general $m$ and $\cE$ we must take a new approach based on Hitchin stacks (Definition \ref{defn: circ case}). The dissimilarity to the number field situation comes from the fact that not every vector bundle on a proper curve splits as a sum of line bundles, while every vector bundle over the ring of integers of a number field splits as a direct sum of line bundles.

For singular $a$, the construction of $[\Cal{Z}_{\cE}^r(a)]$ is more complicated. The cycle $\Cal{Z}_{\cE}^r(a)$ admits an open-closed decomposition according to the possible kernels of the map $a$, and the contribution from each stratum is the product of a virtual fundamental class constructed from a Hitchin stack and the top Chern class of a certain tautological bundle. The construction is completed in Definition \ref{def:special cycle classes}. It may be a useful guide for the number field case, where no definition of special cycle classes in the arithmetic Chow group is currently known, for singular Fourier coefficients, at the time of this writing. 


Having defined $[\cZ_{\cE}^r(a)]$ for each $a$, we then assemble them into higher theta series. More precisely, if $\rank \cE = m$, then we consider the quasi-split unitary group (with respect to the double cover $X'/X$) of rank $2m$ over $X$, abbreviated $U(2m)$, and the standard Siegel parabolic $P_m$. (In the main body of the paper, starting in \S \ref{ss:Eis}, we use the notation $H_m$ for $U(2m)$.) We write down a function on $U(2m)(\BA)$ valued in $\Ch_{r(n-m)}(\Sht_{U(n)}^r)$:
\begin{align*}
\xymatrix{\wt Z^{r}_{m}:U(2m)(\BA)\ar[r]& \Ch_{r(n-m)}(\Sht^{r}_{U(n)})}
\end{align*}
characterized by the following properties:
\begin{enumerate}
\item $\wt Z^{r}_{m}$ is left invariant under $P_{m}(F)$ and right invariant under $K=U(2m)(\wh\cO)$;  
\item for any point in $P_m(F)\bs P_m(\BA)/K\cap P_m(\BA)\simeq P_m(F)\bs U(2m)(\BA)/K$ represented by $(\cG,\cE)$, where $\cG$ is a rank $2m$ vector bundle on $X'$ with a skew-Hermitian structure $h:\cG\simeq \s^* \cG^*$ and $\cE$ is a Lagrangian sub-bundle of $\cG$, we have a ``Fourier expansion"  (in the sense of \cite[\S 2.6]{FYZ})
\begin{equation}\label{eq: intro generating series} \wt Z^{r}_{m}(\cG,\cE) =\chi(\det\cE)q^{n(\deg \cE-\deg\om_{X}) /2}\sum_{a\in \cA_{\cE}(k)}\psi_{0}(\j{e_{\cG,\cE}, a})\z_{*}[\cZ^{ r}_{\cE}(a)].
\end{equation}
 Here we refer to \S\ref{ss:mod} for the undefined notation in the right hand side. We note that, in the special case $\cE=\cO_{X'}^{\oplus m}$ the trivial bundle of rank $m$,  the set of all such $(\cG,\cE)$ is naturally isomorphic to $N_m(F)\bs N_m(\BA)/K\cap N_m(\BA)$, where $N_m$ denotes the unipotent radical of $P_m$.  Then $\cA_\cE(k)$ is naturally isomorphic to the Pontryagin dual of $N_m(F)\bs N_m(\BA)/K\cap N_m(\BA)$ (depending on the choice of a non-trivial character $\psi_0:k\to\BC^\times$). For this $\cE$, \eqref{eq: intro generating series} more closely resembles the expressions for arithmetic theta series on Shimura varieties, as one finds for example in \cite[(5.4)]{Kud04}. 

\end{enumerate}

\begin{conj}[Modularity conjecture]The function $\wt Z^{r}_{m}$ descends to a function
\begin{align*}
\xymatrix{ Z^{r}_{m}:U(2m)(F)\bs U(2m)(\BA)\ar[r]& \Ch_{r(n-m)}(\Sht^{r}_{U(n)}),}
\end{align*}
i.e., $\wt Z^{r}_{m}$ is left $ U(2m)(F)$-invariant. 
\end{conj}

In other words, the  class $\wt Z^{r}_{m}(\cG,\cE)\in \Ch_{r(n-m)}(\Sht^{r}_{U(n)})$ should depend only on the Hermitian bundle $\cG$ and not on its Lagrangian sub-bundle $\cE$.

When $r=0$, $\Ch_{0}(\Sht^{0}_{U(n)})$ is simply the space of $\Q$-valued functions on $\Bun_{U(n)}(k)$ and the evaluation map turns $\wt Z^{r}_{m}$ into a two-variable function
$$\xymatrix{U(2m)(\BA)\times U(n)(\BA) \ar[r]& \Q.}
$$In this case, we obtain the classical theta function and the modularity conjecture essentially follows from the Poisson summation formula.

\begin{remark}A conjecture can also be formulated in the case $X'=X\coprod X$. The special cycles then live on the more familiar moduli stack of $\GL(n)$-shtukas, and we refer to  \S \ref{ssec: split case conj} for the details.
\end{remark}


\subsection{Main results}

Our main results give some evidence towards the modularity conjecture.

One type of evidence, considered in Part III, is of numerical nature: we prove modularity of the functions obtained by intersecting our arithmetic series with classes analogous to what would be called {\em CM (Complex Multiplication) cycles} for unitary Shimura varieties. In particular, this entails proving the modularity of our arithmetic series for rank 1 unitary groups. 

A second type of evidence, studied in Part II, concerns more abstract ``coherence properties'' of the special cycles. For example, we prove that the product of special cycle classes in the Chow ring behaves as predicted in \cite{Kud04}. Perhaps surprisingly, the proofs rely crucially on the methods of \emph{derived algebraic geometry}, and in particular on a construction of \emph{derived special cycles} which yield our virtual fundamental classes. This will be discussed more in \S \ref{ssec: intro DAG}. This is a novelty of the singular terms; derived algebraic geometry has not played a role so far in studying the non-singular terms. It leads us to suspect that derived algebraic geometry may also prove useful in the more classical Shimura variety context of the Kudla program.

\subsubsection{Linear invariance}

We establish compatibility properties of the special cycles under various natural operations. Here we state an example (Theorem \ref{thm: linear invariance}), which we call the \emph{linear invariance} following the analog in the number field case considered by Howard in \cite{How12}.
\begin{thm}
Given a decomposition $\cE \approx  \cE_1 \oplus \cE_2 \oplus \ldots  \oplus \cE_j$, and $a_i \in \cA_{\cE_i}(k)$, the intersection product 
\[
 [\cZ_{\cE_1}^r(a_{1})] \cdot_{\Sht_{U(n)}^r} [\cZ_{\cE_2}^r(a_{2})]  \cdot_{\Sht_{U(n)}^r}  \cdots  \cdot_{\Sht_{U(n)}^r}  [\cZ_{\cE_j}^r(a_{j})]  
\]
coincides with the sum of $[\cZ_{\cE}^r(a)]$ over all $a \co \cE \rightarrow \sigma^*\cE^{\vee}$ satisfying the condition that
\begin{equation}
\text{ the composition $\cE_i \rightarrow \cE \xrightarrow{a} \sigma^* \cE^{\vee} \rightarrow \sigma^*  \cE_i^{\vee}$ is $a_i$ for each $1 \leq i \leq j$.}
\end{equation}\end{thm}

 Although in principle both sides of the equality may be expressed in terms of elementary constructions, our proof relies on the derived algebraic geometry interpretation of the special cycle classes and we do not know a proof without derived methods.

\subsubsection{A refinement of the main result of \cite{FYZ}} 

The stack $\Sht^r_{U(n)}$ is a disjoint union of two open-closed substacks and the modularity conjecture predicts that the restriction of the generating series $\wt Z_{m}^r$ to each of them should also be modular. In \cite{FYZ} we identified the degree of the $[\cZ_\cE^r(a)]\in\Ch_0(\cZ_\cE^r(a))$ for non-singular $a$ with the $r^{\mrm{th}}$ central derivative of the (suitably normalized) $a^{\mrm{th}}$ Fourier coefficient of Siegel--Eisenstein series. In Theorem \ref{thm: deg equal on components}, we refine this result and show that the restriction of $[\cZ_\cE^r(a)]$  for non-singular $a$ to each of the two open-closed substacks has equal degree. The proof turns out to be non-elementary.

\subsubsection{The case $n=m=1$} 
\begin{thm}
The modularity conjecture  holds when $n=m=1$.
\end{thm}
 In this case the higher theta series are valued in the Chow group of proper zero-cycles, and are therefore essentially determined by their degrees. We show that the degrees are given by explicit automorphic functions, namely higher derivatives of a suitably normalized Eisenstein series. In fact this was already established for non-singular Fourier coefficients in \cite{FYZ}, so the remaining work is to calculate the singular term, which turns out to be the Chern class of a certain tautological bundle, and to relate it to the Taylor expansion of the corresponding $L$-function. This computation is carried out in \S \ref{sec: constant term}. Analogous results over number fields (for $r=1$) were established by Kudla-Rapoport-Yang \cite{KRY99}.

\subsubsection{Intersection with ``CM cycles"}
For $\theta \co Y \rightarrow X$ a degree $n$ cover (possibly ramified), we have a ``CM cycle'' $\Sht_{U(1)_Y}^r$ of dimension $r$ and a finite morphism  $\Th:  \Sht^r_{U(1)_{Y}}\to \Sht_{U(n)}^r$. We consider the intersection number of the resulting cycle class $ \Th_*[\Sht^r_{U(1)_{Y}} ] $ with the arithmetic theta series $\wt Z^{r}_{m=1}(g)$ in codimension $r$ (i.e. the generating series of  corank $m=1$ special cycles), for $g \in U(2)(\BA)$. 

\begin{thm}
For any $n$, the function $U(2)(\BA)  \ni g\mapsto  \Big\langle \wt Z^{r}_{m=1}(g) , \Th_*[\Sht^r_{U(1)_{Y}} ] \Big\rangle\in\BC$ is left invariant under $U(2)(F)$.
\end{thm}

In fact, we can identify the intersection number with the $r^{\mrm{th}}$ derivative of an explicit Eisenstein series. For the non-singular terms, this could be thought of as proving a higher-derivative, function-field analogue of \cite{How12}. For the singular terms, it could be thought of as a higher-derivative, function-field analogue of the proof of the ``averaged Colmez conjecture'' in \cite{AGHMP} (also obtained by a different method in \cite{YZ18}).

One reason that we are limited to the corank $m=1$ case is that, in order to intersect with $[\cZ_\cE^r(a)]$ in the corank $m>1$ case,  we need to construct  natural proper cycles of higher (than $r$) dimension on $\Sht_{U(n)}^r$. Some candidates are given by the analogs of basic loci on unitary Shimura varieties, which may reach nearly (but nevertheless strictly smaller than) half of the dimension of $\Sht_{U(n)}^r$. If we demand proper cycles that are surjective to the base $X'^r$, then we only know how to construct examples of dimension $r$ but not higher, see Example \ref{ex:crk n-1}  and Example \ref{ex:CM}.

\subsubsection{Geometric properties of special cycles} In \S \ref{ssec: corank one unitary}, we study the geometric properties of the special cycles $\cZ_{\cE}^r(a)$ in the special case where $\rank \cE = m = 1$. We show that if $a$ is non-singular then it is LCI of the correct dimension, and that the virtual fundamental class $[\cZ_{\cE}^r(a)]$ coincides with the naive fundamental class. This fulfills a result promised in \cite[Remark 7.10]{FYZ}, and allows us to prove that the general constructions of cycle classes considered in this paper recovers the more naive definitions studied in \cite{FYZ}.

\subsection{Some remarks on the derived algebro-geometric method}\label{ssec: intro DAG}

Although we are able to give an explicit formula for the special cycle classes in Part I using only ``classical" algebraic geometry, the key foundation for the structural results proved in Part II is another interpretation of these classes in terms of \emph{derived algebraic geometry}. We emphasize that the formulation of the modularity conjecture itself requires no input from derived algebraic geometry, while the evidence does.

To summarize, in \S \ref{sec: derived} we define derived enhancements of the special cycles and show (Theorem \ref{thm: VFC from derived shtuka}) that their intrinsic derived fundamental cycles coincide with the virtual classes defined earlier. One advantage of this approach is that it does not involve separating the non-singular and singular cases, and so gives a uniform, conceptual derivation of the virtual fundamental classes for special cycles indexed by all Fourier coefficients. We find this to be compelling philosophical evidence for our definition of the singular terms.

Let us elaborate on the role of derived algebraic geometry. A derived scheme/stack has an underlying classical scheme/stack which we call its \emph{classical truncation}, and in this sense the derived object can be thought of as ``enhancing'' the classical object with some kind of ``derived structure''. For example, a \emph{quasi-smooth} (i.e., derived analogue of LCI) derived scheme provides a ``perfect obstruction theory'', in the sense of Behrend-Fantechi, for its classical truncation. Now, the process of classical truncation can lose good geometric properties; for example, \emph{any} (arbitrarily singular) finite type affine scheme can arise as the classical truncation of a derived scheme which is quasi-smooth. The ``hidden smoothness'' philosophy of Deligne, Drinfeld, and Kontsevich holds that many naturally occurring singular moduli spaces are the classical truncations of natural quasi-smooth derived moduli spaces\footnote{In modern terms, ``hidden quasi-smoothness'' would be a more accurate name for this philosophy. As far as we know, the name ``quasi-smooth'' is due to Lurie.}, and this was one of the early motivations to consider derived algebraic geometry. 

In fact, it has been understood since the seminal work of Kudla-Rapoport \cite{KRII} that the special cycles comprising arithmetic theta series need to be defined in a ``derived'' way. The physical special cycles are often not even of the ``correct'' dimension, and may be quite singular, so instead of considering their naive fundamental classes one wants to construct \emph{virtual} fundamental classes. Kudla-Rapoport did this for the non-singular terms on unitary Shimura varieties, by presenting the cycles as a ``derived intersection'' of classical schemes with the correct ``expected dimension''. Then the virtual fundamental class was defined as a refined intersection product in Fulton's sense. Our construction of the non-singular terms on Hermitian shtukas also fits this mold. 

For singular terms, we do not know of a presentation that may be used to carry out a similar strategy. What we shall see, however, is that all special cycles (even for singular coefficients) can be promoted to derived stacks in a natural way, which always have the correct dimension in the derived sense, and are always quasi-smooth. This gives another example of the ``hidden smoothness'' philosophy. Moreover, quasi-smooth derived stacks have an intrinsic notion of fundamental class, which can be viewed as a virtual fundamental class of the underlying classical stack. This gives an intrinsic construction of a virtual fundamental class to each special cycle, which is uniform with respect to the Fourier coefficient (whether singular or not).

From this perspective, the reason that cycles indexed by non-singular Fourier coefficients can be defined more easily is that the derived structure on such cycles can be constructed in an elementary way, by taking the derived intersection of classical stacks. We do not know of such an elementary construction for singular coefficients, nor is it necessary for us. This suggests that derived algebraic geometry may also be relevant for the classical Kudla program (over number fields), where the cycles indexed by singular coefficients had previously been defined in a more ad hoc manner. However, the methods we use to construct the derived special cycles do not have an obvious analogue in the number field situation.

\subsection*{Acknowledgment}

We thank Adeel Khan for discussions on derived intersection theory. We thank Chao Li and Ben Howard for their comments on our draft. TF was supported by an NSF Postdoctoral Fellowship under grant No. 1902927, as well as the Friends of the Institute for Advanced Study. ZY was supported by a Packard Fellowship and a Simons Investigator grant. WZ was supported by the NSF grant DMS \#1901642 and a Simons Investigator grant. 


\subsection{Notation}\label{ssec: notation}

Throughout this paper, $k=\F_{q}$ is a finite field of odd characteristic $p$. Let $\ell\ne p$ be a prime.   Let $\psi_{0}: k\to\Qlbar^{\times}$ be a nontrivial character. For a space $S$ over $\F_q$, we denote by $\Frob_S$ the $q$-th Frobenius endomorphism of $S$; and sometimes we omit the subscript $S$ when it is clear from the context.

\subsubsection{} Let $X$ denote a smooth, projective, geometrically connected curve over $k$, of genus $g_X$. Let $\om_{X}$ be the line bundle of $1$-forms on $X$. 

Let $F=k(X)$ denote the function field of $X$.  Let $|X|$ be the set of closed points of $X$. For $v\in |X|$, let $\cO_{v}$ be the completed local ring of $X$ at $v$ with fraction field $F_{v}$ and residue field $k_{v}$.   Let $\BA=\BA_{F}$ denote the ring of ad\`eles of $F$, and $\wh\cO=\prod_{v\in |X|}\cO_{v}$. Let $\deg(v)=[k_{v}:k]$, and $q_{v}=q^{\deg(v)}=\#k_{v}$. Let $|\cdot|_{v}: F^{\times}_{v}\to q^{\Z}_{v}$ be the absolute value such that  $|\vp_{v}|_{v}=q^{-1}_{v}$ for any uniformizer  $\vp_{v}$ of $\cO_{v}$. Let $|\cdot|_{F}: \BA_{F}^{\times}\to q^{\Z}$ be the absolute value that is $|\cdot|_{v}$ on $F_{v}^{\times}$. 

\subsubsection{} Let $X'$ be another smooth curve over $k$ and let $\nu:X'\to X$ be a finite \'{e}tale map of degree $2$. We denote by $\sigma$ the non-trivial automorphism of $X'$ over $X$. The case where $X'$ is geometrically disconnected is allowed unless stated otherwise; it is usually allowed in Parts 1 and 2 but not in Part 3. Let $F'$ be the ring of rational functions on $X'$, which is either a quadratic extension of $F$ or $F\times F$. We let $k'$ be the ring of constants in $F'$, which may be $\F_q$, $\F_{q^2}$ or $\F_q \times \F_q$. The notations $\om_{X'}, |X'|, F'_{v'}, \cO_{v'}, k_{v'}, \BA_{F'}, |\cdot|_{v'}, |\cdot|_{F'}, q_{v'}$ and $\deg(v')$ (for $v'\in |X'|$) are defined similarly as their counterparts for $X$. Additionally, for $v\in |X|$, we use $\cO'_{v}$ to denote the completion of $\cO_{X'}$ along $\nu^{-1}(v)$, and define $F'_{v}$ to be its total ring of fractions.

For a vector bundle $\cF$ on $X'_S$, we denote $\ft \cF := (\Id_{X'} \times \Frob_S)^* \cF$ for its $\Frob_S$-twist. 

\subsubsection{Notation for cycle classes} For a (derivd) stack $\cY$, $\Ch(\cY)$ denotes its \emph{rational} Chow group in the sense of \cite{KhanI}. We denote by $[\cY]^{\nai} \in \Ch(\cY)$ the fundamental class of $\cY$. Typically we will work with ``virtual fundamental classes'' in $\Ch(\cY)$ which do not (at least a priori) coincide with the naive ones, and we shall denote such by $[\cY] \in \Ch(\cY)$, although they will in fact depend on some auxiliary construction, such as a realization of $\cY$ as a fibered product or as the classical truncation of a derived stack $\sY$.

\subsubsection{Derived notation} In \S \ref{sec: derived} -- \S \ref{sec: compatibility}, we adopt some notational conventions that differ from the rest of the paper. Namely, in those sections we operate within $\infty$-categories, so fibered products mean ``derived fibered products'', limits mean ``homotopy limits'', etc. unless noted otherwise. We refer to \S \ref{sec: derived} for the precise explanation of the notation used in those sections.

\subsubsection{Some notational departures from \cite{FYZ}}  We emphasize that some notation has changed from our first paper \cite{FYZ} regarding Hitchin spaces and Hitchin bases. There we introduced certain Hitchin stacks $\cM \subset \cM^{\all}$ and Hitchin bases $\cA \subset \cA^{\all}$, decorated by indices, but in this paper they would be denoted $\cM^{\circ} \subset \cM$ and $\cA^{\ns} \subset \cA$. This will be explained more precisely when it comes up in the text. 

\part{Formulation of the conjecture}

\section{Some (more) special cycles on moduli of shtukas}

In this section we introduce a variant and a generalization of the special cycles defined in \cite{FYZ}. The variant, which plays a technical role in later definitions and proofs, is obtained by replacing $U(n)$ with the general linear group. For the generalization of special cycles, we consider Hermitian shtukas with similitude line bundles. Later we will formulate the modularity conjecture in this generality.

\subsection{Shtukas for $\GL(n)'$}\label{ssec: shtukas for GL(n)'}

We denote $\GL(n)' := \Res_{X'/X} \GL(n)$, a group scheme over $X$. In this subsection we define stacks $\Sht_{\GL(n)'}^r$ parametrizing certain special types of shtukas for $\GL(n)'$, and establish their basic geometric properties. Their role in the study of Hermitian shtukas is of a somewhat technical nature, stemming from the fact that the Hitchin spaces corresponding to $\GL(n)'$ have better technical properties. They appear in an intermediate stage in the construction of cycle classes labeled by singular Fourier coefficients. 

We begin by explicating the appropriate notion of bundles and Hecke correspondences. Let $\Bun_{\GL(n)'}$ be the moduli stack of $\GL(n)'$-bundles on $X$. By general properties of Weil restriction, there is a canonical equivalence of groupoids
\[
\{\text{$\GL(n)'$-bundles on $X \times S$}\} \cong  \{ \text{$\GL(n)$-bundles on $X' \times S$}\}.
\]
Hence the datum of a $\GL(n)'$-bundle on $X \times S$ is equivalent to the datum of a rank $n$ vector bundle on $X' \times S$, and $\Bun_{\GL(n)'}$ is simply equivalent to the moduli stack of $\GL(n)$-bundles on $X'$.

\begin{defn}\label{defn: Hk GL(n)'}
Let $r \geq 0$ be an integer. The \emph{Hecke stack} $\Hk_{\GL(n)'}^{r}$ has as $S$-points the groupoid of the following data: 
\begin{enumerate}
\item $x_i' \in X'(S)$ for $i = 1, \ldots, r$, with graphs denoted $\Gamma_{x_i'} \subset X' \times S$.
\item A sequence of vector bundles $\Cal{F}_0, \ldots, \Cal{F}_r$ of rank $n$ on $X' \times S$.
\item Isomorphisms $f_i \co \Cal{F}_{i-1}|_{X' \times S - \Gamma_{x'_i} - \Gamma_{\sigma x'_i}} \xrightarrow{\sim} \Cal{F}_i|_{X' \times S - \Gamma_{x'_i}-\Gamma_{\sigma x'_i}}$, for $1 \leq i \leq r$, which are lower of length $1$ at $x'_i$ and upper of length $1$ at $\sigma x'_i$, in the terminology of \cite[Definition 6.5]{FYZ}.
\end{enumerate}
\end{defn}

\begin{warning}\label{warning: Hk}
The stack $\Hk_{\GL(n)'}^{r}$ is different from the usual iterated Hecke stack for rank $n$ vector bundles on $X'$, for example as considered for $n=2$ in \cite{YZ}, because we have demanded modifications to occur over conjugate pairs of points on the curve.
\end{warning}

\begin{lemma}\label{lem: smooth Bun_G}
The (Artin) stack $\Bun_{\GL(n)'}$ is smooth. 
\end{lemma}

\begin{proof}This follows from the standard obstruction theory argument, cf. \cite[Prop. 1]{Hei10}.
\end{proof}



\begin{defn}
Let $r \geq 0$ be an integer. We define $\Sht_{\GL(n)'}^{r}$ by the Cartesian diagram
\[
\begin{tikzcd}
\Sht_{\GL(n)'}^{r} \ar[r] \ar[d] & \Hk_{\GL(n)'}^{r} \ar[d] \\
\Bun_{\GL(n)'} \ar[r, "{(\Id, \Frob)}"]  & \Bun_{\GL(n)'} \times \Bun_{\GL(n)'}
\end{tikzcd}
\]
A point of $\Sht_{\GL(n)'}^{r}$ will be called a ``$\GL(n)'$-shtuka''. (But see Warning \ref{warning: G-shtukas}.) 

Concretely, the $S$-points of $\Sht_{\GL(n)'}^{r}$ are given by the groupoid of the following data: 
\begin{enumerate}
\item $x_i' \in X'(S)$ for $i = 1, \ldots, r$, with graphs denoted $\Gamma_{x_i'} \subset X \times S$. 
\item A sequence of vector bundles $\Cal{F}_0, \ldots, \Cal{F}_r$ of rank $n$ on $X' \times S$.
\item Isomorphisms $f_i \co \Cal{F}_{i-1}|_{X' \times S - \Gamma_{x'_i}-\Gamma_{\sigma x'_i}} \xrightarrow{\sim} \Cal{F}_i|_{X' \times S - \Gamma_{x'_i}-\Gamma_{\sigma x'_i}}$, which are lower of length $1$ at $x'_i$ and upper of length $1$ at $\sigma x'_i$.  

\item An isomorphism of vector bundles $\varphi \co \Cal{F}_r \cong \ft \Cal{F}_0=(\Id_{X'}\times \Frob_{S})^{*}\cF_{0}$.
\end{enumerate}
\end{defn}

\begin{warning}\label{warning: G-shtukas}
For the same reason as Warning \ref{warning: Hk}, the stack $\Sht_{\GL(n)'}^{r}$ is different from the usual iterated stack of rank $n$ shtukas on $X'$, for example as considered for $n=2$ in \cite{YZ}. 
\end{warning}

\begin{lemma}\label{lem: sht dimension} (1) The projection map $(\pr_{X}, \pr_{r}):\Hk_{\GL(n)'}^r \rightarrow (X')^r \times \Bun_{\GL(n)'}$ recording $\{x_{i}\}_{i=1}^r$ and $\cF_{r}$ is smooth of relative dimension $2r(n-1)$. 

(2) $\Sht_{\GL(n)'}^r$ is a smooth Deligne-Mumford stack, locally of finite type, and separated, of pure dimension $r(2n-1)$. 
\end{lemma}

\begin{proof}The proof of (1) is similar to the proof of \cite[Lemma 6.9(1)]{FYZ}, except that in the case $r=1$ the upper and lower modifications are independent, so the map $\Hk_{\GL(n)'}^1 \rightarrow X' \times \Bun_{\GL(n)'}$ is (\'{e}tale locally on target) a $\PP^{n-1}$-fibration over a $\PP^{n-1}$-fibration.

Part (2) follows from (1) by applying \cite[Lemma 2.13]{Laff18} in the analogous way as in \cite[Proposition 2.11]{Laff18}. 
\end{proof}

\subsection{Special cycles} We will define some special cycles on $\Sht_{\GL(n)'}^r$.

\begin{defn}\label{def: Z} Let $\Cal{E}$ be a rank $m$ vector bundle on $X'$.

 We define the stack $\Cal{Z}_{\Cal{E},\GL(n)'}^{r}$ whose $S$-points are the groupoid of the following data: 
 \begin{itemize}
 \item A $\GL(n)'$-shtuka $(\{x'_1, \ldots, x'_r\}, \{\Cal{F}_0, \ldots, \Cal{F}_r\}, \{f_1, \ldots, f_r\}, \varphi) \in  \Sht_{\GL(n)'}^{r}(S)$.
 \item Maps of coherent sheaves  $t_i \co \Cal{E} \boxtimes \cO_{S} \rightarrow \Cal{F}_i$ on $X'\times S$ such that the isomorphism $\varphi \co \Cal{F}_r \cong \ft \Cal{F}_0$ intertwines $t_r$ with $\ft t_0$, and the maps $t_{i-1}, t_{i}$ are intertwined by the modification $f_i \co  \Cal{F}_{i-1} \dashrightarrow \Cal{F}_{i}$ for each $i = 1, \ldots, r$, i.e. the diagram below commutes. 
\[
\begin{tikzcd}
\Cal{E}\boxtimes\cO_{S} \ar[d, "t_0"] \ar[r, equals] &  \Cal{E}\boxtimes\cO_{S} \ar[r, equals]  \ar[d, "t_1"]& \ldots \ar[d] \ar[r, equals] & \Cal{E}\boxtimes\cO_{S} \ar[r, "\sim"] \ar[d, "t_{r}"] & \ft (\Cal{E} \boxtimes\cO_{S})\ar[d, "\ft t_0"] \\
\Cal{F}_0 \ar[r, dashed, "f_0"] & \Cal{F}_1 \ar[r, dashed, "f_1"]  &  \ldots \ar[r, dashed, "f_r"] &  \Cal{F}_{r} \ar[r, "\varphi"] & \ft \Cal{F}_0  
\end{tikzcd}
\]
\end{itemize}
In the sequel, when writing such diagrams we will usually just omit the ``$\boxtimes \cO_S$'' factor from the notation.

We define $\ol{\Cal{Z}}_{\Cal{E},\GL(n)'}^{r}$ to be the stack quotient $[\Cal{Z}_{\Cal{E},\GL(n)'}^{r}/(\Aut(\cE)(\F_q))]$.

We define ${\cZ}_{\Cal{E},\GL(n)'}^{r ,\circ} \subset \Cal{Z}_{\Cal{E},\GL(n)'}^r$ to be the open substack where the maps $\{t_i\}$ are all injective over every geometric point of $S$ (equivalently, any one of $\{t_{i}\}$ is injective).

We will call the $\Cal{Z}_{\Cal{E},\GL(n)'}^r$, ${\cZ}_{\Cal{E},\GL(n)'}^{r ,\circ}$ (or unions of their irreducible components) {\em special cycles of corank $m$} (with $r$ legs) on $\Sht_{\GL(n)'}^r$. 
\end{defn}

\begin{prop}\label{prop: Z finite}
Let $\cE$ be any vector bundle of rank $m$ on $X'$. Then the projection map $\Cal{Z}_{\Cal{E},\GL(n)'}^{r} \rightarrow \Sht_{\GL(n)'}^{r}$ is finite. 
\end{prop}

\begin{proof}This follows from 
similar argument as for \cite[Proposition 7.5]{FYZ}. 
\end{proof}

\subsection{Hermitian shtukas with similitude}\label{ss:Sht sim}

In \cite[\S 6]{FYZ} we worked with Hermitian shtukas based on the notion of a Hermitian bundle, which there was defined as a vector bundle $\cF$ with a Hermitian structure $h \co \cF \xrightarrow{\sim}  \sigma^* \cHom(\cF, \omega_{X'})$.

 In this section we consider a more general situation, where the notion of Hermitian structure is expanded to include maps $h \co \cF \xrightarrow{\sim} \sigma^* \cHom(\cF, \omega_{X'} \otimes \nu^* \LL)$ for any line bundle $\LL$ on $X$. These can be seen as torsors for a similitude unitary group. Most of the arguments of \cite{FYZ} already work at this level and generality, and it encompasses interesting situations not seen in the case $\LL = \cO_X$; for example, when $n$ is odd and $\mf{L}$ is not a norm from $X'$, the moduli space of shtukas with an odd number of legs is non-empty. The methods of \cite{FYZ} then give ``Kudla-Rapoport style'' identities between odd order Taylor coefficients of Siegel-Eisenstein series, whose functional equation has sign $-1$, and special cycles with an odd number of legs; see \S \ref{ssec: similitude higher KR} for the precise statements. 

\begin{defn}\label{defn: twisted similitude shtuka definitions}
Let $\LL$ be a line bundle on $X$. 
\begin{enumerate}
\item We define $\Bun_{U(n),\LL}$ analogously to \cite[Definition 6.1]{FYZ} but with the appearances of ``$\cF^{\vee}$'' ($=\cHom(\cF, \omega_{X'})$) in \emph{loc. cit.} replaced by $\cHom(\cF, \omega_{X'}  \otimes \nu^* \LL)$. Thus the $S$-points of the moduli stack $\Bun_{U(n),\LL}$ is the groupoid of pairs $(\cF, h)$ where $\cF$ is a rank $n$ vector bundle on $X'_S$, $h$ is an isomorphism $\cF \xrightarrow{\sim} \s^*\cHom(\cF, \omega_{X'} \otimes \nu^* \LL)$ satisfying $\sigma^* h^{\vee}= h \otimes \Id_{\nu^* \LL}$ (which we call an $\LL$-\emph{twisted Hermitian structure}), and morphisms $(\cF, h) \xrightarrow{\sim} (\cF', h')$ are isomorphisms $\cF \rightarrow \cF'$ intertwining $h$ with $h'$.

Similarly, for an integer $r \geq 0$, we define $\Hk_{U(n), \LL}^r$ analogously to \cite[Definition 6.3]{FYZ}. It has $S$-points the groupoid of the following data: 
\begin{enumerate}
\item $x_i' \in X'(S)$ for $i = 1, \ldots, r$, with graphs denoted by $\Gamma_{x'_i} \subset X' \times S$.
\item A sequence of vector bundles $\Cal{F}_0, \ldots, \Cal{F}_r$ of rank $n$ on $X' \times S$, each equipped with \emph{$\LL$-twisted Hermitian structures} $h_0, \ldots, h_r$. 
\item Isomorphisms $f_i \co \Cal{F}_{i-1}|_{X' \times S - \Gamma_{x'_i} - \Gamma_{\sigma(x'_i)}} \xrightarrow{\sim} \Cal{F}_i|_{X' \times S - \Gamma_{x'_i}-\Gamma_{\sigma(x'_i)}}$, for $1 \leq i \leq r$, compatible with the $h_i$, which are lower of length $1$ at $x_i'$ and upper of length $1$ at $\s x_i'$ (cf. \cite[Remark 6.4]{FYZ} for the terminology). 
\end{enumerate}

We define $\Sht_{U(n), \LL}^r$ analogously to \cite[Definition 6.6]{FYZ}, by the Cartesian diagram
\[
\begin{tikzcd}
\Sht_{U(n), \LL}^{r} \ar[r] \ar[d] & \Hk_{U(n),\LL}^{r} \ar[d, "{(\pr_{0}, \pr_{r})}"] \\
\Bun_{U(n),\LL} \ar[r, "{(\Id, \Frob)}"]  & \Bun_{U(n),\LL} \times \Bun_{U(n),\LL}
\end{tikzcd}
\]
where $\pr_{i}: \Hk_{U(n),\LL}^{r}\to \Bun_{U(n),\LL}$ records $(\cF_{i},h_{i})$.

Let $\Cal{E}$ be a rank $m$ vector bundle on $X'$. We define $\cZ_{\cE,\LL}^{r}$ analogously to \cite[Definition 7.1]{FYZ}. The $S$-points of the stack $\Cal{Z}_{\Cal{E},\LL}^{r}$ form the groupoid of the following data: 
 \begin{itemize}
 \item An $S$-point $(\{x'_1, \ldots, x'_r\}, \{\Cal{F}_0, \ldots, \Cal{F}_r\}, \{f_1, \ldots, f_r\}, \varphi) \in  \Sht_{U(n), \LL}^{r}(S)$.
 \item Maps of coherent sheaves  $t_i \co \Cal{E} \boxtimes \cO_{S} \rightarrow \Cal{F}_i$ on $X'\times S$ such that the isomorphism $\varphi \co \Cal{F}_r \cong \ft \Cal{F}_0$ intertwines $t_r$ with $\ft t_0$, and the maps $t_{i-1}, t_{i}$ are intertwined by the modification $f_i \co  \Cal{F}_{i-1} \dashrightarrow \Cal{F}_{i}$ for each $i = 1, \ldots, r$, i.e. the diagram below commutes. 
\[
\begin{tikzcd}
\Cal{E} \ar[d, "t_0"] \ar[r, equals] &  \Cal{E} \ar[r, equals]  \ar[d, "t_1"]& \ldots \ar[d] \ar[r, equals] & \Cal{E} \ar[r, "\sim"] \ar[d, "t_{r}"] & \ft \Cal{E} \ar[d, "\ft t_0"] \\
\Cal{F}_0 \ar[r, dashed, "f_0"] & \Cal{F}_1 \ar[r, dashed, "f_1"]  &  \ldots \ar[r, dashed, "f_r"] &  \Cal{F}_{r} \ar[r, "\varphi"] & \ft \Cal{F}_0  
\end{tikzcd}
\]
\end{itemize}

We will call the $\Cal{Z}_{\Cal{E},\LL}^{r}$ (or their connected components) {\em special cycles of corank $m$} (with $r$ legs), where we remind that $m = \rank \cE$.

For each $\cE$, we denote by $\Aut(\cE)(\F_q)$ its (finite) group of automorphisms as a vector bundle over $X'$. We define $\ol{\Cal{Z}}_{\Cal{E},\LL}^{r}$ to be the stack quotient $[\Cal{Z}_{\Cal{E},\LL}^{r}/(\Aut(\cE)(\F_q))]$.

\item The \emph{$\LL$-twisted Hitchin base} $\cA_{\cE, \LL}$ parametrizes maps $a \co \cE \rightarrow \sigma^* \cHom(\cE, \omega_{X'} \otimes \nu^* \LL)$ such that $\sigma^*a^{\vee}= a$ where $a^{\vee}$ is the map obtained by dualizing $a$ and then twisting by $\omega_{X'} \otimes \nu^* \LL$. The open subscheme $\cA_{\cE, \LL}^{\ns} \subset \cA_{\cE, \LL}$ parametrizes $a$ whose restriction to all geometric points of the test scheme are injective as maps of coherent sheaves. 

Note when $\LL=\cO_{X}$, $\cA_{\cE, \LL}(k)$ is what was denoted $\cA^{\all}_{\cE}(k)$ in \cite[Definition 7.2]{FYZ}; $\cA^{\ns}_{\cE, \LL}(k)$ is what was denoted $\cA_{\cE}(k)$ in \emph{loc. cit.}.

\item We have a decomposition $\cZ_{\cE,\LL}^{r}=\coprod_{a\in \cA_{\cE, \LL}(k)}\cZ_{\cE,\LL}^{r}(a)$. We define (analogously to \cite[Definition 7.4]{FYZ}) $\cZ_{\cE, \LL}^{r, \circ} \subset \cZ_{\cE, \LL}^r$ to be the open substack where the $t_i$ are injective, and $\cZ_{\cE, \LL}^{r, *} \subset \cZ_{\cE, \LL}^r$ to be the open substack where the $t_i$ are non-zero. For $a\in \cA_{\cE, \LL}(k)$, define $\cZ_{\cE,\LL}^{r}(a)^{\circ} := \cZ_{\cE, \LL}^r(a) \cap \cZ_{\cE, \LL}^{r, \circ}$ and $\cZ_{\cE,\LL}^{ r}(a)^*  := \cZ_{\cE, \LL}^r(a) \cap \cZ_{\cE, \LL}^{r, *} $, with the intersections formed in $\cZ_{\cE}^r$.
\end{enumerate}

\end{defn}

Let $\eta: \BA^\times/F^\times\rightarrow \{\pm1\}$ be the quadratic character associated to $F'/F$ by class field theory. Since $X'/X$ is \'etale, the character descends to $\y: \Pic_{X}(k)/\Pic_{X'}(k)\to \{\pm1\}$, and for $\LL \in \Pic_X(k)$ we have $\eta(\LL) = 1$ if and only if $\LL$ is a norm from $X'$. 

\begin{lemma}\label{lem: when non-empty}
With notation as above, $\Sht^{r}_{U(n),\LL}$ is non-empty if and only if $(-1)^r = \eta(\LL)^n$. 
\end{lemma}

\begin{proof}
The case $n=1$ is established later in Lemma \ref{lem: U(1) non-empty}. Here we shall assume this case and then establish the general case. 

Note that taking determinants induces a map $\Sht_{U(n), \LL}^r \rightarrow \Sht_{U(1), \LL^{\otimes n} \otimes \omega_X^{\otimes (n-1)}}^r$. By the result for the $n=1$ case, this shows that $\Sht_{U(n), \LL}^r  = \emptyset$ if $(-1)^r \neq \eta(\LL)^n$. It remains to prove that whenever $(-1)^r = \eta(\LL)^n$, then $\Sht_{U(n), \LL}^r$ is non-empty. 

If $\eta(\LL) = 1$, then $\Sht_{U(n), \LL}$ is isomorphic to $\Sht_{U(n)}$ by twisting, so the result follows from \cite[Lemma 6.7]{FYZ}. 

Suppose $\eta(\LL) = -1$. With $\Sht^r_{\GL(1)/X'} $ defined as in \cite[(5.4)]{YZ}, there is a map $\Sht^r_{\GL(1)/X'} \rightarrow \Sht^r_{U(2), \LL}$ sending $\cF_{0} \dashrightarrow  \ldots \dashrightarrow \ft \cF_r \cong \ft \cF_0 $ to $\cF_{0} \oplus  (\sigma^* \cF_{0}^{\vee} \otimes \nu^* \LL) \dashrightarrow \ldots \dashrightarrow \cF_{r} \oplus (\sigma^* \cF_{r}^{\vee} \otimes \nu^* \LL)  \cong \ft (\cF_{0} \oplus  (\sigma^* \cF_{0}^{\vee} \otimes \nu^* \LL))$. Since $\Sht^r_{\GL(1)/X'}$ is non-empty whenever $r$ is even, we find that $\Sht^r_{U(2), \LL}$ is non-empty whenever $r$ is even. Taking direct sums induces a map 
\begin{equation}\label{eq: direct sum map}
\Sht^0_{U(2), \LL}  \times \Sht^r_{U(n-2), \LL} \rightarrow \Sht^r_{U(n), \LL}
\end{equation}
which then inductively shows that $\Sht^r_{U(n), \LL}$ is non-empty whenever $r$ and $n$ are even. 

It remains to show that if $\eta(\LL) = -1$ and $n$ is odd, then $\Sht_{U(1), \LL}^r$ is non-empty whenever $r$ is odd. Since we are assuming the $n=1$ case, we know that $\Sht_{U(1), \LL}^r$ is non-empty for all odd $r$. Then iterating \eqref{eq: direct sum map} shows that $\Sht_{U(n), \LL}^r$ is non-empty for all odd $n$ and $r$. 
\end{proof}

For properties of the objects in Definition \ref{defn: twisted similitude shtuka definitions} whose proofs are the same for general $\LL$ as written in the case $\LL = \cO_X$ in \cite{FYZ}, we will just cite the statements from \cite{FYZ}. For example, by the same proofs as for \cite[Lemma 6.8, Lemma 6.9]{FYZ}, we have the following geometric properties.

\begin{lemma}\label{lem: shtuka geometry} 
Let $\mf{L}$ be any line bundle on $X$. 

\begin{enumerate}
\item The stack $\Bun_{U(n), \mf{L}}$ is smooth and equidimensional of dimension $n^{2}(g_{X}-1)$. 

\item The projection map $(\pr_{X}, \pr_{r}):\Hk_{U(n), \mf{L}}^r \rightarrow (X')^r \times \Bun_{U(n), \mf{L}}$ recording $\{x_{i}\}$ and $(\cF_{r},h_{r})$ is smooth of relative dimension $r(n-1)$. 

\item $\Sht_{U(n), \mf{L}}^r$ is a Deligne-Mumford stack locally of finite type. The map $\Sht_{U(n), \mf{L}}^r \rightarrow (X')^r$ is smooth, separated, equidimensional of relative dimension $r(n-1)$. 
\end{enumerate}
\end{lemma}

Forgetting the Hermitian structures give maps $\Bun_{U(n),\LL}\to \Bun_{\GL(n)'}$ and $\Hk^{r}_{U(n),\LL}\to \Hk^{r}_{\GL(n)'}$, which induce a map over $(X')^{r}$
\begin{equation*}
\Sht^{r}_{U(n),\LL}\to \Sht^{r}_{\GL(n)'}
\end{equation*}

\begin{lemma}\label{lem: special cycle GL(n)' to U(n)}
Let $\cE$ be any vector bundle of rank $m$ on $X'$. Then we have 
\[
\Cal{Z}_{\Cal{E},\LL}^{r}  \cong 
\Cal{Z}_{\cE, \GL(n)'}^{r}  \times_{\Sht_{\GL(n)'}^r} \Sht_{U(n), \LL}^r,
\]
\[
\Cal{Z}_{\Cal{E},\LL}^{r, \circ}  \cong 
\Cal{Z}_{\cE,\GL(n)'}^{r, \circ}  \times_{\Sht_{\GL(n)'}^r} \Sht_{U(n), \LL}^r,
\]
as stacks over $\Sht_{U(n), \LL}^r$. 
\end{lemma}

\begin{proof}
Immediate from the definitions. 
\end{proof}

\subsection{The case $n=1$}\label{ss: sht for n=1} We now undertake a closer analysis of $\Sht_{U(n), \LL}^r$ for $n=1$. We first set up some notation. Let $\Prym_{\LN}=\Bun_{U(1),\LN}$ be $\Nm^{-1}(\LN)$ where $\Nm: \Pic_{X'}\to \Pic_{X}$ is the norm map. When $\LN=\cO_X$ we omit the subscript $\LN$. In this case, there are two connected components and we let $\Prym^{0}$ denote its neutral component, and $\Prym^{1}$ for the other component (both are defined over $k$). Explicitly, there is a map $\Pic(X') \rightarrow \Prym$ taking a line bundle $L$ to $L \otimes \sigma^* L^{-1}$, and $\Prym^0$ is the image of the components of $\Pic(X')$ where $L$ has even degree, while $\Prym^1$ is the image of the components where $L$ has odd degree. This is explained in \cite[p. 186-188]{Mum71}.

Recall that $\Prym^0$ and $\Prym^1$ are also geometrically connected. Since $\Prym_{\frN}$ is a torsor under $\Prym$, $\Prym_{\frN}$ also has two \emph{geometric} connected components. However, its number of ($k$-rational) connected components depends on $\y(\frN)$, as explained below. 

\begin{lemma}\label{lem: prym ratl components}
If $\eta(\frN) = 1$ then $\Prym_{\frN}$ has two connected components. If $\eta(\frN) = -1$ then $\Prym_{\frN}$ has one connected component.
\end{lemma}

\begin{proof}When $\eta(\frN)=1$, i.e., $\frN$ is a norm,  $\Prym_{\frN}$ has a $k$-point hence is a trivial $\Prym$-torsor, therefore both geometric components of $\Prym_{\frN}$ are defined over $k$. When $\y(\frN)=-1$, $\Prym_{\frN}$ has no $k$-point, which implies that the two geometric components of $\Prym_{\frN}$ are permuted by $\Gal(\ov k/k)$ (for otherwise a Frobenius stable geometric component, being a torsor under the connected group $\Prym^{0}$, would contain a $k$-point by Lang's theorem), hence $\Prym_{\frN}$ is connected.
\end{proof}

\begin{lemma}\label{lem: Prym component}
Let $\epsilon(\LN)\in\{0,1\}$ be such that $\eta(\LN)=(-1)^{\epsilon(\LN)}$. Then the Lang map 
\begin{equation*}
\begin{aligned}\Lang\colon\Prym_{\LN}&\to \Prym\\
\cF&\mapsto {}^{\tau}\cF\ot\cF^{-1}
\end{aligned}
\end{equation*}
lands in $\Prym^{\epsilon(\LN)}$. 

\end{lemma}

\begin{proof}

Given $\LN_1, \LN_2 \in \Pic_X(k)$ such that $\LN_2 \otimes (\LN_1)^{-1} = \Nm(\LN')$ for some $\LN' \in \Pic_{X'}(k)$, twisting by $\LN'$ induces an isomorphism $\Prym_{\LN_1} \xrightarrow{\sim} \Prym_{\LN_2}$. Hence if $\LN$ is a norm, then $\Prym_{\LN} \xrightarrow{\sim} \Prym_{\cO_X}$, in which the claim is a result of Wirtinger explained in \cite[\S 2]{Mum71}. 

If $\LN$ is not a norm, by the twisting argument above, it suffices to show the statement for a single choice of $\LN$. We take $\LN = \cO(x)$ for a closed point $x \in |X|$ which is inert in $X'$. 

We claim that it suffices to check that the statement for a single geometric point $\cF\in \Prym_{\LN}$. Indeed, since $\Prym_{\LN}$ is a $\Prym$-torsor, any geometric point of $\Prym_{\LN}$ is of the form $\cF\ot \cF'$ for some $\cF'\in \Prym$, and $\Lang (\cF\ot \cF')\cong\Lang(\cF)\ot\Lang (\cF')$ lies in the same component of $\Prym$ as $\Lang (\cF)$ since $\Lang(\Prym)\subset \Prym^{0}$.


To describe a geometric point $\cF\in \Prym_{\LN}(\ov k)$, write $x \times_{\Spec k} \Spec \ol{k} = \{x_1, \ldots, x_d\}$ such that $\Frob(x_i) = x_{i+1 \pmod{d}}$, etc. Denoting $x'$ the point of $X'$ lying above $x$, we have $x' \times_{\Spec k} \Spec \ol{k} = \{x_1', x_2', \ldots, x_{2d}'\}$ where $\Frob (x_i') = x_{i+1 \pmod{2d}}'$ and $\sigma x_i' = x_{i+d \pmod{2d}}'$, etc. Then $\cF := \cO(x_1' + x_2' + \ldots + x_d')$ lies in $\Prym_{\LN}(\ol{k})$, and $\ft \cF \otimes \cF^{-1} = \cO(x_{d+1}'-x_1') = \cO(\sigma x_1'-x_1')$, which  lies in the non-neutral component of $\Prym$. 
\end{proof}

Let $r$ be even (resp. odd) if $\epsilon(\LN)=0$ (resp. $\epsilon(\LN)=1$). By unwinding definitions one sees directly that the diagram below is Cartesian: 
\begin{equation}\label{eq: prym cartesian square}
\xymatrix{\Sht^{r}_{U(1), \LN}\ar[r]^-{p} \ar[d]_{p_{[1,r]}:=(p_{1},\cdots, p_{r})}& \Prym_{ \LN}\ar[d]^{\Lang}\\
X'^{r}\ar[r]^{\AJ^{r}} & \Prym}
\end{equation}
Here $\AJ^{r}: X'^{r}\to \Prym$ is the map $(x_{1},\cdots,x_{r})\mapsto \cO(\sum_{i=1}^{r}(\s x_{i}-x_{i}))$. 
The map $p_{i}: \Sht^{r}_{U(1)}\to X'$ records the $i$-th leg ($1\le i\le r$), $p_{[1,r]}:=(p_{1},\cdots, p_{r}): \Sht^{r}_{U(1)}\to (X')^{r}$ and $p: \Sht^{r}_{U(1)}\to \Prym_{ \LN}$ records $\cF_{0}$.

\begin{lemma}\label{lem: U(1) non-empty}
 $\Sht^{r}_{U(1),\LL}$ is non-empty if and only if $(-1)^r=\eta(\LL)$. 
 \end{lemma}
 
\begin{proof}Combine Lemma \ref{lem: Prym component} and \eqref{eq: prym cartesian square}.
\end{proof}

\begin{lemma}\label{lem: 2 connected components}
If $r>0$ and $\Sht^r_{U(1), \LL}$ is non-empty, then $\Sht^r_{U(1), \LL}$ has two geometric connected components. Under these same assumptions, $\Sht^r_{U(1), \LL}$ is connected if and only if $r$ is odd. 
\end{lemma}

\begin{proof}

Let $\LN=\om_{X}\ot \LL$. We know that $\Sht^{r}_{U(1),\frL}\ne\vn$ if and only if $\y(\frN)=\y(\LL)=(-1)^{r}$.  We assume this in the following.

 First we establish that there are two geometric connected components. Consider the Cartesian square \eqref{eq: prym cartesian square}.  For any $\ep\in \Irr(\Prym_{\frN, \ol{k}})$ (a torsor for $\Z/2\Z$), let $\Sht^{r, \ep}_{U(1),\frL, \ol{k}}$ be the preimage of $\Prym^{\ep}_{\frN, \ol{k}}$ under $p$. We need to show that $\Sht^{r,\ep}_{U(1),\frL, \ol{k}}$ is connected. 

As a $\Prym^{0}(k)$-torsor over $X'^{r}_{\ov k}$ (cf. \eqref{eq: prym cartesian square}), $\Sht^{r,\ep}_{U(1),\frL, \ov k}$  is given by the homomorphism
\begin{equation}\label{eq: lang composition}
\pi_1(X'^{r}_{\ov k}) \xr{\AJ^{r}_{*}}\pi_1(\Prym^{\un r}_{\ov k}) \xr{\L_{\Prym}} \Prym^0(k)
\end{equation}
where the first map is induced by $\AJ^{r}$ (where $\un r=r\mod 2\in \{0,1\}$), and the second map is given by the Lang torsor $\Lang: \Prym^{\ep}_{\frN}\to \Prym^{\un r}$. It suffices to show that  \eqref{eq: lang composition} is surjective. Since $\Prym^{\ep}_{\frN}$ is geometrically connected, $\L_{\Prym}$ is surjective.  It remains to show that $\AJ^{r}_{*}$ is surjective.

Fixing $z=(z_{1},\cdots, z_{r-1})\in X'^{r-1}(\ov k)$ and letting $\D_{z}=\AJ^{r-1}(z)$, we have a commutative diagram
\begin{equation}
\xymatrix{ X'_{\ov k}\ar[r]^{\AJ_{X'}}\ar[d]^{i_{z}} & \Pic^{1}_{X',\ov k}\ar[r]^{\s-1}  & \Prym^{1}_{\ov k}\ar[d]^{\ot \D_{z}}\\
X'^{r}_{\ov k}\ar[rr]^{\AJ^{r}} && \Prym^{\un r}_{\ov k}}
\end{equation} 
Here $\AJ_{X'}:X'\to \Pic^{1}_{X'}$ is the Abel-Jacobi map for $X'$, $i_{z}(x)=(x,z_{1},\cdots, z_{r-1})$. It induces a commutative diagram on fundamental groups
\begin{equation}
\xymatrix{\pi_{1}(X'_{\ov k})\ar[r]^{\AJ_{X',*}}\ar[d]^{i_{z*}} & \pi_{1}(\Pic^{1}_{X',\ov k})\ar[r]^{(\s-1)_{*}}  & \pi_{1}(\Prym^{1}_{\ov k})\ar[d]^{\cong}\\
\pi_{1}(X'^{r}_{\ov k})\ar[rr]^{\AJ^{r}_{*}} && \pi_{1}(\Prym^{\un r}_{\ov k})
}
\end{equation}
By geometric class field theory, $\AJ_{X',*}: \pi_{1}(X'_{\ov k})\to \pi_{1}(\Pic^{1}_{X',\ov k})$ is surjective, realizing the latter as the abelianization of the former. On the other hand,  $\s-1: \Pic^{1}_{X',\ov k}\to \Prym^{1}_{\ov k}$ is a torsor under $\Pic^{0}_{X,\ov k}$ which is connected, it induces a surjection on $\pi_{1}$. These then imply that the top row of the above diagram is surjective. Therefore the bottom row is surjective as well, i.e.,  $\AJ^{r}_{*}$ is surjective.

To prove the last assertion in the Lemma, we show that $\Frob$ swaps the two geometric connected components of $\Sht^r_{U(1), \LL}$ if and only if $r$ is odd. For $\cF \in \Sht^r_{U(1), \LL}(\ol{k})$, $\ft \cF \otimes \cF^{-1}$ is the tensor product of $r$ line bundles of the form $\cO(x - \sigma x)$, each of which lies in $\Prym^1$, so the tensor product lies in $\Prym^0$ if and only if $r$ is even. 	

\end{proof}

\section{Hitchin stacks}\label{sec: hitchin spaces}

In this section we introduce certain stacks which will be used to analyze special cycles, generalizing the constructions in \cite[\S 8]{FYZ}.

\subsection{Moduli of sections of gerbes} 

In order to encompass the moduli stacks $\Bun_{U(n)}$ and $\Bun_{U(n), \LL}$ in a common framework, it will be advantageous to adopt a more general perspective of moduli stacks of sections of gerbes. 

\begin{example}\label{ex: gerbe 1}
Let $G$ be a group scheme over any scheme $S$. Then the relative classifying stack $BG$ is equipped with the structure of a gerbe over $S$, and the groupoid of sections of $BG$ over $S$ is equivalent to the groupoid over $G$-torsors over $S$. In particular, for a group scheme $G$ over the curve $X$, $\Bun_G$ can be interpreted as a moduli stack of sections of the gerbe $BG$ over $X$. 

In particular, for the group scheme $\GL(n)'=\Res_{X'/X}(\GL(n)\times X')$ over $X$, the gerbe $B\GL(n)'$ represents the following moduli problem: for any $k$-scheme $S$, $(B\GL(n)')(S)$ is the groupoid of pairs $(s,\cF)$, where $s:S\to X$ and $\cF$ is a rank $n$ vector bundle over $S'=S\times_{X}X'$.
\end{example}

In the context of this paper, the moduli stack of Hermitian bundles $\Bun_{U(n)}$ over $X$ play a more fundamental role than the group scheme $U(n)$ itself. Indeed, to recover $U(n)$ from $\Bun_{U(n)}$ we need to choose a base point $(\cF,h)\in \Bun_{U(n)}(k)$ and define $U(n)$ to be the group scheme of automorphisms of $(\cF,h)$. Better yet, we should consider the gerbe $BU(n)$ over $X$ rather than the group scheme $U(n)$ over $X$. Then sections of the gerbe $BU(n)$ are equivalent to $U(n)$-torsors. This point of view generalizes better to include spaces like $\Bun_{U(n), \LL}$, which are not moduli stacks of torsors for a group scheme, but can be seen as moduli stacks of sections of a gerbe $BU(n)_{\frL}$, which will be defined next.

\begin{defn}\label{def: Bun gerbe} Let $\sG$ be a gerbe over $X$. We define the stack $\Bun_{\sG}$ over $k$ to be
\begin{equation}
\Bun_{\sG}:=\Sect(X,\sG)=R_{X/k}\sG.
\end{equation}
In other words, the $S$-points of $\Bun_{\sG}$ form the groupoid of maps $X\times S\to \sG$ over $X$.
\end{defn}

In view of Example \ref{ex: gerbe 1}, we have $\Bun_{BG}=\Bun_{G}$ for a group scheme $G$  over $X$.

\subsubsection{Unitary gerbes}
Fix a line bundle $\frL$ over $X$. We define the gerbe $BU(n)_{\frL}$ over $X$ to represent the following moduli problem: for any scheme $S$ with a map  $s: S\to X$, liftings of $s$ to $BU(n)_{\frL}$ form the groupoid of Hermitian vector bundles $(\cF,h)$ over $S' :=S\times_{X}X'$ valued in $s^{*}\frL$, i.e., $h$ is an isomorphism $\cF\isom \s^{*}_{S}\cHom(\cF,\nu_{S}^{*}s^{*}(\om_{X}\ot\frL))$ satisfying $h=\s^{*}h^{\vee}$ (here $\s_{S}:S'\to S'$ and $\nu_{S}:S'\to S$ are induced from $\s$ and $\nu$). Forgetting the datum of $h$ defines the \emph{standard map} $BU(n)_{\LL} \rightarrow B\GL(n)'$.

We call $BU(n)_{\frL}$ the {\em unitary gerbe over $X$ of rank $n$ and similitude line bundle $\frL$}. With this definition and Definition \ref{def: Bun gerbe}, we have
\begin{equation}
\Bun_{BU(n)_{\LL}}=\Bun_{U(n),\frL}.
\end{equation}


For most of the paper, the only gerbes that will concern us are $BU(n)_{\frL}$ or $B \GL(n)'$. However, in \S \ref{ssec: H_2 functoriality} and \S \ref{sec: CM} it will be necessary to deal with a more general class of gerbes, which we introduce next.

\sss{Gerbes of unitary type}\label{sssec: unitary type}
We define a class of gerbes  over $X$ that we call {\em gerbes of unitary type over $X$}. 

Let $Y$ be another smooth projective curve over $k$ (not assumed to be geometrically connected) and $\th: Y\to X$ be a finite morphism (possibly ramified). Let $\Irr(Y)$ be the set of irreducible components of $Y$ and $\un n: \Irr(Y)\to \Z_{>0}$ be a function.  For $Y_{\a}\in \Irr(Y)$ we denote $\un n(Y_{\a})$ by $n_{\a}$. Let 
\begin{equation}\label{rk sum}
n=\sum_{Y_{\a}\in \Irr(Y)}n_{\a}[Y_{\a}:X]
\end{equation}
where $[Y_{\a}:X]$ is the degree of $\th_{\a}:=\th|_{Y_{\a}}:Y_{\a}\to X$.

Let $\frL\in\Pic(X)$ and  $\frL_{\a}=\th^{*}_{\a}\frL$. Consider the unitary gerbe $BU(n_{\a})_{\frL_{\a}}$ over $Y_{\a}$ with similitude line bundle $\frL_{\a}$ defined using the double cover $Y'_{\a}=Y_{\a}\times_{X}X'$. 

We claim there is  a canonical map of gerbes over $X$
\begin{equation}\label{map from prod BU}
\prod_{Y_{\a}}R_{Y_{\a}/X}BU(n_{\a})_{\frL_{\a}}\to BU(n)_{\frL}.
\end{equation}
We describe the map on the level of $S$-points. For $s:S\to X$,  $(R_{Y_{\a}/X}BU(n_{\a})_{\frL_{\a}})(S)$ is the groupoid of Hermitian bundles $(\cF_{\a}, h_{\a})$ over $S'\times_{X}Y_{\a}$ with similitude line bundle the pullback of $s^{*}\frL$ to $S_{\a}:=S\times_{X}Y_{\a}$. Given such $S$-points $(\cF_{\a}, h_{\a})$ for each $Y_{\a}$, \eqref{map from prod BU} sends them to the direct sum $\cF=\op_{\a}\th'_{S_{\a}*}\cF_{\a}$ (where $\th'_{S_{\a}}: S'_{\a}=S'\times_{X}Y_{\a}\to S'$ is the projection). The pushforward of $h_{\a}$ induces a map
\begin{equation}\label{push herm}
\th'_{S_{\a}*}h_{\a}: \th'_{S_{\a}*}\cF_{\a}\isom \th'_{S_{\a}*} \sigma_S^*\cHom(\cF_{\a}, (\om_{Y_{\a}}\ot  \frL)|_{S'}).
\end{equation}
The relative dualizing sheaves satisfy $\om_{S'_{\a}/S'}\cong \om_{Y'_{\a}/X'}|_{S'}\cong \om_{Y_{\a}/X}|_{S'}$. Grothendieck-Serre duality gives
\begin{equation}
\th'_{S_{\a}*}\cHom(\cF_{\a},\om_{S'_{\a}/S'})=\cHom(\th'_{S_{\a}*}\cF_{\a}, \cO_{S'}).
\end{equation}
Therefore the right side of \eqref{push herm} is isomorphic to
\begin{eqnarray}
&&\sigma_S^*\th'_{S_{\a}*}\cHom(\cF_{\a}, \om_{Y'_{\a}/X'}|_{S'}\ot \th'^{*}_{S_{\a}}( \om_{X}\ot\frL)|_{S'})\\
& \cong& \sigma_S^*\th'_{S_{\a}*}\cHom(\cF_{\a}, \om_{S'_{\a}/S'})\ot (\om_{X}\ot\frL)|_{S'}\\
&\cong& \sigma_S^* \cHom(\th'_{S_{\a}*}\cF_{\a}, \cO_{S'})\ot (\om_{X}\ot\frL)|_{S'}\\
& \cong& \sigma_S^* \cHom(\th'_{S_{\a}*}\cF_{\a}, \nu_{S}^{*}s^{*} (\om_{X}\ot\frL)).
\end{eqnarray}
In other words, $\th'_{S_{\a}*}h_{\a}$ is a Hermitian form on $\th'_{S_{\a}*}\cF_{\a}$ with similitude line bundle $\frL$. Then the direct sum of $\th'_{S_{\a}*}h_{\a}$ gives a Hermitian form $h$ on $\cF$ with similitude line bundle $\frL$.

\begin{defn}\label{def: gerbe unitary type} An \emph{$n$-framed gerbe} is a gerbe $\sG$ over $X$ together with a map $\sG \rightarrow B\GL(n)'$. We say that an $n$-framed gerbe is \emph{smooth} if $\sG$ is smooth. In notation we suppress the map to $B\GL(n)'$ if it is clear from context.

Let $\frL\in\Pic(X)$. A gerbe $\sG$ over $X$ together with a map $i: \sG\to BU(n)_{\frL}$ over $X$ is called {\em a gerbe of unitary type of rank $n$ and similitude line bundle $\frL$}, if there exists the data $\th:Y\to X$ and $\un n: \Irr(Y)\to \Z_{>0}$ as above (satisfying \eqref{rk sum}) such that $(\sG,i)$ is isomorphic to $\prod_{Y_{\a}}R_{Y_{\a}/X}BU(n_{\a})_{\frL_{\a}}$ (product over $X$) with the canonical map to $BU(n)_{\frL}$ defined in \eqref{map from prod BU}. In notation we suppress $i$ if it is clear from context. 

The \emph{standard map} $\sG \rightarrow B\GL(n)'$ is inflated via $i$ from the standard map for $BU(n)_{\LL}$. Given $i: \sG\to BU(n)_{\frL}$ a gerbe of unitary type of rank $n$ (and similitude line bundle $\frL$), we say that the \emph{standard $n$-framed structure} is the $n$-framed gerbe obtained by composing $i$ with the standard map for $BU(n)_{\LL}$. 

The \emph{tautological $n$-framed gerbe} is $\sG \xrightarrow{=} B\GL(n)'$.

\end{defn}

\subsubsection{Hecke stacks and shtukas for  gerbes of unitary type}\label{sssec: hecke for unitary gerbe} 

Let $\sG \cong \prod_{Y_{\alpha}} R_{Y_{\alpha}/X} BU(n_{\alpha})_{\LL_{\alpha}}$ be a gerbe of unitary type. We set $Y' := Y \times_X X' \cong  \coprod Y_{\alpha}'$ with involution $\s_{Y}=\Id_{Y}\times \s$. We have
\begin{equation}
\Bun_{\sG}\cong \prod_{\a}\Bun_{U(n_{\a})/Y_{\a},\LL_{\a}}
\end{equation}
where $\Bun_{U(n_{\a})/Y_{\a},\LL_{\a}}$ is the moduli of  $\LL_{\a}$-twisted Hermitian bundles of rank $n_{\a}$ over $Y'_{\a}$.

Then we define $\Hk_{\sG}^r$ to be the moduli stack with $S$-points being the groupoid of the following data:
\begin{itemize}
\item $(y_1',\cdots, y_{r}') \in Y'(S)^{r}$, 
\item Hermitian bundles $(\cF_{i}, h_{i})_{i=0}^r$, with each $\cF_i$ a vector bundle on $Y'\times S$, of rank $n_{\alpha}$ on $Y_{\alpha}'$, and $h_i$ is an $\LL_{\alpha}=\th_{\a}^{*}\LL$-twisted Hermitian structure on $\cF_i$, and 
\item Isomorphisms $f_i \co \Cal{F}_{i-1}|_{Y' \times S - \Gamma_{y'_i} - \Gamma_{\sigma_{Y} y'_i}} \xrightarrow{\sim} \Cal{F}_i|_{Y' \times S - \Gamma_{y'_i}-\Gamma_{\sigma_{Y} y'_i}}$, for $1 \leq i \leq r$, which are lower of length $1$ at $y'_i$ and upper of length $1$ at $\sigma y'_i$.
\end{itemize}

By recording how many of $y'_{i}$ are lying over each component of $Y$, we have a decomposition
\begin{equation}
\Hk^{r}_{\sG}\cong \coprod_{\un r} \prod_{\a\in \Irr(Y)}\Hk^{r_{\a}}_{U(n_{\a})/Y_{\a}, \LL_{\a}}
\end{equation}
where $\un r$ runs over the set of functions $\un r: \Irr(Y)\to \Z_{\ge0}$ such that $|\un r|:=\sum_{\a}r_{\a}$ is equal to $r$.

We define $\Sht_{\sG}^r$ by the Cartesian diagram 
\[
\begin{tikzcd}
\Sht_{\sG}^r \ar[r] \ar[d] & \Hk_{\sG}^r \ar[d,"{(\pr_{0}, \pr_{r})}"] \\
\Bun_{\sG} \ar[r, "{(\Id, \Frob)}"] & \Bun_{\sG}  \times \Bun_{\sG} 
\end{tikzcd}
\]
Similarly we have
\begin{equation}
\Sht^{r}_{\sG}\cong \coprod_{|\un r|=r} \prod_{\a\in \Irr(Y)}\Sht^{r_{\a}}_{U(n_{\a})/Y_{\a}, \LL_{\a}}.
\end{equation}

In order to make the notation more uniform, we will denote a gerbe of unitary type by $BG \rightarrow X$ (even if it does not arise as the classifying stack of a group scheme $G$). We will write $\Hk_G^r$ for $\Hk_{BG}^r$ (in the unitary type case) or if $G = \GL(n)'$.

\subsection{Hitchin stacks}\label{ssec: Hitchin stacks}
We introduce Hitchin stacks $\cM_{H_{1},H_{2}}$ for certain gerbes $BH_{1}$ and $BH_{2}$, generalizing the construction in \cite[\S8]{FYZ}.

There is an equivalence of categories between the groupoid of $\GL(n)'$-torsors over $X$, and the groupoid of vector bundles of rank $n$ on $X'$ (with maps being isomorphisms). If $\cE$ is a $\GL(n)'$-torsor, we denote by $\Vect(\Cal{E})$ the vector bundle associated by this equivalence. We introduce this notation because we shall frequently need to talk about maps between vector bundles which are not isomorphisms (and so do not come from maps of torsors). 

Because of Example \ref{ex: gerbe 1}, we will use the notation $BH$ for a gerbe over a base $S$, and refer to a global section of $BH$ over $S$ as an ``$H$-torsor over $S$'', even when the gerbe does not actually come as the classifying space of a group scheme $H$. More generally, given an $n$-framed gerbe $BH \rightarrow B\GL(n)'$ over $X$, and an $H$-torsor $\cE$ on $X \times S$, we will denote by $\Vect(\cE)$ the associated rank $n$ vector bundle on $X' \times S$.

\begin{defn}\label{defn: hitchin}
Let $BH_1$ be an $m$-framed gerbe over $X$ and $BH_2$ be an $n$-framed gerbe over $X$. We define the ``Hitchin-type space'' $\Cal{M}_{H_1,H_2}$ whose $S$-points are the groupoid of data: 
\begin{itemize}
\item $\Cal{E}_{H_1}$, an $H_1$-torsor over $X\times S$.
\item $\Cal{F}_{H_2}$, an $H_2$-torsor over $X\times S$. 
\item A map of vector bundles $t \co \Vect(\Cal{E}_{H_1}) \rightarrow \Vect(\Cal{F}_{H_2})$ over $X' \times S$.
\end{itemize}
We define $\Cal{M}_{H_1, H_2}^{\c} \subset \Cal{M}_{H_1, H_2}$ to be the open substack where the map $t$ is injective as a map of coherent sheaves. Note that the definition is in terms of the gerbes $BH_1, BH_2$ and their maps to $B\GL(m)', B\GL(n)'$, but in the notation we only put $H_1, H_2$ as a shorthand for $BH_1, BH_2$ and these maps (this is just notational shorthand -- there may not be an actual group scheme $H_i$ from which $BH_i$ comes). 
\end{defn}

\begin{remark} Let us comment on what generality of gerbes will appear. In all examples of interest, $BH_1 \rightarrow B\GL(m)'$ comes from a map of smooth group schemes over $X$, and $BH_2$ will be either a gerbe of unitary type with the standard map to $B\GL(n)'$, or simply $B \GL(n)'$ (with the identity map). The reader may focus on the cases where the gerbes arise as classifying stacks of smooth group schemes over $X$, without missing the main ideas.
\end{remark}

\begin{example}\label{ex: main hitchin example} Let $BH_1 \xrightarrow{=} B\GL(m)'$ and $BH_2 = BU(n) \rightarrow B\GL(n)'$ be the standard map. In this case, the stack $\cM_{H_1, H_2}$ (resp. $\cM_{H_1, H_2}^{\c}$) is the Hitchin stack denoted $\cM^{\all}(m,n)$ (resp. $\cM(m,n)$) in \cite{FYZ}. (Note the notational inconsistency with \cite{FYZ}: in this paper we do not use the superscript ``$\all$'' to indicate all maps are allowed, and we use the superscript $\c$ to indicate the substack where the map $t$ is injective.) 

More generally, for $\LL$ a line bundle on $X$ we may take $BH_2 = BU(n)_{\LL} \rightarrow B\GL(n)'$ to be the standard map. We also denote the corresponding Hitchin stack $\cM_{H_{1},H_{2}}$ by  $ \cM_{H_1, U(n), \LL}$. It parametrizes $\cE_{H_1} \in \Bun_{H_1}$, $(\cF,h) \in \Bun_{U(n), \LL}$, and a map of vector bundles
\[
t \co \Vect(\Cal{E}_{H_1} ) \rightarrow \cF
\]
over $X' \times S$. Its open substack $\cM_{H_1, H_{2}}^{\circ}=\cM_{H_{1},U(n),\LL}^{\c}$ is the locus where $t$ is injective as a map of coherent sheaves (fiberwise over the test scheme $S$). \emph{Henceforth when $\LL$ is understood, we may suppress it from the notation.} 
\end{example}

\begin{example}
In this paper we shall also be interested in the case $H_1 = \GL(m_1)' \times \ldots \times \GL(m_j)'$ where $m= m_1 + \ldots + m_j$, the map $BH_1 \rightarrow B\GL(m)'$ is induced by the standard block diagonal inclusion, and $BH_2 = BU(n)_{\LL} \rightarrow B\GL(n)'$ induced by the standard embedding. This comes up, for example, in \S \ref{ss:pf LI}.
\end{example}

\subsection{Hitchin base}\label{ssec: hitchin base} We construct Hitchin bases for our Hitchin stacks, generalizing \cite[\S 8.2]{FYZ}.

\begin{defn}\label{def: Hitchin base} Let $BH_1$ be an $m$-framed gerbe over $X$ and fix a line bundle $\LL$ on $X$. The \emph{$\LL$-twisted Hitchin base} $\cA_{H_1, \LL}$ is the stack whose $S$-points are the groupoid of the following data: 
\begin{itemize}
\item $\cE$ an $H_1$-torsor on $X \times S$. 
\item $a \co \Vect(\cE) \rightarrow \sigma^* \cHom(\Vect(\Cal{E}), \omega_{X'} \otimes \nu^{*}\LL)=\s^{*}\Vect(\cE)^{\vee}\ot \nu^{*}\LL$ is a map of coherent sheaves on $X'\times S$ such that $\sigma^*(a^{\vee}) = a$.
\end{itemize}
 We define the \emph{non-degenerate locus} $\cA^{\ns}_{H_1, \LL} \subset \cA_{H_1, \LL}$ to be the open substack where $a$ is injective  fiberwise over the test scheme $S$. \emph{When $\LL$ is understood, we will omit it from the notation in the future.}
\end{defn}

\begin{defn}\label{def: hitchin fibration} Let $BH_1$ be an $m$-framed gerbe over $X$ and fix a line bundle $\LL$ on $X$. Take $BH_2 = BU(n)_{\LL}$  with the standard $n$-framed structure. We define the \emph{Hitchin fibration} $f \co \cM_{H_1, H_{2}} \rightarrow \Cal{A}_{H_1,\LL}$ sending $(\Cal{E}, (\Cal{F},h), t)$ to
the composition 
\[
a\co \Vect(\Cal{E}) \xrightarrow{t} \cF \xrightarrow{h} \sigma^* \cF^{\vee}\ot \nu^{*}\LL\xrightarrow{\sigma^* t^{\vee}} \s^{*}\Vect(\cE)^{\vee}\ot\nu^{*}\LL.
\]
We define $\cM_{H_1, H_2}^{\ns} = \cM_{H_1, H_{2}} |_{ \Cal{A}_{H_1,\LL}^{\ns}}$; note that $\cM_{H_1, H_2}^{\ns}$ is an open substack of $\cM_{H_1, H_{2}}^{\c}$. 
\end{defn}

\begin{example} For $BH_1=B\GL(m)'$ the tautological $m$-framed gerbe, and $\LL = \cO_X$, $\cA_{H_1,\LL}$ (resp. $\cA_{H_1, \LL}^{\ns}$) coincides with the Hitchin base denoted $\cA^{\all}(m)$ (resp. $\cA(m)$) in \cite[Definition 8.2]{FYZ}.  Note the notational inconsistency with \cite{FYZ}: in this paper we do not use the superscript ``$\all$'' to indicate that all maps are allowed, and we use the superscript $\c$ to indicate the substack where the map $t$ must be injective.
\end{example}

\subsection{Smoothness of some Hitchin stacks} \label{ssec: tangent complex}
We will use the description of the tangent complex for the following general situation. Suppose that $G \rightarrow X$ is a smooth group scheme acting linearly on a vector bundle $V \rightarrow X$. Then the relative Lie algebra $\Lie(G/X)$ acts on $V$, and the relative tangent complex for $V/G \rightarrow X$ at a point $(x,v)$  (where $v\in V_{x}$) is represented by the complex
\begin{equation}\label{eq: tangent complex 1}
\begin{tikzcd}[row sep = tiny]
 & \a_{v}:\underbrace{\Lie G_{x}}_{\deg -1}  \ar[r] &  \underbrace{V_{x}}_{\deg 0} \\
& Y \ar[r] & Y \cdot v
\end{tikzcd}
\end{equation}

Let $\kappa$ be a field containing $k$. A $\kappa$-point of $ \Sect(X, V/G)$ can be identified with the data of a $G$-bundle $\Cal{E}$ over $X_{\kappa}$ plus a $G$-equivariant map $s \co \Cal{E} \rightarrow V$ lying over the identity map on $X_{\kappa}$. 
It is explained in \cite[\S 4.14]{Ngo10} that the tangent space to $ \Sect(X, V/G)$ at this $\kappa$-point is 
\begin{equation}\label{eq: Hitchin tangent space}
\upH^0(X_{\kappa}, \underbrace{\Cal{E} \times^{G} \Lie(G/X)}_{\deg -1} \xrightarrow{\a_{s}} \underbrace{\Cal{E} \times^{G}  V}_{\deg 0})
\end{equation}
where the map $\a_{s}: \Cal{E} \times^{G} \Lie(G/X) \rightarrow \Cal{E} \times^{G}  V$ is given by the action of $\Lie(G/X)$ on $s$ (so that its fiber over $x \in X$ is identified with \eqref{eq: tangent complex 1} upon choosing a trivialization of $\Cal{E}$ at $x$), and the obstructions to deformation lie in 
\begin{equation}\label{eq: Hitchin obstruction space}
\upH^1(X_{\kappa}, \underbrace{\Cal{E} \times^{G} \Lie(G/X)}_{\deg -1} \xrightarrow{\a_{s}} \underbrace{\Cal{E} \times^{G}  V}_{\deg 0}).
\end{equation}
In particular, $ \Sect(X, V/G)$ is smooth at $\kappa$-points where \eqref{eq: Hitchin obstruction space} vanishes. 

\begin{prop}\label{prop: Hitchin smooth}  
(1)  Let $BH_1$ be a smooth $m$-framed gerbe induced from any homomorphism of smooth group schemes $H_{1}\to \GL(m)'$ over $X$, and $BH_2$ the tautological $n$-framed gerbe. Then the stack $\Cal{M}_{H_1, H_2}^{\c}$ is smooth. 

(2) Let $BH_1=B\GL(m)'$ be the tautologoical $m$-framed gerbe and $BH_2 = BU(n)_{\LL}$. Then the stack $\Cal{M}_{H_1, H_2}^{\ns}$ is smooth. 

\end{prop}

\begin{proof} 
It is immediate from the definitions that $\cM_{H_1, H_2}$ is a special case of $\Sect(X, V/G)$ where $G = H_1 \times H_2$, and $V$ is the vector bundle of homomorphisms from the standard representation of $\GL(m)'$ inflated to $H_1$ via the given map $H_1 \rightarrow \GL(m)'$, to the standard representation of $\GL(n)'$ inflated to $H_2$ similarly. 

We will show that the obstruction group to $\cM_{H_1, H_2}$ vanishes at any geometric point of $\Cal{M}_{H_1, H_2}^{\c}$. Consider a geometric point $\Spec \ol{\kappa} \rightarrow \Cal{M}_{H_1, H_2}^{\c}$, which is identified with the data of an $H_1$-torsor $\cE$, an $H_2$-torsor $\cF$ (here the notation differs slightly from \S\ref{ss:Sht sim}, where $\cF$ denoted the associated vector bundle), and an injective map of the associated vector bundles $t \co \Vect (\cE) \rightarrow \Vect(\cF)$. Specializing \eqref{eq: Hitchin obstruction space} to this situation, the obstruction group is 
\begin{equation}\label{eq: obs group}
\upH^1 (X_{\ol{\kappa}}, \underbrace{\cE \times^{H_{1}} \Lie(H_1/X)\oplus \cF\times^{H_2}  \Lie(H_2/X)}_{\deg -1} \xrightarrow{\a_{t}} \underbrace{\cHom(\Vect(\cE), \Vect(\cF))}_{\deg 0} ).
\end{equation}
When $BH_{2}=BU(n)_{\LL}$, $\Lie(H_{2}/X)$ is not a priori defined. In this case, $\Vect(\cF)$ is equipped with an $\LL$-twisted Hermitian form $h$, and we understand $\cF\times^{H_2}  \Lie(H_2/X)$ as the vector bundle $\cEnd^{\mrm{asa}}(\Vect(\cF))$ (over $X_{\ol{\kappa}}$) of anti-self-adjoint endomorphisms of $\Vect(\cF)$, i.e., locally $B:\Vect(\cF)\to \Vect(\cF)$ such that $h(Bx,y)+h(x,By)=0$ for $x,y\in\Vect(\cF)$. 

Returning the generality of (1) and (2), the differential ``$\a_{t}$'' is given by $(A,B)\mapsto -tA+Bt$, where $A\in \cE \times^{H_{1}}\Lie(H_{1}/X), B\in  \cF \times^{H_{2}}\Lie(H_{2}/X)$. Since the coherent cohomology of a torsion sheaf on a curve $X$ vanishes in positive cohomological degrees, it therefore suffices to show that the cokernel of the differential ``$\a_{t}$'' is torsion, or in other words $\a_{t}$ is generically surjective. 


Let $V$ (resp. $U$) be the generic fiber of $\Vect(\cF)$ (resp. $\Vect(\cE)$), a vector space of rank $n$ (resp. $m$) over $K'=F'\ot_{k}\ol{\kappa}$. Let $K=F\ot_{k}\ol{\kappa}$. Let $T$ be the generic fiber of $t$. By the assumption that the $\ol{\kappa}$-point lies in $\cM_{H_1,H_2}^{\c}$, $T: U\to V$ is a $K'$-linear injective map.

In case (1),  the generic fiber of $\cF\times^{H_{2}}\Lie(H_{2}/X)$ is $\End_{K'}(V)$. The map $\End_{K'}(V)\to \Hom_{K'}(U,V)$ given by $B\mapsto BT$ is already surjective since $T$ is injective. This shows that $\a_{t}$ is generically surjective, and the obstructions vanish,  as desired.

In case (2), we argue as follows. In this case,  upon trivializing the generic fiber of $\om_{X}\ot \LL$, $V$ carries a Hermitian form $h: V\ot_{K'} \s^{*}V\to K'$. The generic fiber of $\cE\times^{H_{1}}\Lie(H_{1}/X)$ is $\End_{K'}(U)$ and   the generic fiber of $\cF\times^{H_{2}}\Lie(H_{2}/X)$ can be identified with the $K$-vector space $\End^{\mrm{asa}}_{K'}(V)$ of anti-self-adjoint endomorphisms $B:V\to V$. By the assumption that $\ol{\kappa}$-point lies over the non-degenerate locus $\cA_{H_1}^{\ns} \inj \cA_{H_1}$, $T$ is injective and $h|T(U)$ is non-degenerate. Therefore we may assume $(V,h)=(U,h_{U})\op (W,h_{W})$ is a direct sum of two non-degenerate Hermitian spaces, and $T$ is the inclusion of $U$ in $V$. We have 
$$
 \Hom(U,U\op W)\cong\End(U)\op \Hom(U,W),
 $$
 and 
$$
\End^{\mrm{asa}}(V)\cong \End^{asa}(U)\op \Hom(U,W)\op \End^{\mrm{asa}}(W),
$$ where the last isomorphism is given by $B\mapsto (\pr_{U}(B|_{U}), \pr_{W}(B|_{U}), \pr_{W}(B|_{W}))$. Under these identifications,  the generic fiber of $\a_{t}$ then takes the form
\begin{eqnarray*}
\End(U)\op \End^{\mrm{asa}}(U)\op \Hom(U,W)\op \End^{\mrm{asa}}(W)&\to&\End(U)\op \Hom(U,W)\\
(A,B_{1}, B_{2}, B_{3})&\mapsto &(-A+B_{1},-B_{2})
\end{eqnarray*}
 from which we see that $\a_{t}$ is generically surjective.

\end{proof}


The following variant will be used below in Lemma \ref{lem: Hk^ns smooth}. Following the proof of \cite[Lemma 8.14]{FYZ}, we define an ``$\LL$-twisted almost-Hermitian bundle with defect at $(x', \sigma(x'))$'' to be the data of a vector bundle $\cF^{\flat}$ on $X' \times S$ equipped with a Hermitian map $h \co \cF^{\flat} \inj \sigma^* (\cF^{\flat})^{\vee} \otimes \nu^* \LL$ with cokernel an invertible sheaf on the union of the graphs of $x'$ and $\sigma(x')$. Let $\cM^{\flat}_{H_1, U(n), \LL}$ be the Hitchin stack parametrizing $x'\in X', \cE \in \Sect(X, BH_1)$, and $\LL$-twisted almost-Hermitian bundle $\cF^{\flat}$ with defect at $(x', \sigma(x'))$ and $t \co \Vect(\cE) \rightarrow \cF^{\flat}$. There is a Hitchin fibration $\cM^{\flat}_{H_1, U(n), \LL} \rightarrow \cA_{H_1, \LL}$ defined completely analogously to Definition \ref{def: hitchin fibration}.

\begin{lemma}\label{lem: almost Hermitian Hitchin smooth}  
The map $\cM^{\flat}_{H_1, U(n), \LL}|_{\cA_{H_1}^{\ns}}\to X'$ is smooth if $BH_1$ is the classifying stack of a smooth group scheme $H_1/X$. 
\end{lemma}

\begin{proof} The proof is similar to that of Proposition \ref{prop: Hitchin smooth}. The obstruction group for the map $\cM^{\flat}_{H_1, U(n), \LL}\to X'$ at $(\cF^{\flat},x,t) \in \cM^{\flat}_{H_1, U(n), \LL}|_{\cA_{H_1}^{\ns}}(\ol{\kappa})$ lying over $x \in X'(\ol{\kappa})$ is 
\begin{equation}\label{eq: obs group Hk}
\upH^1 (X_{\ol{\kappa}}, \underbrace{\cE \times^{H_{1}} \Lie(H_1/X)\oplus \cEnd^{\mrm{asa}}(\cF^\flat,h)}_{\deg -1} \xrightarrow{\a_{t}} \underbrace{\cHom(\Vect(\cE), \cF^{\flat})}_{\deg 0} )
\end{equation}
where $\cEnd^{\mrm{asa}}(\cF^\flat,h)$ is the sheaf of endomorphisms of $\cF^\flat$ compatible with $h$. As in the proof of Proposition \ref{prop: Hitchin smooth}(2), it suffices to show that $\alpha_t$ is generically surjective when $(\cF^{\flat}, x,t)$ lies over the non-singular locus of $\cA_{H_1}^{\ns}$. This can be checked on the generic fiber, where $h$ is an isomorphism and we may apply the same argument as in proof of Proposition \ref{prop: Hitchin smooth}(2).
\end{proof}

\subsection{Hecke stacks for Hitchin stacks} 
\begin{defn}[Hecke stacks for Hitchin spaces]\label{defn: hitchin hecke}
Let $BH_1$ be an $m$-framed gerbe over $X$ and $BH_2 \rightarrow B\GL(n)'$ be as in Definition \ref{defn: hitchin}. Further assume that $BH_2$ is of unitary type or $B\GL(n)'$, so that $\Hk_{H_2}^r$ has been defined (cf. \S \ref{sssec: hecke for unitary gerbe} for the first case, and \S \ref{defn: Hk GL(n)'} for the second case). For $r \geq 0$, we define $\Hk_{\Cal{M}_{H_1, H_2}}^{r}$ to be the stack whose $S$-points are given by the groupoid of the following data: 	
\begin{enumerate}
\item $(\{x'_{i}\}, \Cal{F}_0 \dashrightarrow \Cal{F}_1 \dashrightarrow \ldots \dashrightarrow \Cal{F}_r) \in \Hk^{r}_{H_2}(S)$. 
\item $\Cal{E}$ an $H_1$-torsor on $X\times S$. 
\item Maps $t_i \co \Vect(\Cal{E} )  \rightarrow \Vect(\Cal{F}_i )$, fitting into the commutative diagram below. 
\end{enumerate}
\[
\begin{tikzcd}
\Vect(\Cal{E} )  \ar[r, equals] \ar[d, "t_0"] & \Vect(\Cal{E} ) \ar[r, equals]  \ar[d, "t_1"] & \ldots \ar[r, equals] \ar[d] & \Vect(\Cal{E} )  \ar[d, "t_r"]  \\
\Vect(\Cal{F}_0 )  \ar[r, dashed] & 
\Vect(\Cal{F}_1 ) \ar[r, dashed] & \ldots \ar[r, dashed]  &  \Vect(\Cal{F}_r )
\end{tikzcd}
\]
(The dashed notation follows \cite[Definition 6.5]{FYZ}.) Let $\pr_i \co \Hk_{\Cal{M}_{H_1, H_2}}^{r}  \rightarrow \Cal{M}_{H_1, H_2}$ be the map recording $(\cE, \cF_{i}, t_{i})$, for $0\le i\le r$. 

We define $\Hk_{\Cal{M}_{H_1, H_2}}^{r ,\circ } \subset \Hk_{\Cal{M}_{H_1, H_2}}^{r}$ to be the fibered product of $\pr_0$ (equivalently, any $\pr_i$) with $\Cal{M}_{H_1, H_2}^{\c} \inj \Cal{M}_{H_1, H_2}$, and $\pr_i^{\c} \co \Hk_{\Cal{M}_{H_1, H_2}}^{r, \c} \to \Cal{M}_{H_1, H_2}^{\c} $ to be the restriction of $\pr_i$.

\end{defn}

Let $(x', \cE, \Cal{F}_0 \dashrightarrow \Cal{F}_1, t_0, t_1)$ be an $S$-point of $\Hk^{1}_{\cM_{H_1, H_2}}$. By means of the given rational isomorphism between $\Cal{F}_0$ and $\Cal{F}_1$, we may form the intersection $\cF^{\flat}_{1/2}=\Cal{F}_0 \cap \Cal{F}_1$, which is an $S$-point of $\Bun_{\GL(n)'}$. By definition, the maps $t_0, t_1$ factor through a unique map $t^{\flat} \co \Vect(\cE) \rightarrow \cF^{\flat}_{1/2}$. The data of $(\cE, \cF^{\flat}_{1/2}, t)$ determines an $S$-point of $\cM_{H_1, H_2}$. We define a map $\pr_{1/2}: \Hk^{1}_{\cM_{H_1, H_2}} \rightarrow \cM_{H_1, H_2} $ by $\pr_{1/2}(x', \cE, \Cal{F}_0 \dashrightarrow \Cal{F}_1, t_0, t_1)=(\cE, \cF^{\flat}_{1/2}, t)$.  

\begin{lemma}\label{l:HkM sm}Let $BH_1$ be a smooth $m$-framed gerbe induced from any homomorphism of smooth group schemes $H_{1}\to \GL(m)'$ over $X$, and $BH_2$ be the tautological $n$-framed gerbe over $X$.
\begin{enumerate}
\item The map  $(\pr_{1/2},\pr_{X'}): \Hk^{1}_{\cM_{H_1, H_2}} \rightarrow \cM_{H_1, H_2} \times X'$ is smooth and of relative (equi)dimension $2(n-1)$. In particular, $\Hk^{1}_{\cM_{H_1, H_2}^{\c}} $ is smooth. 

\item For any geometric point $\xi\in \Hk^{1}_{\cM^{\c}_{H_1, H_2}}$, the local dimension of $\Hk^{1}_{\cM^{\c}_{H_1, H_2}}$ at $\xi$ satisfies
\begin{equation}\label{dim Hk minus dim Ma}
\dim_{\xi} \Hk^{1}_{\cM^{\c}_{H_1, H_2}}-\dim_{\pr_{i}(\xi)} \cM^{\c}_{H_1, H_2}=2n-1-m, \quad i=0,1.
\end{equation}

\end{enumerate}

\end{lemma}
\begin{proof}
(1) To recover $(x',\cE, \cF_0 \dashrightarrow \cF_1, t_0, t_1)$ from its image $(\cE, \cF^{\flat}_{1/2}, t,x')$ under $(\pr_{1/2}, \pr_{X'})$ is equivalent to giving the datum of a line in each of the fibers of $\cF^{\flat}_{1/2}$ at $x'$ and $\sigma(x')$. Hence $(\pr_{1/2}, \pr_{X'})$ is 
a $\PP^{n-1} \times \PP^{n-1}$-fiber bundle. In particular  it is smooth of relative dimension $2(n-1)$.



Now $\Hk^{1}_{\cM^{\c}_{H_{1},H_{2}}}$ is also the preimage of $\cM^{\c}_{H_{1},H_{2}}$ under $\pr_{1/2}$.  The smoothness of $\Hk^{1}_{\cM^{\c}_{H_{1},H_{2}}}$ follows by combining the relative smoothness with Proposition \ref{prop: Hitchin smooth}. 
 
(2) Let $\xi=(x'_{1}, \cE, \cF_{0}\dashrightarrow \cF_{1}, t_{0}, t_{1})$ be a geometric point of $\Hk^{1}_{\cM^{\c}_{H_{1},H_{2}}}$. Comparing the tangent complexes of $\cM_{H_{1},H_{2}}^{\c}$ at $\pr_{i}(\xi)=(\cE, \cF_{i}, t_{i})\in \cM_{H_{1},H_{2}}^{\c}$ and at $\pr_{1/2}(\xi)=(\cE, \cF^{\flat}_{1/2}, t)\in \cM_{H_{1},H_{2}}^{\c}$ given in the proof of Proposition \ref{prop: Hitchin smooth}, we see that (for $i=0,1$) 
\begin{align}\label{dim diff Mab}
&\dim_{\pr_{i}(\xi)}\cM_{H_{1},H_{2}}^{\c}-\dim_{\pr_{1/2}(\xi)} \cM_{H_{1},H_{2}}^{\c}\notag
\\=&\deg\cHom(\Vect(\cE),\cF_{i})-\deg\cHom(\Vect(\cE),\cF^{\flat}_{1/2})=m.
\end{align}
On the other hand, by (1) we know that 
\begin{align}\label{dim Hk minus dim Mb}
\dim_{\xi}\Hk^{1}_{\cM^{\c}_{H_{1},H_{2}}}=\dim_{\pr_{1/2}(\xi)} \cM_{H_{1},H_{2}}^{\c}+2n-1.
\end{align}
Combining \eqref{dim diff Mab} and \eqref{dim Hk minus dim Mb} we get \eqref{dim Hk minus dim Ma}.

\end{proof}

When $BH_{2}=BU(n)_{\LL}$, the composition $\Hk_{\cM_{H_{1},H_{2}}}^{r}\xr{\pr_{i}}\cM_{H_{1},H_{2}}\xr{f}\cA_{H_{1},\LL}$ is independent of $i$, so that  $\Hk_{\cM_{H_{1},H_{2}}}^{r}$ has a well-defined map to $\cA_{H_{1},\LL}$.  

 \begin{lemma}\label{lem: Hk^ns smooth}
Let $BH_1$ be the tautological $m$-framed gerbe and $BH_2 = BU(n)_{\LL}$ with the standard $n$-framed structure. Define $\Hk_{\cM_{H_{1},H_{2}}^{\ns}}^{1}:=\Hk_{\cM_{H_{1},H_{2}}}^{1}|_{\cA^{\ns}_{H_{1},\LL}}$. 
Then: 
\begin{enumerate}
\item The projection map $\Hk^{1}_{\cM_{H_1, H_2}^{\ns}} \rightarrow X'$ is smooth. In particular, $\Hk^{1}_{\cM_{H_1, H_2}^{\ns}}$ is smooth.
\item For any geometric point $\xi\in \Hk^{1}_{\cM_{H_1, H_2}^{\ns}}$, the local dimension of $\Hk^{1}_{\cM_{H_1, H_2}^{\ns}}$ at $\xi$ satisfies
\begin{equation}\label{dim Hk minus dim M ns}
\dim_{\xi} \Hk^{1}_{\cM_{H_1, H_2}^{\ns}}-\dim_{\pr_{i}(\xi)} \cM_{H_1, H_2}^{\ns}=n-m, \quad i=0,1.
\end{equation}
\end{enumerate}
\end{lemma}

\begin{proof}
(1) We claim that $\Hk^1_{\cM_{H_1, H_2}}\cong \cM^{\flat}_{H_1, U(n),\LL}$, which is defined in Lemma \ref{lem: almost Hermitian Hitchin smooth}. If we admit this then the assertion  follows from the smoothness of $\cM^{\flat}_{H_1, U(n),\LL} |_{\cA^{\ns}_{H_{1},\LL}} \to X'$ established in Lemma \ref{lem: almost Hermitian Hitchin smooth}. So it suffices to establish the claim. 

We define a map $ \Hk^1_{\cM_{H_1, H_2}} \rightarrow \cM^{\flat}_{H_1, U(n),\LL} $. Take $(x', \Cal{F}_0 \dashrightarrow \Cal{F}_1) \in \Hk_{U(n), \LL}^1$ and set $\cF_{1/2}^{\flat} := \cF_0 \cap \cF_1$. The generically compatible Hermitian structures on $\Cal{F}_0$ and $\Cal{F}_1$ equip $\cF_{1/2}^{\flat}$ with an $(\LL$-twisted) almost Hermitian structure (cf. \cite[Proof of Lemma 8.14]{FYZ} for the definition of ``almost Hermitian'') with defect at $(x', \sigma(x'))$.  

 Given $(\Cal{F}^{\flat},h^{\flat})$ almost Hermitian with defect at $(x', \s(x'))$, define $\cF_{0}$ (resp. $\cF_{1}$) as the upper modification of $\cF^{\flat}$ at $x'$ (resp. $\s(x')$) inside $\s^{*}(\cF^{\flat})^{\vee}\ot \nu^{*}\LL$.  It is easy to see that this defines the inverse map.

(2) Let $\cM^{\flat}_{H_1, U(n),\LL,x'}$ be the fiber of $\cM^{\flat}_{H_1, U(n),\LL}$ over $x' \in X'$. Consider a geometric point
$$
\xi= (\cE, (\cF_{1/2}, h), t_{1/2} \co \Vect(\cE) \rightarrow \cF_{1/2}^\flat)
$$
of $\cM^{\flat}_{H_1, U(n),\LL,x}|_{\cA^{\ns}_{H_{1},\LL}}$. By the smoothness established in (1) and Proposition \ref{prop: Hitchin smooth}, the local dimensions at $\xi$ and $\pr_i(\xi)$ may be computed as the Euler characteristic of the respective tangent complexes. Comparing the tangent complexes of $\cM^{\flat}_{H_1, U(n),\LL,x'}$ at $\xi$ and of $\cM_{H_{1},H_{2}}^{\ns}$ at $\pr_{i}(\xi)=(\cE, \cF_{i}, t_{i})\in \cM_{H_{1},H_{2}}^{\ns}$ using the proof of Proposition \ref{prop: Hitchin smooth}, we see that (for $i=0,1$) 
\begin{align*}
&\dim_{\xi} \cM^{\flat}_{H_1, U(n),\LL,x'} - \dim_{\pr_{i}(\xi)} \cM_{H_{1},H_{2}}^{\ns}\\ =&  - \deg \cEnd^{\mrm{asa}}(\cF_{1/2}^{\flat})  + \deg\cHom(\Vect(\cE),\cF^{\flat}_{1/2}) -  \deg\cHom(\Vect(\cE),\cF_{i}) 
\end{align*}
where $\cEnd^{\mrm{asa}}(\cF_{1/2}^{\flat})$ is the space of anti-self-adjoint morphisms with respect to the Hermitian map $h \co \cF_{1/2}^{\flat} \inj \sigma^* (\cF_{1/2}^{\flat} )^{\vee} \otimes \LL$. Here we have used $ \deg \cEnd^{\mrm{asa}}(\cF_{i}) =0$.

We have $\deg\cHom(\Vect(\cE),\cF^{\flat}_{1/2}) -  \deg\cHom(\Vect(\cE),\cF_{i}) = - m$, as in the proof of Lemma \ref{l:HkM sm}. To compute $\deg \cEnd^{\mrm{asa}}(\cF_{1/2}^{\flat})$, we reduce to the case where the double cover is split, by base changing along $X' \rightarrow X$. In that case, $X' = X \sqcup X$ and we may assume $x'$ lies in the first copy of $X$ and its image in $X$ is denoted by $x$. Then the datum of $\cF_{1/2}^{\flat}$ may be identified with a pair of vector bundles $\cF^{(1)}, \cF^{(2)}$ on $X$ and a map $h \co \cF^{(1)} \rightarrow (\cF^{(2)})^{\vee} \otimes \LL$ whose cokernel is flat of length $1$ along the graph of $x$.  Then $\cEnd^{\mrm{asa}}(\cF_{1/2}^{\flat})$ consists of endomorphisms $(B_1,B_2)$ of $\cF^{(1)} \boxplus \cF^{(2)}$ such that for every local section $v^{(1)}\in\cF^{(1)}(U)$ and $v^{(2)} \in \cF^{(2)}(U)$  on  an open subset $U \subset X$, we have $h(B_1 v^{(1)} ) = - B_2^{\vee} ( h v^{(1)})$. Hence any such endomorphism is determined by its restriction to $\cF^{(2)}$, giving an injection $\cEnd^{\mrm{asa}}(\cF_{1/2}^{\flat}) \inj \cEnd( (\cF^{(2)})^{\vee} \otimes \LL).$ Let us abbreviate $\wt{\cF}^{(2)} := (\cF^{(2)})^{\vee} \otimes \LL$, which we remind is a rank $n$ vector bundle on $X$. The image of the preceding injection consists of those maps in $\cEnd(\wt{\cF}^{(2)} ) $ preserving $\cF^{(1)}$, viewed as a subsheaf of $\wt{\cF}^{(2)}$, hence the image is equal to  the kernel of the composition of arrows below 
\[
\begin{tikzcd}
\cEnd(\wt{\cF}^{(2)} )  \ar[r, twoheadrightarrow] & \cHom(\wt{\cF}^{(2)}  ,\, \wt{\cF}^{(2)} / \cF^{(1)}) \ar[d] \\
&\cHom(\cF^{(1)},\, \wt{\cF}^{(2)}  / \cF^{(1)})
\end{tikzcd}
\]
The sheaf $\cHom(\cF^{(1)}, \wt{\cF}^{(2)}  / \cF^{(1)})$ is torsion of degree $n$ on $X$. The long exact sequence for $\cHom(-, \wt{\cF}^{(2)} / \cF^{(1)})$ shows that the image of the vertical map is the kernel of the surjection 
$$
\cHom(\cF^{(1)}, \wt{\cF}^{(2)}  / \cF^{(1)}) \rightarrow \cExt^1(\wt{\cF}^{(2)}  / \cF^{(1)}, \wt{\cF}^{(2)}  / \cF^{(1)})
$$
 whose codomain is an invertible sheaf along the graph of $x$. As $\cEnd(\wt{\cF}^{(2)} )$ has degree $0$, we conclude that $\cEnd^{\mrm{asa}}(\cF_{1/2}^{\flat})$ has degree $-(n-1)$. 

\end{proof}

\subsection{Hitchin shtukas}\label{ssec: Hitchin shtukas}
We now discuss a notion of shtukas for Hitchin-type spaces $\Cal{M}_{H_1, H_2}$.

\begin{defn}[Shtukas for Hitchin spaces]\label{defn: hitchin shtuka}
Let $BH_1$ be an $m$-framed gerbe and $BH_2$ an $n$-framed gerbe as in Definition \ref{defn: hitchin hecke}. For $r \geq 0$, we define $\Sht_{\Cal{M}_{H_1,H_2}}^{r}$ as the fibered product 
\begin{equation}\label{eq: shtuka all as intersection}
\begin{tikzcd}
\Sht_{\Cal{M}_{H_1,H_2}}^{r}  \ar[r] \ar[d] & \Hk_{\Cal{M}_{H_1,H_2}}^{r} \ar[d, "{(\pr_0, \pr_r)}"] \\
\Cal{M}_{H_1,H_2} \ar[r, "{(\Id, \Frob)}"] & \Cal{M}_{H_1,H_2}  \times \Cal{M}_{H_1,H_2}
\end{tikzcd}
\end{equation}
We define the open substack $\Sht_{\Cal{M}_{H_1,H_2}^{\c}}^{r} \subset \Sht_{\Cal{M}_{H_1,H_2}}^{r}$ as the fibered product 
\begin{equation}\label{eq: shtuka as intersection}
\begin{tikzcd}
\Sht_{\Cal{M}_{H_1,H_2}^{\c}}^{r}  \ar[r] \ar[d] & \Hk_{\Cal{M}_{H_1,H_2}^{\c}}^{r} \ar[d, "{(\pr_0, \pr_r)}"] \\
\Cal{M}_{H_1,H_2}^{\c} \ar[r, "{(\Id, \Frob)}"] & \Cal{M}_{H_1,H_2}^{\c}  \times \Cal{M}_{H_1,H_2}^{\c}
\end{tikzcd}
\end{equation}
Note that $\Sht_{\Cal{M}_{H_1,H_2}^{\c}}^{r} \inj \Sht_{\Cal{M}_{H_1,H_2}}^r$ can be equivalently described as the base change of $\Cal{M}_{H_1,H_2}^{\c} \inj  \Cal{M}_{H_1,H_2} $ against any of the projection maps $\pr_i \co \Sht_{\Cal{M}_{H_1,H_2}}^r \rightarrow \Cal{M}_{H_1,H_2}$. 
\end{defn}

\begin{example}\label{ex: hitchin shtuka U(n)} Let $BH_1$ be the tautological $m$-framed gerbe and $BH_2  = BU(n)_{\LL}$ with the standard $n$-framed structure. Recall that we defined $\ol{\cZ}_{\cE, \LL}^r := [\cZ_{\cE, \LL}^r/(\Aut(\cE)(\F_q))]$, where $\cZ_{\cE, \LL}^r$ was defined in Definition \ref{defn: twisted similitude shtuka definitions}. Then we have
\[
\Sht_{\cM_{H_1, H_2}}^r = \coprod_{\cE \in \Bun_{\GL(m)'}(k)} \ol{\cZ}^{r}_{\cE,\LL}. 
\]
and for $\ol{\cZ}^{ r,\circ}_{\cE,\LL} := [\cZ^{ r,\circ}_{\cE,\LL}/(\Aut(\cE)(\F_q))]$, we have
\begin{equation}
\Sht_{\cM_{H_1, H_2}^{\c}}^r = \coprod_{\cE \in \Bun_{\GL(m)'}(k)} \ol{\cZ}^{ r,\circ}_{\cE,\LL}. 
\end{equation}

When $\LL = \cO_X$, $\Hk_{\Cal{M}_{H_1, H_2}}^{r, \circ }$ (resp. $\Hk_{\Cal{M}_{H_1, H_2}}^{r}$) is the stack denoted by $\Hk_{\cM(m,n)}^r$ (resp. $\Hk_{\cM^{\all}(m,n)}^r$) in \cite[\S 8]{FYZ}, and $\Sht_{\Cal{M}_{H_1,H_2}^{\c}}^{r}$ (resp. $\Sht_{\Cal{M}_{H_1,H_2}}^{r}$) is the stack denoted by $\Sht_{\Cal{M}(m,n)}^{r}$ (resp. $\Sht_{\Cal{M}^{\all}(m,n)}^{r}$) in \cite[\S 8]{FYZ}. 
\end{example}

\begin{example}\label{ex: hitchin shtuka GL(n)'}
Let $BH_1$ be the tautological $m$-framed gerbe and $BH_2$ be the tautological $n$-framed gerbe. Recall that we defined $\ol{\cZ}_{\cE, \GL(n)'}^r := [\cZ_{\cE, \GL(n)'}^r/(\Aut(\cE)(\F_q))]$, where $\cZ_{\cE, \GL(n)'}^r$ was defined in Definition \ref{def: Z}. Then we have
\[
\Sht_{\cM_{H_1, H_2}}^r = \coprod_{\cE \in \Bun_{\GL(m)'}(k)} \ol{\cZ}^r_{\cE,\GL(n)'}, 
\]
and for $\ol{\cZ}_{\cE, \GL(n)'}^{r, \circ} := [\cZ_{\cE, \GL(n)'}^{r, \circ}/(\Aut(\cE)(\F_q))]$ we have 
\begin{equation}\label{eq: sht_M decomposition}
\Sht_{\cM_{H_1, H_2}^{\c}}^r = \coprod_{\cE \in \Bun_{\GL(m)'}(k)} \ol{\cZ}^{r, \circ}_{\cE, \GL(n)'}.
\end{equation}
\end{example}

\begin{remark}\label{rem: Sht H_2}
Note that if we take $H_1$ to be the trivial $0$-framed gerbe, then $\cM_{H_1, H_2} = \Bun_{H_2} = \Sect(X, BH_2)$. Furthermore, if $BH_2$ is of unitary type  or $B\GL(n)'$ then the definition of $\Hk^r_{\cM_{H_1, H_2}}$ (resp. $\Sht^r_{\cM_{H_1, H_2}}$) above specializes to $\Hk^r_{H_2}$ (resp. $\Sht^r_{H_2}$) as defined in \S \ref{sssec: hecke for unitary gerbe}.
\end{remark}

\subsection{Cycle classes from Hitchin shtukas}\label{ss:ShtM cycle}

\begin{defn}\label{def:Phi} For any stack $S$ over $k$ we define a morphism
\begin{equation}
\Phi^{r}_{S}: S^{r+1}\to S^{2r+2}
\end{equation}
by the formula $\Phi^{r}_{S}(\xi_{0},\cdots, \xi_{r})=(\xi_{0},\xi_{1},\xi_{1},\xi_{2},\xi_{2},\cdots, \xi_{r-1}, \xi_r, \xi_r, \Frob(\xi_{0}))$. When $r$ is fixed in the context, we simply write $\Phi_{S}$.
\end{defn}

We rewrite $\Sht^{r}_{\cM_{H_1,H_2}^{\c}}$ as the fiber product
\begin{equation}\label{rewrite ShtM}
\xymatrix{\Sht^{r}_{\cM_{H_1,H_2}^{\c}} \ar[r]\ar[d] & (\Hk^{1}_{\cM_{H_1,H_2}^{\c}})^{r} \times \cM_{H_1,H_2}^{\c} \ar[d]^{(\pr_{0}^{\c},\pr_{1}^{\c})^{r} \times \Delta} \\
(\cM_{H_1,H_2}^{\c})^{r+1}\ar[r]^{\Phi^{r}_{\cM_{H_1,H_2}^{\c}}} & (\cM_{H_1,H_2}^{\c})^{2r+2}}
\end{equation}

\begin{defn}\label{defn: sht_M class for GL(n)'} Let $BH_1 \rightarrow B\GL(m)'$ be induced from any homomorphism of smooth group schemes $H_{1}\to \GL(m)'$ over $X$, and $BH_2$ the tautological $n$-framed gerbe. By Lemma \ref{l:HkM sm}, the fundamental class of $(\Hk^{1}_{\cM_{H_1,\GL(n)'}^{\c}})^{r}$ is defined, which we denote by $[(\Hk^{1}_{\cM_{H_1,\GL(n)'}^{\c}})^{r}]^{\nai}\in \Ch_{*}((\Hk^{1}_{\cM_{H_1,\GL(n)'}^{\c}})^{r})$. Then we define the cycle class $[\Sht^{r}_{\cM_{H_1, \GL(n)'}^{\c}}]\in \Ch_{*}(\Sht^{r}_{\cM_{H_1,\GL(n)'}^{\c}})$ as the the image of $[(\Hk^{1}_{\cM_{H_1,\GL(n)'}^{\c}})^{r} \times \cM_{H_1,\GL(n)'}^{\c}]^{\nai}$  under the refined  Gysin map along $\Phi^{r}:(\cM_{H_1,\GL(n)'}^{\c})^{r+1}\to (\cM_{H_1,\GL(n)'}^{\c})^{2r+2}$ (which is defined since $\cM_{H_1,\GL(n)'}^{\c}$ is smooth by Proposition \ref{prop: Hitchin smooth} -- see \cite[\S A.1.4]{YZ})  
\begin{equation}
[\Sht^{r}_{\cM_{H_1,\GL(n)'}^{\c}}]:=(\Phi^{r}_{\cM_{H_1,\GL(n)'}^{\c}})^{!}[(\Hk^{1}_{\cM_{H_1,\GL(n)'}^{\c}})^{r} \times \cM_{H_1,\GL(n)'}^{\c} ]^{\nai} \in \Ch_{*}(\Sht^{r}_{\cM_{H_1,\GL(n)'}^{\c}}).
\end{equation}
In particular, when $BH_1 \xrightarrow{=} B\GL(m)'$, the dimension formula in Lemma \ref{l:HkM sm} implies that $[\Sht^{r}_{\cM_{H_1,\GL(n)'}^{\c}}] \in \Ch_{r(2n-1)-rm}(\Sht^{r}_{\cM_{H_1,\GL(n)'}^{\c}})$. 
\end{defn}

\begin{remark}
Definition \ref{defn: sht_M class for GL(n)'} will be used in the next section to define cycle classes $[\cZ_{\cE}^r(a)]$. Even though we are in some sense more interested in the case $BH_2 = BU(n)_{\LL}$, for the purpose of constructing cycle classes corresponding to singular $a$, it was crucial to take $BH_2 = B\GL(n)'$ in Definition \ref{defn: sht_M class for GL(n)'}, because Proposition \ref{prop: Hitchin smooth} gives smoothness of the $\cM_{H_1, \GL(n)'}^{\circ}$ even over the singular part of the Hitchin base. Because we lack such control when $BH_2 = BU(n)_{\LL}$, we cannot make an analogous definition in that case. 
\end{remark}






\section{Formulation of the modularity conjecture}\label{sec: formulation of the conjecture}
  
Let $\cE$ be a vector bundle on $X'$ of rank $m$, and let $\LL$ be a line bundle on $X$. For any $a\in \cA_{\cE,\LL}(k)$, we have defined a special cycle $\cZ^{r}_{\cE,\LL}(a)\to \Sht^{r}_{U(n),\LL}$, cf. Definition \ref{defn: twisted similitude shtuka definitions}. The goal of this section is to construct a virtual fundamental class $[\cZ^{r}_{\cE,\LL}(a)]\in \Ch_{r(n-m)}(\cZ^{r}_{\cE,\LL}(a))$ for every $a$, and formulate a conjecture that a generating series of such cycle classes is modular. We note that $\dim \cZ^r_{\cE,\LL}(a)$ can differ significantly from $r(n-m)$ in general, so we really need a virtual fundamental class. 

It turns out that when $a$ is non-singular, $[\cZ^{r}_{\cE,\LL}(a)]$ can be defined directly using Hitchin stacks. For possibly singular $a$, we define $[\cZ^{r}_{\cE,\LL}(a)]$ in two steps. First, we define the cycle class on the open-closed substack $\cZ^{r}_{\cE,\LL}(a)^{\c}$ consisting of generically injective maps from $\cE$.  Next, on the rest of the connected components of $\cZ^{r}_{\cE,\LL}(a)$, we reduce to the case of an already-defined cycle class (of smaller corank), and cap it with an appropriate Chern polynomial coming from tautological bundles over $\Sht^{r}_{U(n),\LL}$. (Later in \S \ref{sec: VFC}, specifically Theorem \ref{thm: VFC from derived shtuka}, we will see how this recipe arises from a natural derived enhancement of $\cZ^{r}_{\cE,\LL}(a)$.) 

In this section, we fix a similitude line bundle $\LL$ on $X$ and consider $\LL$-twisted Hermitian bundles. When there is no confusion we will omit $\LL$ from the notation, e.g., we write $\cA_{\cE}$ and $\cZ^{r}_{\cE}(a)$ for $\cA_{\cE,\LL}$ and $\cZ^{r}_{\cE,\LL}(a)$.

In Sections \S \ref{ssec: kernel cycles} -- \S \ref{ssec: test intersection numbers}, we will assume that $\nu \co X' \rightarrow X$ is a \emph{non-split} double cover. Then in \S \ref{ssec: split case conj} we will treat the split case. Although the two cases could in principle be treated uniformly, this would make the language complicated (the issue is that for $X' = X \coprod X$, we need to consider vector bundles $\cE'$ on $X'$ which have different ranks on the two connected components; this is necessary in order to define the virtual fundamental class even for a starting $\cE$ which has the same rank $m$ on the two components). Furthermore, in the split case we can formulate \emph{more refined} conjectures.  

\subsection{Decomposition according to kernel}\label{ssec: kernel cycles}

Let $\cK \subset \subset \cE$ be a sub-bundle of $\cE$ (the notation $\cK \subset\subset \cE$ means that $\cK$ is a sub-bundle of $\cE$, i.e., the quotient $\cE  / \cK$ is a vector bundle) and $\ov\cE=\cE/\cK$ be the quotient bundle. We define the closed substack $\cZ_{\cE}^r[\cK] \subset \cZ_{\cE}^r$ to parametrize those $(\cE \xrightarrow{t_i} \cF_i)$ such that $\cK \subset \ker (t_i)$. We define $\cZ_{\cE}^r[\cK]^{\circ} \subset  \cZ_{\cE}^r[\cK]$ to be the open substack where $\ker(t_i) = \cK$. 

Each $\cZ_{\cE}^r[\cK]^{\circ}$ is locally closed in $\cZ_{\cE}^r$. It is clear that $\cZ_{\cE}^r[\cK]^{\circ}$ for varying $\cK$ form a partition of $\cZ_{\cE}^{r}$. 
In particular,
\begin{eqnarray}
\cZ_{\cE}^r[0]=\cZ^{r}_{\cE}; \quad \cZ_{\cE}^r[0]^{\circ}=\cZ_{\cE}^{r,\c}\\\label{ZK circ}\cZ_{\cE}^{r}[\cK]^{\c}=\cZ_{\cE}^r[\cK]\setminus (\bigcup_{\cK\subsetneq\cK'}\cZ^{r}_{\cE}[\cK']).
\end{eqnarray}
We will show that $\cZ_{\cE}^r[\cK]^{\c}$ are in fact open-closed in $\cZ_{\cE}^r$.

\begin{lemma}\label{l:iE open-closed}
The substack $\cZ^r_{\cE}[\cK] \subset \cZ^r_{\cE}$ is open-closed. 
\end{lemma}
\begin{proof}
Consider the natural map
\begin{equation}
r_{\cK}: \cZ^{r}_{\cE}\to \cZ^{r}_{\cK}
\end{equation}
by restricting $t_{i}: \cE\to \cF_{i}$ to $\cK$. Let  $z:\Sht^{r}_{U(n)} \cong \cZ^r_{\cK}[\cK] \incl \cZ^{r}_{\cK}$ be the locus of zero maps $\cK \to \cF_{i}$ (for varying $\{\cF_{i}\}\in \Sht^{r}_{U(n)}$). Its complement is the union of $\cZ^{r}_{\cK}(a)$ for non-zero $a\in \cA_{\cK}(k) $  and $\cZ^{r}_{\cK}(0)^{*}$. Note that $\cZ^{r}_{\cK}(a)$ is open-closed, and $\cZ^{r}_{\cK}(0)^{*}$ is proper over $\Sht^{r}_{U(n)}$ by the same proof as for \cite[Proposition 7.5]{FYZ}. 
Therefore $\cZ^{r}_{\cK}(0)^{*}\incl \cZ^{r}_{\cK}(0)$ is open-closed and $z$ is open-closed. The inclusion $\cZ_{\cE}^r[\cK] \inj \cZ_{\cE}^r$ is the base change of $z$ along $r_{\cK
}$, hence also open-closed.
\end{proof}

\begin{lemma}\label{lem: kernel circ open-closed} The substack $\cZ^{r}_{\cE}[\cK]^{\c} \subset \cZ^{r}_{\cE}$ is open-closed. In particular, $\cZ^{r,\c}_{\cE}$ is open-closed in $\cZ^{r}_{\cE}$.
\end{lemma}
\begin{proof}
Combine Lemma \ref{l:iE open-closed} with \eqref{ZK circ}. 
\end{proof}


Thus we have a decomposition of $\cZ^{r}_{\cE}$ in open-closed substacks
\begin{equation}\label{ZEK decomp}
\cZ_{\cE}^r = \coprod_{\text{sub-bundles }\cK \subset\subset \cE} \cZ_{\cE}^r[\cK]^{ \circ}.
\end{equation}

\begin{remark}\label{rem: identification with quotient special cycle}
For a sub-bundle $\cK \subset \subset\cE$, there is an identification over $\Sht^{r}_{U(n)}$
\begin{equation}\label{Z quot}
\cZ^{r}_{\cE/\cK} \xrightarrow{\sim} \cZ^{r}_{\cE}[\cK]
\end{equation}
given by inflating $t_i \co \cE/\cK  \rightarrow \cF_i$ along $\cE \surj \cE/\cK$. It restricts to an isomorphism
\begin{equation}\label{Z quot circ}
\cZ^{r,\c}_{\cE/\cK} \xrightarrow{\sim} \cZ^{r}_{\cE}[\cK]^{\c}.
\end{equation}
\end{remark}

\subsection{The cycle class $[\cZ_{\cE}^{r,\c}]$}\label{ssec: circ cycle class}


Consider $\cM_{\GL(m)', \GL(n)'}^{\c}$, which is smooth by Proposition \ref{prop: Hitchin smooth}. Form the stack of Hitchin-shtukas $\Sht^r_{\cM_{\GL(m)', \GL(n)'}^{\c}}$. In Definition \ref{defn: sht_M class for GL(n)'} we have defined a cycle class
\begin{equation}
[\Sht^r_{\cM_{\GL(m)', \GL(n)'}^{\c}}]\in \Ch_{r(2n-1-m)}(\Sht^r_{\cM_{\GL(m)', \GL(n)'}^{\c}}).
\end{equation}
Note we have a decomposition
\begin{equation}
\Sht^r_{\cM_{\GL(m)', \GL(n)'}^{\c}}=\coprod_{\cE\in \Bun_{\GL(m)'}(k)} \ol{\cZ}^{r, \c}_{\cE, \GL(n)'}.
\end{equation}
We define $[\ol{\cZ}^{r,\c}_{\cE, \GL(n)'}]\in \Ch_{r(2n-1-m)}(\ol{\cZ}^{r,\c}_{\cE, \GL(n)'})$ to be the projection of $[\Sht^r_{\cM_{\GL(m)', \GL(n)'}^{\c}}]$ to the summand indexed by $\cE$, and we define $[\cZ^{r,\c}_{\cE, \GL(n)'}]\in \Ch_{r(2n-1-m)}(\cZ^{r,\c}_{\cE, \GL(n)'})$ to be the pullback of $[\ol{\cZ}^{r,\c}_{\cE, \GL(n)'}]$ via the finite \'{e}tale map $\cZ^{r,\c}_{\cE, \GL(n)'}  \rightarrow \ol{\cZ}^{r,\c}_{\cE, \GL(n)'}$. 


We have a Cartesian diagram from Lemma \ref{lem: special cycle GL(n)' to U(n)}
\[
\xymatrix{\cZ^{r,\c}_{\cE}
\ar[r]\ar[d] & \cZ^{r,\c}_{\cE,\GL(n)'} \ar[d]  \\
\Sht_{U(n)}^r \ar[r]^{u} & \Sht_{\GL(n)'}^r
}
\]
Note that $u$ is a regular local immersion, so that the refined Gysin pullback $u^!$ is defined.

\begin{defn}\label{defn: circ case}
We define
\begin{equation}
[\cZ_{\cE}^{r,\c} ]:=u^{!}[\cZ^{r,\c}_{\cE,\GL(n)'}]\in \Ch_{r(n-m)}(\cZ^{r,\c}_{\cE}).
\end{equation}
\end{defn}

Here we are using the equality
\begin{equation}
r(2n-1-m)-(\dim \Sht^{r}_{\GL(n)'}-\dim \Sht^{r}_{U(n)})=r(n-m)
\end{equation}
to determine the (virtual) dimension of the resulting cycle.

\begin{defn}\label{def: Z quotient version} 
From the Cartesian diagram 
\[
\begin{tikzcd} 
\Sht_{\cM^{\circ}_{\GL(m)', U(n)}}^r \ar[r] \ar[d] & \Sht_{\cM^{\circ}_{\GL(m)', \GL(n)'}}^r \ar[d] \\
\Bun_{\GL(m)'}(k) \times \Sht_{U(n)'}^r \ar[r, "u"] & \Bun_{\GL(m)'}(k) \times \Sht_{\GL(n)'}^r
\end{tikzcd}
\]
we define $[\Sht_{\cM^{\circ}_{\GL(m)', U(n)}}^r]  := u^! [\Sht_{\cM^{\circ}_{\GL(m)', \GL(n)'}}^r ]$, where the RHS was defined in Definition \ref{defn: sht_M class for GL(n)'}. Then we define $[\ol{\cZ}_{\cE}^{r,\circ}] \in \Ch_{r(n-m)}(\ol{\cZ}_{\cE}^{r, \circ})$ to be the projection of $[\Sht_{\cM^{\circ}_{\GL(m)', U(n)}}^r] $ to $\ol{\cZ}_{\cE}^{r, \circ}$ using the decomposition of Example \ref{ex: hitchin shtuka U(n)}. This class $[\ol{\cZ}_{\cE}^{r,\circ}]$ pulls back to the $[\cZ_{\cE}^{r, \circ}]$ from Definition \ref{defn: circ case} under the finite \'{e}tale map $\cZ_{\cE}^{r, \circ} \rightarrow \ol{\cZ}_{\cE}^{r,\circ}$. We introduced it later because it will only be used as an auxiliary device to compare Definition \ref{defn: circ case} with another construction of cycle classes. 
\end{defn}

\subsection{Tautological line bundles}\label{ssec: tautological bundles}
For $i=1,\cdots, r$ we have a line bundle $\ell_{i}$ on $\Hk^{r}_{U(n)}$ whose fiber at $(\{x'_{j}\}, \{\cF_{j},h_{j}\})$ is the line $\cF_{i}/\cF^{\flat}_{i-1/2}$ (supported at $\s(x'_{i})$), where $\cF^{\flat}_{i-1/2}=\cF_{i-1}\cap \cF_{i}$. We use the same notation $\ell_{i}$ to denote its pullback to $\Sht^{r}_{U(n)}$. We call them {\em tautological line bundles} on $\Sht^{r}_{U(n)}$.

\begin{defn}\label{def:cap Chern} Let $\cK \subset \cE$ be a sub-bundle. In Definition \ref{defn: circ case} we have defined a cycle class $[\cZ^{r ,\c}_{\cE/\cK}]\in \Ch_{r(n-m+m_0)}(\cZ^{r, \circ }_{\cE/\cK})$, where $m_0$ is the rank of $\cK$. Using \eqref{Z quot circ} 
we view $[\cZ^{r, \c}_{\cE/\cK}] \in \Ch_{r(n-m + m_0)}(\cZ^{r }_{\cE}[\cK]^{\circ})$  We define
\begin{equation}
[\cZ^{r}_{\cE}[\cK]^{\circ}]:=(\prod_{i=1}^{r}c_{\mrm{top}}(p_{i}^{*}\s^{*}\cK^{*}\ot\ell_{i}))\cap[\cZ^{r, \c}_{\cE/\cK}]\in \Ch_{r(n-m)}(\cZ^{r }_{\cE}[\cK]^{\circ}).
\end{equation}
Here  $\cK^{*} := \cHom(\cK,\cO_{X'})$ is the linear dual of $\cK$, and recall that $p_{i}: \Sht^{r}_{U(n)}\to X'$ records the leg $x'_{i}$. The notation $c_{\mrm{top}}(\ldots)$ denotes the ``top Chern class''. 

\end{defn}

\begin{remark}
More generally, if $BH_2$ is any gerbe of unitary type as in Definition \ref{def: gerbe unitary type}, then the same formula defines a tautological bundle $\ell_i$ on $\Sht_{H_2}^r$. We may then define an analogous class $[\cZ_{\cE, H_2}^r[\cK]^{\circ}]$. We will not have much need for this extra generality, so we prefer to focus on the case $BH_2 = BU(n)_{\LL}$ for concreteness. The general unitary gerbe case is only invoked in Example \ref{ex: CM pullback special cycle} and \S \ref{sec: CM}.
\end{remark}


\subsection{Virtual fundamental classes for special cycles}\label{ssec: VFC for special cycle}
Finally we have the definition of the cycle class $[\cZ^{r}_{\cE}]$.

\begin{defn}[Definition of special cycle classes]\label{def:special cycle classes}\hfill
\begin{enumerate}
\item Under the decomposition \eqref{ZEK decomp}, let $[\cZ^{r}_{\cE}]\in \Ch_{r(n-m)}(\cZ^{r}_{\cE})$ be the cycle class whose restriction to the open-closed substack $\cZ^{r}_{\cE}[\cK]^{\c}$ is the class $[\cZ^{r}_{\cE}[\cK]^{\c}]$ from Definition \ref{def:cap Chern}, for all sub-bundles $\cK$ of $\cE$. 
\item Let $a\in \cA_{\cE}(k)$. Define $[\cZ^{r}_{\cE}(a)]\in \Ch_{r(n-m)}(\cZ^{r}_{\cE}(a))$ to be the projection of $[\cZ_{\cE}^r]$ to the summand $\Ch_{r(n-m)}(\cZ^{r}_{\cE}(a))$.
\end{enumerate}
Define $[\ol{\cZ}_{\cE}^r] \in \Ch_{r(n-m)} (\ol{\cZ}_{\cE}^r)$ analogously. We note that $[\ol{\cZ}_{\cE}^r]$ pulls back to $[\cZ_{\cE}^r]$ under the finite \'{e}tale map $\cZ_{\cE}^r \rightarrow \ol{\cZ}_{\cE}^r$; conversely, under this same map $[\cZ_{\cE}^r]$ pushes forward to $\# \Aut(\cE)(\F_q) \cdot [\ol{\cZ}_{\cE}^r]$, so the two classes are essentially interchangeable. Again, $[\ol{\cZ}_{\cE}^r]$ will only be used as an auxiliary device to compare the virtual class $[\cZ^{r}_{\cE}]$ with another construction that will appear later. 
\end{defn}

Let $a \in \cA_{\cE}(k)$. We define substacks of $\cZ^{r}_{\cE}(a)$:
\begin{equation}
\cZ_{\cE}^r [\cK](a):=\cZ_{\cE}^r [\cK]\cap \cZ_{\cE}^r(a), \quad \cZ_{\cE}^r [\cK](a)^{\c}:=\cZ_{\cE}^r [\cK]^{\c}\cap \cZ_{\cE}^r(a).
\end{equation}
It is clear that $\cZ_{\cE}^r [\cK](a)$ is non-empty only when $\cK\subset \ker(a)$. 

The open-closed decomposition \eqref{ZEK decomp} restricts to an open-closed decomposition of $\cZ_{\cE}^r(a)$, 
\begin{equation}\label{decomp Z kernel}
\cZ_{\cE}^r(a) = \coprod_{\text{sub-bundles }\cK \subset\subset \ker(a)} \cZ_{\cE}^r[\cK](a)^\circ.
\end{equation}
We define $[\cZ_{\cE}^r[\cK](a)^\circ]\in \Ch_{r(n-m)}(\cZ_{\cE}^r[\cK](a)^\circ)$ to be the restriction of $[\cZ_{\cE}^r(a)]$.

\begin{remark}
Note that a different definition of $[\cZ_{\cE}^{r}(a)^{\circ}]$ has already been given in \cite[\S 7]{FYZ} when $a$ is non-singular, at least in special cases where $\cE$ is a direct sum of line bundles or $\rank \cE = n$. We will establish later (Proposition \ref{prop: compatibilitiy}) that the definitions are consistent. 
\end{remark}

We denote the natural projection from special cycles to $\Sht^{r}_{U(n)}$ by 
\begin{equation}
\z: \cZ^{r}_{\cE}(a)\to \Sht^{r}_{U(n)}.
\end{equation}
Recall from \cite[Proposition 7.5]{FYZ} that $\z$ is finite, the map $\z_{*}$ on Chow groups is therefore defined. In particular we have the Chow class
\begin{equation}
\z_{*}[\cZ^{r}_{\cE}(a)]\in \Ch_{r(n-m)}(\Sht^{r}_{U(n)})
\end{equation}
for any $a\in \cA_{\cE}(k)$.

Recall from the decomposition \eqref{decomp Z kernel} that for a singular $a$, $\cZ^{r}_{\cE}(a)$ may have infinitely many components $\cZ^{r}_{\cE}[\cK](a)^{\c}$ indexed by sub-bundles $\cK\subset \ker(a)$. The cycle $\z_{*}[\cZ^{r}_{\cE}(a)]$ is still well-defined because $\z$ is finite on the whole $\cZ^{r}_{\cE}(a)$ and not just on each $\cZ^{r}_{\cE/\cK}(a)^{\c}$.  Although not logically needed, we give an independent proof of the following fact that assures us that $\z_{*}[\cZ^{r}_{\cE}(a)]$ is a locally finite union of algebraic cycles.

\begin{lemma}\label{l:loc fin} Fix $(\cE,a)$ as above. For each sub-bundle $\cK \subset \ker(a)$, let $\frZ^{r}_{\cE/\cK }(\ov a)^{\c}\subset \Sht^{r}_{U(n)}$ be the image of  $\cZ^{r}_{\cE/\cK }(\ov a)^{\c}\cong \cZ^{r}_{\cE}[\cK](a)^{\c}$ under $\z$.  Then the collection of closed substacks $\{\frZ^{r}_{\cE/\cK }(\ov a)^{\c}\}_{\cK \subset \ker(a)}$ of $\Sht^{r}_{U(n)}$ is locally finite.
\end{lemma}
\begin{proof}
For $d\in\Q$, let $\Sht^{r, \le d}_{U(n)}$ be the open substack of those Hermitian Shtukas $\cF_{\bu}$ such that all slopes of $\cF_{0}$ (as vector bundles over $X'$) are $\le d$. It suffices to show that the intersection $\Sht^{r, \le d}_{U(n)}\cap \frZ^{r}_{\cE/\cK }(\ov a)^{\c}$ is non-empty only for finitely many sub-bundles $\cK \subset \ker(a)$. Now suppose $\Sht^{r, \le d}_{U(n)}\cap \frZ^{r}_{\cE/\cK }(\ov a)^{\c}\ne \vn$, and let $\ov\cE=\cE/\cK $. For any vector bundle $\cV$ on $X'$  let $\mu_{\max}(\cV)$ and $\mu_{\min}(\cV)$ be the maximal and minimal slopes of $\cV$. On one hand, a $\ov k$-point in $\Sht^{r, \le d}_{U(n)}\cap \frZ^{r}_{\cE/\cK }(\ov a)^{\c}$  gives an injective map $\ov\cE_{\ov k}\to \cF_{0}$, which implies $\mu_{\max}(\ov \cE)\le \mu_{\max}(\cF_{0})\le d$. On the other hand, $\ov \cE$ being a quotient of $\cE$ implies that $\mu_{\min}(\ov\cE)\ge\mu_{\min}(\cE)$.  Thus all slopes of $\ov\cE$ are within the range $[\mu_{\min}(\cE),d]$. This leaves finitely many possibilities for vector bundles $\ov \cE$ over $X'$ of rank bounded by the rank of $\ker(a)$.
\end{proof}

\subsection{Hermitian and skew-Hermitian bundles} 
In the classical formalism of reductive dual pairs, unitary groups for a Hermitian space and a skew-Hermitian space form a reductive dual pair. Similarly, our formulation of the modularity conjecture will involve both Hermitian and skew-Hermitian bundles. Points in the Hitchin base $\cA_{\cE}(k)$ will appear as Fourier parameters for an automorphic form on a unitary group attached to a skew-Hermitian form.  We spell out the distinction between Hermitian and skew-Hermitian bundles and the duality responsible for the Fourier expansion.

When we assemble the special cycles into a Fourier series, the dual of the Hitchin space $\cA_{\cE}(k)$ will be the Fourier parameter. For a fixed similitude line bundle $\LL$, the Hitchin space $\cA_{\cE}(k)$ is a subspace of $\Hom(\cE, \sigma^* \cE^\vee \otimes \nu^* \LL)$, which by Serre duality is dual to $\Ext^1(\sigma^* \cE^*  \otimes \nu^* \LL, \cE)$. Elements of this $\Ext$ group may be viewed as bundles with a conjugate dual structure,  so in preparation for the formulation of the modularity conjecture, we spell out the distinction between Hermitian and skew-Hermitian bundles. 

\sss{Hermitian sheaf}
For a vector bundle $\cE$ on $X'$ and a line bundle $\LL$ on $X$, let 
$$\cHerm(\cE;\LL)\subset \nu_{*}\cHom(\cE\ot \s^{*}\cE,\nu^{*}\LL)\cong \nu_{*}\cHom(\cE, \s^{*}\cE^{*}\ot\nu^{*}\LL)$$ 
be the subsheaf (over $X$) of Hermitian maps, i.e., local sections of  $\cHerm(\cE;\LL)$ are maps $a \co \cE \rightarrow \sigma^* \cE^* \otimes \nu^* \LL$ such that $\sigma^* a^\vee \co \cE \otimes  \sigma^* \nu^* \LL^{*} \rightarrow \sigma^* \cE^{*}$ agrees with $a$ after using the tautological descent datum $\sigma^* \nu^* (\LL^{*}) \cong \nu^* \LL^{*}$ and then twisting by $\nu^* \LL$.  

Note by definition that
\begin{equation}\label{eq: A as H0}
\cA_{\cE,\LL}(k)=H^{0}(X, \cHerm(\cE;\LL)\ot \om_{X}).
\end{equation}

\sss{Skew-Hermitian bundles}
Now we recall the notion of skew-Hermitian bundles.  Let $\frM$ be a line bundle on $X$. For a vector bundle $\cG$ on $X'$, a \emph{skew-Hermitian structure (with similitude $\frM$)} on $\cG$ is an isomorphism $h: \cG\isom \s^{*}\cG^{\vee}\ot \nu^{*}\frM=\s^{*}\cG^{*}\ot\nu^{*}(\om_{X}\ot \frM)$ such that $\s^* h^\vee \co \cG  \otimes \sigma^* (\nu^* \frM)^{*}\rightarrow \s^* \cG^\vee$ is identified with $-h$ upon twisting by $\nu^* \frM$ and using descent datum $\sigma^* \nu^* \frM   \cong \nu^* \frM$. (Recall that a Hermitian structure is similar to such $h$ except that we instead ask $\s^* h^{\vee} = h$ under the same identifications.) The datum $(\cG, h)$ is referred to as a skew-Hermitian bundle on $X$.

Let $L_{X'/X}$ be the line bundle on $X$ associated to the double covering $\nu \co X' \rightarrow X$, i.e., $\nu_* \cO_{X'} \cong \cO_X \oplus L_{X'/X}$. The map $L_{X'/X} \rightarrow \nu_* \cO_{X'}$ induces an isomorphism $\phi \co \nu^* L_{X'/X} \cong \cO_{X'}$. Under $\phi$, the descent datum $\sigma^* \nu^* L_{X'/X} \cong \nu^* L_{X'/X}$ is intertwined with multiplication by $-1$ on $\cO_{X'}$. 

\begin{lemma}
There is an equivalence of categories between Hermitian bundles on $X$ with similitude factor $\LL$ and skew-Hermitian bundles on $X$ with similitude factor $\LL \otimes L_{X'/X}$. 
\end{lemma}

\begin{proof}
Let $\LL' := \LL \otimes_{\cO_X} L_{X'/X}$. Given a Hermitian bundle $(\cG,h)$, the map
\[
 h \otimes \phi^{-1} \co \cG \otimes_{\cO_{X'}} \cO_{X'} \xrightarrow{\sim} (\s^* (\nu^* \cG)^\vee \otimes \nu^* \LL )\otimes_{\cO_{X'}} \nu^* L_{X'/X}
 \]
 is skew-Hermitian. Obviously $\cG \otimes_{\cO_{X'}} \cO_{X'} \cong \cG$, so the functor $(\cG, h) \mapsto (\cG \otimes_{\cO_{X'}} \cO_{X'},  h \otimes \phi^{-1})$ defines an equivalence between Hermitian bundles on $X$ with similitude factor $\LL$ and skew-Hermitian bundles on $X$ with similitude factor $\LL \otimes L_{X'/X}$.  
\end{proof}

Now let $\LL$ be a line bundle on $X$. We form the sheaf 
$$
\cHerm(\s^{*}\cE^{*}\ot \nu^{*}\LL; \LL)\subset \nu_{*}\cHom((\s^{*}\cE^{*}\ot\nu^{*}\LL)\ot (\cE^{*}\ot \s^{*}\nu^{*}\LL),\nu^{*}\LL)\\
\cong \nu_{*}\cHom(\s^{*}\cE^{*}\ot \nu^{*}\LL, \cE).
$$
Define
\begin{equation}\label{eq: def Ext Herm}
\Ext^{1}_{\Herm}(\s^{*}\cE^{*}\ot\nu^{*}\LL, \cE):=H^{1}(X, \cHerm(\s^{*}\cE^{*}\ot \nu^{*}\LL; \LL))\subset \Ext^{1}_{X'}(\s^{*}\cE^{*}\ot \nu^{*}\LL, \cE).
\end{equation}


\begin{prop}\label{prop: Lag in G}
The subspace $\Ext^1_{\Herm}(\sigma^* \cE^*  \otimes \nu^* \LL, \cE) \subset \Ext^1_{X'}(\sigma^* \cE^*  \otimes \nu^* \LL, \cE)$ is canonically isomorphic to the group of equivalence classes of extensions
\[
0 \rightarrow \cE \rightarrow \cG \rightarrow \sigma^* \cE^*\ot\nu^{*}\LL  \rightarrow 0
\]
equipped with a skew-Hermitian structure $h$ on $\cG$ (with similitude bundle $\frM=\om_{X}^{-1}\ot \LL$), with respect to which $\cE$ is Lagrangian, such that $h$ restricts to the tautological pairing 
$$\cE\otimes \s^{*}(\cG/\cE)\cong \cE\ot(\cE^{*}\ot \nu^{*}\LL)\to \nu^{*}\LL.$$
\end{prop}

\begin{proof}

If $(\cG, h)$ is a skew-Hermitian bundle on $X'$ with similitude $\frM=\om_{X}^{-1}\ot \LL$, we view $h$ as a pairing 
$$h: \cG\ot\s^{*}\cG\to \nu^{*}\LL.$$ 

For a Lagrangian-subbundle $\cE \subset \cG$, the datum of $(\cG, h, \cE)$ has automorphism sheaf $\cAut^1(\cG, h, \cE)$, whose local sections consist of local $\cO_{X'}$-linear automorphisms $\alpha \co \cG \rightarrow \cG$ preserving $h$ such that $\alpha|_{\cE} = \Id$. Sending $\alpha$ to $\alpha-\Id$ induces a bijection between such $\alpha$ and $\cO_{X'}$-linear maps $\varphi \co \sigma^* \cE^* \otimes \nu^* \LL \rightarrow \cE$ such that (in order to preserve $h$)
\begin{equation}\label{eq: varphi condition}
h(\varphi v, w) + h(v, \varphi w) = 0
\end{equation}
for all local sections $v,w$ of $\sigma^* \cE^* \otimes \nu^* \LL$. Using that $h(v, \varphi w) = - \sigma^{*} h (\varphi w,v)$, we may rewrite \eqref{eq: varphi condition} as $h(\varphi v,w) = \sigma^{*} h(\varphi w), v)$. Then $(v, w) \mapsto h(\varphi v,w)$ locally defines a Hermitian form on $\sigma^* \cE^* \otimes \nu^* \LL$ valued in $\nu^* \LL$. This exhibits an isomorphism $\cAut^1(\cG,h, \cE) \cong \cHerm(\sigma^* \cE^* \otimes \nu^* \LL; \nu^* \LL)$ of sheaves of abelian groups on $X$. Hence equivalence classes of $(\cG, h)$ (for fixed $\cE$) as above are classified by elements $e_{\cG, \cE} \in H^1(X, \cHerm(\sigma^* \cE^* \otimes \nu^* \LL;	 \nu^* \LL))$, which is $\Ext^1_{\Herm}(\sigma^* \cE^*  \otimes \nu^* \LL, \cE)$ by definition. 
\end{proof}

\subsubsection{Duality}\label{sss: duality A and Ext}
Linear duality between $\cHom(\cE, \s^{*}\cE^{*}\ot\nu^{*}\LL)$ and $\cHom(\s^{*}\cE^{*}\ot\nu^{*}\LL, \cE)$ restricts a perfect pairing between vector bundles over $X$
\begin{equation*}
\cHerm(\cE;\LL)\ot_{\cO_{X}}\cHerm(\s^{*}\cE^{*}\ot \nu^{*}\LL; \LL)\to \cO_{X}.
\end{equation*}
Serre duality gives a perfect pairing between the $k$-vector spaces
\begin{equation*}
\cA_{\cE,\LL}(k)\stackrel{\eqref{eq: A as H0}}{=}H^{0}(X, \cHerm(\cE;\LL)\ot \om_{X})
\end{equation*}
and
\begin{equation*}
\Ext^{1}_{\Herm}(\s^{*}\cE^{*}\ot \nu^{*}\LL, \cE)\stackrel{\eqref{eq: def Ext Herm}}{=}H^{1}(X, \cHerm(\s^{*}\cE^{*}\ot \nu^{*}\LL; \LL)).
\end{equation*}

\subsection{The modularity conjecture}\label{ss:mod}

Let $\Bun_{GU^-(2m)}$ be the moduli stack of triples $(\cG, \frM, h)$ where $\cG$ is a vector bundle of rank $2m$ over  $X'$, $\frM$ is a line bundle over $X$, and $h$ is a \emph{skew-}Hermitian isomorphism $h: \cG\isom \s^{*}\cG^{\vee}\ot \nu^{*}\frM=\s^{*}\cG^{*}\ot\nu^{*}(\om_{X}\ot \frM)$. Let $c:\Bun_{GU^-(2m)}\to \Pic_{X}$ be the map recording $\om_{X}\ot \frM$. Then for any $\frL\in\Pic_{X}(k)$, $c^{-1}(\frL)=\Bun_{U^-(2m),\om_{X}^{-1}\ot\frL}$ as in \S\ref{ss:Sht sim}.

 A priori $\Bun_{GU^-(2m)}(k)$ has a decomposition
\begin{equation}
\Bun_{GU^-(2m)}(k)=\coprod_{\xi} \wt H_{\xi}(F)\bs \wt H_{\xi}(\BA)/\wt H_{\xi}(\wh\cO)
\end{equation}
where $\xi$ runs through $2m$-dimensional skew-Hermitian spaces over $F'$ that are locally split at all places, and $\wt H_{\xi}$ is the corresponding unitary similitude group. By the Hasse principle for Hermitian spaces \cite[Theorem 6.2]{Sch85}
$\xi$ must be globally split. Let $\wt H_{m}=GU^-(2m)$ be the unitary similitude group for a fixed split $2m$-dimensional skew-Hermitian space over $F'$.  Then
\begin{equation}
\Bun_{GU^-(2m)}(k)=\wt H_{m}(F)\bs \wt H_{m}(\BA)/\wt H_{m}(\wh\cO).
\end{equation}

We can similarly define the moduli $\Sht^{r}_{GU(n)}$ of shtukas for $GU(n)$. It simply adds the similitude line bundle $\frL$ as part of the data which is invariant under Frobenius pullback, and it is the disjoint  union
\begin{equation}\label{eq: GU shtuka decomposition}
\Sht^{r}_{GU(n)}=\coprod_{\frL\in\Pic_{X}(k)}\Sht^{r}_{U(n),\frL} \times B(\Aut(\frL)(k)).
\end{equation}
We emphasize that in our definition, the modifications cannot occur on the similitude factor $\frL$.


Let $\Bun_{\wt P_{m}}$ be the moduli stack of quadruples $(\cG,\frM, h,\cE)$ where $(\cG,\frM,h)\in \Bun_{GU^-(2m)}$, and $\cE\subset \cG$ is a Lagrangian sub-bundle (of rank $m$). Let $\Bun_{P_{m},\frM}$ be the substack with the fixed similitude line bundle $\frM$. We usually omit $h$ and write a point in  $\Bun_{P_{m},\frM}$ as $(\cG,\cE)$.

The map $\Bun_{\wt P_{m}}\to \Bun_{GU^-(2m)}$ forgetting the Lagrangian sub-bundle is surjective as map of stacks, and it is also surjective on $k$-points. Indeed,  since the generic fiber of any  $(\cG,\frM,h)\in \Bun_{GU^-(2m)}(k)$ is a split skew-Hermitian space over $F'$ of dimension $2m$, it  has a Lagrangian sub-bundle at the generic point, hence a Lagrangian sub-bundle over $X'$ by saturation.  If we write $\wt P_{m}\subset \wt H_{m}$ for the Siegel parabolic subgroup stabilizing a Lagrangian subspace, then
\begin{equation}\label{adelic points Pm}
\Bun_{\wt P_{m}}(k)=\wt P_{m}(F)\bs \wt H_{m}(\BA)/\wt H_{m}(\wh\cO).
\end{equation}



Now fix $\frL\in\Pic_{X}(k)$. For $(\cG,\cE)\in \Bun_{P_{m},\om_{X}^{-1}\ot \frL}(k)$, we have a short exact sequence
\begin{equation}
0\to \cE\to \cG\to \s^{*}\cE^{*}\ot \nu^{*}\LL\to 0
\end{equation}
which by Proposition \ref{prop: Lag in G} gives an extension class
\begin{equation}
e_{\cG,\cE}\in \Ext^{1}_{\Herm}(\s^{*}\cE^{*}\ot\nu^{*}\LL,\cE).
\end{equation}
Recall the perfect pairing
\begin{equation}
\j{\cdot,\cdot}: \cA_{\cE,\LL}(k)\times \Ext^{1}_{\Herm}(\s^{*}\cE^{*}\ot\nu^{*}\LL,\cE)\to k
\end{equation}
given in \S\ref{sss: duality A and Ext}.
It is the restriction of the Serre duality pairing
\begin{equation}
\j{\cdot,\cdot}: \Ext^{1}_{X'}(\s^{*}\cHom(\cE,\nu^{*}\frL),\cE)\times \Hom_{X'}(\cE, \s^{*}\cHom(\cE,\nu^{*}\frL)\ot \om_{X'})\to k.
\end{equation}
For $a\in \cA_{\cE,\LL}(k)$, the element $\j{e_{\cG,\cE}, a}\in k$ is therefore defined.


Recall that we have fixed a nontrivial character $\psi_{0}: k\to \Qlbar^{\times}$. 
Finally, recall that $\y:\Pic_X(k)\to \{\pm1\}$ is the character with kernel $\Nm(\Pic_{X'}(k))$. Let $\chi: \Pic_{X'}(k)\to  \Qlbar^{\times}$ be a character such that $\chi|_{\Pic_X(k)}=\y^{n}$. (The existence of such $\chi$ is justified in \cite[Remark 2.1]{FYZ}.)

\begin{defn}\label{def: gen series before descent}
Define a map
\begin{eqnarray}
\wt Z^{r}_{m,\frL}: \Bun_{P_{m}, \om_{X}^{-1}\ot \frL}(k) &\to& \Ch_{r(n-m)}(\Sht^{r}_{U(n), \LL}) \otimes_{\Q} \ol{\Q}_{\ell} \\
 (\cG,\cE) & \mapsto  & \chi(\det\cE)q^{n(\deg \cE-\deg\frL-\deg\om_{X}) /2}\sum_{a\in \cA_{\cE,\frL}(k)}\psi_{0}(\j{e_{\cG,\cE}, a})\z_{*}[\cZ^{ r}_{\cE,\LL}(a)].\notag
\end{eqnarray}
Taking the union over $\frL\in\Pic_{X}(k)$, we get a map
\begin{equation}
\wt Z^{r}_{m}: \Bun_{\wt P_{m}}(k) \to \Ch_{r(n-m)}(\Sht^{r}_{GU(n)})  \otimes_{\Q} \ol{\Q}_{\ell} .
\end{equation}
\end{defn}

\begin{remark}
Using \eqref{adelic points Pm}, we may identify $\wt Z^{r}_{m}$ as a function
\begin{equation}\label{Zm function}
\wt H_{m}(\BA)\ni g\mapsto \wt Z^{r}_{m}(g)\in  \Ch_{r(n-m)}(\Sht^{r}_{GU(n)})  \otimes_{\Q} \ol{\Q}_{\ell}
\end{equation}
such that
\begin{itemize}
\item $\wt Z^{r}_{m}$ is left invariant under the Siegel parabolic $\wt P_{m}(F)$ and right invariant under $\wt H_{m}(\wh\cO)$ (everywhere unramified);
\item if $g\in \wt H_{m}(\BA)$ has similitude factor $c(g)\in \BA^{\times}$ that projects  to the line bundle $\frL\in \Pic_{X}(k)=F^{\times}\bs  \BA^{\times}/\wh\cO^{\times}$, then $\wt Z^{r}_{m}(g)$ is supported on $\Sht^{r}_{U(n),\frL}\subset \Sht^{r}_{GU(n)}$.
\end{itemize}
\end{remark}

The following is the main conjecture of the paper.
\begin{conj}[Modularity conjecture]\label{c:mod} The map $\wt{Z}^{r}_{m}$ descends to a map
\begin{equation}
Z^{r}_{m}: \Bun_{GU^-(2m)}(k)\to \Ch_{r(n-m)}(\Sht^{r}_{GU(n)})  \otimes_{\Q} \ol{\Q}_{\ell}.
\end{equation}
i.e., the function \eqref{Zm function} is left $\wt H_{m}(F)$-invariant. 
\end{conj}

In other words, the Chow class $\wt Z^{r}_{m}(\cG,\cE)\in \Ch_{r(n-m)}(\Sht^{r}_{GU(n)})  \otimes_{\Q} \ol{\Q}_{\ell}$ should depend only on the skew-Hermitian bundle $\cG$ and not on its Lagrangian sub-bundle $\cE$.

\begin{remark} When $r=0$, $\Ch_{0}(\Sht^{0}_{GU(n)})$ is simply the space of $\Q$-valued functions on $\Bun_{GU(n)}(k)$.  The conjecture in this case follows from the automorphy of the theta series constructed from the Weil representation for the dual pair $(GU^-(2m), GU(n))$.
\end{remark}

\begin{remark}
Suppose $r>0$ and $n>1$. We expect based on \S \ref{sssec: chow groups of zero-cycles} that $\Ch_{0}(\Sht^{r}_{GU(n)})$ vanishes, making the conjecture vacuous for $m=n$ in this situation. In \cite{FYZ}, for the non-singular terms we constructed cycle classes in the Chow group of \emph{proper} cycles on $\Sht_{U(n)}^r$, and proved a higher Siegel-Weil formula for those terms. It remains an open problem to formulate a more refined version of the generating series where the singular terms also have a meaningful notion of degree. 
\end{remark}

\subsection{Special cases}

Let $\cE$ be a rank $m$ vector bundle on $X'$. Let $\cE'=\s^{*}\cHom(\cE,\nu^{*}\frL)$. Consider the Hermitian vector bundle $\cG=\cE\op \cE'$ with the natural Hermitian form isotropic on each summand and induces the natural pairing between the two summands.  In this case, both $(\cG,\cE)$ and $(\cG,\cE')$ are points of $\Bun_{P_{m},\om_{X}^{-1}\ot\frL}(k)$ over $\cG\in \Bun_{U^-(2m),\om_{X}^{-1}\ot\frL}(k)$. Conjecture \ref{c:mod} specializes to the following identity.

\begin{conj} In the above situation, we have an identity in $\Ch_{r(n-m)}(\Sht^{r}_{U(n),\frL})$:
\begin{equation}
\chi(\det\cE)q^{n\deg \cE /2}\sum_{a\in \cA_{\cE}(k)}\z_{*}[\cZ^{r}_{\cE}(a)]=\chi(\det\cE')q^{n\deg \cE' /2}\sum_{a'\in \cA_{\cE'}(k)}\z_{*}[\cZ^{r}_{\cE'}(a')].
\end{equation}
Equivalently,
\begin{equation}
\y(\frL)^{mn}q^{n(\deg \cE -m\deg\frL)}\sum_{a\in \cA_{\cE}(k)}\z_{*}[\cZ^{r}_{\cE}(a)]=\sum_{a'\in \cA_{\cE'}(k)}\z_{*}[\cZ^{r}_{\cE'}(a')].
\end{equation}
\end{conj}
In the equivalent formulation above, we use that
\begin{equation}
\deg\cE'=-\deg\cE+2m\deg\frL, \quad \det(\cE')\cong\s^{*}(\det\cE)^{-1}\ot \nu^{*}\frL^{\ot m}.
\end{equation}

We may further specialize to the case where $\cE'$ has large slopes, or equivalently $\cE$ has small slopes, so that $\cA_{\cE'}(k)$ only contains the zero Hermitian map.

\begin{conj} Suppose the maximal slope $\mu_{\max}(\cE)$ satisfies 
\begin{equation}
\mu_{\max}(\cE)<\deg\frL-\deg\om_{X}.
\end{equation}
Then we have an identity in $\Ch_{r(n-m)}(\Sht^{r}_{U(n),\frL})$:
\begin{equation}
\y(\frL)^{mn}q^{n(\deg \cE -m\deg\frL)}\sum_{a\in \cA_{\cE}(k)}\z_{*}[\cZ^{r}_{\cE}(a)]=\z_{*}[\cZ^{r}_{\cE'}(0)].
\end{equation}
\end{conj}

\subsection{Test intersection numbers}\label{ssec: test intersection numbers}

To give evidence for Conjecture \ref{c:mod}, we may start with any cycle {\em with compact support} $\xi\in \Ch_{rm,c}(\Sht^{r}_{U(n),\frL})$, and form the numerical function by intersecting $\wt Z^{r}_{m}$ with $\xi$:
\begin{equation}
\wt I^{r}_{m,\xi}:=\j{\wt Z^{r}_{m}(-), \xi}_{\Sht^{r}_{U(n),\frL}}: \Bun_{P_{m},\om_{X}^{-1}\ot\frL}(k)\to \Q.
\end{equation}
Conjecture \ref{c:mod} predicts that $\wt I^{r}_{m,\xi}(\cG,\cE)$ is independent of $\cE$, hence  descends to a function on $\Bun_{U^-(2m),\om_{X}^{-1}\ot\frL}(k)$.
We give two families of examples compact $r$-dimensional cycles $\xi$ on $\Sht^{r}_{U(n),\frL}$, hence giving test grounds for Conjecture \ref{c:mod} in the case $m=1$.


\begin{example}[Corank $n-1$ special cycles]\label{ex:crk n-1} Let  $\cE$ be a rank $n-1$ vector bundle over $X'$, and $a\in \cA_{\cE}^{\ns}(k)$ be a non-singular Hermitian map. Then the special cycle $\cZ^{r}_{\cE}(a)$ is proper (we omit the proof here). We have the cycle class $[\cZ^{r}_{\cE}(a)]\in \Ch_{r}(\cZ^{r}_{\cE}(a))$ by Definition \ref{defn: circ case}. Its direct image in $\Sht^{r}_{U(n),\frL}$ is then a compact cycle 
\begin{equation}
\xi:=\z_{*}[\cZ^{r}_{\cE}(a)]\in \Ch_{r,c}(\Sht^{r}_{U(n),\frL}).
\end{equation}
\end{example}

\begin{example}[CM cycles]\label{ex:CM} Let  $Y$ be another smooth projective curve over $\F_{q}$, and $\th: Y\to X$ be a map of degree $n$, possibly ramified. Let $\nu_{Y}: Y'=X'\times_{X}Y\to Y$, and assume this double covering is nonsplit over each connected component of $Y$. Let $\Sht^{r}_{U(1)/Y, \th^{*}\frL}$ be the moduli stack of rank $1$ $\th^{*}\frL$-twisted Hermitian shtukas (cf. \S\ref{sssec: hecke for unitary gerbe} for the definition) on $Y'$ (with respect to the double cover $\nu_{Y}$). 
Then push-forward along $\nu_{Y}$ gives a map $\Th: \Sht^{r}_{U(1)/Y, \th^{*}\frL}\to \Sht^{r}_{U(n), \frL}$. Now $\Sht^{r}_{U(1)/Y, \th^{*}\frL}$ is smooth and proper of pure dimension $r$, we have the compact cycle class 
\begin{equation}
\xi=\Th_{*}[\Sht^{r}_{U(1)/Y, \th^{*}\frL}]\in \Ch_{r,c}(\Sht^{r}_{U(n), \frL}).
\end{equation}
The intersection number of the generating series of corank $1$ and this cycle will be calculated in \S\ref{sec: CM}.  In particular, we will verify the modularity of such intersection numbers. 
\end{example}

\begin{remark}It is possible to give a general construction that includes both examples as special cases, but the details will not be included here. 
\end{remark}

\subsection{The split case}\label{ssec: split case conj}
In the case where $X' = X^{(1)} \coprod X^{(2)}$ is the split double cover of $X$ (so each $X^{(i)}\cong X$), the definition of the cycle classes $[\cZ^{r}_{\cE}]$ can be spelled out more explicitly as follows. In this case, an $\LL$-twisted Hermitian bundle $\cF$ on $X'$ identifies with a pair of vector bundles $(\cF^{(1)},\cF^{(2)})$, each living on one copy of $X$, equipped with an isomorphism $\cF^{(2)}\cong \cF^{(1), \vee} \otimes \LL$.  Therefore we have $\Bun_{U(n),\LL}\cong \Bun_{\GL(n)}$ by recording only $\cF^{(1)}$. Since every $\LL$ is a norm, without loss of generality we can and will assume  $\LL=\cO_X$. Then we have a disjoint union 
 $$\Sht^r_{U(n),\LL}=\coprod_{\mu\in\{\pm 1\}^{r}}\Sht_{\GL(n)}^{\mu} $$ where the $\mu=(\mu_1,\cdots,\mu_r)$-th component is empty unless $\sum_{i=1}^r\mu_i=0$;
 see \cite[\S12.3]{FYZ} which also recalled the definition of $\Sht_{\GL(n)}^{\mu}$. In particular, this implies that $r$ is even, so that $r/2$ is an integer. 

A vector bundle $\cE$ on $X'$ of rank $m$ corresponds to two rank $m$ vector bundles $(\cE^{(1)},\cE^{(2)})$, each living on one copy of $X$. Now $\cA_{\cE}(k)=\cA_{\cE^{(1)}, \cE^{(2)}}(k)$ may be identified with the set of maps $a: \cE^{(1)}\to \cE^{(2), \vee}$. 

We now fix a $\mu=(\mu_1,\cdots,\mu_r)$ such that $\sum_{i=1}^r\mu_i=0$. The special cycle $\cZ_{\cE}^{\mu}=\cZ_{\cE^{(1)}, \cE^{(2)}}^{\mu}$ in the split case parametrizes
\begin{equation}\label{pt ZE split}
\{ (\{x_{i}\}_{1\le i\le r}, \cF_0 \dashrightarrow \ldots \dashrightarrow \cF_r \cong \ft \cF_0, \cE^{(1)}\xrightarrow{t^{(1)}_{i}} \cF_i,   \cE^{(2)} \xrightarrow{t^{(2)}_{i}} \cF_i^{\vee})\}
\end{equation}
where  $x_{i}\in X$, $\cF_{i}$ are vector bundles of rank $n$ on $X$, the dashed arrow $\cF_{i-1}\dashrightarrow \cF_{i}$ is a lower modification of length $1$ at $x_{i}$  if $\mu_i=-1$, and an upper modification of length $1$ at $x_{i}$ if $\mu_{i}=+1$. The maps  $t^{(1)}_{i}$ and $t^{(2)}_{i}$ are required to be compatible with the chain of modifications.

The kernel decomposition of $\cZ_{\cE}^{\mu}$ in this case is indexed by $\cK = (\cK^{(1)}, \cK^{(2)}) \subset (\cE^{(1)}, \cE^{(2)})$ where we note that the ranks of $\cK^{(1)}$ and $\cK^{(2)}$ may be different. We have an open-closed decomposition 
\[
\cZ_{\cE^{(1)},\cE^{(2)}}^\mu = \coprod_{\cK^{(1)}\subset \subset \cE^{(1)}} \coprod_{\cK^{(2)} \subset \subset \cE^{(2)}} \cZ_{\cE^{(1)},\cE^{(2)}}^\mu[\cK^{(1)}, \cK^{(2)}]^{\circ}
\]
where $\cZ_{\cE^{(1)},\cE^{(2)}}^\mu[\cK^{(1)}, \cK^{(2)}]^{\circ} $ is the substack of those points in \eqref{pt ZE split} 
where $\ker  t_i^{(1)} = \cK^{(1)}$ and $\ker t_i^{(2)} = \cK^{(2)}$ for any (equivalently, all) $0\le i\le r$. With $\ol{\cE}^{(1)} = \cE^{(1)}/\cK^{(1)}$ and $\ol{\cE}^{(2)} = \cE^{(2)}/\cK^{(2)}$, we have $\cZ_{\cE^{(1)},\cE^{(2)}}^\mu[\cK^{(1)}, \cK^{(2)}]^{\circ} \cong \cZ_{\ol{\cE}^{(1)}, \ol{\cE}^{(2)}}^{\mu, \circ}$. The virtual classes are then defined by summing over all $(\cK^{(1)},\cK^{(2)})$ the product of $[\cZ_{\ol{\cE}^{(1)}, \ol{\cE}^{(2)}}^{\mu, \circ}]$ with $\prod_{i=1}^{r}c_{\mrm{top}}(p_{i}^{*}\cK^{(?_i)}\ot\ell_{i})$ where $?_i = 1$ if $\mu = 1$ and $?_i = 2$ if $\mu = -1$.

In this case, $ \Bun_{U(2m), \om_{X}^{-1}}$ from \S \ref{ss:mod} is isomorphic to $\Bun_{\GL(2m)}$. However, we do not restrict ourselves to even rank $2m$ in the split case; we will be able to formulate a modularity conjecture for a higher theta function on $\Bun_{\GL(m)}$. For $a+b  = m$, let $P_{(a,b)}$ be the corresponding maximal parabolic subgroup of $\GL(m)$. Then $\Bun_{P(a,b)}$ classifies pairs $(\cG, \cE^{(1)})$ where $\cE^{(1)}$ is a rank $a$ sub-bundle of a rank $m$ vector bundle $\cG$ on $X'$. From the pair $(\cG, \cE^{(1)})$ we obtain a rank $b$ bundle $\cE^{(2)}$ by the exact sequence
\[
0 \rightarrow \cE^{(1)} \rightarrow \cG \rightarrow \cE^{(2),*}  \rightarrow 0. 
\]
Given $(\cE^{(1)},\cE^{(2)})$, the space of such extensions is $\Ext^1(\cE^{(2), *} , \cE^{(1)})$, which is dual to $\cA_{\cE^{(1)}, \cE^{(2)}}(k)= \Hom(\cE^{(1)}, \cE^{(2), \vee} )$. The class of $\cG$ defines $e_{\cG,\cE^{(1)}}\in \Ext^1(\cE^{(2), *} , \cE^{(1)})$ and we denote $\langle e_{\cG,  \cE^{(1)}}, - \rangle$ the induced $k$-linear functional on  $\cA_{\cE^{(1)}, \cE^{(2)}}(k)$.

\begin{remark}
The case $a=b=m$ resembles the discussion in the case where $X'/X$ is non-split, while the $a \neq b$ case has no counterpart there. We note however that even in the case $a=b=m$ the definition of the virtual class $[\cZ^{\mu}_{\cE^{(1)}, \cE^{(2)}}]$ involves considering $\cZ^{\mu}_{\ol{\cE}^{(1)}, \ol{\cE}^{(2)}}$ for $\ol{\cE}^{(1)}$ and $\ol{\cE}^{(2)}$ that have unequal rank. 
\end{remark}

We define the \emph{higher theta function} $\wt{Z}_{a,b}^\mu \co \Bun_{P_{(a,b)}}(k)\to \Ch_{\frac{r}{2}(2n-m)}(\Sht^\mu_{\GL(n)})$  as 
\[
(\cG, \cE^{(1)}) \mapsto q^{n \deg \cE^{(2)}}
\sum_{a \in \cA_{\cE^{(1)}, \cE^{(2)}}(k) }  \psi_0 (\langle e_{\cG,  \cE^{(1)}}, a \rangle) \zeta_* [\cZ_{\cE^{(1)}, \cE^{(2)}}^\mu(a) ]
\]
where $\zeta: \cZ_{\cE^{(1)}, \cE^{(2)}}^\mu(a) \to \Sht^\mu_{\GL(n)}$ is the natural projection map. The following is our modularity Conjecture for this split case. 


\begin{conj}
For each $\mu$, the map $\wt{Z}_m^\mu:\Bun_{P_{(a,b)}}(k)\to \Ch_{\frac{r}{2}(2n-m)}(\Sht^\mu_{\GL(n)})$ descends to a map $$\xymatrix{Z_{m}^\mu: \Bun_{\GL(m)}(k)\ar[r]& \Ch_{\frac{r}{2}(2n-m)}(\Sht^\mu_{\GL(n)}).}$$

\end{conj}

\part{Properties of the special cycles}

\section{Derived Hitchin stacks}\label{sec: derived} 

\subsection{Overview}

In the next two sections, we explain the special cycle classes of Definition \ref{def:special cycle classes} from the perspective of derived algebraic geometry. To motivate this, we recall that in \cite{FYZ}, certain ``Hitchin stacks'' $\cM$ were introduced and it was proved that the virtual fundamental class $[\cZ_{\cE}^r(a)]$ for \emph{non-singular} $a$ could be obtained from $\cM$ by taking the derived intersection of a Hecke correspondence $\Hk_{\cM}^r$ for $\cM$ with the graph of Frobenius on $\cM$. This interpretation was key to the proof of the Higher Siegel-Weil formula \cite[Theorem 1.1]{FYZ}. 

The restriction to non-singular $a$ can be explained thus: for such $a$, the intersection involves only the smooth part of the Hitchin stack $\cM$. But if we try to repeat such a construction to obtain the cycles indexed by singular $a$, we necessarily run into loci in $\cM$ whose geometry is too poorly behaved (more precisely, we cannot control the singularities nor the dimension) to carry it out. 

It turns out that these problems can be resolved with \emph{derived algebraic geometry}. In this section we will introduce \emph{derived} Hitchin stacks $\sM$, which are always quasi-smooth (the derived analogue of LCI) and have the ``correct'' virtual dimension, whose classical truncation is $\cM$. By taking the derived intersection of derived Hecke correspondences $\sHk_{\sM}^r$ for $\sM$ with the graph of Frobenius on $\sM$, we then obtain derived enhancements of the special cycles which we call $\sZ_{\cE}^r$. These are similarly always quasi-smooth and of the correct dimension; derived algebraic geometry then associates to them certain virtual fundamental classes in the Chow group of the underlying classical special cycles $\cZ_{\cE}^r$. On general grounds it is non-trivial to compute these virtual fundamental classes ``explicitly'' in terms of classical objects. Nevertheless, we will be able to prove that they coincide with the explicit constructions introduced earlier in Definition \ref{def:special cycle classes}. This gives a pleasing derivation of the cycle classes for singular terms, which is on the same conceptual footing as for the non-singular terms.

The fruits of this labor are not merely philosophical: in \S \ref{sec: linear invariance} we use this derived algebraic geometry interpretation of the cycle classes to prove the \emph{linear invariance} property of our special cycles. The number field analogue of this property is a well-known conjectural property of arithmetic theta series \cite[Problem 5]{Kud04}. The statement can be formulated in purely classical terms, but \emph{we do not know a proof without derived algebraic geometry}. In turn, \S \ref{sec: linear invariance} will also be used later in \S \ref{sec: CM} to provide numerical evidence for modularity conjecture.

\subsection{Derived stacks}

\subsubsection{Orientation on derived algebraic geometry} 

We give an introductory discussion on derived algebraic geometry, in order to help orient readers not accustomed to this formalism. We confine ourselves to informal and sometimes vague remarks, referring to \cite{Lur04, TV08} for more complete treatments. Some relevant introductory references are \cite{To09, To10}. 

Just as Grothendieck's schemes are spaces built locally from ``spectra'' of commutative rings, \emph{derived} schemes are built locally from ``spectra'' of \emph{simplicial commutative rings}. Roughly speaking, one can think of simplicial commutative rings as a model for the concept of topological rings.

We use the adjective ``classical'' (ring, scheme, stack...) to refer to the usual notions of non-derived algebraic geometry. From a formal perspective, the relationship between derived schemes and classical schemes is analogous to the relationship between classical schemes and reduced classical schemes. A derived scheme has an underlying classical scheme, and intuitively one thinks of a derived scheme as an ``infinitesimal thickening'' of its underlying classical scheme. Formation of the underlying classical scheme (an operation called ``classical truncation'') defines a functor which is right adjoint to a fully faithful embedding from classical schemes to derived schemes. 
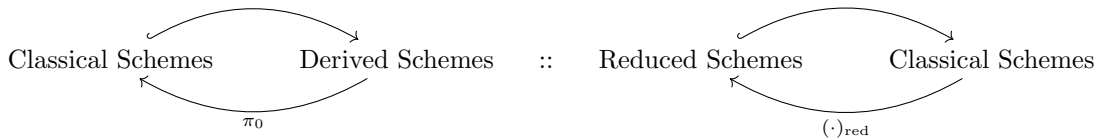
\[
\begin{tikzcd}
\text{Classical Schemes}  \ar[r, hook, bend left]   &  \ar[l, "\pi_0", bend left]  \text{Derived Schemes}
\end{tikzcd}  \quad  :: \quad 
\begin{tikzcd}
\text{Reduced Schemes}  \ar[r, hook, bend left]   &  \ar[l, "(\cdot)_{\mrm{red}}", bend left]  \text{Classical Schemes}
\end{tikzcd}
\]

One advantage of considering non-reduced schemes is that it gives a natural interpretation of the tangent space, as maps from the spectrum of the dual numbers. Analogously, derived algebraic geometry gives a very natural interpretation of the \emph{cotangent complex} (which governs deformation theory), even for a classical scheme. Indeed, the higher cohomology groups of the tangent complex can be viewed in terms of maps from certain ``derived infinitesimal schemes'', which are derived generalizations of dual numbers. The theory of the cotangent complex plays a crucial technical role in this section.

The passage from classical schemes to classical stacks goes through the ``functor of points'' perspective: a scheme can be viewed as a functor from commutative rings to sets, and a stack can be viewed as a functor from commutative rings to groupoids. Generalizing this perspective, a derived scheme can be interpreted as a functor from \emph{simplicial} commutative rings to \emph{simplicial} sets. However, when working with simplicial objects, the notion of equivalence should be homotopy-theoretic. Consequently, the test and target categories should be the ``non-abelian derived categories'' of simplicial commutative rings and of simplicial sets, respectively, which are called the \emph{$\infty$-category of simplicial rings} and the \emph{$\infty$-category of simplicial sets}, respectively. We therefore define a \emph{derived stack} to be a functor from the $\infty$-category of simplicial rings to the $\infty$-category of simplicial sets, satisfying certain descent conditions.  A derived stack has a classical truncation by restricting its domain to classical rings, and this defines a functor which has a fully faithful left adjoint, embedding classical stacks into derived stacks. 
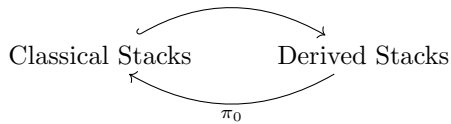
\[
\begin{tikzcd}
\text{Classical Stacks}  \ar[r, hook, bend left]   &  \ar[l, "\pi_0", bend left]  \text{Derived Stacks}
\end{tikzcd}  
\]

For us, derived algebraic geometry will be used in the following way. We wish to attach fundamental classes to our special cycles $\cZ_{\cE}^r(a)$, but these spaces are poorly behaved in general, so we need to construct \emph{virtual} fundamental cycles. For example, the virtual fundamental cycle of $\cZ_{\cE}^r(a)$ should have dimension $(n-\rank \cE) r$, but $\dim \cZ_{\cE}^r(a)$ may be much larger, e.g when $a=0$ then $ \cZ_{\cE}^r(a)$ has components of the maximal dimension $nr$. However, it turns out that there is a natural \emph{derived special cycle} $\sZ_{\cE}^r(a)$ with the correct \emph{virtual} dimension, which is moreover quasi-smooth (the analogue of LCI in derived algebraic geometry), and whose classical truncation recovers $\cZ_{\cE}^r(a)$. The property that $\sZ_{\cE}^r(a)$ is quasi-smooth implies that it has an intrinsic virtual fundamental class $[\sZ_{\cE}^r(a)]$, and derived invariance  of Chow groups (see \S\ref{sssec: derived invariance}) allows to view it as an element of $\Ch_{(n-\rank \cE)r}(\cZ_{\cE}^r(a))$.

\begin{example}Even when working with classical schemes, derived structure often shows up implicitly because of derived intersections. Locally, this is based on the ``derived tensor product'' operation, which can produce a non-classical simplicial commutative ring even when the inputs are classical. From this optic, the virtual fundamental cycles of \cite{FYZ} come from derived stacks obtained by taking derived intersections in a particular presentation of $\cZ_{\cE}^r(a)$ as a fibered product of smooth classical stacks. In particular, for non-singular $a$ the derived stack $\sZ_{\cE}^r(a)$ is a global complete intersection in the derived sense; more generally, derived algebraic geometry provides an \emph{intrinsic} construction of a virtual fundamental class to any derived stack which \emph{locally} looks like a derived fibered product of smooth classical schemes (this is one formulation of quasi-smoothness). Crucially this is a local property and we do not require any \emph{global} presentation as a derived intersection of smooth stacks, which we do not have in the case of singular coefficients. 
\end{example}

\subsubsection{Notational conventions} We will use script letters such as $\sX, \sY$ for derived stacks, and calligraphic letters such as $\cX, \cY$ for classical stacks. We will often use $\cX$ to denote the classical truncation of $\sX$ (defined later in \S \ref{sssec: classical truncation}). 

\subsubsection{Derived (Artin) stacks} For the framework of derived stacks, we follow \cite[\S 1.1]{KhanI}. To summarize, derived stacks are defined as functors from a test category to a target category, satisfying a sheaf condition, where: 
\begin{itemize}
\item The test category is the $\infty$-category of simplicial commutative rings. This can be constructed as in \cite[Definition 4.1.1]{DAGV}; an intrinsic characterization can be found in \cite[\S 5.1]{CS}. Following Clausen-Scholze we call it the category of \emph{animated rings}, and use the phrase ``animated ring'' to indicate an object of this category. 
\item The target category is the $\infty$-category of simplicial sets. Similar remarks apply as above. Following Clausen-Scholze we call it the category of \emph{anima} (also called ``$\infty$-groupoid'', or ``space''), and use the phrase ``anima" to indicate an object of this category. 
\end{itemize}
Thus, derived stacks $\sY$ over $k$ are functors from the category of animated rings to the category of anima, denoted $R_{\bu} \mapsto \sY(R_{\bu})$, satisfying \'{e}tale hyperdescent. 

We define \emph{$n$-geometric derived stacks} as in \cite[\S 1.3.3]{TV08}\footnote{There are differing conventions on $n$-stacks -- for example the above notion differs from the ``$n$-algebraic stacks'' of \cite[\S 5.2]{To10} -- but they all produce the same notion of Artin stack, which is the only one of importance to us.}, and \emph{derived Artin stacks} to be derived stacks which are $n$-geometric for some $n$. 

Functors between $\infty$-categories cannot be constructed by specifying the images of objects and 1-morphisms, although it is common practice in some parts of the literature to describe them in this way (often even omitting the descriptions of 1-morphisms when they are obvious), with an implicit understanding that the reader can fill in the formal construction. Our practice will be to instead construct functors by formal operations bootstrapping off of elemental constructions established in \cite{TV05, TV08}; we will then describe their effects on 0-cells (i.e., objects) for informal intuition but this should not be mistaken for a formal definition.


\subsubsection{Representable morphisms}\label{sssec: representable morphism} \emph{Affine derived schemes} are the representable derived stacks. Derived schemes are the derived stacks that have a Zariski open cover by affine derived schemes. Following \cite[Definition 1.3.3.1, Definition 1.3.3.7]{TV08}, we say that a morphism of derived stacks $f \co \sX \rightarrow \sY$ is \emph{$n$-representable} if for any derived scheme $S$ and any map $S \rightarrow \sY$, the fibered product $\sX \times_{\sY} S$ is $n$-geometric. We say $f$ is \emph{representable} if it is $n$-representable for some $n$. (Note that this condition is much broader than representability for morphisms of classical Artin stacks, the latter of which is analogous to ``$(-1)$-representable'' in our sense.) 

By \cite[Proposition 1.3.3.3]{TV08}, the class of representable morphisms is closed under isomorphisms, (homotopy) pullbacks, and compositions. 

\subsubsection{Derived terminology} We remind the reader that all operations in $\infty$-categories are ``homotopical'', so that tensor products of animated rings correspond to ``derived tensor products'', fiber products of derived stacks correspond to ``homotopy fiber products'', the ``fiber'' of a map of complexes $\cK \xrightarrow{f} \cK'$ in the derived category means the ``derived fiber'' $\mrm{Cone}(f)[-1]$, etc. (If we need to refer to a classical fibered product of classical stacks $\cX$ and $\cY$ over $\cT$, we will denote it by $\cX \stackrel{\mrm{cl}}\times_{\cT} \cS$.) At some points we include the adjectives ``homotopy'' or ``derived'' to emphasize this, but it applies everywhere in this section. 


\subsubsection{Classical truncation}\label{sssec: classical truncation} We shall frequently invoke the notion of the ``underlying classical stack'', i.e. ``classical truncation'', of a derived stack. Here we recall what this means. If $R_{\bu}$ is a simplicial commutative ring, then its ``underlying classical ring'' is $\pi_0(R_{\bu})$. (In topological terminology this is the ``first Postnikov truncation'' of $R_{\bu}$, which explains the synonymous terminology ``classical truncation''.) This descends to a functor on animated rings, which is left adjoint to the inclusion of classical (i.e. discrete) commutative rings into animated rings. On the opposite categories, we get a fully faithful functor from affine schemes to derived affine schemes which is left adjoint to the classical truncation. 

This operation then glues in the Zariski topology to give a functor $T \mapsto \pi_0(T)$ from derived schemes to classical (discrete) schemes, which is right adjoint to a fully faithful inclusion functor from classical schemes into derived schemes. By abuse of notation we may regard $\pi_0(T)$ as a derived scheme via this inclusion; then the unit of the adjunction is a map $\pi_0(T) \rightarrow T$, natural in $T$. 

Finally, if $\sX$ is a derived stack, then its underlying classical stack $\sX_{\mrm{cl}}$ is the restriction of $\sX$ along the embedding $\{ \text{Classical affine schemes}\} \inj \{\text{Derived affine schemes}\}$. The classical truncation functor $\sX \mapsto \sX_{\mrm{cl}}$ has a left adjoint, which can be described as the sheafification of the left Kan extension on the underlying prestacks, and is fully faithful \cite[\S I.2.6]{GRI}. The left adjoint gives an inclusion $\{ \text{Classical stacks}\} \inj \{\text{Derived stacks}\}$, and the unit of the adjunction is the ``classical truncation map'' 
$$\io_{\sX} \co \sX_{\mrm{cl}} \rightarrow \sX$$ 
functorial in $\sX$, which we call the inclusion of the underlying classical stack. We say that $\sX$ is \emph{isomorphic to its underlying classical stack} (or just \emph{classical} for brevity) if $\io_{\sX}$ is an isomorphism.

\subsubsection{Derived mapping stacks} We give some examples of derived Artin stacks which are of particular relevance to this paper. For an animated $k$-algebra $R_{\bu}$ and a scheme $X$ over $k$, we abbreviate $X_{R_{\bu}} := X \times_{\Spec k} \Spec R_{\bu}$. 

\begin{example}\label{ex: derived mapping stack}
Let $X$ be a proper scheme over $k$ and $\sY$ a derived Artin stack locally of finite presentation over $k$. The \emph{derived mapping stack} $\sMap(X, \sY)$ is constructed in \cite[\S 3.6]{TV05} and \cite[\S 2.2.6.3]{TV08} (the first reference constructs an internal hom on the model category of stacks on a site, while the second reference establishes its geometricity properties). At the level of $0$-cells, $\sMap(X, \sY)$ sends an animated $k$-algebra $R_{\bu}$ to the anima of morphisms
\[
X_{ R_{\bu}} \rightarrow \sY
\]
over $k$. According to \cite[Corollary 3.3]{To14}, $\sMap(X, \sY)$ is a derived Artin stack locally of finite presentation over $k$. 

More generally, in the above situation, if both $X$ and $\sY$ are over a scheme $S$ over $k$, we can define the derived mapping stack $\sMap_{S}(X, \sY)$ as the (homotopy) fiber of $\sMap(X, \sY) \rightarrow \sMap(X, S)$, induced by $\sY \rightarrow S$, over the given map $X \rightarrow S$. Assume that $S$ is of finite type over $k$. By \cite[Corollary 3.3]{To14}, if $X$ is proper over $k$ and $\sY$ is locally of finite presentation over $S$, then $\sMap_{S}(X, \sY)$ is a derived Artin stack locally of finite presentation over $k$. At the level of $0$-cells, $\sMap_{S}(X, \sY)$ sends $R_{\bu}$ to the anima of morphisms $X_{ R_{\bu}} \rightarrow \sY$ over $S$. 

When $S=X$ we write $\sSect(X,\sY)$ for $\sMap_{X}(X,\sY)$.
\end{example}

\begin{example}\label{ex: derived Bun_G}
Let  $X$ be a scheme over $k$, and let $G$ be a smooth algebraic group (which for us means by definition that it is of finite type) over $X$. Regard the classical classifying stack $BG = [X/G]$ as a derived stack over $X$ via the embedding discussed above. At the level of 0-cells, the derived mapping stack $\sSect(X, BG)$ sends $R_{\bu}$ to the anima of $G$-bundles on $
X_{ R_{\bu}}$. When $X$ is a smooth projective curve we will see in Corollary \ref{cor: T for Bun_G} that $\sSect(X, BG)$ is isomorphic to its underlying classical stack, which is $\Bun_G$.  
\end{example}

\begin{example}\label{ex: derived standard/gerbe}
Let $X$ be a proper scheme over $k$. Let $G \rightarrow X$ be a smooth algebraic group and $V \rightarrow X$ a vector bundle that is a representation of $G$. We apply Example \ref{ex: derived mapping stack} with $\sY = V/G$ (a classical stack). This results in the derived stack of sections $\sSect(X, V/G)$, which at the level of 0-cells sends an animated ring $R_{\bu}$ to the $\infty$-groupoid of $(\cF, s)$ where 
\begin{itemize}
\item $\cF \xrightarrow{\pi} X_{ R_{\bu}} $ is a $G$-bundle. 
\item $s$ is in the (homotopy) fibered product $\Map_G(\cF, V_{R_{\bu}}) \stackrel{h}\times_{\Map(\cF, X)} \{\pi\}$ where the map $\Map_G(\cF, V) \rightarrow \Map(\cF, X_{R_{\bu}})$ is induced by composition with the tautological map $V \rightarrow X$.
\end{itemize}
The (derived) fiber of the map $\sSect(X, V/G) \rightarrow \Bun_G$ (here we are using Example \ref{ex: derived Bun_G} to identify $\sSect(X, BG)$, which is a priori a ``derived version'' of $\Bun_G$, with $\Bun_G$) over a field-valued point $\cF \in \Bun_G(\kappa)$ is the derived scheme $R\Gamma(X_{\kappa}, \cF \times^G V)$. Note for contrast that the classical fiber of the map of classical stacks $\mrm{Sect}(X, V/G) \rightarrow \Bun_G$ is $H^0(X_\kappa, \cF \times^G V)$. We spell out how $R\Gamma(X_{\kappa}, \cF \times^G V)$ is viewed as a derived scheme:
\begin{enumerate}
\item $R\Gamma(X_\kappa, \cF \times^G V)$ is a connective perfect cochain complex (i.e., cohomology groups vanish in negative degrees) in the derived category of $\kappa$-modules.  
\item Its dual $R\Gamma(X_\kappa, \cF \times^G V)^*$ is a connective  perfect \emph{chain} complex (i.e., homology groups vanish in negative degrees) in the derived category of $\kappa$-modules, which by the Dold-Kan correspondence may be viewed as an animated $\kappa$-module. 
\item The forgetful functor from animated $\kappa$-algebras to animated $\kappa$-modules admits a left adjoint, the derived symmetric algebra functor $\Sym^{\bu}_\kappa$. 
\item The derived scheme $R\Gamma(X_\kappa, \cF \times^G V)$ is the spectrum of $\Sym_\kappa^{\bu} \left(R\Gamma(X_\kappa, \cF \times^G V)^*\right)$. 
\end{enumerate}
\end{example}

\subsubsection{Cotangent complexes}\label{sssec: cotangent complex}
We refer to \cite{TV08, To10} for the theory of the cotangent complex to a morphism $f \co \sX \rightarrow \sY$ of derived stacks, denoted $\bL_f$. The \emph{tangent complex} to $f$ is $\bT_f := \cRHom
_{\cO_{\sX}}(\bL_f, \cO_{\sX})$. Sometimes these will be denoted $\bL_{\sX/\sY}$ and $\bT_{\sX/\sY}$ when the map is clear. When $f$ is the structure morphism $f \co \sX \rightarrow \Spec k$, we abbreviate $\bT_{\sX} := \bT_{\sX/\Spec k}$ and $\bL_{\sX} := \bL_{\sX/\Spec k}$.

A useful characterization of the cotangent complex of $f \co \sX \rightarrow \Spec k$ is as follows \cite[p.37]{To09}. Let $R_{\bu}$ be an animated $k$-algebra and recall that for any animated $R_{\bu}$-module $M_{\bu}$ there is an animated $R_{\bu}$-algebra $R_{\bu} \oplus M_{\bu}$, which on homotopy groups is the square-zero extension of $\pi_*(R_{\bu})$ by $\pi_*(M_{\bu})$. Then for any map $a \co \Spec R_{\bu}\rightarrow \sX$ and any animated $R_{\bu}$-module $M_{\bu}$, there is a natural equivalence between $\RHom_{R_{\bu}-\mrm{Mod}}(a^* \bL_f, M_{\bu})$ and the homotopy fiber of $\sX(R_{\bu} \oplus M_{\bu}) \rightarrow \sX(R_{\bu})$ over $a \in \sX(R_{\bu})$. 

The following fundamental facts will be used frequently: 
\begin{itemize}
\item For a sequence of morphisms $\sX \xrightarrow{f} \sY \xrightarrow{g} \sZ$, there is an exact triangle in $\mrm{QCoh}(\sX)$: 
\[
f^* \bL_g \rightarrow \bL_{g \circ f} \rightarrow  \bL_f.
\]
\item For a Cartesian square 
\[
\begin{tikzcd}
\sX' \ar[r, "g'"] \ar[d, "f'"] & \sX \ar[d,"f"] \\
\sY' \ar[r,"g"] & \sY
\end{tikzcd}
\]
we have $(g')^* \bL_f \xrightarrow{\sim} \bL_{f'}$. Given compatible maps to a base derived stack $\sS$, we then deduce an exact triangle 
\[
\bL_{\sX'/\sS} \rightarrow (g')^* \bL_{\sX/\sS} \oplus (f')^* \bL_{\sY'/\sS} \rightarrow  (f \circ g')^* \bL_{\sY/\sS}.
\]
\end{itemize}

\begin{lemma}\label{lem: classicality of smooth derived stacks}
Let $\sY$ be a locally finite type derived Artin stack over $k$. Suppose that the cotangent complex $\bL_{\sY/k}$ has tor-amplitude in $[0, \infty)$. Then $\iota_{\sY} \co (\sY)_{\mrm{cl}} \rightarrow \sY$ is an equivalence.
\end{lemma}

\begin{proof}
This is well-known, but at the referee's suggestion we sketch a proof. Since the assertion can be checked locally in the smooth topology, we may assume that $\sY$ is a connected derived affine scheme of finite type over $k$. Then $\bL_{\sY/k}$ has tor-amplitude in $[-\infty, 0]$ on general grounds, so the assumption forces $\bL_{\sY/k}$ to be represented by a vector bundle (in degree $0$). The map on cotangent complexes induced by $\iota_{\sY}$ is always an isomorphism on $H^0$ and a surjection on $H^{-1}$ (for example dualize \cite[(7.11)]{GV18}). Then by obstruction theory $\sY_{\mrm{cl}}$ is smooth of dimension equal to the virtual dimension of $\sY$, so we deduce that $\iota_{\sY}$ induces an isomorphism on cotangent complexes. As $\iota_{\sY}$ induces an isomorphism on $\pi_0$ by definition, it then induces an isomorphism on all homotopy groups by \cite[Corollary 3.2.27]{Lur04} (an analogue of the ``Hurewicz Theorem'' for simplicial commutative rings), hence is an equivalence. 
\end{proof}

\subsubsection{Quasi-smoothness} A key role is played by the notion of \emph{quasi-smooth} derived Artin stacks, and more generally quasi-smooth morphisms. Recall that a morphism $f \co \sX \rightarrow \sY$ of derived Artin stacks is \emph{quasi-smooth} if it is locally of finite presentation and the relative cotangent complex $\bL_{f}$ is perfect of Tor-amplitude $[-1,\infty)$. (Here we are using \emph{cohomological} grading, as opposed to the homological grading of \cite{KhanI}, so this means that $H^i(\bL_f \otimes_{\cO_{\sX}} \cE)$ vanishes for $i<-1$ for every discrete quasi-coherent sheaf $\cE$ on $\sX$.) Given $f$ locally of finite presentation with $\bL_f$ a perfect complex, $f$ is quasi-smooth if and only if the fiber of $\bL_f$ at all geometric points is acyclic in (cohomological) degrees $<-1$ \cite[\S 2.1]{AG15}. This is the derived analogue of being LCI, and for that reason is also sometimes referred to as ``derived LCI''. In particular, a classical LCI morphism between classical stacks, regarded as derived stacks, is quasi-smooth.

The following facts are immediate from basic properties of the cotangent complex:
\begin{itemize}
\item The composition of quasi-smooth morphisms is quasi-smooth. 
\item The (derived) base change of any quasi-smooth morphism is quasi-smooth. Note that the classical analogue is completely false for classical LCI morphisms! 
\end{itemize}
If $\sX \rightarrow \Spec k$ is a quasi-smooth morphism, then we simply say that \emph{$\sX$ is quasi-smooth} in particular, a classical LCI stack over $k$ is quasi-smooth when regarded as a derived stack. As we shall see later, quasi-smooth derived Artin stacks are those to which we can naturally associate a virtual fundamental class, which is why this notion is important for us. 

We recall for comparison that a morphism of derived stacks $f \co \sX \rightarrow \sY$ is \emph{smooth} if it is locally of finite presentation and $\bL_f$ is perfect of Tor-amplitude $[0, \infty)$. In particular, this includes smooth morphisms of classical stacks.

\begin{example}
Suppose $\sX \rightarrow \Spec k$ is smooth. Then the diagonal map $\sX \rightarrow \sX \times_k \sX$ is quasi-smooth. 
\end{example}

\subsection{Tangent complexes to derived mapping stacks}\label{ss:tangent} 

\begin{lemma}\label{lem: tangent to derived mapping stack} 

Let $S$ be a derived stack over $k$ with perfect cotangent complex and $X$ be a smooth proper scheme over $k$ with a map to $S$. Let $\sY$ be a finite type derived stack over $S$ with perfect relative cotangent complex $\bL_{\sY/S}$. Then the cotangent complex $\bL_{\sMap_{S}(X,\sY)}$ is perfect, and its pullback to any $R_{\bu}$-point $f \co X_{R_{\bu}} \rightarrow \sY$, for any animated ring $R_{\bu}$, is naturally in $R_{\bu}$ isomorphic to $R\pr_*( f^* \bL_{\sY/S} \otimes \omega_{X/k})$, where $\pr$ is the projection map $X_{R_{\bu}} \rightarrow \Spec R_{\bu}$ and $\omega_{X/k}$ is the dualizing complex of $X$. 

In particular, $\bT_{\sMap_{S}(X,\sY)}|_f$ is naturally in $R_{\bu}$ isomorphic to $R\pr_*(f^* \bT_{\sY/S})$.
\end{lemma} 

In the statement of the Lemma, ``naturally isomorphic in $R_{\bu}$'' means that there is natural transformation between the two functors; informally speaking, that the isomorphisms base change coherently along $R'_{\bu} \rightarrow R_{\bu}$.

\begin{proof}We apply \cite[Proposition 5.1.10]{HLP}, which implies that $\bL_{\sMap_{S}(X,\sY)}|_f$ is isomorphic to $\pr_+ (f^* \bL_{\sY/S})$, where $\pr_+$ is the left adjoint to $\pr^*$. Since $X/k$ is smooth and proper, we have $\pr^! (-) = \omega_{X/k} \otimes \pr^* (-)$, and that the left adjoint of $\pr^!$ exists and is $R\pr_*$, so the left adjoint of $\pr^*$ is $R\pr_*(- \otimes \omega_{X/k})$. 

The last sentence follows from applying Serre duality to the description of $\bL_{\sMap_{S}(X,\sY)}|_f$. 

\end{proof}



\begin{cor}\label{cor: T for Bun_G} Let $X$ be a smooth proper curve over $k$.  Let $\sG$ be a smooth (classical) gerbe over $X$. Then $\sSect(X, \sG)$ is isomorphic to its classical truncation $\Sect(X, \sG)$, which is smooth.
\end{cor}

\begin{proof}
Because a smooth gerbe is locally in the smooth topology isomorphic to the classifying stack of a group scheme, the relative tangent complex $\bT_{\sG/X}$ is concentrated in degree $-1$, hence for any section $f \co X_{R_{\bu}} \rightarrow \sG$ the cohomology groups of $R\pr_*(f^* \bT_{\sG/X})$ are non-vanishing only in degrees $-1, 0$. We conclude by applying Lemma \ref{lem: classicality of smooth derived stacks}. 
\end{proof}

\begin{example}\label{ex: T for Bun_G} Let $X$ be a smooth proper scheme over $k$ and $G \rightarrow X$ a smooth group scheme. Then $\bT_{BG/X} \cong \Lie(G/X)[1]$. Lemma \ref{lem: tangent to derived mapping stack} implies that $\bL_{\sSect(X, BG)}$ is perfect, and $\bT_{\sSect(X,BG)}$ pulled back to $\Spec R_{\bu}$ via a $G$-torsor $\cF$ over $X_{R_{\bu}}$ is isomorphic to $R\pr_*(\cF\times^{G} \Lie(G/X))[1]$ naturally in $R_{\bu}$ .

\end{example}

\begin{cor}\label{cor: TC Sect to Bun_G} Let $S$ be a derived stack over $k$ with perfect cotangent complex.  Let $X$ be a smooth proper scheme over $k$ with a map to $S$. Suppose $\mu \co \sY \rightarrow \sY'$ is a morphism of finite type derived stacks over $S$ such that $\bL_{\mu}$ is perfect.  Then the induced map $\mu_* \co \sMap_{S}(X,\sY)  \rightarrow \sMap_{S}(X, \sY')$ has perfect relative cotangent complex, and for any $R_{\bu}$-point $f \co X_{R_{\bu}} \rightarrow \sY$ of $\sMap_{S}(X,\sY)$, $\bT_{\mu_{*}}|_f$ is isomorphic to $R\pr_* (f^{*}\bT_{\mu})$ naturally in $R_{\bu}$. 

In particular, if $X$ is a smooth projective curve and $\mu$ is smooth, then $\mu_*$ is quasi-smooth. 
\end{cor}

\begin{proof}
The perfectness of the cotangent complex for $\mu_*$ and description of the tangent complex follow from the functoriality of Lemma \ref{lem: tangent to derived mapping stack} with respect to $\sY$. The last sentence follows because in this situation, $R\pr_* $ has cohomological amplitude $1$ and  $\bT_{\mu}$ is concentrated in degree $\le 0$, so $R\pr_* ( f^* \bT_{\mu})$ has cohomology concentrated in degrees $\le 1$. 
\end{proof}

\begin{example}\label{ex: TC Sect to Bun_G} Let $X$ be a smooth proper scheme over $k$. Let $G$ be a smooth group scheme over $X$ and $V\to X$ a vector bundle that is a representation of $G$ over $X$. Consider an $R_{\bu}$-point of $\sSect(X, V/G)$, represented by a $G$-torsor $\cF$ on $X_{R_{\bu}}$ and $s \in R\Gamma(X_{R_{\bu}}, \cF \times^G V)$.

We give a more concrete description of various tangent complexes in this situation. 
\begin{enumerate}
\item The tangent complex to $\sSect(X, V/G)$ at the $R_{\bu}$-point $(\cF, s)$ is naturally in $R_{\bu}$ isomorphic to 
\[
R\pr_*(\underbrace{\cF \times^{G} \Lie(G/X)}_{\deg -1} \xrightarrow{\cdot s} \underbrace{\cF \times^{G}  V}_{\deg 0} )
\]
where the meaning of the differential $\cdot s$ is as in \S \ref{ssec: tangent complex}, and $\pr \co X_{R_{\bu}} \rightarrow \Spec R_{\bu}$ is the projection map. 

\item The map of tangent complexes induced by $\sSect(X, V/G) \xrightarrow{\pi} \Bun_G$ pulled back to $\Spec R_{\bu}$ via $(\cF,s)$ is naturally in $R_{\bu}$ isomorphic to 
\begin{equation}\label{eq: TC Sect to Bun_G}
R\pr_*(\underbrace{\cF \times^{G} \Lie(G/X)}_{\deg -1} \xrightarrow{\cdot s} \underbrace{\cF \times^{G}  V}_{\deg 0} ) \xrightarrow{t} R\pr_*  (\underbrace{\cF \times^{G} \Lie(G/X)}_{\deg -1})
\end{equation}
where the map $t$ is induced by the truncation of complexes. 
\item If $X$ is a smooth projective curve, $\sSect(X, V/G) \xrightarrow{\pi} \Bun_G$ is quasi-smooth, and $\sSect(X, V/G) $ is quasi-smooth. 
\end{enumerate}
\end{example}

We consider the following ``relative" variant of mapping stacks.

\begin{defn}\label{def: relative mapping stack} Let $B$ be a scheme, let $X$ be a smooth proper scheme over $B$ and $\sY$ be a finite type derived stack over $B$ with perfect relative cotangent complex $\bL_{\sY/B}$. The \emph{relative mapping stack} (cf. \cite[Corollary 3.3]{To14}) $\sMap(X/B, \sY/B)$ is the internal Hom between $X$ and $\sY$ viewed as derived stacks over $B$; at the level of $0$-cells it sends $\Spec R_{\bu} \rightarrow B$ to the anima of morphisms $X \times_B \Spec R_{\bu} \rightarrow \sY$ over $B$.

As a variant, if there is a map $\pi \co \sY \rightarrow X$ over $B$ then we write $\sSect(X/B, \sY/B)$ for the homotopy fiber of $\sMap(X/B, \sY/B)$ over the identity element of $\sMap(X/B, X/B)$. At the level $0$-cells, $\sSect(X/B, \sY/B)$ sends $\Spec R_{\bu} \rightarrow B$ to the anima of commutative diagrams
\[
\begin{tikzcd} &  \sY \ar[d, "\pi"] \\
X \times_B \Spec R_{\bu} \ar[r, "\pr_1"] \ar[ur, dashed]  & X
\end{tikzcd}
\]
Let $\cY$ be the classical truncation of $\sY$. Then the classical truncation of $\sMap(X/B, \sY/B)$ is the classical mapping stack $\Map(X/B, \cY/B)$ and the classical truncation of $\sSect(X/B, \sY/B)$ is $\Sect(X/B, \cY/B)$.
\end{defn}

\begin{lemma}\label{lem: tangent to relative mapping stack} 
(1) Let $S$ be a derived stack, let $X$ be a smooth proper scheme over $S$, and let $\sY$ be a finite type derived stack over $S$ with perfect relative cotangent complex $\bL_{\sY/S}$. Then the relative cotangent complex $\bL_{\sMap(X/S,\sY/S) / S}$ is perfect, and its pullback to any $R_{\bu}$-point $f \co X \times_S \Spec(R_{\bu}) \rightarrow \sY$, for any animated ring $R_{\bu}$, is naturally in $R_{\bu}$ isomorphic to $R\pr_*( f^* \bL_{\sY/S} \otimes \omega_{X/S})$, where $\pr$ is the projection map $X \times_S \Spec R_{\bu} \rightarrow \Spec R_{\bu}$ and $\omega_{X/S}$ is the relative dualizing complex of $X \rightarrow S$. 

In particular, the relative tangent complex $\bT_{\sMap(X/S,\sY/S)/S}|_f$ is naturally in $R_{\bu}$ isomorphic to $R\pr_*(f^* \bT_{\sY/S})$.

(2) Suppose that in addition to the setup above, there is a map $\pi \co \sY \rightarrow X$ over $S$ such that $\bL_{\pi}$ is perfect. Then the relative cotangent complex $\bL_{\sSect(X/S,\sY/S) / S}$ is perfect, the relative tangent complex $\bT_{\sSect(X/S,\sY/S)/S}|_f$ is naturally in $R_{\bu}$ isomorphic to $R\pr_*(f^* \bT_{\pi})$.
\end{lemma}

\begin{proof}
(1) Similar calculation as for Lemma \ref{lem: tangent to derived mapping stack}. Then (2) follows from (1) as in Corollary \ref{cor: TC Sect to Bun_G}.
\end{proof}

\subsection{(Un)derived $\Hk^{r}_{G}$ and $\Sht^{r}_{G}$}

We now resume our convention that $X$ is a smooth projective curve over $k$.  For a smooth gerbe $\sG \rightarrow X$, Corollary \ref{cor: T for Bun_G} implies that the derived stack $\Bun_{\sG} = \sSect(X, \sG)$ is isomorphic to its classical truncation. Next we define and analyze derived Hecke stacks for gerbes of unitary type and $B\GL(n)'$. They will also turn out to be isomorphic to their classical truncations, a fact that is needed later to compute cycle classes in explicit terms. 

We will rewrite Definition \ref{defn: hitchin hecke} in a way that is amenable to a derived construction.
\begin{defn}\label{defn: hecke gerbe unitary}
We define a \emph{rank n (unitary) Hecke gerbe} $\sH\rightarrow X' \times X$ that sends a $k$-scheme $S$ to the groupoid of $(x', \xi, \cF, h)$ where 
\begin{itemize}
\item $x' \in X'(S)$. 
\item $\xi \in X(S)$. 
\item $\cF$ is a vector bundle of rank $n$ on $S'$, which is the fibered product 
\[
\begin{tikzcd}
S' \ar[r,"\xi'"] \ar[d, "\nu_S"] & X' \ar[d, "\nu"] \\
S \ar[r, "\xi"] & X
\end{tikzcd}
\]
Let $\s_{S}: \Aut_{S}(S')$ be the base change of $\s\in \Aut_{X}(X')$ along $\xi$. Let $\D(X'), \Gamma_{\s}\subset X'\times X'$ be the diagonal and the graph of $\s$ respectively. Using the map $(x'\c\nu_{S},\xi'):S'\to X'\times X'$, we form the closed subschemes
\begin{equation*}
S'_{\xi'=x'}:=S'\times_{X'\times X'}\D(X'), \quad \mbox{ and } S'_{\xi'=\s(x')}:=S'\times_{X'\times X'}\Gamma_{\s}
\end{equation*}
of $S'$ that are disjoint.
\item $h \co \cF \rightarrow \sigma_{S}^* \cF^* \otimes \nu_S^*\xi^* (\omega_X \otimes \LL)$ is a map such that $\sigma^* h^* = h$, and $\coker(h)$ is locally free of rank 1 supported on $S'_{\xi'=x'}\coprod S'_{\xi'=\s(x')}$.
\end{itemize}
Note that upper modifying\footnote{This refers to the unique bundle $\cF^{\sharp}$ such that $\cF^{\sharp}/\cF$ is isomorphic to the given torsion sheaf; see \cite[Definition 6.5]{FYZ} for more details.} $\cF$ along either $\coker(h)|_{S'_{\xi'=x'}}$ or $\coker(h)|_{S'_{\xi'=\s(x')}}$ produces two maps $\sH$ to $BU(n)_{\LL}\times X'$ as gerbes over $X \times X'$, 
\[
\begin{tikzcd}
& \sH \ar[dl, "h_0"'] \ar[dr, "h_1"] \\
BU(n)_{\LL} \times X' & & BU(n)_{\LL} \times X' 
\end{tikzcd}
\]
\end{defn}

Now regarding $\sH \rightarrow X' \times X$ as a map over $X'$, we may consider the derived relative mapping stack $\sSect(X' \times X/X', \sH/X')$ and (noting that $\sH$ is a classical stack) its classical truncation $\Sect(X' \times X/X', \sH/X')$. Then $\sSect(X' \times X/X', \sH/X')$ is a derived stack over $X'$, and the two maps $h_{0},h_{1}: \sH \rightarrow BU(n)_{\LL} \times X'$ of gerbes over $X \times X'$ induce maps 
\begin{equation}\label{eq: gerbe hecke map}
\begin{tikzcd}
& \sSect(X' \times X/X', \sH/X')  \ar[dl, "h_0"'] \ar[dr, "h_1"]  \\
\Bun_{U(n), \LL} \times X'  & & \Bun_{U(n), \LL} \times X'
\end{tikzcd}
\end{equation}
with classical truncation 
\begin{equation}\label{eq: gerbe hecke map classical}
\begin{tikzcd}
& \Sect(X' \times X/X', \sH/X') \ar[dl, "h_0"'] \ar[dr, "h_1"]  \\
\Bun_{U(n), \LL} \times X' & & \Bun_{U(n), \LL} \times X'
\end{tikzcd}
\end{equation}

\begin{lemma}\label{lem: hecke as sect}
There is a canonical isomorphism $\Sect(X' \times X/X', \sH/X') \cong \Hk_{U(n), \LL}^1$ intertwining \eqref{eq: gerbe hecke map classical} with 
\[
\begin{tikzcd}
& \Hk_{U(n), \LL}^1 \ar[dl, "\pr_0 \times \pr_X"']  \ar[dr, "\pr_1  \times \pr_X"] \\
 \Bun_{U(n), \LL}  \times X' & & \Bun_{U(n), \LL} \times X'
\end{tikzcd}
\]
\end{lemma}

\begin{proof}
Let $t:T\to X'$ be a scheme over $X'$. Let us examine the $T$-points of $\Sect(X' \times X/X', \sH/X')$. By definition the groupoid of $T$-points is the groupoid of diagrams
\[
\begin{tikzcd}
& \sH \ar[d]  \\
T \times X \ar[ur, dashed] \ar[r, "t \times \Id_X"] &  X' \times X
\end{tikzcd}
\]
We unpack the meaning of such a diagram. From the definition of $\sH$, this is the groupoid of pairs $(\cF, h)$ where $\cF$ is a vector bundle on $T \times X'$ and $h \co \cF \rightarrow (\Id_{T}\times\sigma)^* \cF^* \otimes \pr_{X'}^{*}\nu^* (\omega_X \otimes \LL)$, such that $\coker(h)$ is a line bundle on 
\[
(T \times X') \times_{X' \times X'}( \Delta(X') \coprod \Gamma_{\s})\cong \Gamma_t \cup \Gamma_{\sigma (t)}
\]
where $\Gamma_t$ and $\Gamma_{\sigma(t)}$ are the graphs of $t \in X'(T)$ and $\sigma(t) \in X'(T)$, respectively. 

From this data we construct $(\cF_0 \leftarrow \cF^\flat_{1/2} \rightarrow \cF_1) \in \Hk^1_{U(n), \LL}(T)$ as follows. Let $\cF_{1/2}^\flat$ be $\cF$ above, $\cF_0$ the upper modification of $\cF_{1/2}^\flat$ along $\coker(h)|_{\Gamma_t}$, and $\cF_1$ the upper modification of $\cF_{1/2}^\flat$ along $\coker(h)|_{\Gamma_{\sigma(t)}}$. Then $h$ extends to Hermitian structures on $\cF_0$ and $\cF_1$. 

Conversely, given $(t, \cF_0 \leftarrow \cF^\flat_{1/2} \rightarrow \cF_1) \in \Hk^1_{U(n), \LL}(T)$, the vector bundle  $\cF:=\cF^{\flat}_{1/2}$ on $T\times X'$ together with the Hermitian map $h$ induced from that of $\cF_{0}$ gives a $T$-point $(\cF,h)$ of $\Sect(X' \times X/X', \sH/X')$. It is clear that this construction gives an inverse to the construction in the previous paragraph, hence giving a canonical isomorphism $\Sect(X' \times X/X', \sH/X') \cong \Hk_{U(n), \LL}^1$.

Furthermore, under this identification the map projecting $(\cF,h)$ to $(\cF_0, t)$ intertwines with $(\pr_0 \times \pr_X)$, and the map projecting to the data of $(\cF_1, \sigma(t))$ intertwines with $(\pr_1 \times \pr_X)$. 

\end{proof}

More generally, let $\sG \cong \prod_{Y_{\alpha}} R_{Y_{\alpha}/X} BU(n_{\alpha})_{\LL_{\alpha}}$ be a gerbe of unitary type on $X$ in the notation of \S \ref{sssec: hecke for unitary gerbe}, with the standard $n$-framed structure, so $Y' = Y \times_X X' = \coprod Y_{\alpha}$, Then there is a \emph{unitary Hecke gerbe} $\sH_{\sG}$ over $Y' \times X$ defined analogously to $\sH$, which in the case $\sG = BU(n)_{\LL}$ recovers the $\sH$ defined in Definition \ref{defn: hecke gerbe unitary}.  Moreover, $\sH_{\sG}$ is equipped with two maps 
\begin{equation}\label{eq: general hecke gerbe unitary}
\begin{tikzcd}
& \sH_{\sG} \ar[dl, "h_0"'] \ar[dr, "h_1"] \\
\sG \times Y' & & \sG \times Y'
\end{tikzcd}
\end{equation}
of gerbes over $X \times Y'$, and has an isomorphism (by a similar argument as that for Lemma \ref{lem: hecke as sect})
\begin{equation}\label{eq: sect as hecke}
\Sect(Y' \times X/Y', \sH_{\sG}/Y') \cong \Hk_{\sG}^1
\end{equation}
with the RHS as defined in \S \ref{sssec: hecke for unitary gerbe}, intertwining 
\begin{equation}\label{eq: general gerbe hecke map classical}
\begin{tikzcd}
& \Sect(Y' \times X/Y', \sH_{\sG}/Y') \ar[dl, "h_0"'] \ar[dr, "h_1"]  \\
\Sect(Y' \times X/Y', \sG \times Y'/Y') & & \Sect(X' \times X/Y', \sG \times Y'/Y')
\end{tikzcd}
\end{equation}
with 
\[
\begin{tikzcd}
& \Hk_{\sG}^1 \ar[dl, "\pr_0 \times \pr_{Y'}"']  \ar[dr, "\pr_1  \times \pr_{Y'}"] \\
 \Bun_{\sG}\times Y' & & \Bun_{\sG} \times Y'
\end{tikzcd}
\]
 We have not spelled out the full definitions because the notation would be so heavy that it would be unenlightening. 

\begin{defn}[Derived Hecke stacks for gerbes of unitary type]\label{def: derived Hk_U}
Let $\sG \cong \prod_{Y_{\alpha}} R_{Y_{\alpha}/X} BU(n_{\alpha})_{\LL_{\alpha}}$ be a gerbe of unitary type, with the standard $n$-framed structure. We define the derived Hecke stack
\[
\sHk_{\sG}^1 := \sSect(Y' \times X/Y', \sH_{\sG}/Y')
\]
It is equipped with two maps induced by \eqref{eq: general hecke gerbe unitary}
\[
\begin{tikzcd}
& \sHk_{\sG}^1 \ar[dl, "\pr_0 \times \pr_{Y'}"'] \ar[dr, "\pr_1 \times \pr_{Y'}"] \\
\Bun_{\sG} \times Y' & & \Bun_{\sG}  \times Y'
\end{tikzcd}
\]


We define $\sHk_{\sG}^r$ to be  the $r$-fold (derived) fibered product 
\begin{equation}\label{eq: sHk^r derived fibered product 1}
\sHk_{\sG}^r := \sHk_{\sG}^1 \times_{\Bun_{\sG}} \sHk_{\sG}^1  \times_{\Bun_{\sG}}  \dots \times_{\Bun_{\sG}} \sHk_{\sG}^1
\end{equation}
where on the $i^{\mrm{th}}$ factor of $\sHk_{\sG}^1$, parametrizing $\cF_{i-1} \leftarrow \cF_{i-1/2}^{\flat} \rightarrow \cF_i$, the left and right maps to $\Bun_{\sG}$ project to $\cF_{i-1}$ and $\cF_i$ respectively. A point of $\sHk_{\sG}^r$ will be denoted
\[
\begin{tikzcd}
 &   \cF_{1/2}^{\flat} \ar[dl] \ar[dr]   & & \ldots  & & \cF_{r-1/2}^{\flat} \ar[dl] \ar[dr] \\ 
\cF_0 \ar[rr, dashed, "f_0"]  &  & \cF_{1} \ar[r, dashed] & \ldots \ar[r, dashed] & \cF_{r-1} \ar[rr, dashed, "f_{r-1}"] & & \cF_r 
\end{tikzcd}   
\]
and the projection map to $\cF_i$ denoted $\pr_i \co \sHk_{\sG}^r \rightarrow \Bun_{\sG}$. 
\end{defn}

\begin{remark}
By Lemma \ref{lem: hecke as sect} (and its generalization to arbitrary gerbes of unitary type), the classical truncation of $\sHk_{\sG}^r$ is $\Hk_{\sG}^r$ as defined in \S \ref{sssec: hecke for unitary gerbe}. 
\end{remark}

Next we make the analogous construction for $B\GL(n)'$.  In the definition below, we will use the notations introduced in Definition \ref{defn: hecke gerbe unitary}.

\begin{defn}\label{defn: hecke gerbe GL(n)'}
We define a \emph{(linear) Hecke gerbe} $\sH_{\GL(n)'} \rightarrow X' \times X$ that sends a $k$-scheme $S$ to the groupoid of $(x', \xi, \cF,L)$ where 
\begin{itemize}
\item $x' \in X'(S)$. 
\item $\xi \in X(S)$. 
\item $\cF$ is a vector bundle of rank $n$ on $S'=S\times_{\xi, X}X'$.
\item A line sub-bundle $L$ of the restriction of $\cF$ to $S'_{\xi'=x'}\coprod S'_{\xi'=\s(x')}$.
\end{itemize}
Note that upper modifying $\cF$ along either $L|_{S'_{\xi'=x'}}$ or $L|_{S'_{\xi'=\s(x')}}$ produces two maps $\sH_{\GL(n)'}$ to $B\GL(n)' \times X'$ as gerbes over $X \times X'$, 
\begin{equation}\label{eq: hecke gerbe maps GL(n)'}
\begin{tikzcd}
& \sH_{\GL(n)'} \ar[dl, "h_0"'] \ar[dr, "h_1"] \\
B\GL(n)' \times X' & & B\GL(n)' \times X' 
\end{tikzcd}
\end{equation}
\end{defn}
We denote by 
\begin{equation}\label{eq: linear hecke gerbe}
\begin{tikzcd}
\cF_0  & \cF \ar[l, hook'] \ar[r, hook] & \cF_1
\end{tikzcd}
\end{equation}
these two modifications so that $\cF_{0}/\cF$ (resp. $\cF_{1}/\cF$) is the line bundle $L|_{S'_{\xi'=x'}}\ot \cO_{S'}(S'_{\xi'=x'})|_{S'_{\xi'=x'}}$ (resp. $L|_{S'_{\xi'=\s(x')}}\ot \cO_{S'}(S'_{\xi'=\s(x')})|_{S'_{\xi'=\s(x')}}$). Given $(x',\xi', \cF)$, the the datum of $L$ is equivalent to the datum of the diagram \eqref{eq: linear hecke gerbe}. 

Consider the diagram 
\[
\begin{tikzcd}
& \Sect(X \times X'/X', \sH_{\GL(n)'}/X') \ar[dl, "h_0"'] \ar[dr, "h_1"] \\
\Sect(X \times X'/X', B\GL(n)' \times X'/X') & & \Sect(X \times X'/X', B\GL(n)' \times X' /X')
\end{tikzcd}
\]
By a similar argument as for Lemma \ref{lem: hecke as sect}, it is identified with 
\[
\begin{tikzcd}
& \Hk_{\GL(n)'}^1  \ar[dl, "\pr_0 \times \pr_{X'}"']  \ar[dr, "\pr_1  \times \pr_{X'}"] \\
\Bun_{\GL(n)'} \times X' & & \Bun_{\GL(n)'} \times X'
\end{tikzcd}
\]

\begin{defn}[Derived Hecke stacks for $\GL(n)'$] 
We define the derived Hecke stack
\[
\sHk_{\GL(n)'}^1 := \sSect(X' \times X/X', \sH_{\GL(n)'}/X')
\]
It is equipped with two maps induced by \eqref{eq: general hecke gerbe unitary}
\[
\begin{tikzcd}
& \sHk_{\GL(n)'}^1  \ar[dl, "\pr_0 \times \pr_{X'}"'] \ar[dr, "\pr_1 \times \pr_{X'}"] \\
\Bun_{\GL(n)'} \times X' & & \Bun_{\GL(n)'}  \times X'
\end{tikzcd}
\]

At the level of zero-cells, it sends $x' \in X'(R_{\bu})$ to the anima of $\cF_0, \cF_1 \in \Bun_{\GL(n)'}(R_{\bu})$, and a rank $n$ vector bundle $\cF_{1/2}^{\flat}$ on $X'_{R_{\bu}}$ plus a diagram 
\[
\begin{tikzcd}
 &   \cF_{1/2}^{\flat} \ar[dl, hook, "h^{\leftarrow}"'] \ar[dr, hook, "h^{\rightarrow}"]  \\ 
 \cF_0  \ar[rr, dashed, "f_0"]  &  & \cF_{1}
\end{tikzcd}   
\]
such that $\mrm{cone}(h^{\leftarrow})$ (resp. $\mrm{cone}(h^{\rightarrow})$) is supported on $\Gamma_{x'}$ (resp. $\Gamma_{\s(x')}$) and locally isomorphic to $R_{\bu}$ (as $R_{\bu}$-modules). 


We define $\sHk_{\GL(n)'}^r$ to be  the $r$-fold
 (derived) fibered product 
\begin{equation}\label{eq: sHk^r derived fibered product 2}
\sHk_{\GL(n)'}^r := \sHk_{\GL(n)'}^1 \times_{\Bun_{\GL(n)'}} \sHk_{\GL(n)'}^1  \times_{\Bun_{\GL(n)'}}  \dots \times_{\Bun_{\GL(n)'}} \sHk_{\GL(n)'}^1
\end{equation}
where on the $i^{\mrm{th}}$ factor of $\sHk_{\GL(n)'}^1$, parametrizing $\Vect(\cF_{i-1}) \leftarrow \cF_{i-1/2}^{\flat} \rightarrow \Vect(\cF_i)$, the left and right maps to $\Bun_{\GL(n)'}$ project to $\cF_{i-1}$ and $\cF_i$ respectively. 
\end{defn}

In order to make the notation more uniform, we will denote by $\sG$ either a gerbe of unitary type over $X$, or $\sG=B\GL(n)'$, in which case $\sHk_{\sG}^{r}=\sHk_{\GL(n)'}^{r}$. 
We have seen that the classical truncation of $\sHk_{\sG}^r$ is the $\Hk_{\sG}^r$ from \S \ref{sssec: hecke for unitary gerbe} (for $\sG$ of unitary type) or Definition \ref{defn: Hk GL(n)'} (for $\sG =B \GL(n)'$). We prove below that the canonical map $\Hk_{\sG}^r \rightarrow \sHk_{\sG}^r$ is an isomorphism in both cases. 

\begin{lemma}\label{lem: shtuka derived cartesian} Let $\sG$ be either a gerbe of unitary type or $B\GL(n)'$ over $X$. Then: 
\begin{enumerate}
\item $\iota \co \Hk^r_{\sG} \rightarrow \sHk^r_{\sG}$ is an isomorphism of derived stacks. 
\item The following diagram of classical stacks is derived Cartesian
\[
\begin{tikzcd}
\Sht_{\sG}^r  \ar[r] \ar[d] & \Hk_{\sG}^r \ar[d, "{(\pr_0, \pr_r)}"] \\
\Bun_{\sG} \ar[r, "{(\Id, \Frob)}"] & \Bun_{\sG} \times \Bun_{\sG}
\end{tikzcd}
\]
\end{enumerate}
\end{lemma}

\begin{proof}

(1) By the same argument as in \cite[Lemma 6.9]{FYZ} (for $\sG$ of unitary type) or Lemma \ref{lem: sht dimension} (for $\sG=B\GL(n)'$), the map $\sHk^1_{\sG} \xrightarrow{(\pr_0, \pr_{X'})} \sMap(X, \sG) \times X'$ is a projective space bundle, and $\sMap(X, \sG) = \Bun_{\sG}$ is classical by Corollary \ref{cor: T for Bun_G}, so $\sHk^1_{\sG}$ is also classical. Consider the derived Cartesian square in the diagram 
\[
\begin{tikzcd}
\sHk^r_{\sG} \ar[r, "\Pi^{r}"] \ar[d] 
& \sHk^{r-1}_{\sG} \ar[d, "\pr_{r-1}"]  \\ 
\sHk^1_{\sG} \ar[r, "\pr_0"] & \Bun_{\sG} 
\end{tikzcd}
\]
where $\Pi^{r}$ is the projection to the first $r-1$ factors of $\sHk^1_{\sG}$. This diagram implies that $(\Pi^{r}, \pr_{X',r}):\sHk_{\sG}^r \rightarrow \sHk^{r-1}_{\sG} \times X'$ is a projective space bundle, so by induction $\sHk_{\sG}^r$ is classical for each $r$.

(2) Since $\Bun_{\sG}$ is smooth, it suffices to show that the maps $(\pr_0, \pr_r)$ and ${(\Id, \Frob)}$ are transversal. The differential of $\Frob$ is zero, so this follows from the smoothness of $\pr_r$ (see the argument of \cite[Lemma 6.9(1)]{FYZ} for $\sG$ of unitary type, or Lemma \ref{lem: sht dimension} for $\sG = B\GL(n)'$). 
\end{proof}

Thanks to Lemma \ref{lem: shtuka derived cartesian}, we may and do write $\Hk_{\sG}^r$ instead of $\sHk_{\sG}^r$ in the sequel.


\subsection{Derived Hitchin stacks}
We now define derived versions of the Hitchin stacks introduced in \S \ref{sec: hitchin spaces}. Unlike $\Bun_{\sG}$ and $\Hk_{\sG}^r$, these will be genuinely non-classical in general. We will adopt the conventions in \S\ref{ssec: Hitchin stacks} to denote a gerbe over $X$ by $BH$, even when the group scheme $H$ is not defined, and call a section of $BH$ over $S\to X$ an ``$H$-torsor over $S$''.

\begin{defn}
\begin{enumerate}
\item The \emph{standard representation} of the gerbe $B\GL(m)' \times_X B\GL(n)' \rightarrow X$ is the stack $\Std'(m,n)$ over $X$ assigning to $s \co S \rightarrow X$ the groupoid of $(\cE, \cF, t)$ where 
\begin{itemize}
\item $\cE$ is a vector bundle of rank $m$ on $S \times_X X'$, 
\item $\cF$ is a vector bundle of rank $n$ on $S\times_X X'$, 
\item $t  \in \Hom(\cE, \cF)$.
\end{itemize}
It is equipped with a forgetful map
\begin{equation*}
\om_{m,n}: \Std'(m,n)\to B\GL(m)' \times_X B\GL(n)'. 
\end{equation*}

\item If $BH_1$ is an $m$-framed gerbe and $BH_2$ is an $n$-framed gerbe, the \emph{standard representation} of $BH_1 \times_X BH_2$, denoted $\Std'(BH_{1},BH_{2})$, is defined using the Cartesian diagram
\begin{equation*}
\xymatrix{\Std'(BH_{1},BH_{2}) \ar[d]\ar[r] & \Std'(m,n)\ar[d]^{\om_{m,n}}\\
BH_1 \times_X BH_2\ar[r] & B\GL(m)'\times_{X}B\GL(n)'}
\end{equation*}
Here the bottom arrow is given by the framings of $BH_{1}$ and $BH_{2}$.
\end{enumerate}
\end{defn}

\begin{example}
The relative tangent complex of $\Std'(m,n)\rightarrow X$ at 
\[
(\cE, \cF, t) \in \Std'(m,n)(\Spec R)
\]
is the complex of $R$-modules,
\[
\underbrace{\cEnd(\cE) \oplus \cEnd(\cF)}_{\deg -1} \xrightarrow{\alpha_t} \underbrace{\cHom(\cE, \cF}_{\deg 0})
\]
where the map $\alpha_t$ is given by the action on $t \in \Hom(\cE, \cF)$, 
\[
(A \in \cEnd(\cE), B \in \cEnd(\cF)) \mapsto Bt - tA.
\]
More generally, the tangent complex of $\Std'(BH_{1}, BH_2)$ is of a similar form.
\end{example}


\begin{defn}\label{defn: derived hitchin} 
Let $BH_1$ be a smooth $m$-framed gerbe over $X$ and $BH_2$ be a smooth $n$-framed gerbe over $X$. We define the \emph{derived Hitchin stack} $\sM_{H_1,H_2} $ to be the derived mapping stack $\sSect(X, \Std'(BH_{1},BH_{2}))$. At the level of $0$-cells, it sends an animated ring $R_{\bu}$ to the anima of $(\cE, \cF, t)$ where 
\begin{itemize}
\item $\Cal{E}\in \sSect(X, BH_1)(R_{\bu})$,
\item $\Cal{F} \in \sSect(X, BH_2)(R_{\bu}) $, and 
\item $t \in R\Gamma(X'_{R_{\bu}}, \cHom(\Vect(\Cal{E}),  \Vect(\Cal{F})))$.
\end{itemize}
The space of such triple $(\cE, \cF, t)$ is an anima as discussed in Example \ref{ex: derived standard/gerbe}.
As the map $\sM_{H_1, H_2} \rightarrow \sSect(X, BH_1 \times BH_2)$ is evidently representable, this is a derived Artin stack. It is immediate from the definition that the classical truncation of $\sM_{H_1, H_2} $ is $\cM_{H_1, H_2}$. 

We define $\sM_{H_1, H_2}^{\c} \subset \sM_{H_1, H_2}$ to be the open derived substack whose classical truncation is $\cM_{H_1, H_2}^{\c} \subset \cM_{H_1, H_2}$ as in Definition \ref{defn: hitchin} (cf. \cite[\S 2.2.2]{TV08} for the notion of Zariski open immersion in derived algebraic geometry). 
\end{defn}

\begin{remark}
The definition makes sense even when $BH_1$ and $BH_2$ are not smooth, but we never consider such an example. The smoothness assumption guarantees that the mapping stacks $\sMap(X, BH_1)$ and $\sMap(X, BH_2)$ are smooth and classical, and this is necessary for much of what we say below. 
\end{remark}

\begin{remark}
In \S \ref{ssec: derived vector bundles} we will give an alternative description of $\sM_{H_1, H_2}$, which might be more enlightening. 
\end{remark}




From the general description of tangent complexes for derived mapping stacks in \S \ref{ss:tangent}, we will deduce the following. 

\begin{cor}\label{cor: Hitchin to Bun tangent} Let $BH_{1}$ and $BH_{2}$ be as in Definition \ref{defn: derived hitchin}. For any animated ring $R_{\bu}$, the tangent complex of the morphism $\scr{M}_{H_1, H_2}\xrightarrow{\pi} \Bun_{H_1} \times \Bun_{H_2}$ at $(\cE, \cF, t \in R\Gamma(X'_{R_{\bu}}, \cHom(\Vect(\cE) , \Vect(\cF)))) \in  \scr{M}_{H_1, H_2} (R_{\bu})$ is naturally in $R_{\bu}$ isomorphic to 
\[
R\pr_*(\cHom_{X_{R_{\bu}}'}(\Vect(\cE) , \Vect(\cF))).
\]
where $\pr \co X_{R_{\bu}}' \rightarrow \Spec R_{\bu}$ is the projection map. In particular, $\pi$ is quasi-smooth, hence $\sM_{H_1, H_2}$ is quasi-smooth. 
\end{cor}
\begin{proof} This is a special case of Corollary \ref{cor: TC Sect to Bun_G} applied with $\sY=\Std'(BH_{1}, BH_2)$, $\sY' = BH_1 \times_{X} BH_2$ and $\mu: \sY \rightarrow \sY'$ being the forgetful map.
\end{proof}

\begin{cor}\label{cor: M already classical} 
\begin{enumerate}
\item Let $BH_1$ be a smooth $m$-framed gerbe over $X$. Then the classical truncation map $\cM_{H_1, \GL(n)'}^{\c} \rightarrow \sM_{H_1, \GL(n)'}^{\c}$ is an isomorphism, and both stacks are smooth.
\item Let $BH_{1}=B\GL(m)'$ be the tautological $m$-framed gerbe over $X$. Let $\cM^{\ns}_{H_{1},U(n),\LL}\subset \cM_{H_{1},U(n),\LL}$ be the preimage of $\cA^{\ns}_{H_{1},\LL}$ under the Hitchin fibration (see \S \ref{ssec: hitchin base}), and $\sM^{\ns}_{H_{1},U(n),\LL}\subset \sM_{H_{1},U(n),\LL}$ be the corresponding open derived substack.  Then the classical truncation map $\cM_{H_{1}, U(n),\LL}^{\ns} \rightarrow \sM_{H_{1}, U(n),\LL}^{\ns}$ is an isomorphism, and both stacks are smooth.
\end{enumerate}
\end{cor}
\begin{proof}

In both cases we will verify the classicality using Lemma \ref{lem: classicality of smooth derived stacks}. By Corollary \ref{cor: Hitchin to Bun tangent}, the tangent complex in each case is perfect and has cohomology concentrated in degrees $\leq 1$. It remains to show that the tangent complex in each case is concentrated in degrees $\leq 0$, i.e., that the respective first cohomology groups vanish. These first cohomology groups are the obstruction groups \eqref{eq: obs group}, and in Proposition \ref{prop: Hitchin smooth} we calculated that they vanish. 
\end{proof}

\begin{remark}
It is important that we restrict to the injective locus for Corollary \ref{cor: M already classical}(1). The statement would not be true for $\sM_{H_1, \GL(n)'}$ in place of $\sM_{H_1, \GL(n)'}^{\c}$. Furthermore, in part (2), 
we would not have been able to make the same argument with $\sM_{H_{1},U(n),\LL}^{\c}$ in place of $\sM_{H_{1},U(n),\LL}^{\ns}$. 
\end{remark}

\subsection{Derived Hecke stacks for derived Hitchin stacks}
We shall define derived Hecke stacks $\Hk_{\sM_{H_1, H_2}}^r$ and $\Hk_{\sM_{H_1, H_2}^{\c}}^r$ whose classical truncation recovers the classical stacks $\Hk_{\cM_{H_1, H_2}}^r$ and 
$\Hk_{\cM_{H_1, H_2}^{\c}}^r$.

Recall the gerbe $\sH_{\GL(n)'}$ from Definition \ref{defn: hecke gerbe GL(n)'}.
\begin{defn}
\begin{enumerate}
\item The \emph{standard representation} of the gerbe $B\GL(m)' \times_X \sH_{\GL(n)'} \rightarrow X \times X'$ is the stack $\Std'(B\GL(m)',\sH_{\GL(n)'})$ over $X$ assigning to a $k$-scheme $S$ the groupoid of $(\cE, x', \xi, \cF, L, t)$ where 
\begin{itemize}
\item $(x', \xi, \cF, L) \in \sH_{\GL(n)'}(S)$. In particular, $\xi: S\to X$.
\item $\cE$ is a vector bundle of rank $m$ on $S'=S\times_{\xi, X}X'$.
\item $t \in \Hom(\cE, \cF)$.
\end{itemize}

\item More generally, if $BH_1$ is an $m$-framed gerbe, the \emph{standard representation} of $BH_1 \times_X \sH_{\GL(n)'}$, denoted $\Std'(BH_{1},\sH_{\GL(n)'})$, is defined by the Cartesian diagram
\begin{equation*}
\xymatrix{\Std'(BH_{1},\sH_{\GL(n)'})\ar[d]\ar[r] & \Std'(m,\sH_{\GL(n)'})\ar[d]^{\om_{m}}\\
BH_1 \ar[r]& B\GL(m)' }
\end{equation*}
where $\om_{m}$ records the rank $m$ bundle and the bottom map is the $m$-framing of $BH_{1}$.

\item If $BH_1$ is an $m$-framed gerbe and $BH_2$ is a smooth $n$-framed gerbe of unitary type, then we define the \emph{standard representation} of $BH_1 \times_X \sH_{H_2}$, denoted $\Std'(BH_{1},\sH_{H_2})$, by the Cartesian diagram
\[
\begin{tikzcd}
\Std'(BH_{1},\sH_{H_2})\ar[d]\ar[r] & \Std'(BH_1,\sH_{\GL(n)'})\ar[d] \\
\sH_{H_2} \ar[r]& \sH_{\GL(n)'}
\end{tikzcd}
\]
\end{enumerate}
\end{defn}

\begin{defn}\label{def: derived Hk} Let $BH_{1}$ be a smooth $m$-framed gerbe and $BH_2$ a smooth $n$-framed gerbe of unitary type or $\GL(n)'$. We define the derived stack $\Hk^1_{\sM_{H_1, H_2}}$ to be the derived mapping stack 
\[
\sSect(X \times X'/X', \Std'(BH_1,\sH_{H_2} )).
\]
We define $\Hk^1_{ \sM_{H_1, H_2}^{\c}}$ to be the pre-image of the open substack $ \sM_{H_1, H_2}^{\c} \inj  \sM_{H_1, H_2}$ under either projection. 


At the level of $0$-cells, $\Hk_{\sM_{H_1, H_2}}^1$ sends an animated ring $R_{\bu}$ to the anima of an $H_{1}$-torsor $\cE$ over $X'_{R_{\bu}}$, a diagram
\[
\begin{tikzcd}
 &   \cF_{1/2}^{\flat} \ar[dl, hook] \ar[dr, hook]  \\ 
 \Vect(\cF_0) \ar[rr, dashed, "f_0"] & & \Vect(\cF_{1})
\end{tikzcd}   \in \Hk_{H_2}^1(R_{\bu})
\]
and  $t_{1/2}^\flat \in  R\Gamma(X_{R_{\bu}}, \cHom(\Vect(\cE), \cF_{1/2}^{\flat}))$, and $\Hk_{\sM_{H_1, H_2}^\circ}^1$ is the open substack where $t_{1/2}^\flat$ is injective. For $\sM = \sM_{H_1, H_2}$ or $\sM_{H_1, H_2}^{\circ}$, we denote such an $R_{\bu}$-point of $\Hk^1_{\sM}$ by the diagram 
\[
\begin{tikzcd}
& \Vect(\cE) \ar[d, "t_{1/2}^{\flat}"]  \ar[ddl, "t_0"', bend right] \ar[ddr, "t_{1}", bend left] \\
 &   \cF_{1/2}^{\flat} \ar[dl, hook] \ar[dr, hook]  \\ 
\Vect( \cF_0)  \ar[rr, dashed, "f_0"] & & \Vect(\cF_{1})
\end{tikzcd}
\]

Then we define $\Hk_{\sM}^r$ as the $r$-fold derived fibered product 
\begin{equation}\label{eq: sHk_M derived fibered produt}
\Hk_{\sM}^r := \underbrace{\Hk_{\sM}^1 \times_{\sM}\Hk_{\sM}^1 \times_{\sM} \ldots \times_{\sM} \Hk_{\sM}^1}_{\text{$r$ times}}
\end{equation}
where the maps are like in Definition \ref{def: derived Hk_U}. 
\end{defn}

Again it is clear that the classical truncation of $\Hk_{\sM}^r$ is $\Hk_{\cM}^r$ where $\cM$ is the classical truncation of $\sM$, i.e. $\cM = \cM_{H_1, H_2}$ or $\cM_{H_1, H_2}^{\c}$.

\begin{lemma}\label{lem: TC Hk_M to Hk}Let $BH_{1}$ be a smooth $m$-framed gerbe over $X$ and $BH_{2}$ a gerbe of unitary type or $B \GL(n)'$, equipped with the standard $n$-framed structure. Then the morphism $\Hk_{\sM_{H_1,H_2}}^r \xrightarrow{\pi} \Bun_{H_{1}}\times\Hk^r_{H_2}$ is quasi-smooth. For any animated ring $R_{\bu}$, the tangent complex of $\pi$ at any $(\cE, \{x'_{i}\}, \{\cF_i\}, \{t_{i}\}) \in \Hk_{\sM_{H_1, H_2}}^r (R_{\bu})$ is naturally in $R_{\bu}$ isomorphic to (using notation from Definition \ref{def: derived Hk})
\[
R\pr_* \left(\underbrace{ \bigoplus_{i=1}^{r} \cHom(\Vect(\cE), \cF_{i-1/2}^{\flat} ) }_{\deg 0}\rightarrow  \underbrace{ \bigoplus_{i=1}^{r-1}  \cHom(\Vect(\cE), \Vect(\cF_{i})) }_{\deg 1} \right) 
\]
where $\pr \co X_{R_{\bu}}' \rightarrow \Spec R_{\bu}$ is the projection, and the differential is induced by the map 
\begin{equation}\label{eq: Hk cokernel bundle}
\bigoplus_{i=1}^{r} \cF_{i-1/2}^{\flat} \rightarrow  \bigoplus_{i=1}^{r-1}\Vect(\cF_{i}) 
\end{equation}
sending $(v_{1/2}, \ldots, v_{r-1/2})  \mapsto (v_{1/2}-v_{3/2}, v_{3/2}-v_{5/2}, \ldots, v_{r-3/2} - v_{r-1/2})$, using the given maps $\cF_{i-1/2}^{\flat} \inj \Vect(\cF_{i-1})$ and $\cF_{i-1/2}^{\flat} \inj \Vect(\cF_{i})$. In particular, $\Hk_{\sM_{H_1, H_2}}^r$ is quasi-smooth. 
\end{lemma}

\begin{proof}
The case $r=1$ is a direct consequence of the definition and Lemma \ref{lem: tangent to derived mapping stack}. Abbreviate $\sM = \sM_{H_1, H_2}$. Consider an $R_{\bu}$-point of $\Hk^r_{\sM}$, which is represented by a diagram
\begin{equation}
\begin{tikzcd}
 & & &  \Vect(\cE)   \ar[dll, "t_0"] \ar[drr, "t_{r-1}"'] \\
& \cF_{1/2}^{\flat} \ar[dl, hook] \ar[dr, hook]  & & \ldots   & &   \cF_{r-1/2}^{\flat}\ar[dl, hook] \ar[dr, hook]  \\ 
 \Vect(\cF_0) & &  \Vect(\cF_1)  & \ldots &  \Vect(\cF_{r-1}) & &  \Vect(\cF_r )
\end{tikzcd}
\end{equation}
The presentation of $\Hk_{\sM}^r $ from \eqref{eq: sHk_M derived fibered produt} induces an exact triangle
\[
\begin{tikzcd}
\bT_{\Hk_{\sM}^r } \ar[r] & \bT_{\Hk_{\sM}^1}^{\oplus r}|_{\Hk_{\sM}^r}\ar[r] & \bT_{\sM}^{\oplus (r-1)}|_{\Hk_{\sM}^r}
\end{tikzcd}
\]
Similarly, the presentations \eqref{eq: sHk^r derived fibered product 1},\eqref{eq: sHk^r derived fibered product 2} induce an exact triangle 
\[
\begin{tikzcd}
\bT_{\Bun_{H_{1}}\times \Hk_{H_2}^r } \ar[r] & \bT_{\Bun_{H_{1}}\times \Hk_{H_2}^1}^{\oplus r}|_{\Bun_{H_{1}}\times \Hk_{H_2}^r } \ar[r] & \bT_{\Bun_{H_{1}}\times \Bun_{H_2}}^{\oplus (r-1)}|_{\Bun_{H_{1}}\times \Hk_{H_2}^r }
\end{tikzcd}
\]
So the (derived) fiber of $\bT_{\Hk_{\sM}^r } \rightarrow \bT_{\Hk_{H_2}^r }|_{\Hk_{\sM}^r}$ is the limit of the diagram 
\[
\begin{tikzcd}
 \bT_{\Hk_{\sM}^1}^{\oplus r} |_{\Hk_{\sM}^r} \ar[r] \ar[d]  & \bT_{\sM}^{\oplus (r-1)} |_{\Hk_{\sM}^r} \ar[d]  \\ 
  \bT_{\Bun_{H_{1}}\times \Hk_{H_2}^1}^{\oplus r} |_{\Hk_{\sM}^r} \ar[r]& \bT_{\Bun_{H_{1}}\times\Bun_{H_2}}^{\oplus (r-1)} |_{\Hk_{\sM}^r}
\end{tikzcd}
\]
This limit may be calculated by first forming vertical fibers along the columns, and then taking the horizontal fiber. The vertical fibers restricted to the given $R_{\bu}$-point are described by Corollary \ref{cor: Hitchin to Bun tangent}: 
\[
\begin{tikzcd}
R\pr_* \left(\bigoplus_{i=1}^{r} \cHom(\Vect(\cE), \cF_{i-1/2}^{\flat} ) \right)  \ar[r] & 
R\pr_*\left(\bigoplus_{i=1}^{r-1}  \cHom(\Vect(\cE), \Vect(\cF_{i})) \right)
\end{tikzcd}
\]
In particular, as the map \eqref{eq: Hk cokernel bundle} is supported on a subscheme of $X'_{R_{\bu}}$ which is finite over $\Spec R_{\bu}$, the fiber of this map is concentrated in degrees $\leq 1$. 

\end{proof}

\begin{cor}\label{cor: HkM already classical} 
\begin{enumerate}
\item Let $BH_1$ be a smooth $m$-framed gerbe over $X$. Then the classical truncation map $\Hk_{\cM_{H_1, \GL(n)'}^{\c}}^1 \rightarrow \Hk_{\sM_{H_1, \GL(n)'}^{\c}}^1 $ is an isomorphism, and both stacks are smooth.

\item Let $BH_{1}=B\GL(m)'$ be the tautological $m$-framed gerbe over $X$. Let $\Hk^{1}_{\cM_{H_1, U(n),\LL}^{\ns}}\subset \Hk^{1}_{\cM_{H_{1},U(n),\LL}}$ be the preimage of $\cA^{\ns}_{H_{1},\LL}$ (under the Hitchin fibration composed with $\pr_1$)	, and $\Hk^{1}_{\sM_{H_1, U(n),\LL}^{\ns}}\subset \Hk^{1}_{\sM_{H_{1},U(n),\LL}}$ be the corresponding open derived substack.  Then the classical truncation map $\Hk_{\cM_{H_1, U(n),\LL}^{\ns}}^1 \rightarrow \Hk_{\sM_{H_1, U(n), \LL}^{\ns}}^1 $ is an isomorphism, and both stacks are smooth.
\end{enumerate}
\end{cor}
\begin{proof}

(1) There is a map $p_{1/2} \co \Hk_{\sM_{H_1, \GL(n)'}}^1 \rightarrow \sM_{H_1, H_2} \times X'$ defined as in Lemma \ref{l:HkM sm}, which is a $\PP^{n-1} \times \PP^{n-1}$ fiber bundle. Since $\sM_{H_1, H_2}^{\circ}$ is smooth and classical by Corollary \ref{cor: M already classical}, $\Hk_{\sM_{H_1, \GL(n)'}^{\c}}^1$ is smooth and classical.  

(2) We will verify the classicality using Lemma \ref{lem: classicality of smooth derived stacks}. By Lemma \ref{lem: TC Hk_M to Hk}, the tangent complex in each case is perfect and has cohomology concentrated in degrees $\leq 1$. It remains to show that the tangent complex in each case is concentrated in degrees $\leq 0$, i.e., that the respective first cohomology groups vanish. By Lemma \ref{lem: tangent to relative mapping stack}, the relative tangent complex of $\Hk_{\sM_{H_1, U(n),\LL}}^1 \rightarrow X'$ at $(\cE, \cF^\flat, h,t) \in \Hk_{\sM_{H_1, U(n),\LL}}^1$ 
\[
R\Gamma(X_{R_{\bu}}, \underbrace{\cEnd(\Vect(\cE)) \oplus \cEnd^{\mrm{asa}}(\cF^\flat,h)}_{\deg -1} \xrightarrow{\a_{t}} \underbrace{\cHom(\Vect(\cE), \cF^{\flat})}_{\deg 0} )
\]
The obstruction group is \eqref{eq: obs group Hk}, and we showed there that it vanishes if $(\cE, \cF^\flat, h,t)  \in \Hk_{\sM_{H_1, U(n),\LL}^{\ns}}^1(R_{\bu})$. 
\end{proof}


\subsection{Description as derived vector bundles}\label{ssec: derived vector bundles}
There is a derived symmetric algebra functor \cite[\S 1]{KhanII} from connective perfect complexes on a base (derived) scheme $S$ to animated commutative $\cO_S$-algebras, denoted $\cK \mapsto \Sym(\cK)$. Given $\cK \in \Perf(S)$ which is \emph{co-connective}, meaning locally represented by a complex of vector bundles in non-negative degrees, $\cK^{*}$ is connective and we denote by $\Tot_S(\cK)$ the relative spectrum of $\Sym_{\cO_S}(\cK^*)$. We call $\Tot_S(\cK)$ the \emph{derived vector bundle associated to $\cK$}; if $\cK$ a locally free sheaf concentrated in degree $0$ then $\Tot_S(\cK)$ is the usual vector bundle associated to a locally free sheaf.

\begin{remark}\label{rem: dvb tangent complex} It is immediate upon unwinding the definitions that the relative tangent complex of $\Tot_S(\cK) \rightarrow S$ is the pullback of $\cK$ from $S$ to $\Tot_S(\cK)$. 
\end{remark}

Let $BH_1$ be a smooth $m$-framed gerbe and $BH_2$ be a smooth $n$-framed gerbe over $X$. We indicate a more concrete description of the derived structure on $\sM:=\sM_{H_1,H_2}$ in Definition \ref{defn: derived hitchin}. We have a tautological bundle $\cH$ over $\Bun_{H_1} \times \Bun_{H_2}\times X'$ whose restriction to $\{(\cE,\cF)\}\times X'$ is $\cHom(\Vect(\cE),\Vect(\cF))$. Let $p_{\Bun}: \Bun_{H_1}  \times \Bun_{H_2}\times X'\to  \Bun_{H_1}  \times \Bun_{H_2}$ be the projection. Let $\cK =Rp_{\Bun, *}\cH$, a co-connective perfect complex on $\Bun_{H_1} \times \Bun_{H_2}$. The pullback of $\cK$ to $\Spec R_{\bu} \rightarrow \Bun_{H_1} \times \Bun_{H_2}$ is $R\Gamma(X'_{R_{\bu}}, \cHom(\Vect(\cE), \Vect(\cF)))$.

The local structure of $\sM$ can now be made fairly concrete. Let $U \rightarrow \Bun_{H_1}  \times \Bun_{H_2}$ be an open substack on which $(Rp_{\Bun, *}\cH)|_{U}$ can be represented by a two-step perfect complex $\cK^{0}\xr{d}\cK^{1}$ over $U$ in degrees $0$ and $1$. Then $\sM|_{U}$ is isomorphic to the {\em derived} fiber product
\begin{equation}
\begin{tikzcd}
\sM|_{U}\ar[d]\ar[r] & \Tot_{U}(\cK^{0})\ar[d, "d"]\\
0_{U}\ar[r, hook] & \Tot_{U}(\cK^{1})
\end{tikzcd}
\end{equation}
where $0_{U}\cong U$ denotes the zero section in $\Tot_{U}(\cK^{1})$. 

\begin{prop}
The derived Hitchin stack $\sM_{H_1, H_2}$ is isomorphic to $\Tot_{\Bun_{H_1} \times \Bun_{H_2}} (\cK)$ over $\Bun_{H_1} \times \Bun_{H_2}$. 
\end{prop}

\begin{proof}

By construction, $\Std'(H_{1}, H_{2})$ is the total stack of a vector bundle $\cV_{\Std}$ over $BH_{1}\times_{X}BH_{2}$. We claim that for a proper scheme $X$ over $k$, and a vector bundle $\cV$ over a derived stack $S$ over $X$,  there is a canonical equivalence of derived stacks
\begin{equation}\label{eq: sect vb as total space}
\sSect(X, \Tot_{S}(\cV))\isom \Tot_{T}(p_{T*}e^{*}\cV)
\end{equation}
where $T=\sSect(X,S)$, $e: X\times T\to S$ is the evaluation map, and $p_{T}: X\times T\to T$ is the projection. The proposition follows by applying the above claim to $S=BH_{1}\times_{X}BH_{2}$ and $\cV=\cV_{\Std}$ (noting that $e^{*}\cV_{\Std}\cong \cH$ hence $p_{T*}e^{*}\cV_{\Std}\cong \cK$).

To prove \eqref{eq: sect vb as total space}, we construct maps in both directions. First we construct a map
\begin{equation*}
\a: \sM:=\sSect(X, \Tot_{S}(\cV))\to\Tot_{T}(p_{T*}e^{*}\cV).
\end{equation*}
There is a natural map $\om: \sM\to T$ induced from the projection $\pi: \Tot_{S}(\cV)\to S$. Giving such a map $\a$ covering $\om$ is the same as giving a global section of $\om^{*}p_{T*}e^{*}\cV$ over $\sM$. By base change we have $\om^{*}p_{T*}\cong p_{\sM*}(\Id_{X}\times \om)^{*}$, where $p_{\sM}: X\times\sM\to \sM$ is the projection. Therefore a global section of $\om^{*}p_{T*}e^{*}\cV$ is the same as a global section of $p_{\sM*}(\Id_{X}\times \om)^{*}e^{*}\cV$, which amounts to the same thing as a global section of $(\Id_{X}\times \om)^{*}e^{*}\cV$ on $X\times \sM$. By the commutative diagram (where $\wt e$ is the evaluation map)
\begin{equation*}
\xymatrix{X\times \sM\ar[r]^{\wt e}\ar[d]^{\Id_{X}\times \om} & \Tot_{S}(\cV)\ar[d]^{\pi}\\
X\times T\ar[r]^{e} & S
}
\end{equation*}
we have canonically
\begin{equation*}
(\Id_{X}\times \om)^{*}e^{*}\cV\cong \wt e^{*}\pi^{*}\cV.
\end{equation*}
Now $\pi^{*}\cV$ has a tautological section $\s_{\taut}$ over $\Tot_{S}(\cV)$. Then $\wt e^{*}\s_{\taut}$ gives a global section of $\wt e^{*}\pi^{*}\cV\cong (\Id_{X}\times \om)^{*}e^{*}\cV$ over $X\times \sM$, which then gives the desired map $\a$.

Next we give a map in the other direction
\begin{equation*}
\b: \Tot_{T}(p_{T*}e^{*}\cV)\to \sM=\sSect(X, \Tot_{S}(\cV)).
\end{equation*}
This is equivalent to giving a map
\begin{equation*}
X\times \Tot_{T}(p_{T*}e^{*}\cV)\to \Tot_{S}(\cV),
\end{equation*}
Note that  $X\times \Tot_{T}(p_{T*}e^{*}\cV)=\Tot_{X\times T}(p_{T}^{*}p_{T*}e^{*}\cV)$. The counit map $p_{T}^{*}p_{T*}\to \Id$ induces a map 
\begin{equation*}
\Tot_{X\times T}(p_{T}^{*}p_{T*}e^{*}\cV)\to \Tot_{X\times T}(e^{*}\cV).
\end{equation*}
Composing with the natural map $\Tot_{X\times T}(e^{*}\cV)\to \Tot_{S}(\cV)$ (base change of $e: X\times T\to S$), we get the desired map $\b$. We omit the details of checking that $\a$ and $\b$ are inverse to each other.

\end{proof}

Let $BH_1$ be a smooth $m$-framed gerbe and $BH_2$ be a gerbe of unitary type or $B\GL(n)'$ with the standard $n$-framing. Then there is a perfect complex $\cK^\flat$ on $\Bun_{H_1} \times \Hk_{H_2}^1$, whose pullback to $\Spec R_{\bu}$ via $(\cE, \cF^{\flat}) \in (\Bun_{H_1} \times \Hk_{H_2}^1)(R_{\bu})$ is $R\Gamma(X'_{R_{\bu}}, \cHom(\Vect(\cE),\Vect( \cF^\flat)))$. Similarly, one has a description of the derived Hecke stack for derived Hitchin spaces as a derived vector bundle: 

\begin{prop}
The derived Hitchin stack $\Hk_{\sM_{H_1, H_2}}^1$ is isomorphic to $\Tot_{\Bun_{H_1} \times \Hk_{H_2}^1} (\cK^\flat)$ over $\Bun_{H_1} \times \Hk_{H_2}^1$. 
\end{prop}

\begin{remark}
There is an analogous generalization to $\Hk_{\sM_{H_1, H_2}}^r$ for any $r \geq 0$, although we omit the statement because the corresponding complex $\cK^\flat$ becomes more complicated to describe. 
\end{remark}
\subsection{Derived shtukas for derived Hitchin spaces}\label{sssec: derived shtukas}
We now introduce derived stacks of shtukas $\Sht_{\sM_{H_1, H_2}}^r$ and $\Sht_{\sM_{H_1, H_2}^{\c}}^r$, whose classical truncation recovers the classical stacks $\Sht_{\cM_{H_1, H_2}}^r$ and $\Sht_{\cM_{H_1, H_2}^{\c}}^r$ from \S \ref{sec: hitchin spaces}. 

\begin{defn}\label{def: derived shtuka} Let $BH_{1}$ and $BH_{2}$ be as in Definition \ref{defn: derived hitchin}.   Let $\sM=\sM_{H_{1},H_{2}}$ or $\sM_{H_{1},H_{2}}^{\c}$. We define $\Sht_{\sM}^r$ by the (homotopy) Cartesian diagram 
\begin{equation}\label{eq: Sht_M cartesian square}
\begin{tikzcd}
\Sht_{\sM}^r \ar[r] \ar[d] & \Hk_{\sM}^r \ar[d,  "\pr_{0} \times \pr_{r}"]  \\
\sM \ar[r, "{(\Id, \Frob)}"] & \sM\times\sM
\end{tikzcd}
\end{equation}
\end{defn}

We are primarily interested in the case where $BH_2 = BU(n)_{\LL}$ (although we will make some remarks on the unitary type case below). In order to study the tangent complex of $\Sht_{\sM}^r$ in this case, we introduce a vector bundle on $\Hk^{r}_{U(n), \LL}$.

\subsubsection{Excess bundle}\label{sssec: excess bundle} 
Let $\cE$ be a vector bundle on $X'$. We denote by $\cV^r_{\cE}$ the rank $mr$ tautological vector bundle over $\Hk^r_{U(n), \LL}$ whose fiber at $(\{x'_{j}\}, \{\cF_{j}, h_{j}\}) \in \Hk^r_{U(n), \LL}(R)$ is the direct image under $X_R \rightarrow \Spec R$ of the cokernel of the map
\begin{equation}\label{eq: excess bundle map}
\bigoplus_{i=1}^{r}  \cHom(\cE, \cF^{\flat}_{i-1/2}) \to \bigoplus_{i=1}^{r}  \cHom(\cE, \cF_{i})
\end{equation}
given by
\begin{equation*}
(v_{1/2},v_{3/2}, \cdots, v_{r-1/2})\mapsto (v_{1/2}-v_{3/2}, v_{3/2}-v_{5/2},\cdots, v_{r-3/2}-v_{r-1/2}, v_{r-1/2}).
\end{equation*}
Here we use the natural inclusions $\cF^{\flat}_{i-1/2}\incl \cF_{i-1}$ and $\cF^{\flat}_{i-1/2}\incl \cF_{i}$. (Note that as \eqref{eq: excess bundle map} is injective, the cone coincides with the cokernel, which is supported on a subscheme of $X_R$ finite over $\Spec R$ by the definition of the modification types in $\Hk^r_{U(n), \LL}$.) 

We use the same notation $\cV^{r}_{\cE}$ to denote the pullback of $\cV^{r}_{\cE}$ to  $\Sht^{r}_{U(n), \LL}$. In the future we will typically consider the latter object. We define $\cV^r$ to be the bundle on $\Bun_{\GL(m)'}(k) \times \Sht^{r}_{U(n), \LL}$ whose restriction to $\{\cE\} \times \Sht^{r}_{U(n), \LL}$ is $\cV_{\cE}^r$.

\begin{lemma}\label{lem: excess bundle graded} The bundle $\cV^{r}_{\cE}$ on $\Hk^{r}_{U(n), \LL}$ carries a filtration with associated graded
\[
\bigoplus_{i=1}^r p_i^* \s^{*}\cE^{*} \otimes \ell_i .
\]
Here $\ell_{i}$ are the tautological bundles over $\Hk^{r}_{U(n), \LL}$ introduced in \S \ref{ssec: tautological bundles}. In particular, we have
\begin{equation}
\prod_{i=1}^{r}c_{m}(p_{i}^{*}\s^{*}\cE^{*}\ot\ell_{i})=c_{mr}(\cV^{r}_{\cE})\in \Ch^{mr}(\Hk^{r}_{U(n), \LL}).
\end{equation}
\end{lemma}

\begin{proof}
Filter the cokernel of 
\[
\bigoplus_{i=1}^r \cF^{\flat}_{i-1/2} \to \bigoplus_{i=1}^{r}  \cF_{i}
\]
by the grading on $i$. Then its associated graded is $\bigoplus_{i=1}^r \ell_i$. This induces a filtration on $\cV^r_{\cE}$ with associated graded 
\[
\bigoplus_{i=1}^r p_i^* \s^{*}\cE^{*} \otimes \ell_i .
\]
Here the $\s^{*}\cE^{*}$ (rather than $\cE^{*}$) is due to the fact that the line bundle $\ell_i$  is supported at $\sigma x'_i$.
\end{proof}

\begin{lemma}\label{lem: cotangent complex for pi}

Let $BH_{1}$ be a smooth $m$-framed gerbe over $X$ and $BH_2 = BU(n)_{\LL}$ with the standard $n$-framed structure. Abbreviate $\sM := \sM_{H_1, H_2}$. Then the relative tangent complex for the map $\Sht_{\sM}^r \xrightarrow{\pi} \Bun_{H_{1}}(k)\times \Sht_{U(n), \LL}^r$ is perfect, and for $\cE\in \Bun_{H_{1}}(k)$ we have $$\bT_{\pi}|_{\pi^{-1}(\{\cE\}\times \Sht_{U(n), \LL}^r)}\cong \pi^{*}\cV_{\Vect(\cE)}^{r}[-1].$$

In particular, $\Sht_{\sM}^r \rightarrow \Bun_{H_{1}}(k) \times \Sht_{U(n), \LL}^r$ is quasi-smooth, so $\Sht_{\sM}^r $ is quasi-smooth. 
\end{lemma}

\begin{proof}
We calculate the tangent complex of $\Sht_{\sM}^r$ using the presentation \eqref{eq: Sht_M cartesian square}. Consider an $R_{\bu}$-point of $\Sht_{\sM}^r$, represented by the data $(\ul{x}', \cE, \cF_0, \cF_1, \ldots, \cF_r)$ and a diagram 
\begin{equation}
\begin{tikzcd}
 & & &  \Vect(\cE)   \ar[dll, "t_0"] \ar[drr, "t_{r-1}"']  \ar[rrrr, equals] &  & & & \Vect(\ft \cE )\ar[dd, "\ft  (i_0 \circ t_{0})"] \\
& \cF_{1/2}^{\flat} \ar[dl, hook, "i_0"] \ar[dr, hook]  & & \ldots   & &  \cF_{r-1/2}^{\flat} \ar[dl, hook] \ar[dr, hook]   \\ 
\cF_0 & & \cF_1  & \ldots & \cF_{r-1} & & \cF_r \ar[r, "\sim"] & \ft \cF_0 
\end{tikzcd}
\end{equation}
By the behavior of cotangent complexes in Cartesian squares, we see from Corollary \ref{cor: Hitchin to Bun tangent} and Lemma \ref{lem: TC Hk_M to Hk} that $\bL_{\pi}$ is perfect, and the tangent complex of $\Sht_{\sM}^r$ is the derived fiber (i.e. cone shifted by $1$) of the map 
\[
\bT_{\Hk_{\sM}^r}|_{\Sht_{\sM}^r } \times \bT_{\sM} |_{\Sht_{\sM}^r } \rightarrow  \bT_{\sM^{2}}|_{\Sht_{\sM}^r }
\]

The Cartesian square \eqref{eq: Sht_M cartesian square} fits into a commutative diagram where the back and front faces are Cartesian 
\[
\begin{tikzcd}[column sep = huge]
\Sht_{\sM}^r \ar[r] \ar[d] \ar[ddr] & \Hk^{r}_{\sM}  \ar[d, "\pr_{0}\times \pr_{r}"] \ar[ddr]  \\
\sM \ar[r, "{(\Id, \Frob)}"]  \ar[ddr]  & \sM^2  \ar[ddr] \\
& \Bun_{H_1}(k) \times \Sht_{U(n), \LL}^r \ar[r] \ar[d] & \Bun_{H_1} \times \Hk^r_{U(n), \LL} \ar[d, "{(\Delta, \pr_0 \times \pr_r)}"] \\
& \Bun_{H_1} \times \Bun_{U(n), \LL} \ar[r, "{({(\Id, \Frob)}, {(\Id, \Frob)})}"] & \Bun_{H_1}^2  \times \Bun_{U(n), \LL}^2
\end{tikzcd}
\]

To shorten notation, we write $S:=\Sht^{r}_{\sM}$. By the same argument as in the proof of Lemma \ref{lem: TC Hk_M to Hk}, $\bT_{\pi}$ is the (homotopy) limit of the diagram 
\[
\begin{tikzcd}
\bT_{\Hk_{\sM}^r} |_{S}  \oplus \bT_{\sM} |_{S}  \ar[r]  \ar[d] & \bT_{\sM}^{\oplus 2}|_{S}  \ar[d] \\
(\bT_{\Bun_{H_1}} |_{S}\op \bT_{\Hk_{U(n), \LL}^r} |_{S}  )\oplus (\bT_{\Bun_{H_1}} |_{S}  \oplus \bT_{\Bun_{U(n), \LL}}|_{S} )\ar[r] & \bT_{\Bun_{H_1}}^{\oplus 2}|_{S}  \oplus  \bT_{\Bun_{U(n), \LL}}^{\oplus 2}|_{S} 
\end{tikzcd}
\]
To compute this we take fiber of the vertical morphisms, using Corollary \ref{cor: Hitchin to Bun tangent} and Lemma \ref{lem: TC Hk_M to Hk}. This says that for any $R_{\bu}$-point of $S$, the pullback of the above diagram to $R_{\bu}$ is (naturally in $R_{\bu}$) isomorphic to $R\pr_*(-)$ applied to the (homotopy) limit of the diagram of complexes on $X_{R_{\bu}}'$ below: 

\[
\begin{tikzcd}
\left( \bigoplus_{i=1}^{r} \cHom(\Vect(\cE), \cF_{i-1/2}^{\flat}) \right)   \oplus \cHom(\Vect(\cE), \cF_0) \ar[r]  \ar[d] &  \cHom(\Vect(\cE), \cF_0) \oplus \cHom(\Vect(\cE), \cF_r)   \\
\bigoplus_{i=1}^{r-1} \cHom(\Vect(\cE), \cF_i) 
\end{tikzcd}
\]
where the maps are:
\begin{equation*}
\xymatrix{ (v_{1/2},\cdots, v_{r-1/2}; v_{0})\ar[r]\ar[d] & (v_{0}+v_{1/2}; v_{r-1/2})\\
 (v_{1/2}- v_{3/2}, \cdots, v_{r-3/2}-v_{r-1/2})}
\end{equation*}
We may rewrite the (homotopy) limit as the complex
\[
\begin{tikzcd}
\left(  \bigoplus_{i=1}^{r} \cHom(\Vect(\cE), \cF_{i-1/2}^{\flat})   \right)  \oplus \cHom(\Vect(\cE), \cF_0) \ar[d] \\
 \left( \bigoplus_{i=1}^{r-1} \cHom(\Vect(\cE), \cF_i)  \right)  \oplus \cHom(\Vect(\cE), \cF_0) \oplus \cHom(\Vect(\cE), \cF_r )
 \end{tikzcd}
 \]
Unraveling down the definitions of the maps, this is seen to be quasi-isomorphic to \eqref{eq: excess bundle map} after cancelling the summand $\cHom(\Vect(\cE), \cF_0)$.

\end{proof}

\begin{remark}
More generally, for any gerbe of unitary type $BH_2$ we can define the excess bundle $\cV_{\cE}^r$ on $\Sht_{H_2}^r$ by the same formulas. It carries a filtration with associated graded given by the same formulas as in Lemma \ref{lem: excess bundle graded}. The analogue of Lemma \ref{lem: cotangent complex for pi} holds, by the same argument. 
\end{remark}

\section{Fundamental classes of derived special cycles}\label{sec: VFC}

\subsection{Summary of derived intersection theory}

For the framework of intersection on derived stacks, we will use the work \cite{KhanI} of Khan. In order to make this paper as self-contained as possible, we give a quick summary of the basic facts from \cite[\S 2,3]{KhanI} that we will need, simplified to our situation of interest. 

\subsubsection{Motivic Borel-Moore homology} 

The role of Chow groups of a locally finite type derived Artin stack $\sX$ will be played by its \emph{motivic Borel-Moore homology} groups $\mrm{H}_s^{\BM}(\sX/\Spec k, \Q(r))$ as  defined in \cite[Definition 2.1, Example 2.10]{KhanI}. (Only the case $s=2r$ will be of interest to us.) Henceforth we omit the ``$/\Spec k$'' when the base is $\Spec k$. 

According to \cite[Example 2.10]{KhanI}, for $\cX$ a classical Artin stack locally of finite type over $k$, $\mrm{H}_{2r}^{\mrm{BM}}(\cX , \Q(r) )$ identifies with the Chow groups (with $\Q$-coefficients) of Joshua \cite{Josh02}; when $\cX$ is of finite type they are identified with the Chow groups (with $\Q$-coefficients) of Kresch \cite{Kr99}. We shall see shortly in \S \ref{sssec: derived invariance} that for a locally finite type derived Artin stack $\sX $ over $\Spec k$, $\mrm{H}_{2r}^{\mrm{BM}}(\sX, \Q(r) )$ can be identified with the motivic Borel-Moore homology of the underlying classical stack $\cX := \sX_{\mrm{cl}}$, and thereby interpreted in terms of Chow groups. 

More generally, if $\sX \rightarrow \sS$ is a locally finite type morphism of derived Artin stacks over $k$, then there is a theory of \emph{relative motivic Borel-Moore homology} groups $\mrm{H}_s^{\BM}(\sX/\sS, \Q(r))$. In this paper we are mainly concerned with the absolute groups; the relative groups play a technical role in some intermediate statements. 

We next discuss the basic functorialities enjoyed by $\mBM(\sX/\sS)$. 

\subsubsection{Proper pushforward} (\cite[\S 2.2.1]{KhanI}) 
If $f \co \sX \rightarrow \sY$ is a representable (cf. \S \ref{sssec: representable morphism}) proper morphism of derived Artin stacks, locally of finite type over $\sS$, then there are functorial direct image morphisms
\[
f_* \co \mrm{H}_{s}^{\mrm{BM}}(\sX / \sS , \Q(r) ) \rightarrow \mrm{H}_{s}^{\mrm{BM}}(\sY/ \sS , \Q(r) ).
\]

\subsubsection{Smooth pullback} (\cite[\S 2.2.2]{KhanI}) If $f \co \sX \rightarrow \sY$ is a representable smooth morphism of derived Artin stacks, locally of finite type over $\sS$, of relative dimension $d$, then there is a functorial pullback 
\[
f^! \co \mrm{H}_{s}^{\mrm{BM}}(\sY/ \sS  , \Q(r) ) \rightarrow \mrm{H}_{s+2d}^{\mrm{BM}}(\sX/ \sS  , \Q(r+d) ).
\]

\subsubsection{Derived invariance}\label{sssec: derived invariance} For any derived Artin stack $\sX$ over $k$, we denote by $i_{\sX} \co \cX  \rightarrow \sX$ the inclusion of the underlying classical stack (cf. \S \ref{sssec: classical truncation}). According to \cite[Theorem 2.19(ii)]{KhanI}, if $\sX$ is locally finite type over $\sS$ then the direct image 
\[
(\iota_{\sX})_* \co  \mrm{H}_{s}^{\mrm{BM}}(\cX/\sS , \Q(r) ) \rightarrow \mrm{H}_{s}^{\mrm{BM}}(\sX /\sS , \Q(r))
\]
is an isomorphism. 

\subsubsection{Base change}\label{sssec: base change} Consider a commutative square of derived Artin stacks
\[
\begin{tikzcd}
\sY \ar[r] \ar[d] & \sT \ar[d, "f"] \\
\sX \ar[r] & \sS
\end{tikzcd}
\]
which is Cartesian on the underlying classical stacks. There is a base change homomorphism \cite[\S 2.2.3]{KhanI}
\[
f^* \co \mrm{H}_{s}^{\BM}(\sX/\sS; \Q(r)) \rightarrow \mrm{H}_{s}^{\BM}(\sY/\sT, \Q(r))
\]

\begin{example}
We note that in the special case where $\sT = \sS$, $f$ is the identity map, and $\sY = \cX$ (the underlying classical truncation of $\sX$) with its canonical map to $\sX$, chasing through the definitions reveals $f^* = (i_{\sX})_{*}^{-1}$ to be the isomorphism of derived invariance \S \ref{sssec: derived invariance}.
\end{example}

\subsubsection{Quasi-smooth pullback}\label{sssec: qs pullback} If $f \co \sX \rightarrow \sY$ is a quasi-smooth morphism of derived Artin stacks locally finite type over $\sS$, then we may define the \emph{relative virtual dimension} of $f$ at $x \in \sX$ to be the Euler characteristic of $\bL_{f}$ at $x$ (which could be negative). 

Letting $d$ be the relative virtual dimension of $f \co \sX \rightarrow \sY$, there is a Gysin map \cite[Construction 3.4]{KhanI}
\[
f^! \co \mrm{H}_{s}^{\BM}(\sY/\sS, \Q(r)) \rightarrow \mrm{H}_{s+2d}^{\BM}(\sX/\sS, \Q(r+d))
\]
By \cite[\S 3.3]{KhanI}, if $\sX$ and $\sY$ are classical and $f$ is representable in (classical) Deligne-Mumford stacks, then the resulting $f^!$ agrees with the Gysin pullback of classical stacks \cite{Man12}.

\subsubsection{Compatibility with the refined Gysin homomorphism}We shall need the following compatibility of the quasi-smooth pullback with the classical refined Gysin homomorphism. Suppose $f \co \cS \rightarrow \cT$ is a quasi-smooth morphism between classical Artin stacks of relative dimension $d$ representable by Deligne-Mumford stacks, such that $f$ satisfies the hypotheses of \cite[Construction 3.6]{Man12}, $\sY$ is a quasi-smooth derived Artin stack, and $g \co \sY \rightarrow \cT$ is locally of finite type. Let $\cY$ be the classical truncation of $\sY$, and suppose that the classical fiber product $\cX := \cY \stackrel{\mrm{cl}}\times_{\cT} \cS \rightarrow \cY$ satisfies the hypotheses of \cite[Construction 3.6]{Man12}. Note that $\cX$ is the classical truncation of $\sX  := \sY \times_{\cT} \cS$. Consider the diagram with the bottom square being derived Cartesian and the outer square being Cartesian as classical stacks: 
\[
\begin{tikzcd}
\cX \ar[d, "\iota_{\sX}"] \ar[r] & \cY \ar[d,  "\iota_{\sY}"] \\
\sX \ar[r, "f'"] \ar[d] & \sY \ar[d, "g"]  \\
\cS \ar[r, "f"] & \cT 
\end{tikzcd}
\]
The hypotheses ensure that the refined virtual pullback $f_{\mrm{ref}}^! \co \Ch_*(\cY) \rightarrow \Ch_{*+d}(\cX)$ is defined \cite[Construction 3.6]{Man12}. On the other hand, we have the identification $\iota_{\sY*} \co \Ch_s(\cY) \cong \mBM_{2s}(\cY, \Q(s)) \xrightarrow{\sim} \mBM_{2s}(\sY, \Q(s) )$ from derived invariance.

\begin{lemma}\label{lem: refined gysin compatibility}
Following the notation above, the diagram below commutes.
\[
\begin{tikzcd}
 \Ch_s(\cY) \ar[d, "i_{\sY*}", "\sim"'] \ar[r, "f_{\mrm{ref}}^!"]   & 
\Ch_{s+d}(\cX ) \ar[d, "i_{\sX*}", "\sim"'] \\ 
  \mBM_{2s}(\sY, \Q(s) )  \ar[r, "{(f')^!}"]    &  \mBM_{2s+2d}(\sX, \Q(s+d) ) 
\end{tikzcd}
\]
\end{lemma}

\begin{proof}
The argument is very similar to that of \cite[\S 3.3]{KhanI}, which handles the case where $\sY \rightarrow \cT$ is the identity map (and in particular $\sX$ is classical). We explain the necessary adjustments in the present situation. Let $\mrm{C}_{\cX/\cY}$ be the intrinsic normal cone for $\cX \rightarrow \cY$, and $\mrm{D}_{\cX/\cY}$ be Kresch's deformation, so we have a diagram 
\[
\begin{tikzcd}
\mrm{C}_{\cX/\cY} \ar[r] \ar[d]& \mrm{D}_{\cX/\cY} \ar[d] & \cY \times \G_m \ar[l] \ar[d] \\ \
 0  \ar[r, hook] & \A^1 & \G_m \ar[l, hook'] 
 \end{tikzcd}
\]
Since $f$ is quasi-smooth, it has a normal bundle stack $N_f$,  which is the stack associated to the (co-connective) two-term complex $\bT_f[1]$ \cite[\S 1.3]{KhanI}. 
The intrinsic normal cone $\mrm{C}_{\cX/\cY}$ admits an embedding into $f^* \mrm{N}_{f} \cong  \iota_{\sX}^* \mrm{N}_{f'}$, which fits into a commutative diagram 
\[
\begin{tikzcd}
\mrm{C}_{\cX/\cY} \ar[r] \ar[d, "a"] & \mrm{D}_{\cX/\cY} \ar[d] & \cY \times \G_m \ar[l] \ar[d] \\
f^* \mrm{N}_{f} \ar[r] & \mrm{D}_{\sX/\sY} & \sY \times \G_m \ar[l] 
\end{tikzcd}
\]
where $\mrm{D}_{\sX/\sY} $ is the deformation to the normal bundle stack for the quasi-smooth morphism $f': \sX \rightarrow \sY$ \cite[\S 1.4]{KhanI}. The rest of the argument concludes as in \cite[\S 3.3]{KhanI}. 
\end{proof}

\begin{remark} Lemma \ref{lem: refined gysin compatibility} implies that the intersection product of \S \ref{sssec: intersection product} is compatible with that of \cite[\S A.1.4]{YZ}, a fact that we will repeatedly use without further comment.
\end{remark}

\subsubsection{Top Chern class} (\cite[\S 2.2.4]{KhanI}) If $\cE$ is a finite locally free sheaf of rank $r$ on a derived Artin stack $\sX$ of finite type over $k$, then there is a \emph{top Chern class} $c_r(\cE) \in \mrm{H}^{\BM}_{-2r}(\sX/\sX , \Q(-r))$. To compare this with the usual formulation of Chern classes, we observe that $\mBM_{-s}(\sX/\sX, \Q(-r))$ is naturally isomorphic to the motivic \emph{cohomology} groups $\mrm{H}^s(\sX, \Q(r))$.

Next we will discuss some operations on these motivic Borel-Moore homology groups.

\subsubsection{Composition product}(\cite[\S 2.2.5]{KhanI}) Given a derived Artin stack $\sT$ locally of finite type over $\sS$, and a derived Artin stack $\sX$ locally of finite type over $\sT$, there is a composition product 
\begin{equation}\label{eq: composition product}
\circ \co \mBM_s(\sX/\sT, \Q(r)) \otimes \mBM_{s'}(\sT/\sS, \Q(r')) \rightarrow \mBM_{s+s'}(\sX/\sS, \Q(r+r')).
\end{equation}

\subsubsection{Virtual fundamental classes}\label{sssec: VFC}
We next discuss one of the key features provided by derived algebraic geometry, namely the intrinsic construction of virtual fundamental classes. 

 Let $f \co \sX \rightarrow \sS$ be a quasi-smooth morphism of derived Artin stacks, of relative virtual dimension $d$. Write $1_{\sS}$ for the unit of $\mBM_0(\sS/\sS; \Q(0))$. Then the \emph{relative fundamental class} of $f$ is \cite[Construction 3.6]{KhanI} 
$$[\sX/\sS] := f^!(1_{\sS})  \in \mrm{H}_{2d}^{\BM}(\sX/\sS; \Q(d)).$$ 
Of particular importance is the case $\sS = \Spec k$, in which case we write $[\sX] := [\sX/ \Spec k]$ and call it the \emph{virtual fundamental class of $\sX$}. Note that by \S \ref{sssec: derived invariance}, we may view $[\sX] \in \mBM_{2d}(\cX; \Q(d)) \cong \Ch_d(\cX)$ where $\cX$ is the underlying classical stack of $\sX$, and we will frequently do so. 

When $\iota_{\sX}: \cX\to \sX$ is an isomorphism  and $\cX$ is smooth, then $[\sX]$ is the usual fundamental class $[\cX]^{\nai}$. 

We next establish some basic properties of these virtual fundamental classes. 



\subsubsection{Intersection product of virtual fundamental classes}\label{sssec: intersection product}

Let $\sX, \sY$ and $\sY'$ be derived Artin stacks locally finite type and equidimensional over $k$, and suppose furthermore that $\sX$ is smooth and $\sY, \sY'$ are quasi-smoooth over $k$. Suppose we have maps (not necessarily quasi-smooth) $f \co \sY \rightarrow \sX$ and $f' \co \sY' \rightarrow \sX$. Consider the Cartesian square 
\[
\begin{tikzcd}
\sY \times_{\sX} \sY' \ar[d] \ar[r, "\Delta'_{\sX}"] & \sY \times \sY' \ar[d] \\
\sX \ar[r, "\Delta_{\sX}"] & \sX \times \sX 
\end{tikzcd}
\]
The morphism $\Delta'_{\sX}$ is quasi-smooth as it is the base change of the quasi-smooth morphism $\D_{\sX}$. In particular $\sY \times_{\sX} \sY'$ is also quasi-smooth over $k$ of virtual dimension $r=\dim \sY+\dim \sY'-\dim \sX$.
We write 
\[
[\sY] \cdot_{\sX} [\sY' ]:= (\Delta'_{\sX})^! ([\sY \times \sY']) \in \mrm{H}_{2r}^{\BM}(\sY \times_{\sX} \sY', \Q(r)).
\]

\begin{lemma}\label{lem: VFC fibered product} In the situation above, we have
\begin{equation}\label{eq: refined intersection product}
[\sY] \cdot_{\sX} [\sY' ] =  [\sY \times_{\sX} \sY'].
\end{equation}
\end{lemma}

\begin{proof}
By definition $[\sY \times \sY'] = \mrm{pr}^! [\Spec k]$ where $\mrm{pr} \co \sY \times \sY' \rightarrow \Spec k$ is the structure map, and $\mrm{pr}^!$ is defined because $\sY \times \sY'$ is quasi-smooth. Hence we have 
\[
[\sY] \cdot_{\sX} [\sY' ]  = ( \Delta'_{\sX})^!  \pr^! [\Spec k] = (\pr \circ \Delta'_{\sX})^! [\Spec k]  =  [\sY \times_{\sX} \sY'].
\]
\end{proof}	


\subsubsection{Excess intersection formula}\label{sssec: excess intersection} We shall make crucial use of the following \emph{excess intersection formula}, which is \cite[Proposition 3.15]{KhanI}\footnote{We note that the conventions of \cite{KhanI} are off from ours by a dualization, e.g. the virtual fundamental class of a self-intersection in \cite[Corollary 3.17]{KhanI} is the top Chern class of what is called the ``conormal bundle'' in \emph{loc. cit.}, whereas we would call it the normal bundle.}. Suppose we have a commutative (but \emph{not} necessarily Cartesian) square of derived Artin stacks over $k$, 
\begin{equation}\label{eq: excess square}
\begin{tikzcd}
\sX' \ar[r, "g"] \ar[d, "p"] & \sY' \ar[d, "q"] \\
\sX \ar[r, "f"] & \sY 
\end{tikzcd}
\end{equation}
where $f$ and $g$ are quasi-smooth, and equidimensional. We say that \eqref{eq: excess square} is an \emph{excess intersection square} if it is Cartesian on underlying classical stacks, and the homotopy fiber of the canonical map $p^* \bL_{\sX/\sY} [-1] \rightarrow \bL_{\sX'/\sY'} [-1]$ is a locally free $\cO_{\sX'}$-module of finite rank $r$, whose dual we call the \emph{excess bundle} $\cE$. Then we have the top Chern class $c_r(\cE) \in \mBM_{-2r}(\sX'/\sX', \Q(-r))$. The excess intersection formula asserts that
\[
q^*[\sX/\sY] = c_r(\cE) \circ [\sX'/\sY'] \in \mrm{H}_{2d}^{\BM} (\sX'/ \sY' , \Q(d)) ,
\]
where $d$ is the virtual dimension of $f$ and $q^*$ is the base change map of \S \ref{sssec: base change}. 

\begin{lemma}\label{lem: excess VFC}
Let $p: \sX'\to \sX$ be a map of quasi-smooth derived Artin stacks locally finite type over $k$ that induces an isomorphism on their classical truncations $p_{\mrm{cl}}: \cX'\isom \cX$. Assume $\bL_{p}[-2]$ is a locally free $\cO_{\sX'}$-module of finite rank $r$. Then 
\begin{equation*}
[\sX]=c_{r}(\bT_{p}[2]) \c [\sX'] \in \Ch_{d}(\cX),
\end{equation*}
where $d$ is the virtual dimension of $\sX$ (note here $\bT_{p}[2]$ is a locally free $\cO_{\sX'}$-module of finite rank $r$). 
\end{lemma}
\begin{proof}
Apply the excess intersection formula to the square
\begin{equation*}
\xymatrix{ \sX' \ar[d]^{p}\ar[r] & \Spec k\ar@{=}[d]\\
\sX\ar[r] & \Spec k
}
\end{equation*} 
\end{proof}

\subsection{Calculation of virtual fundamental classes}
We now return to the (derived) Hitchin stacks. Fix $m \leq n$ and let $\sM  = \sM_{\GL(m)', U(n), \LL}$ and $\cM =  \cM_{\GL(m)', U(n), \LL}$, which is the classical truncation of $\sM$. In the future we will suppress $\LL$ for notational simplicity. As $\Sht_{\sM}^r$ is quasi-smooth by Lemma \ref{lem: cotangent complex for pi}, the virtual fundamental class $[\Sht^{r}_{\sM}]  \in \Ch_*(\Sht^{r}_{\cM})$ is defined by \S \ref{sssec: VFC}. We can now confirm that the virtual fundamental classes of $\cZ_{\cE}^{r}$ constructed earlier in \S \ref{ssec: VFC for special cycle} agree with the (pullbacks of) components of $[\Sht^{r}_{\sM}]$.

\begin{thm}\label{thm: VFC from derived shtuka}  Recall $\sM  = \sM_{\GL(m)', U(n)}$ and $\cM =  \cM_{\GL(m)', U(n)}$. We have 
\[
[\Sht_{\sM}^r]|_{\cZ_{\cE}^r}  = [\cZ_{\cE}^r] \in \Ch_{(n-m)r}(\cZ_{\cE}^r),
\]
where the latter is as in Definition \ref{def:special cycle classes}(1). Here, the notation $( \cdot )|_{\cZ_{\cE}^r}$ means pullback along the map $\cZ_{\cE}^r \rightarrow \Sht_{\cM}^r$ which is finite \'etale over an open-closed substack. 
\end{thm}



The rest of this subsection is devoted to the proof of Theorem \ref{thm: VFC from derived shtuka}. Recall the open-closed decomposition \eqref{ZEK decomp} of $\cZ^{r}_{\cE}$ into $\cZ^{r}_{\cE}[\cK]^{\c}$ according to the kernel of $t_{i}: \cE\to \cF_{i}$. The proof will proceed by describing a stratification of $\Sht_{\sM}^r$ that induces a decomposition of $[\Sht_{\sM}^r]$ into many pieces. We will identify a piece of $[\Sht_{\sM}^r]$ that matches the non-degenerate term $\cZ^{r,\c}_{\cE}$, then another piece that matches the most degenerate term $\cZ^{r}_{\cE}[\cE]$, and then calculating the contribution of all pieces by reducing them to the two extreme cases.

\subsubsection{Kernel decomposition} Recall that $\Sht^{r}_{\cM}$ is the disjoint union of open-closed substacks $\ol{\cZ}^{r}_{\cE}$ for $\cE\in \Bun_{\GL(m)'}(k)$. Each $\cZ^{r}_{\cE}$ is the disjoint union of open-closed substacks $\cZ^{r}_{\cE}[\cK]^{\c}$ indexed by sub-bundles $\cK\subset \subset\cE$ (see \S \ref{ssec: kernel cycles} and Lemma \ref{lem: kernel circ open-closed}). The substacks where $\cK=  0$ and $\cK = \cE$ are stable under $\Aut(\cE)(k)$, we denote $\ol{\cZ}^r_{\cE}[\cK]^{\c} := [\cZ^{r}_{\cE}[\cK]^{\c} /( \Aut(\cE)(\F_q))]$ in these cases. 

Since $\Sht^{r}_{\sM}$ has underlying classical stack $\Sht^{r}_{\cM}$, it similarly decomposes into open-closed derived substacks  $\ol{\sZ}^{r}_{\cE}$ (whose underlying classical stack is $\ol{\cZ}^{r}_{\cE}$), whose pullbacks to $\cZ_{\cE}^r$ further decompose into open-closed derived substacks $\sZ^{r}_{\cE}[\cK]^{\c}$ (whose underlying classical stack is $\cZ^{r}_{\cE}[\cK]^{\c}$).  Similarly, we have the open-closed derived substack $\sZ^{r}_{\cE}[\cK]\subset \sZ^{r}_{\cE}$ whose classical truncation is $\cZ^{r}_{\cE}[\cK]$. 

To summarize, we have a decomposition into open-closed derived substacks
\begin{equation}
\Sht^{r}_{\sM} =\coprod_{\substack{\cE\in \Bun_{\GL(m)'}(k), \\ \cK\subset\subset \cE}} \ol{\sZ}^{r}_{\cE}
\end{equation}
whose classical truncation recovers Example \ref{ex: hitchin shtuka U(n)}, and this is refined by a finite \'{e}tale covering
\[
\coprod_{\substack{\cE\in \Bun_{\GL(m)'}(k), \\ \cK\subset\subset \cE}} \sZ^{r}_{\cE} [\cK]^{\c} \surj \Sht^{r}_{\sM} 
\]

The virtual fundamental classes $[\sZ^{r}_{\cE}[\cK]^{\c}]\in \Ch_{r(n-m)}(\cZ^{r}_{\cE}[\cK]^{\c})$ and $[\sZ^{r}_{\cE}[\cK]]\in \Ch_{r(n-m)}(\cZ^{r}_{\cE}[\cK])$ are defined as the restriction of the virtual fundamental class $[\Sht^{r}_{\sM}]$. We define $\sZ_{\cE}^r$ as the (derived) fibered product 
\[
\begin{tikzcd}
\sZ^{r}_{\cE} \ar[r] \ar[d] &  \ol{\sZ}^r_{\cE} \ar[d] \\
\Spec k \ar[r] & {[\Spec k / (\Aut(\cE)(\F_q))]}
\end{tikzcd}
\]
and we define $\sZ^{r}_{\cE}[\cK], \sZ^{r}_{\cE}[\cK]^{\c}$, etc. similarly.

%
%

\begin{remark}By pulling back along the finite \'{e}tale map $\cZ_{\cE}^r \rightarrow \ol{\cZ}_{\cE}^r$, we see that Theorem \ref{thm: VFC from derived shtuka} implies (and is equivalent to) 
\[
[\sZ_{\cE}^r] = [\cZ_{\cE}^r] \in \Ch_{(n-m)r} (\cZ_{\cE}^r). 
\]
\end{remark}

\subsubsection{Non-degenerate terms} 
We consider $\ol{\sZ}^{r,\c}_{\cE}:= \ol{\sZ}^{r}_{\cE}[0]^{\c}$ whose underlying classical stack is $ \ol{\cZ}_{\cE}^{r , \circ}$. We will show:


\begin{prop}\label{lem: non-degen terms} We have
\begin{equation}\label{eq: non-degen terms}
[\ol{\sZ}_{\cE}^{r,\c}] = [\ol{\cZ}_{\cE}^{r, \circ}]\in \Ch_{r(n-m)}(\ol{\cZ}_{\cE}^{r , \circ}). 
\end{equation}
\end{prop}

Let $H_{1}\to \GL(m)'$ be any homomorphism of smooth group schemes over $X$ (although we shall only need the case where this map is the identity). 

\begin{lemma}\label{lem: GL(n)' to U(n) cartesian} \hfill
\begin{enumerate}
\item The following square is (derived) Cartesian: 
\[
\begin{tikzcd}
\sM_{H_1, U(n)}  \ar[d] \ar[r] & \sM_{H_1, \GL(n)'} \ar[d] \\
\Bun_{U(n)}  \ar[r] & \Bun_{\GL(n)'} 
\end{tikzcd}
\]

\item The following square is (derived) Cartesian: 
\[
\begin{tikzcd}
\Hk_{\sM_{H_1, U(n)}}^r \ar[d] \ar[r] & \Hk_{\sM_{H_1, \GL(n)'} }^r \ar[d] \\
\Hk_{U(n)}^r \ar[r] & \Hk_{\GL(n)'}^r 
\end{tikzcd}
\]
\item The following square is (derived) Cartesian: 
\[
\begin{tikzcd}
\Sht_{\sM_{H_1, U(n)}}^r  \ar[d] \ar[r] & \Sht_{\sM_{H_1, \GL(n)'}}^r \ar[d] \\
\Sht_{U(n)}^r  \ar[r] & \Sht_{\GL(n)'}^r
\end{tikzcd}
\]
\end{enumerate}
\end{lemma}

\begin{proof}
Parts (1) and (2) are immediate from the definitions. We focus on (3). 

Abbreviate $\sM := \sM_{H_1,U(n)}$ and $\sM' := \sM_{H_1, \GL(n)'}$. Consider the commutative diagram below.
\begin{equation}\label{eq: s6 9-term diag 2}
\begin{tikzcd}
\sM' \ar[r, "{(\Id, \Frob)}"] \ar[d]  & \sM' \times \sM'  \ar[d] & \Hk_{\sM'}^r \ar[l,"{(\pr_{0},\pr_{r})}"'] \ar[d] \\
\Bun_{\GL(n)'} \ar[r, "{(\Id, \Frob)}"]  &\Bun_{\GL(n)'}  \times \Bun_{\GL(n)'} & \Hk^r_{\GL(n)'} \ar[l,"{(\pr_{0},\pr_{r})}"'] \\
\Bun_{U(n)} \ar[r, "{(\Id, \Frob)}"]  \ar[u] & \Bun_{U(n)}   \times \Bun_{U(n)}  \ar[u] & \Hk^r_{U(n)} \ar[l, "{(\pr_{0},\pr_{r})}"'] \ar[u]   
\end{tikzcd}
\end{equation}
The derived fibered products along the rows of \eqref{eq: s6 9-term diag 2} are
\begin{equation}\label{eq: s6 H_2 rows}
\begin{tikzcd}
\Sht_{\sM'}^{r} \ar[d]  \\
\Sht_{\GL(n)'}^r  \\
\Sht_{U(n)}^r \ar[u]
\end{tikzcd}
\end{equation}

By parts (1) and (2), the derived fibered products along the columns of \eqref{eq: s6 9-term diag 2} are 
\begin{equation}\label{eq: s6 H_2 columns}
\begin{tikzcd}
\sM  \ar[r] \ar[r, "{(\Id, \Frob)}"] & \sM  \times \sM   & \Hk_{\sM}^{r} \ar[l,"{(\pr_{0},\pr_{r})}"'] 
\end{tikzcd}
\end{equation}

The same proof as for \cite[Lemma A.9]{YZ} gives canonical isomorphisms of derived stacks between the derived fibered products of \eqref{eq: s6 H_2 rows} and \eqref{eq: s6 H_2 columns}. The derived fibered product of \eqref{eq: s6 H_2 columns} is $\Sht_{\sM}^r$. We then conclude by applying \eqref{eq: refined intersection product} to \eqref{eq: s6 H_2 rows}.

\end{proof}

Applying Lemma \ref{lem: VFC fibered product}, we obtain: 


\begin{lemma}We have $[\sM_{H_1, U(n)}] = [\sM_{H_1, \GL(n)'}] \cdot_{\Bun_{\GL(n)'}} [\Bun_{U(n)}] \in \Ch_*(\sM_{H_1, U(n)})$. 
\end{lemma}

\begin{lemma}\label{cor: VFC for GL(n)' circ} We have 
\[
[\Sht_{\sM_{H_1, \GL(n)'}^{\c}}^r ] =  [\Sht_{\cM_{H_1, \GL(n)'}^{\c}}^r] \in \Ch_*(\Sht_{\cM_{H_1, \GL(n)'}^{\c}}^r)
\]
where the right side is defined in Definition \ref{defn: sht_M class for GL(n)'}.
\end{lemma}

\begin{proof}
We abbreviate $\sM'^{ \circ} := \sM_{H_1, \GL(n)'}^{\c}$ and $\cM'^{ \circ} := \cM_{H_1, \GL(n)'}^{\c}$, which is the classical truncation of $\sM'^{ \circ}$. Consider the Cartesian square
\[
\begin{tikzcd}
 \Sht_{\sM'^{ \circ}}^r \ar[r] \ar[d] & (\Hk_{\sM'^{ \circ}}^1)^r  \times \sM'^{ \circ} \ar[d, "q \times \Delta"]\\
(\sM'^{ \circ})^{r+1} \ar[r, "\Phi"] &(\sM'^{ \circ})^{2r+2}
\end{tikzcd}
\]
By Lemma \ref{lem: VFC fibered product}, we have 
\begin{equation}\label{eq: cartesian pullback}
[ \Sht_{\sM'^{ \circ}}^r ] = [(\Hk_{\sM'^{ \circ}}^1)^r \times \sM'^{ \circ}] \cdot_{(\sM'^{ \circ})^{2r+2}} [(\sM'^{ \circ})^{r+1}].
\end{equation}
According to Corollaries \ref{cor: M already classical}(1) and \ref{cor: HkM already classical}(1),  the three corners of the above diagram (except $\Sht_{\sM'^{ \circ}}^r$) are smooth and isomorphic to their classical truncations. By \S\ref{sssec: VFC}, we then have 
\[
[(\Hk_{\sM'^{ \circ}}^1)^r] = [(\Hk_{\cM'^{ \circ}}^1)^r]^{\nai} \in \Ch_*(\Hk_{\cM'^{ \circ}}^r)
\]
 and 
 \[
 [(\sM'^{ \circ})^r ]= [(\cM'^{ \circ})^r]^{\nai}\in \Ch_*((\cM'^{ \circ})^r).
 \]
 Inserting these into \eqref{eq: cartesian pullback} gives
\[
 [(\Hk_{\sM'^{ \circ}}^1)^r \times \sM'^{ \circ}] \cdot_{(\sM'^{ \circ})^{2r+2}} [(\sM'^{ \circ})^{r+1}] = [(\Hk_{\cM'^{ \circ}}^1)^r \times \cM'^{ \circ} ]^{\nai} \cdot_{(\cM'^{ \circ})^{2r+2}} [(\cM'^{ \circ})^{r+1}]^{\nai} \in \Ch_*(\Sht_{\cM'^{ \circ}}^r).
\]
The right hand side is precisely the definition of $[\Sht_{\cM'^{ \circ}}^r]$ in Definition \ref{defn: sht_M class for GL(n)'}.
\end{proof}

\begin{proof}[Proof of Proposition \ref{lem: non-degen terms}]Applying Lemma \ref{lem: GL(n)' to U(n) cartesian}  and Lemma \ref{lem: VFC fibered product}
, we have 
\[
[\Sht_{\sM_{H_1, U(n)}^{\c}}^r ] = [\Sht_{U(n)}^r] \cdot_{\Sht_{\GL(n)'}^r} [\Sht_{\sM_{H_1, \GL(n)'}^{\c}}^r] \in \Ch_*(\Sht_{\cM_{H_1, U(n)}^{\c}}^r ).
\]
By Lemma \ref{cor: VFC for GL(n)' circ}, the RHS above identifies with $[\Sht_{U(n)}^r] \cdot_{\Sht_{\GL(n)'}^r} [\Sht_{\cM_{H_1, \GL(n)'}^{\c}}^r] \in \Ch_*(\Sht_{\cM_{H_1, U(n)}^{\c}
}^r )$. Now specializing to the case $H_{1}=\GL(m')$, and decomposing both sides of the resulting equality according to $\cE \in \Bun_{H_1}(k)$ yields \eqref{eq: non-degen terms}.
\end{proof}

\subsubsection{The most degenerate term}\label{sssec: most degenerate}

We will next handle the most degenerate term $\cZ^{r}_{\cE}[\cE]$. Let ${}_0 \Sht_{\cM}^r$ be the substack of $\Sht^{r}_{\cM}$ where $t_{i}=0$. Then ${}_0 \Sht_{\cM}^r$ is the disjoint union of $\ol{\cZ}^{r}_{\cE}[\cE]$ over $\cE\in \Bun_{\GL(m)'}(k)$, hence open-closed in $\Sht^{r}_{\cM}$. Let ${}_0 \Sht_{\sM}^r=\coprod \sZ^{r}_{\cE}[\cE]\subset \Sht^{r}_{\sM}$ be the corresponding open-closed derived substack. Note that the underlying classical stack of ${}_0 \Sht_{\sM}^r$ is
\begin{equation}\label{cl 0ShtM}
({}_0 \Sht_{\sM}^r)_{\mrm{cl}} = {}_0 \Sht_{\cM}^r \cong  \Bun_{\GL(m)'}(k) \times \Sht_{U(n)}^r.
\end{equation}

In \S \ref{sssec: excess bundle} we defined a bundle $\cV^r$ on $\Bun_{\GL(m)'}(k) \times \Sht_{U(n)}^r$. Below we write $\cV^r|_{{}_0\Sht_{\sM}^r}$ to denote the restriction of $\cV^r$ to ${}_0\Sht_{\sM}^r$ via the composition 
\[
{}_0\Sht_{\sM}^r \rightarrow \Sht_{\sM}^r \rightarrow \Bun_{\GL(m)'}(k) \times \Sht_{U(n)}^r
\]
which we shall denote $\pi$.


\begin{lemma}\label{lem: most degenerate}
We have 
\begin{equation}\label{eq: most degenerate}
[{}_0 \Sht_{\sM}^r] = c_{mr}(\cV^r|_{{}_0\Sht_{\sM}^r}) \cdot [{}_0\Sht_{\cM}^r]^{\nai} \in \Ch_{(n-m)r}({}_{0}\Sht_{\cM}^r).
\end{equation}
\end{lemma}
\begin{proof} We apply Lemma \ref{lem: excess VFC} to the map $\io: {}_0 \Sht_{\cM}^r \to {}_0 \Sht_{\sM}^r $. Note that ${}_0 \Sht_{\sM}^r $ is quasi-smooth by Lemma \ref{lem: cotangent complex for pi}, and ${}_{0}\Sht^{r}_{\cM}$ is smooth by \eqref{cl 0ShtM} (using \cite[Lemma 6.9]{FYZ}). To apply the excess intersection formula, we claim that $\bL_{\iota}$ is concentrated in degree $-2$, and $H^{-2}\bL_{\iota}\cong (\cV^r)^{*}|_{{}_0\Sht_{\sM}^r}$.

We have an exact triangle 
\begin{equation}\label{eq: VFC eq 1}
\iota^* \bL_{{}_0 \Sht_{\sM}^r} \rightarrow \bL_{{}_0 \Sht_{\cM}^r} \rightarrow \bL_\iota.
\end{equation}
Consider the composition 
\[
{}_0 \Sht_{\cM}^r \xrightarrow{\iota} {}_0 \Sht_{\sM}^r \xrightarrow{\pi} \Bun_{\GL(m)'}(k) \times \Sht_{U(n)}^r.
\]
This induces an exact triangle 
\[
\iota^* \bL_{\pi}  \rightarrow  \bL_{\pi\c\io} \rightarrow \bL_{\iota} .
\]
Note that $\pi\c\io = \Id$, so that $\bL_{\pi \circ \iota}  = 0$. Hence $\bL_{\iota} = \iota^* \bL_{\pi} [1]$, which is  $(\cV^r)^*[2]|_{{}_0\Sht_{\sM}^r}$  by Lemma \ref{lem: cotangent complex for pi}.
\end{proof}

By Lemma \ref{lem: excess bundle graded}, the pullback of $c_{mr}(\cV^r|_{{}_0\Sht_{\sM}^r})$ to $\cZ_{\cE}^r$ agrees with the product of Chern classes used to define $[\cZ_{\cE}^r[\cE]]$. This matches the contributions of the ``most degenerate terms''. 

%

\subsubsection{Intermediate terms} In order to simplify notation, we will conflate $\GL(m)$ torsors with rank $m$ vector bundles in this section.   Also, for ease of language we will give the argument in the case where $X'$ is connected. At the end in Remark \ref{rem: geom split case}, we will summarize the adjustments that need to be made if $X'$ is disconnected.

\begin{lemma}\label{lem: kernel K}
For any sub-bundle $\cK\subset \cE$ with quotient $\ov\cE=\cE/\cK$ of rank $i$, we have
\begin{equation}\label{ctop VK}
[\sZ^{r}_{\cE}[\cK]]=c_{(m-i)r}(\cV^{r}_{\cK}|_{\sZ^{r}_{\cE}})\cdot [\sZ^{r}_{\ov\cE}]\in \Ch_{r(n-m)}(\cZ^{r}_{\cE}[\cK]).
\end{equation}
Here $\cV^{r}_{\cK}|_{\sZ^{r}_{\cE}}$ denotes the pullback of $\cV^{r}_{\cK}$ along $\sZ^{r}_{\cE} \rightarrow  \Bun_{\GL(m)'}(k) \times \Sht_{U(n)}^r$, and $[\sZ^{r}_{\ov\cE}]\in \Ch_{r(n-i)}(\cZ^{r}_{\ov\cE})$ is viewed as an element in $\Ch_{r(n-i)}(\cZ^{r}_{\cE}[\cK])$ via the isomorphism $\cZ^{r}_{\cE}[\cK]\cong \cZ^{r}_{\ov \cE}$.
\end{lemma}

\begin{proof}[Proof of Theorem \ref{thm: VFC from derived shtuka} assuming Lemma \ref{lem: kernel K}] Restricting \eqref{ctop VK} to the open-closed $\cZ^{r}_{\cE}[\cK]^{\c}\cong \cZ^{r ,\c}_{\ov \cE}$, we get
\begin{equation}\label{ZEKc}
[\Sht^{r}_{\sM}]|_{\cZ^{r}_{\cE}[\cK]^{\c}}=c_{(m-i)r}(\cV^{r}_{\cK}|_{\sZ^{r}_{\cE}})\cdot [\sZ^{r ,\c}_{\ov\cE}]\in \Ch_{r(n-m)}(\cZ^{r}_{\cE}[\cK]^{\c}).
\end{equation}
By Proposition \ref{lem: non-degen terms} we have $[\sZ_{\ol{\cE}}^{r, \c} ] = [\cZ_{\ol{\cE}}^{r , \circ}]$. Inserting this into \eqref{ZEKc}, we get exactly the expression for $[\cZ_{\cE}^r]$ from Definition \ref{def:special cycle classes}.
\end{proof}

It remains to establish Lemma \ref{lem: kernel K}. Suppose $0\leq  i \leq m \leq n $. We abbreviate $\sM(m,n) = \sM_{\GL(m)', U(n)}$. We define two auxiliary variants of derived Hitchin stacks. 

\begin{itemize}
\item $\sM'$ classifies $\cE \in \Bun_{\GL(m)'}(R_{\bu})$, $\cF \in \Bun_{U(n)}(R_{\bu})$, a vector sub-bundle $\cK \subset \cE$ of rank $m-i$ (so $\cE/\cK \in \Bun_{\GL(i)'}$) and a derived section $t \in R\Gamma(X'_{R_{\bu}}, \cHom(\cE/\cK, \cF))$. Projecting such data to $(\cE/\cK, \cF, t)$ induces a map $\sM' \rightarrow \sM(i,n) $.

\item $\sM''$ classifies $\cE \in \Bun_{\GL(m)'}(R_{\bu})$,$\cF \in \Bun_{U(n)}(R_{\bu})$, a vector sub-bundle $\cK \subset \cE$ of rank $m-i$ (so $\cE/\cK \in \Bun_{\GL(i)'}$) and a derived section $t \in R\Gamma(X'_{R_{\bu}}, \cHom(\cE, \cF))$. Projecting such data to $(\cE, \cF, t)$ induces a map $\sM'' \rightarrow \sM(m,n)$, while sending it to $(\cK, \cF, t|_{\cK} \in R\Gamma(X'_{R_{\bu}}, \cHom(\cK, \cF)))$ induces a map $\sM'' \rightarrow \sM(m-i, n)$. 

\end{itemize}
From the constructions we get a canonical map $\sM'\to \sM''$ sending $(\cK\subset \cE,\cF, t)$ to $(\cK\subset \cE, \cF, t')$ where $t'$ is the image of $t$ under the natural map $R\Gamma(X'_{R_{\bu}}, \cHom(\cE/\cK, \cF))\to R\Gamma(X'_{R_{\bu}}, \cHom(\cE, \cF))$. So we have a diagram 
\begin{equation}\label{eq: diagram of Ms}
\begin{tikzcd}
\sM' \ar[d] \ar[r] & \sM'' \ar[d] \ar[r] & \sM(m,n)\\
\sM(i,n) & \sM(m-i,n)
\end{tikzcd}
\end{equation}

We define $\Hk_{\sM'}^r := \Hk_{\sM(i,n)}^r \times_{\sM(i,n)} \sM' $ and $\Sht_{\sM'}^r$ by the Cartesian square 
\[
\begin{tikzcd}
\Sht_{\sM'}^r \ar[r] \ar[d] & \Hk_{\sM'}^r \ar[d, "\pr_{0} \times \pr_{r}"] \\
\sM' \ar[r, "{(\Id, \Frob)}"] & \sM' \times \sM'
\end{tikzcd}
\]

We have an open-closed decomposition 
\[
\Sht_{\sM(i,n)}^r = \coprod_{\ol{\cE} \in \Bun_{\GL(i)'}(k)} \ol{\sZ}_{\ol{\cE}}^r 
\]
and an open-closed decomposition of $\Sht_{\sM'}^r $ according to the discrete data $(\cK\subset\cE)$, or equivalently according to $\cE$ and $\ov\cE=\cE/\cK$:
\begin{equation}\label{Sht M'}
\Sht_{\sM'}^r = \coprod_{\substack{\ol{\cE} \in \Bun_{\GL(i)'}(k),\\ \cE \in \Bun_{\GL(m)'}(k),\\ \cE \surj \ol{\cE}}}  \Sht_{\sM'}^r(\cE \surj \ol{\cE}).
\end{equation}

Let $\sA(i,n)$ be the Hitchin base for $\sM(i,n)$, classifying $\ol{\cE} \in \Bun_{\GL(i)'}$ and a derived section $a$ of $\cHom(\ol{\cE} , \s^{*}\ol{\cE} ^{\vee} \otimes \nu^* \LL)$ such that $\s^{*}a^{\vee}=a$. Let $\sA'$ be the Hitchin base for $\sM'$, classifying $\cE \in \Bun_{\GL(m)'}$, a vector sub-bundle $\cK \subset \cE$ of rank $m-i$, with quotient bundle $\ol{\cE}$ of rank $i$, and a derived section $a$ of $\cHom(\ol{\cE} , \s^{*}\ol{\cE}^{\vee} \otimes \nu^* \LL)$ such that $\s^{*}a^{\vee}=a$.

\begin{lemma}\label{lem: frob fixed points}
Let $\sY$ be a locally finite type derived stack over $k$. Viewing the discrete groupoid $\sY(k)$ as a constant stack, the diagram 
\[
\begin{tikzcd}
\sY(k) \ar[r] \ar[d] & \sY \ar[d, "{(\Id, \Frob)}"] \\
\sY \ar[r, "\Delta"] & \sY \times \sY
\end{tikzcd}
\]
is derived Cartesian. 
\end{lemma}

\begin{proof}
Let $\sY^{h \Frob}$ be the derived fibered product 
\[
\begin{tikzcd}
\sY^{h\Frob} \ar[r] \ar[d] & \sY \ar[d, "{(\Id, \Frob)}"] \\
\sY \ar[r, "\Delta"] & \sY \times \sY
\end{tikzcd}
\]
Clearly $\sY^{h\Frob}$ receives a canonical map from $\sY(k)$ (regarded as a constant stack), and we want to show that this map is an isomorphism. It suffices to show that $\sY^{h\Frob}$ is isomorphic to its classical truncation, in which case it follows from the analogous statement for finite type classical stacks. To this end, let us examine the tangent complex of the derived fibered product: since $\Frob$ induces the zero map on tangent complexes, it is the derived fibered product of the diagram of complexes
\[
\begin{tikzcd} 
&  \bT_{\sY}|_{\sY^{h\Frob}} \ar[d, "{(\Id, 0)}"] \\  
\bT_{\sY}|_{\sY^{h\Frob}} \ar[r, "{(\Id, \Id)}"] &  \bT_{\sY}|_{\sY^{h\Frob}} \oplus  \bT_{\sY}|_{\sY^{h\Frob}}
\end{tikzcd}
\]
which is evidently zero. Then we conclude using Lemma \ref{lem: classicality of smooth derived stacks}.

\end{proof}

\begin{lemma}\label{lem: left vertical}
For any choice of $\cE \surj \ol{\cE}$ as in \eqref{Sht M'}, the natural map $ \Sht_{\sM'}^r(\cE \surj \ol{\cE}) \rightarrow \Sht_{\sM(i,n)}^r(\ol{\cE}) = \ol{\sZ}_{\ol{\cE}}^r$ is an isomorphism. 
\end{lemma}

\begin{proof}

We have isomorphisms $\sM' \rightarrow \sA' \times_{\sA(i,n)} \sM(i,n)$ and $\Hk^r_{\sM'} \rightarrow \sA' \times_{\sA(i,n)} \Hk^r_{\sM(i,n)}$, which induce $\Sht_{\sM'}^r \xrightarrow{\sim}  \sA'(k) \times_{\sA(i,n)(k)} \Sht_{\sM(i,n)}^r$ by the diagram below (where we have used Lemma \ref{lem: frob fixed points}). 
\[
\begin{tikzcd} 
\Sht_{\sM'}^r \ar[r] \ar[d] \ar[dddr] & \Hk^r_{\sM'} \ar[dddr] \ar[d] \\
\sM' \ar[dddr] \ar[r, "{(\Id, \Frob)}"] & \sM' \times \sM' \ar[dddr] \\
\\
  & \sA'(k) \times_{\sA(i,n)(k)} \Sht_{\sM(i,n)}^r \ar[r] \ar[d]  & \sA' \times_{\sA(i,n)} \Hk^r_{\sM(i,n)} \ar[d] \\
  & \sA' \times_{\sA(i,n)}  \sM(i,n) \ar[r, "{(\Id, \Frob)}"] & (\sA' \times_{\sA(i,n)}  \sM(i,n)) \times (\sA' \times_{\sA(i,n)}  \sM(i,n)) 
\end{tikzcd}
\]
Decomposing this last isomorphism over $\cA'(k)$ gives the result.
\end{proof}

Similarly, we define $\Hk_{\sM''} := \Hk_{\sM(m,n)}^r \times_{\sM(m,n)} \sM''$ and $\Sht_{\sM''}$ by the Cartesian square
\[
\begin{tikzcd}
\Sht_{\sM''}^r \ar[r] \ar[d] & \Hk_{\sM''}^r \ar[d, "\pr_{0} \times \pr_{r}"] \\
\sM'' \ar[r, "{(\Id, \Frob)}"] & \sM'' \times \sM''
\end{tikzcd}
\]
We have open-closed decompositions 
\[
\Sht_{\sM(m,n)}^r = \coprod_{\cE \in \Bun_{\GL(m)'}(k)} \ol{\sZ}_{\cE}^r ,
\]
\[
\Sht_{\sM(m-i,n)}^r = \coprod_{\cK \in \Bun_{\GL(m-i)'}(k)} \ol{\sZ}_{\cK}^r,
\]
and
\[
\Sht_{\sM''}^r = \coprod_{\substack{\cK \in \Bun_{\GL(m-i)'}(k),\\ \cE \in \Bun_{\GL(m)'}(k),\\ \cK \subset \subset\cE}} \Sht_{\sM''}^r(\cK \subset \cE).
\]
We remind that the notation $\cK \subset \subset \cE$ means that $\cK \subset \cE$ is a vector sub-bundle of $\cE$, i.e. $\cE/\cK$ is a vector bundle (as opposed to merely a sub coherent sheaf). 

\begin{lemma}\label{lem right vertical}
For any rank $m$ vector bundle $\cE$ over $X'$ and any vector sub-bundle $\cK \subset \cE$ of rank $m-i$, the map $\Sht_{\sM''}^r\to \Sht_{\sM(m,n)}^r$ restricts to an isomorphism $\Sht_{\sM''}^r(\cK \subset \cE) \isom  \ol{\sZ}_{\cE}^r$.
\end{lemma}

\begin{proof}
The argument is similar to that for Lemma \ref{lem: left vertical}.
\end{proof}

We have a map $z: \Bun_{\GL(m-i)'}  \times \Bun_{U(n)} \rightarrow \sM(m-i, n)$ sending $(\cK \in \Bun_{\GL(m-i)'}, \cF \in \Bun_{U(n)})$ to $(\cK, \cF,  \cK \xrightarrow{0} \cF ) \in \sM(m-i, n)$. This map fits into a Cartesian square 
\[
\begin{tikzcd}
\sM' \ar[r] \ar[d] & \sM'' \ar[d] \\
\Bun_{\GL(m-i)'}  \times \Bun_{U(n)} \ar[r,"z"] & \sM(m-i, n)
\end{tikzcd}
\]
This in turn induces a Cartesian square
\begin{equation}\label{eq: intermediate cartesian}
\begin{tikzcd}
\Sht^r_{\sM'}\ar[r]  \ar[d] & \Sht^r_{\sM''} \ar[d] \\
\Bun_{\GL(m-i)'}(k)  \times \Sht_{U(n)}^r \ar[r] & \Sht_{\sM(m-i, n)}^r
\end{tikzcd}
\end{equation}

\begin{proof}[Proof of Lemma \ref{lem: kernel K}]
Thanks to Lemma \ref{lem right vertical} and Lemma \ref{lem: left vertical}, we have open-closed decompositions
\[
\Sht_{\sM(m-i, n)}^r = \coprod_{\cK} \ol{\sZ}_{\cK}^r, \quad  \Sht_{\sM''}^r = \coprod_{\cK \subset\subset \cE} \ol{\sZ}_{\sE}^r,   \quad  \Sht_{\sM'}^r = \coprod_{\cE \surj \ol{\cE}} \ol{\sZ}_{\ol{\cE}}^r.
\]
Inserting these decompositions into \eqref{eq: intermediate cartesian} and then pulling back to $\sZ_{\cE}^{r}$, we obtain a Cartesian square (where $\ov\cE=\cE/\cK$)
\begin{equation}\label{eq: cartesian 2}
\begin{tikzcd}
\sZ_{\ol{\cE}}^{r} \ar[r,"\zeta"] \ar[d, "\pi_{\ov\cE}"] &   \sZ_{\cE}^r [\cK] \ar[d] \\
\Sht_{U(n)}^r \ar[r,"\io"] & \sZ_{\cK}^r
\end{tikzcd}
\end{equation}
Note that the classical truncation of the top arrow $\z$ is the canonical isomorphism $\cZ_{\ol{\cE}}^{r}\cong \cZ^{r}_{\cE}[\cK]$. We then apply Lemma \ref{lem: excess VFC} to $\z$. Note that both $\sZ_{\ol{\cE}}^{r}$ and $\sZ_{\cE}^r [\cK]$ are quasi-smooth by Lemma \ref{lem: cotangent complex for pi} because they are finite \'{e}tale over open-closed substacks in $\Sht^{r}_{\sM(i,n)}$ and $\Sht^{r}_{\sM(m,n)}$ respectively. By the proof of Lemma \ref{lem: most degenerate}, $\bL_{\io}\cong (\cV^{r}_{\cK})^{*}[2]$. By the base change property of cotangent complexes, $\bL_{\z}\cong \pi_{\ov\cE}^{*}(\cV^{r}_{\cK})^{*}[2]$, so $\bT_{\z}[2]\cong \pi_{\ov\cE}^{*}\cV^{r}_{\cK}$. Now the formula \eqref{ctop VK} follows from Lemma \ref{lem: excess VFC}. 
\end{proof}

\begin{remark}\label{rem: geom split case}
In the case where $X'$ is disconnected, the sub-bundles $\cK \subset \cE$ occurring in the ``decomposition according to the kernel'' need not have the same rank on the two components of $X' = X \sqcup X$. Hence, in that case one needs to replace the unions over $\cK \in \Bun_{\GL(m-i)'}(k)$ above by unions over all sub-bundles $\cK \subset \subset \cE$, and similarly replace the unions over quotients $\ol{\cE} = \cE/\cK \in \Bun_{\GL(i)'}(k)$ by unions over all quotients $\cE \surj \ol{\cE}$. With these adjustments, the proof goes through exactly as above.
\end{remark}

\section{Linear Invariance}\label{sec: linear invariance}

In this section we prove various ``functoriality'' results for the virtual fundamental cycles $[\Sht_{\sM_{H_1, H_2}}^r]$, regarding their compatibility with respect to morphisms induced by gerbe maps $BH_{1}\to BH'_{1}$ and $BH_{2}\to BH'_{2}$. 

Throughout this section we fix a line bundle $\LL$ on $X$ and all Hermitian bundles will be $\LL$-twisted, all unitary gerbes will be twisted by $\LL$, etc. For conciseness we suppress this from the notation. 

For example, we will prove the following property, which resolves the function field analogue of \cite[Problem 5]{Kud04}. 

\begin{thm}[Linear Invariance for special cycles]\label{thm: linear invariance}
Let $\cE$ be a rank $m$ vector bundle on $X'$ admitting a decomposition $\cE \cong  \cE_1 \oplus \cE_2 \oplus \ldots \op\cE_j$ with $\cE_i$ having rank $m_i$. Let $a_i \in \cA_{\cE_i}(k)$, and $a \in \cA_{\cE}(k)$ whose restriction to $\cE_{i}$ is $a_{i}$, i.e., \begin{equation}\label{eq: condition on a wrt a_i}
\text{ the composition $\cE_i \rightarrow \cE \xrightarrow{a} \sigma^* \cE^{\vee} \rightarrow \sigma^*  \cE_i^{\vee}$ is $a_i$ for each $1 \leq i \leq j$.}
\end{equation}
Then we have an equality of cycle classes in $\Ch_{r(n-m)}(\cZ_{\cE}^r(a))$,
\[
\left( [\cZ_{\cE_1}^r(a_{1})] \cdot_{\Sht_{U(n)}^r} [\cZ_{\cE_2}^r(a_{2})]  \cdot_{\Sht_{U(n)}^r}  \ldots  \cdot_{\Sht_{U(n)}^r}  [\cZ_{\cE_j}^r(a_{j})]  \right) |_{\cZ_{\cE}^r(a)}= [\cZ_{\cE}^r(a)]
\]
where $(\cdot) |_{\cZ_{\cE}^r(a)}$ denotes the projection to the corresponding component of the open-closed decomposition
\[
\cZ_{\cE_1}^r(a_{1}) \stackrel{\mrm{cl}}\times_{\Sht_{U(n)}^r} \cZ_{\cE_2}^r(a_{2}) \stackrel{\mrm{cl}}\times_{\Sht_{U(n)}^r}  \ldots  \stackrel{\mrm{cl}}\times_{\Sht_{U(n)}^r}  \cZ_{\cE_j}^r(a_{j})  = \coprod_{a \text{ satisfying } \eqref{eq: condition on a wrt a_i}} \cZ_{\cE}^r(a). 
\]
\end{thm}

We call this property ``Linear Invariance'' in analogy to a result of Howard \cite{How19}, which could be viewed as a mixed characteristic local analogue of the special case where $\cE_1, \ldots, \cE_j$ are all line bundles. It is closely related to functoriality in the ``$H_1$-variable'' (cf. \S \ref{ssec: H_1} below). Functoriality in the ``$H_2$-variable'', explained in \S \ref{ssec: H_2 functoriality} below, will also be used in the next part in order to compute numerical evidence for modularity.

\subsection{Functoriality in $H_{1}$}\label{ssec: H_1}

Let $BH_1 \xrightarrow{\phi} BH_1'$ be a map of smooth $m$-framed gerbes over $X$, and $BH_2$ be an $n$-framed gerbes over $X$. Then we have a map of the corresponding derived Hitchin stacks $\sM_{H_1, H_2} \rightarrow \sM_{H_1', H_2}$, whose classical truncation is $\cM_{H_1, H_2} \rightarrow \cM_{H_1', H_2}$. Assume further that $BH_2$ is of unitary type or $B\GL(n)'$, with the standard map to $B \GL(n)'$. Then we get induced maps $\Hk^r_{\sM_{H_1, H_2}} \rightarrow \Hk^r_{\sM_{H_1', H_2}}$, $\Sht^r_{\sM_{H_1, H_2}} \rightarrow	 \Sht_{\sM_{H_1', H_2}}^r$, etc.

\begin{lemma}\label{lem: hitch funct in H_1}
In the situation above, $\phi$ induces isomorphisms (where fibered products are derived)
\begin{enumerate}
\item $\sM_{H_1, H_2} \cong \sM_{H_1', H_2} \times_{\Bun_{H_1'}} \Bun_{H_1}  $, and 
\item $\Hk^r_{\sM_{H_1, H_2}}  \cong \Hk^r_{\sM_{H_1', H_2}} \times_{\Bun_{H_1'}} \Bun_{H_1}$. 
\end{enumerate}
\end{lemma}

\begin{proof}
Immediate from the definitions. 
\end{proof}

\begin{prop}\label{prop: sht_M funct in H_1}
In the situation above, $\phi$ induces an isomorphism of derived stacks
\[
\Sht_{\sM_{H_1, H_2}}^r \cong \Sht_{\sM_{H_1', H_2}}^r \times_{\Bun_{H_1'}(k)} \Bun_{H_1}(k),
\]
so that 
\[
[\Sht_{\sM_{H_1, H_2}}^r] = [ \Sht_{\sM_{H_1', H_2}}^r] \cdot_{\Bun_{H_1'}(k)} [\Bun_{H_1}(k)] \in \Ch_*(\Sht_{\cM_{H_1, H_2}}^r).
\]
\end{prop}

\begin{proof}

Abbreviate $\sM := \sM_{H_1,H_2}$ and $\sM' := \sM_{H_1', H_2}$. Consider the commutative diagram below.
\begin{equation}\label{eq: 9-term diag H_1}
\begin{tikzcd}
\sM' \ar[r, "{(\Id, \Frob)}"] \ar[d]  & \sM' \times \sM'  \ar[d] & \Hk_{\sM'}^r \ar[l,"{(\pr_{0},\pr_{r})}"'] \ar[d] \\
\Bun_{H_1'} \ar[r, "{(\Id, \Frob)}"]  &\Bun_{H_1'}  \times \Bun_{H_1'} & \Bun_{H_1'} \ar[l, "\Delta"'] \\
\Bun_{H_1} \ar[r, "{(\Id, \Frob)}"]  \ar[u] & \Bun_{H_1}   \times \Bun_{H_1}  \ar[u] & \Bun_{H_1}  \ar[l, "\Delta"'] \ar[u]   
\end{tikzcd}
\end{equation}
The derived fibered products along the rows of \eqref{eq: 9-term diag H_1} are (using Lemma \ref{lem: frob fixed points})
\begin{equation}\label{eq: H_1 rows}
\begin{tikzcd}
\Sht_{\sM'}^{r} \ar[d]  \\
\Bun_{H_1'}(k)  \\
\Bun_{H_1}(k) \ar[u]
\end{tikzcd}
\end{equation}
Each term is quasi-smooth by Lemma \ref{lem: cotangent complex for pi}, and moreover $\Bun_{H_1'}(k)$ and $\Bun_{H_1}(k) $ are smooth. 

Using Lemma \ref{lem: hitch funct in H_1}, we compute that the derived fibered products along the columns of \eqref{eq: 9-term diag H_1} are 
\begin{equation}\label{eq: H_1 columns}
\begin{tikzcd}
\sM  \ar[r] \ar[r, "{(\Id, \Frob)}"] & \sM  \times \sM   & \Hk_{\sM}^{r} \ar[l,"{(\pr_{0},\pr_{r})}"'] 
\end{tikzcd}
\end{equation}

The same proof as for \cite[Lemma A.9]{YZ} gives canonical isomorphisms of derived stacks between the derived fibered products of \eqref{eq: H_1 rows} and \eqref{eq: H_1 columns}. The derived fibered product of \eqref{eq: H_1 columns} is $\Sht_{\sM}^r$. We then conclude by applying \eqref{eq: refined intersection product} to \eqref{eq: H_1 rows}.
\end{proof}

\subsection{Functoriality in $H_{2}$}\label{ssec: H_2 functoriality} 
Let $BH_1$ be an $m$-framed gerbe, and $BH_2 \xrightarrow{\phi} BH_2'$ be a map of smooth $n$-framed gerbes over $X$. Then we have a map of the corresponding derived Hitchin stacks $\sM_{H_1, H_2} \rightarrow \sM_{H_1, H_2'}$, whose classical truncation is $\cM_{H_1, H_2} \rightarrow \cM_{H_1, H_2'}$. Assume further that $BH_2$ and $BH_2'$ are unitary type or $B\GL(n)'$, with the standard map to $B \GL(n)'$. Then these induce $\Hk^r_{\sM_{H_1, H_2}} \rightarrow \Hk^r_{\sM_{H_1, H_2'}}$, $\Sht^r_{\sM_{H_1, H_2}} \rightarrow	 \Sht_{\sM_{H_1, H_2'}}^r$, etc.

\begin{lemma}\label{lem: hitch funct in H_2}
In the situation above, $\phi$ induces isomorphisms (where fibered products are derived)
\begin{enumerate}
\item $\sM_{H_1, H_2} \cong \sM_{H_1, H_2'} \times_{\Bun_{H_2'}} \Bun_{H_2}  $, and 
\item $\Hk^r_{\sM_{H_1, H_2}}  \cong \Hk^r_{\sM_{H_1, H_2'} } \times_{\Hk_{H_2'}^r} \Hk^r_{H_2}$. 
\end{enumerate}
\end{lemma}

\begin{proof}
Immediate from the definitions. 
\end{proof}

\begin{prop}\label{prop: sht_M funct in H_2}
Then $\phi$ induces an isomorphism of derived stacks 
\[
\Sht_{\sM_{H_1, H_2}}^r \cong \Sht_{\sM_{H_1, H_2'}}^r \times_{\Sht_{H_2'}^r} \Sht_{H_2}^r,
\]
(with the RHS a derived fibered product), so that 
\[
[\Sht_{\sM_{H_1, H_2}}^r] = [ \Sht_{\sM_{H_1, H_2'}}^r] \cdot_{\Sht_{H_2'}^r} [\Sht_{H_2}^r] \in \Ch_*(\Sht_{\cM_{H_1, H_2}}^r).
\]
\end{prop}

\begin{proof}
Abbreviate $\sM := \sM_{H_1,H_2}$ and $\sM' := \sM_{H_1, H_2'}$. Consider the commutative diagram below.
\begin{equation}\label{eq: 9-term diag 2}
\begin{tikzcd}
\sM' \ar[r, "{(\Id, \Frob)}"] \ar[d]  & \sM' \times \sM'  \ar[d] & \Hk_{\sM'}^r \ar[l,"{(\pr_{0},\pr_{r})}"'] \ar[d] \\
\Bun_{H_2'} \ar[r, "{(\Id, \Frob)}"]  &\Bun_{H_2'}  \times \Bun_{H_2'} & \Hk^r_{H_2'} \ar[l,"{(\pr_{0},\pr_{r})}"'] \\
\Bun_{H_2} \ar[r, "{(\Id, \Frob)}"]  \ar[u] & \Bun_{H_2}   \times \Bun_{H_2}  \ar[u] & \Hk^r_{H_2} \ar[l, "{(\pr_{0},\pr_{r})}"'] \ar[u]   
\end{tikzcd}
\end{equation}
The derived fibered products along the rows of \eqref{eq: 9-term diag 2} are
\begin{equation}\label{eq: H_2 rows}
\begin{tikzcd}
\Sht_{\sM'}^{r} \ar[d]  \\
\Sht_{H_2'}^r  \\
\Sht_{H_2}^r \ar[u]
\end{tikzcd}
\end{equation}
Each term is quasi-smooth by Lemma \ref{lem: cotangent complex for pi}, and moreover $\Sht_{H_2'}^r$ and $\Sht_{H_2}^r$ are smooth. 

Using Lemma \ref{lem: hitch funct in H_1}, we compute that the derived fibered products along the columns of \eqref{eq: 9-term diag 2} are 
\begin{equation}\label{eq: H_2 columns}
\begin{tikzcd}
\sM  \ar[r] \ar[r, "{(\Id, \Frob)}"] & \sM  \times \sM   & \Hk_{\sM}^{r} \ar[l,"{(\pr_{0},\pr_{r})}"'] 
\end{tikzcd}
\end{equation}

The same proof as for \cite[Lemma A.9]{YZ} gives canonical isomorphisms of derived stacks between the derived fibered products of \eqref{eq: H_2 rows} and \eqref{eq: H_2 columns}. The derived fibered product of \eqref{eq: H_2 columns} is $\Sht_{\sM}^r$. We then conclude by applying \eqref{eq: refined intersection product} to \eqref{eq: H_2 rows}.
\end{proof}

\begin{example}\label{ex: CM pullback special cycle}
Consider the situation of Example \ref{ex:CM}, with $Y$ be another smooth projective curve over $\F_{q}$, and $\th: Y\to X$ be a map of degree $n$, possibly ramified. Let $\th': Y'\to X'$ (resp. $\nu':Y'\to Y$) be the base change of $\theta$ (resp. $\nu$).

We define the moduli of shtukas $\Sht_{U(1)/Y, \theta^*\LL}^r$ to be $\Sht_{H_2}^r $ for $BH_2 = BR_{Y/X}U(1)_{\theta^* \LL}$ (defined in \S \ref{sssec: hecke for unitary gerbe}). This definition is consistent with our previous definitions when $Y$ is connected, but previously we did not cover the case where $Y$ is disconnected. 


Take $BH_1$ to be the tautological $m$-framed gerbe and $BH_2  \rightarrow BH_2' = U(n)_{\LL}$ the natural map. Proposition \ref{prop: sht_M funct in H_2} implies that 
\[
[\Sht_{\sM_{H_1, H_2}}^r]  = [\Sht_{\sM_{H_1, H_2'}}^r] \cdot_{\Sht_{H_2'}^r}  [\Sht_{H_2}^r] \in \Ch_*(\Sht_{\cM_{H_1, H_2}}^r).
\]
The right hand side is, by Theorem \ref{thm: VFC from derived shtuka}, 
\[
\bigoplus_{\cE \in \Bun_{\GL(m)'}(k)} [\ol{\cZ}_{\cE,\LL}^r] \cdot_{\Sht_{U(n),\LL}^r} [\Sht_{U(1)/Y,\theta^*\LL}^r]
\]
and the left hand side is 
\[
\bigoplus_{\cE \in \Bun_{\GL(m)'}(k)  } [\ol{\cZ}_{\th'^* \cE,\theta^*\LL}^r],
\]where the summands are the special cycles defined relative to $Y'/Y$. 
In particular, projecting to the component indexed by $(\cE, a)$ (where $a \in \cA_{\cE,\LL}(k)$), and pulling back along $\pt \rightarrow B( \Aut(\cE)(\F_q))$, yields  
\[
 \bigoplus_{\substack{\wt{a} \in \cA_{\th'^* \cE, \theta^*\LL}(k) \\ \tr(\wt{a}) = a}} [\cZ_{\th'^* \cE,\theta^*\LL}^r(\wt{a})] = [\cZ_{\cE,\LL}^r(a)] \cdot_{\Sht_{U(n),\LL}^r} [\Sht_{U(1)/Y,\theta^*\LL}^r].
\]

Here the trace map is defined as follows.
Recall that  $\cA_{\cE,\LL}(k)$ is the set of Hermitian maps $a:\cE\to \s^*\cE^\vee\ot  \nu^*\LL$, i.e., $\sigma$-invariant elements in $\Hom_{\cO_{X'}}(\cE, \s^*\cE^\vee\ot \nu^\ast\LL )$. Having defined $ (\th'^*\cE)^\vee= (\th'^*\cE)^{*}\ot \omega_{Y'}$, we have natural isomorphisms
\begin{eqnarray*}
&& \Hom_{\cO_{Y'}}(\th'^*\cE, \s^* (\th'^*\cE)^\vee\ot \nu'^\ast( \th^*\LL))\\
 &\cong &\Hom_{\cO_{Y'}}(\th'^*\cE, \th'^* (\s^*\cE^*\ot \nu^* \LL) \ot \om_{Y'}) \\
 &\cong &\Hom_{\cO_{X'}}(\cE,\s^*\cE^*\ot \nu^* \LL\ot\th'_\ast\om_{Y'})  \quad \text{(By adjunction)}.
\end{eqnarray*}
Post-composition with the trace map
$$
\tr_{Y'/X'}\co \th_\ast\om_{Y'}\to \om_{X'}
$$
defines a map
$$
\xymatrix{\Hom_{\cO_{X'}}(\cE,\s^*\cE^*\ot \nu^* \LL\ot\th'_\ast\om_{Y'}) \ar[r]&\Hom_{\cO_{X'}}(\cE,\s^*\cE^*\ot \nu^* \LL\ot\om_{X'}) }
$$
and hence a trace map
\begin{align}\label{eq:tr Hom}
\xymatrix{\tr\colon \Hom_{\cO_{Y'}}(\th'^*\cE, \s^* (\th'^*\cE)^\vee\ot \nu'^\ast( \th^*\LL))\ar[r]&\Hom_{\cO_{X'}}(\cE, \s^*\cE^\vee\ot \nu^\ast\LL )  }
\end{align} 
It is easy to see that the map \eqref{eq:tr Hom} preserves Hermitian elements and therefore defines the desired trace map
\begin{align}\label{eq:tr A}
\xymatrix{\tr\colon \cA _{\th'^*\cE,  \th^*\LL}(k)\ar[r]&\cA_{\cE, \LL}(k) }.
\end{align} 

\end{example}

\subsection{Proof of Linear Invariance} 
\label{ss:pf LI}
We will work up to the proof of Theorem \ref{thm: linear invariance} with two intermediate steps. 


\begin{lemma}\label{lem: octahedron equality 2}
Fix $m \leq n$ and $m = m_1 + \ldots + m_j$. Let $BH_1^{(1)} = B\GL(m_1 )',  BH_1^{(2)} = B\GL(m_2)', \ldots, BH_1^{(j)} =  B\GL(m_j)'$, $BH_1 = BH_1^{(1)} \times \cdots \times BH_1^{(j)}$, and $BH_2 = BU(n)$. Define the derived Hitchin stacks $\sM_{H_1^{(i)}, H_2}$ using the identity map $BH_1^{(i)} \xrightarrow{=} B\GL(m_i)'$ and the standard map $BH_2 \rightarrow B\GL(n)'$, and $\sM_{H_1, H_2}$ using the standard block diagonal map $BH_1 \rightarrow B\GL(m)'$ and the standard map $BH_2 \rightarrow B\GL(n)'$. 

Then we have the following equality in $\Ch_{r(n-m)}(\Sht_{\cM_{H_1, H_2}}^{r})$: 
\begin{equation}\label{eq: octahedron equality 2}
 [\Sht_{\sM_{H_1^{(1)}, H_2}}^{r}] \cdot_{\Sht_{U(n)}^{r}} \cdots\cdot_{\Sht_{U(n)}^{r}} [\Sht_{\sM_{H_1^{(j)}, H_2}}^{r}]   = [\Sht_{\sM_{H_1, H_2}}^{r}].
\end{equation}
\end{lemma}

\begin{proof} For $1 \leq i \leq j$ we abbreviate $\sM^{(i)} := \sM_{H_1^{(i)}, H_2}$. Consider the diagram below. 
\begin{equation}\label{eq: 9-term diag 1}
\begin{tikzcd}
\sM^{(1)}  \times \cdots \times \sM^{(j)} \ar[r, "{(\Id, \Frob)^{j}}"] \ar[d]  & (\sM^{(1)}  \times \sM^{(1)})\times \cdots \times (\sM^{(j)} \times \sM^{(j)})\ar[d] & \Hk_{\sM^{(1)}}^r  \times \cdots \times \Hk_{\sM^{(j)}}^r \ar[d]  \ar[l, "{(\pr_{0}, \pr_{r})^{j}}"']  \\
(\Bun_{U(n)})^j \ar[r, "{(\Id, \Frob)^{j}}"]  & (\Bun_{U(n)}\times \Bun_{U(n)})^j & (\Hk_{U(n)}^r)^j  \ar[l, "{(\pr_{0}, \pr_{r})^{j}}"']  \\
\Bun_{U(n)} \ar[r, "{(\Id, \Frob)}"] \ar[u, "\Delta_{\Bun_{U(n)}}"]  & \Bun_{U(n)} \times \Bun_{U(n)}  \ar[u, "{\Delta_{\Bun_{U(n)}\times\Bun_{U(n)}}}"]  & \Hk_{U(n)}^r  \ar[u, "\Delta"] \ar[l, "{(\pr_{0}, \pr_{r})}"'] 
\end{tikzcd}
\end{equation}

Using Lemma \ref{lem: shtuka derived cartesian}, we see that the derived fibered products along the rows of \eqref{eq: 9-term diag 1} are
\begin{equation}\label{eq: lem2 rows}
\begin{tikzcd}
\Sht_{\sM^{(1)}}^r  \times \cdots \times \Sht_{\sM^{(j)}}^r \ar[d]  \\
(\Sht_{U(n)}^{r})^j \\
\Sht_{U(n)}^{r} \ar[u, "\Delta"]  
\end{tikzcd}
\end{equation}
Each term is quasi-smooth by Lemma \ref{lem: cotangent complex for pi}, and $(\Sht_{U(n)}^{r})^j $ and $\Sht_{U(n)}^{r}$ are smooth by \cite[Lemma 6.9(2)]{FYZ}.

The derived fibered products along the columns of \eqref{eq: 9-term diag 1} are 
\begin{equation}\label{eq: lem2 columns}
\begin{tikzcd}
\sM  \ar[r, "{(\Id, \Frob)}"] & 
\sM \times 
\sM & \Hk_{\sM}^{r} \ar[l,"{(\pr_{0}, \pr_{r})}"']  
\end{tikzcd}
\end{equation}
where we abbreviated $\sM := \sM_{H_1, H_2}$. 

The same proof as for \cite[Lemma A.9]{YZ} gives canonical isomorphisms of derived stacks between the derived fibered products of \eqref{eq: lem2 rows} and \eqref{eq: lem2 columns}. The derived fibered product of \eqref{eq: lem2 rows} is $\Sht_{\sM}^r$, which is quasi-smooth by Lemma \ref{lem: cotangent complex for pi}. We then conclude by applying \eqref{eq: refined intersection product} to \eqref{eq: lem2 rows}.

\end{proof}

\begin{proof}[Proof of Theorem \ref{thm: linear invariance}]
We have by definition
\begin{equation}\label{eq: LI 1}
[\Sht_{\sM_{\GL(m)', U(n)}}^r] = \bigoplus_{\cE \in \Bun_{\GL(m)'}(k)} [\ol{\sZ}_{\cE}^r] \in \Ch_{r(n-m)}(\Sht^{r}_{\cM_{\GL(m)',U(n)}}),
\end{equation}
where the virtual fundamental classes are defined because $\Sht_{\sM_{\GL(m)', U(n)}}^r$ is quasi-smooth. 

Similarly, we have for each $i=1, \ldots, j$ that 
\begin{equation}
[\Sht_{\sM_{\GL(m_i)', U(n)}}^r] = \bigoplus_{\cE_i \in \Bun_{\GL(m_i)'}(k)} [\ol{\sZ}_{\cE_i}^r].
\end{equation}

Let $H_{1}$ be the subgroup $\GL(m_{1})'\times\cdots\times\GL(m_{j})'$ of $\GL(m)'$ as in the hypotheses of Lemma \ref{lem: octahedron equality 2}. By Lemma \ref{lem: octahedron equality 2}, we then have 
\begin{equation}\label{eq: LI 2}
[\Sht_{\sM_{H_{1}, U(n)}}^r] = \bigoplus_{\cE_i \in \Bun_{\GL(m_i)'}(k), i=1, \ldots, j} [\ol{\sZ}_{\cE_1}^r] \cdot_{\Sht_{U(n)}^r}   \cdots  \cdot_{\Sht_{U(n)}^r}  [\ol{\sZ}_{\cE_j}^r]  .
\end{equation}


Applying Proposition \ref{prop: sht_M funct in H_1} to the inclusion $H_{1}\inj H_{1}'=\GL(m')$ we get
\begin{equation}\label{eq: octahedron equality 3}
[\Sht_{\sM_{H_1, U(n)}}^{r}]  = [\Sht_{\sM_{\GL(m)', U(n)}}^{r}]  \cdot_{\Bun_{\GL(m)'}(k)} [\Bun_{H_{1}}(k)]\in \Ch_{r(n-m)}(\Sht^{r}_{\cM_{H_{1},U(n)}}).
\end{equation}
Projecting the above equality to the component indexed by $(\cE_1, \ldots, \cE_j)\in \Bun_{H_{1}}(k)$ and using \eqref{eq: LI 1} and \eqref{eq: LI 2} yields
\begin{equation}\label{eq: LI 3}
[\ol{\sZ}_{\cE_1}^r] \cdot_{\Sht_{U(n)}^r}   \cdots  \cdot_{\Sht_{U(n)}^r}  [\ol{\sZ}_{\cE_j}^r] = [\ol{\sZ}_{\cE}^r].
\end{equation}
By Theorem \ref{thm: VFC from derived shtuka}, we have
\[
[\ol{\cZ}_{\cE}^r]  = [\ol{\sZ}_{\cE}^r] \in \Ch_{r(n-m)}(\ol{\cZ}_{\cE}^r),
\]
and similarly (using the compatibility of \S \ref{sssec: intersection product} and the refined intersection product) we have 
\begin{align*}
 [\ol{\cZ}_{\cE_1}^r] \cdot_{\Sht_{U(n)}^r} \cdots  \cdot_{\Sht_{U(n)}^r}  [\ol{\cZ}_{\cE_j}^r]  =  [\ol{\sZ}_{\cE_1}^r] \cdot_{\Sht_{U(n)}^r}  \cdots  \cdot_{\Sht_{U(n)}^r}  [\ol{\sZ}_{\cE_j}^r].
\end{align*}
Putting these equalities into \eqref{eq: LI 3}, pulling back along $\pt \rightarrow B \Aut(\cE_1) \times \cdots \times B \Aut(\cE_j)$, and then projecting to the component indexed by $a$ gives the result. 

\end{proof}

\section{Compatibility with the cycle classes of \cite{FYZ}}\label{sec: compatibility}

Let $\cE$ be a rank $m$ vector bundle on $X'$ and $a \in \cA_{\cE}^{\ns}(k)$ (i.e., $a$ is \emph{non-singular}). In \cite{FYZ} we gave a different definition of the virtual fundamental class $[\cZ_{\cE}^r(a)]$ in the case either $\cE$ is a direct sum of line bundles and or $\rank\cE=n$. The goal of this section is to prove that these cycles defined in \cite{FYZ} agree with the same-named cycles defined in Definition \ref{defn: circ case}. Although \cite{FYZ} was written with the twisting line bundle $\LL$ being trivial, a completely analogous construction applies with any $\LL$. We shall fix a choice of $\LL$ throughout this section and suppress it from the notation.

\subsection{Corank one special cycles}\label{ssec: corank one unitary} In this subsection, we show that the ``least degenerate strata'' of any corank one special cycle, either for $BH_{2}=BU(n)_\LL$ or $BH_{2}=B\GL(n)'$,  is LCI.

Throughout this subsection we abbreviate $\cM':=\cM_{\GL(1)', \GL(n)'}^{\c}$.

\begin{prop}\label{prop: KR is LCI} (1) The stack $\Sht^r_{\cM'}$ is LCI of pure dimension $r(2n-2)$. In particular, for any line bundle $\cL$ on $X'$, $\cZ^{r,\c}_{\cL, \GL(n)'}$ is LCI of pure dimension $r(2n-2)$. 

(2) The class $[\Sht^{r}_{\cM'}] \in \Ch_{r(2n-2)} ( \Sht^{r}_{\cM'})$ from Definition \ref{defn: hitchin shtuka} agrees with 
\[
\sum_{\Cal{L} \in \Bun_{\GL(1)'}(k)} [\ol{\Cal{Z}}_{\Cal{L}}^{r , \circ}]^{\nai} \in \Ch_{r(2n-2)} ( \Sht^{r}_{\cM'}),
\]
where $[\ol{\Cal{Z}}_{\Cal{L}}^{r, \circ}]^{\nai} \in \Ch_{r(2n-2)} ( \Sht^{r}_{\cM'})$ is the fundamental class of that component.

\end{prop}
\begin{proof} (1) We may write $\Sht^{r}_{\cM'}$ by the Cartesian diagram
\begin{equation}\label{eq: shtuka cycle class} 
\begin{tikzcd}
 \Sht^{r}_{\cM'}\ar[r]\ar[d]      &  (\Hk^{1}_{\cM'} \times \cM')^{r}\ar[d, "{(p_{0},p_{1})^r}"] \\
{\cM'}^{r+1} \ar[r, "\Phi^{r}_{\cM'}"] &  {\cM'}^{2r+2}
\end{tikzcd}
\end{equation}

By the smoothness of $\cM'$ and $\Hk^{1}_{\cM'}$ and the relative dimension calculations in Proposition \ref{prop: Hitchin smooth} and Lemma \ref{l:HkM sm}, we see that $\Sht^{r}_{\cM'}$ has local dimension $\ge r(2n-2)$ everywhere. 

On the other hand, we will show in Proposition \ref{prop: Sht_M for GL(n)'}(3) and Corollary \ref{cor: Sht_M for GL(1)'} that $\dim\cZ^{r,\c}_{\cL, \GL(n)'}\le r(2n-2)$ for any line bundle $\cL$, where $\cZ^{r,\c}_{\cL, \GL(n)'}$ is the pullback of $\cZ^{r,\c}_{\cL, \GL(n)'}$ along $\pt \rightarrow B \Aut(\cL)$, hence $\dim \Sht^{r}_{\cM'}\le r(2n-2)$. Combining this with the lower bound of local dimension given above, and the fact that $\Sht^r_{\cM'}$ is a fibered product of smooth stacks, we conclude that $\Sht^{r}_{\cM'}$ is LCI of pure dimension $r(2n-2)$.

(2) We have seen that $\Sht^{r}_{\cM'}$ is LCI  and the fibered product in \eqref{eq: shtuka cycle class} exhibits it as a proper intersection, so the claim follows from \cite[Proposition 7.1]{Ful98}. (Strictly speaking, the statement in \emph{loc. cit.} is for schemes, so we apply it after truncating and adding sufficient level structure, and then taking a limit over truncations. When adding level structure along a finite subscheme $D \subset X$, we ask that the leg maps avoid $D$, so this lies over an open substack of $\Sht_{U(n)}^r$. As $D$ varies, these substacks form an open cover as $D$ varies. The equality in question can be checked on this open cover because it is an equality of \emph{top-dimensional} cycles. We omit the details of this step because they require heavy notation, and yet are not very interesting.) 
\end{proof}

We may now establish a result that was promised in \cite[Remark 7.10]{FYZ}.

\begin{cor}\label{cor: proper intersection} 
The Cartesian square 
\[
\begin{tikzcd} \Cal{Z}_{\Cal{L}}^{r, \circ }    \ar[r] \ar[d] & 
\Cal{Z}_{ \Cal{L}, \GL(n)'}^{r ,\circ }  \ar[d]  \\
\Sht_{U(n)}^r \ar[r] & \Sht_{\GL(n)'}^r
\end{tikzcd}
\]
is a proper intersection. Hence for all $a \in \cA_{\cL}(k)$ \footnote{Note our $\cA_{\cL}(k)$ is denoted $\cA_{\cL}^{\all}(k)$ in \cite[\S 7]{FYZ}, i.e., it includes singular $a$.} we have that $\cZ_{\cL}^r(a)^\circ$ is LCI of dimension $r(n-1)$. 

In particular, $[\cZ_{\cL}^r(a)^\circ]^{\nai} = [\cZ_{\cL}^r(a)^\circ] \in \Ch_*( \Cal{Z}_{\Cal{L}}^{r, \circ } )$, the latter being the projection of Definition \ref{defn: circ case} to the closed-open substack indexed by $a$. 
\end{cor}

\begin{proof}
Lemma \ref{lem: special cycle GL(n)' to U(n)} implies that $\coprod_{a \in \cA(k)} \cZ_{\cL}^r(a)^\circ $ is the fibered product of $ \cZ_{\cL, \GL(n)'}^{r , \circ} $ and $\Sht_{U(n)}^r$ over $\Sht_{\GL(n)'}^r$. These have dimensions $r(2n-2)$, $rn$, and $r(2n-1)$ respectively, as established in Proposition  \ref{prop: KR is LCI}, \cite[Lemma 6.9(2)]{FYZ}, and Lemma \ref{lem: sht dimension} respectively. Since $\Sht_{U(n)}^r \rightarrow \Sht_{\GL(n)'}^r$ is a regular local immersion of smooth Deligne-Mumford stacks, this implies that the fibered product has dimension $\geq r(n-1)$.  

On the other hand, it was already established in \cite[Proposition 9.1, 9.5]{FYZ} that $\dim \cZ_{\cL}^r(a)^\circ \leq r(n-1)$, so equality holds. As this presentation realizes $\cZ_{\cL}^{r,\circ}$ as the pullback of the LCI Deligne-Mumford stack $\cZ_{\cL, \GL(n)'}^{r,\circ}$ against a regular local immersion, we conclude that $\cZ_{\cL}^{r,\circ}$ is also LCI. 

For the last statement, we use that $\Sht_{U(n)}^r \rightarrow \Sht_{\GL(n)'}^r$ is a regular local immersion (as both are smooth) and Proposition \ref{prop: KR is LCI}. (Strictly speaking, we need to truncate and add level structure as discussed in the proof of Proposition \ref{prop: KR is LCI}.) 
\end{proof}




\subsection{Agreement of definitions}
Let $\cE$ be a rank $m$ vector bundle on $X'$. In \cite{FYZ} we gave a different definition of the virtual fundamental class $[\cZ_{\cE}^r(a)]$ when $a \in \cA_{\cE}^{\ns}(k)$ (i.e., $a$ is \emph{non-singular}), in the following cases: 
\begin{enumerate}
\item $\cE \approx \cL_1 \oplus \ldots \oplus \cL_m$ is a direct sum of line bundles on $X'$ \cite[\S 7.8]{FYZ}, or 
\item $\rank \cE = n$ \cite[\S 7.9]{FYZ}. 
\end{enumerate} 
We denote the class defined in \cite{FYZ} by $[\cZ_{\cE}^r(a)]^{\old}$.  In this section we prove that $[\cZ_{\cE}^r(a)]^{\old}$ agrees with the class $[\cZ_{\cE}^r(a)]$ defined in Definition \ref{defn: circ case}.

\begin{prop}\label{prop: compatibilitiy}
For $\cE$ as in Cases (1) or (2) above and $a \in \cA_{\cE}^{\ns}(k)$, we have 
\[
[\cZ_{\cE}^r(a)]^{\old} = [\cZ_{\cE}^r(a)] \in \Ch_{r(n-m)}(\cZ_{\cE}^r(a)).
\]
\end{prop}

\begin{proof} (1) If $\cE = \cL_1 \oplus \ldots \oplus \cL_m$, \cite[\S 7.8]{FYZ} defined $[\cZ_{\cE}^r(a)]^{\old}$ to be the projection of $[\cZ_{\cL_1}^r(a_1)^{\c}]^{\nai} \cdot_{\Sht_{U(n)}^r} \ldots \cdot_{\Sht_{U(n)}^r}  [\cZ_{\cL_m}^r(a_m)^{\c}]^{\nai}$ to the components indexed by $a$ (where $a_{i}\in\cA_{\cL_{i}(k)}$ is the restriction of $a$ to $\cL_{i}$). By Theorem \ref{thm: linear invariance}, $[\cZ_{\cE}^r(a)]$ is described in the same way with respect to $[\cZ_{\cL_1}^r(a_1)] \cdot_{\Sht_{U(n)}^r} \ldots \cdot_{\Sht_{U(n)}^r}  [\cZ_{\cL_m}^r(a_m)]$. Moreover, since $a$ is non-singular, $\cZ_{\cE}^r(a)$ is contained in the fiber product of $\cZ_{\cL_i}^r(a_i)^{\c}$ over $\Sht^{r}_{U(n)}$, hence $[\cZ_{\cE}^r(a)]$ is the projection of $[\cZ_{\cL_1}^r(a_1)^{\c}] \cdot_{\Sht_{U(n)}^r} \ldots \cdot_{\Sht_{U(n)}^r}  [\cZ_{\cL_m}^r(a_m)^{\c}]$ to the components indexed by $a$. So we are reduced to showing that $[\cZ_{\cL_{i}}(a_{i})^{\c}]^{\nai}=[\cZ_{\cL_{i}}^r(a_i)^{\c}]$ (here we allow $a_{i}=0$). This follows from Corollary \ref{cor: proper intersection}.


(2) Suppose $\rank \cE  = n$. Take $BH_1 = B\GL(n)'$, $BH_2 = BU(n)_{\LL}$, and form $\sM^{\ns} := \sM_{H_1, H_2}^{\c}|_{\sA_{H_1}^{\ns}}$ with respect to the identity map $H_1 \xrightarrow{=} \GL(n)'$ and the standard map $BH_2 \rightarrow B \GL(n)'$. Let $\cM^{\ns}$ be the classical truncation of $\sM^{\ns}$. 

Recall that there is a finite \'{e}tale map 
\[
\coprod_{\cE \in \Bun_{\GL(n)'}(k)} \coprod_{a \in \cA_{\cE}^{\ns}(k)} \cZ_{\cE}^r(a) \rightarrow \Sht_{\cM^{\ns}}^r
\]
and \cite[Theorem 10.1]{FYZ} establishes that $[\Sht_{\cM^{\ns}}^r]|_{\cZ_{\cE}^r(a)}$, defined as in \cite[Definition 8.16]{FYZ} using the (classical) Gysin pullback, agrees with
$[\cZ_{\cE}^r(a)]^{\old}$. Hence it suffices to show that $[\Sht_{\sM^{\ns}}^r] =[\Sht_{\cM^{\ns}}^r] \in \Ch_{0}(\Sht_{\cM^{\ns}}^r )$.


Rewriting $\Hk_{\sM^{\ns}}^r$ as the derived fibered product $(\Hk_{\sM^{\ns}}^1 \times \sM^{\ns})^r \times_{\sM^{\ns,2r+2}} \sM^{\ns,r+1}$,  we see that $\Sht_{\sM^{\ns}}^r$ may be rewritten as the fibered product below, where $\Phi_{\sM^{\ns}}^r$ is as in Definition \ref{def:Phi}.
\begin{equation}\label{eq: Sht_M cartesian 2}
\begin{tikzcd}
\Sht_{\sM^{\ns}}^r \ar[r] \ar[d] & (\Hk_{\sM^{\ns}}^1)^r \times \sM^{\ns} \ar[d] \\
\sM^{\ns,r+1} \ar[r, "\Phi_{\sM^{\ns}}^{r}"] & \sM^{\ns,2r+2}
\end{tikzcd}
\end{equation}

By Corollaries \ref{cor: M already classical}(2) and \ref{cor: HkM already classical}(2), the canonical maps $\cM^{\ns} := (\sM^{\ns})_{\mrm{cl}} \rightarrow \sM^{\ns}$ and $\Hk_{\cM^{\ns}}^1 \rightarrow \Hk_{\sM^{\ns}}^1$ are isomorphisms of smooth stacks, so in particular $[\Hk_{\sM^{\ns}}^1] = [\Hk_{\cM^{\ns}}^1]^{\nai}$. Lemma \ref{lem: VFC fibered product} then implies that $[\Sht_{\sM^{\ns}}^r] =[ (\Hk_{\sM^{\ns}}^1)^r \times \sM^{\ns} ] \cdot_{\sM^{\ns,2r+2}} [ \sM^{\ns,r+1}]$, which is the same as $(\Phi^{r}_{\sM^{\ns}})^! [ (\Hk_{\sM^{\ns}}^1)^r \times \sM^{\ns} ]$. By Lemma \ref{lem: refined gysin compatibility}, $(\Phi^{r}_{\sM^{\ns}})^! [ (\Hk_{\sM^{\ns}}^1)^r  \times \sM^{\ns}]$ agrees with $[\Sht_{\cM^{\ns}}^r]$.


\end{proof}

\begin{remark}
Although the paper \cite{FYZ} focused on the case $\rank \cE = n$, the construction of $[\Sht_{\cM^{\ns}}^r]$ found there can be performed equally well for other values of $\rank \cE$, and the above proof shows that the $[\Sht_{\cM^{\ns}}^r]$ so defined agrees with $[\Sht_{\sM^{\ns}}^r ]$. 
\end{remark}

\subsection{Stratification of shtukas for $\Cal{M}_{\GL(1)', \GL(n)'}$}\label{sec: corank one}
The goal of this subsection is to show that $\dim\cZ_{\cL, \GL(n)'}^r\le r(2n-2)$ for line bundles $\cL$ on $X'$, as has been used in the proof of Proposition \ref{prop: KR is LCI}. The idea is similar to that of \cite[\S 9]{FYZ}. 

Fix a line bundle $\cL$ on $X'$. We define $\cZ' := (\cZ_{\cL, \GL(n)'}^r)_{\ol{k}}$. Let $I_{0}\sqcup I_{+}\sqcup I_{-} \sqcup I_{\pm} $ be a partition of $\{1,2,\cdots, r\}$. We denote this partition simply by $I_{\bu}$. For any $N\in \Z_{\ge0}$, define $\frD(N;I_{\bu})$ to be the moduli space of sequences of effective divisors $(D_{i})_{0\le i\le r}$ on $X'_{\ov k}$ such that
\begin{enumerate}
\item $\deg(D_{0})\le N$.
\item For $? \in \{0,+,-, \pm\}$, and $i\in I_{?}$, the pair $(D_{i-1},D_{i})$ belongs to the corresponding Case (?) beow
\begin{enumerate}
\item[(0)] $D_{i}=D_{i-1}$;
\item[(+)] $D_{i}=D_{i-1}+\s x'_{i}$ for some $x'_{i}\in X'_{\ov k}$;
\item[(--)] $D_{i}=D_{i-1}-x'_{i}$ for some $x'_{i}\in X'_{\ov k}$;
\item[($\pm$)] $D_{i}=D_{i-1}-x'_{i}+\s x'_{i}$ for some $x'_{i}\in X'_{\ov k}$.
\end{enumerate}
\item $D_{r}={}^{\t}D_{0}$. 
\end{enumerate}
For $?=+,-$ or $\pm$, $(D_{i-1},D_{i})$ determines a point $x'_{i}\in X_{\ov k}$. This gives a map. 
\[
(\pi_+,\pi_-, \pi_{\pm}) \co \frD(N, I_{\bu}) \rightarrow (X'_{\ol{k}})^{I_+\sqcup I_-\sqcup I_{\pm}}.
\]

\begin{lemma}\label{l:Dqf}
The map $\pi_{+}:\frD(N, I_{\bu}) \rightarrow (X'_{\ol{k}})^{I_+}$ is quasi-finite. In particular, $\dim \frD(N, I_{\bu})\le |I_{+}|$.
\end{lemma}
\begin{proof} Let $D_{\bu}\in \frD(N;I_{\bu})$ and $\ov D_{i}=\nu(D_{i})$ be the image of $D_{i}$ in $X_{\ov k}$. Let $x_{i}=\nu(x'_{i})$, then $\ov D_{i}=\ov D_{i-1}$ if $i\in I_{0}\sqcup I_{\pm}$, and $\ov D_{i}=\ov D_{i-1}+x_{i}$ if $i\in I_{+}$, $\ov D_{i}=\ov D_{i-1}-x_{i}$ if $i\in I_{-}$. By condition (3) above, $\ov D_{0}$ satisfies the equation
\begin{equation}\label{D0bar}
\ov D_{0}+\sum_{i\in I_{+}}x_{i}={}^{\t}\ov D_{0}+\sum_{j\in I_{-}}x_{j}.
\end{equation}
By  \cite[Lemma 9.4]{FYZ}, for fixed $\{x_{i}\}_{i\in I_{+}}$, there are only finitely many $\ov D_{\bu}$ satisfying \eqref{D0bar} and $\deg D_{0}\le N$. If  $\ov D_{\bu}$ is fixed then $D_{\bu}$ has finitely many choices. We conclude that there are finitely many $\ov k$-points in $\frD(N;I_{\bu})$ with fixed image in $(X'_{\ov k})^{I_{+}}$.
\end{proof}

For a partition $I_{\bu} = (I_0,  I_+, I_-, I_{\pm})$ of $\{1, 2, \ldots, r\}$, let $\Cal{Z}'[N, I_{\bu}] $ be the stack classifying $(\{D_{i}\}_{0\le i\le r}, \{x'_{i}\}_{1\le i\le r}, \{\cL\xr{t_{i}}\cF_{i}\}) $ where $ (\{x'_{i}\}_{1\le i\le r}, \{\cL\xr{t_{i}}\cF_{i}\})  \in \cZ'(S)$, and $\{D_{i}\}\in \frD(N;I_{\bu})$ with image $\{x'_{i}\}_{i\in I_{?}}$ under $\pi_{?}$ ($?=+,-,\pm$), and $t_{i}$ extends to a saturated embedding $\cL(D_{i})\inj \cF_{i}$. Since $D_i$ is determined by $t_i$, the natural map $\Cal{Z}'[N, I_{\bu}] \inj \cZ'$ is a locally closed immersion. As in \cite[\S 9.2.2]{FYZ}, we define the map 
\[
\pi'[N; I_{\bu} ] \co \Cal{Z}'[N, I_{\bu}] \rightarrow (X'_{\ol{k}})^{I_0} \times \frD(N; I_{\bu}).
\]

\begin{cor}[of Lemma \ref{l:Dqf}]\label{cor: Sht_M for GL(1)'} When $n=1$, $\dim\cZ'[N;I_{\bu}]=0$.
\end{cor}
\begin{proof}
When $n=1$,  $\cZ'[N;I_{\bu}]$ classifies $(\{D_{i}\}_{0\le i\le r}, \{x'_{i}\}_{1\le i\le r}, \{\cL\xr{t_{i}}\cF_{i}\})$ such that $t_{i}$ extends to an isomorphism $\cL(D_{i})\cong \cF_{i}$. This implies that $I_{\pm}=\{1,2,\cdots, r\}$, and the forgetful map $\cZ'[N;I_{\bu}]\to \frD(N;I_{\bu})$ is an isomorphism. By Lemma \ref{l:Dqf}, $\dim\cZ'[N;I_{\bu}]=\dim\frD(N;I_{\bu})=0$.
\end{proof}

\begin{prop}\label{prop: Sht_M for GL(n)'} Assume $n\ge2$.
\begin{enumerate}
\item For varying $N\in \Z_{\ge0}$ and partitions $I_{\bu}$ of $\{1,2,\cdots, r\}$ with $|I_+|=|I_-|$, the substacks  $\cZ'[N;I_{\bu}]$ give a stratification of $\cZ'$. 
\item The fibers of the map $\pi'[N;I_{\bu}]$ have dimension at most $(2n-3)|I_{0}| + (2n-2)|I_{+}| + (n-2) |I_-| +(n-1) |I_{\pm}|$.
\item We have $\dim \cZ'[N;I_{\bu}]\le r(2n-2)$. Moreover,  the equality is achieved only when $I_{0}=\{1,2,\cdots, r\}$, i.e., all $D_{i}$ are equal to the same divisor of $X'$ defined over $\F_{q}$.
\end{enumerate}
\end{prop}
\begin{proof}
(1) is clear (note that $|I_+| = |I_-|$ is implied by the assumption that $D_r = \ft D_0$). 

(2) The analysis is similar to that of \cite[Proposition 9.1]{FYZ}, although the cases behave differently, so let us explain how they play out. 

Fix $D_{\bu}\in \frD(N; I_{\bu})(\ov k)$, let $\cZ'[D_{\bu}]$ be the fiber of the projection $\cZ'[N;I_{\bu}]\to \frD(N; I_{\bu})$ over $D_{\bu}$. Let $\cM'=\cM_{\GL(1)', \GL(n)'}^{\c}$.

Let $\cH'[D_{\bu}]$ be the substack of $\Hk^{r}_{\cM', \ol{k}}$ classifying data $(x'_{\bu},\cL\xr{t_{i}}\cF_{i})$ such that $t_{i}$ extends to a map $t'_{i}: \cL(D_{i})\to \cF_{i}$. Note that for $i\notin I_{0}$, the $x'_{i}$ are determined by $D_{\bu}$. 
Let $\cM'[D_{i}]$ be the substack of $\cM'_{\ol{k}}$ classifying maps $t:\cL\to \cF$ that extend to a saturated map $t': \cL(D_{i})\to \cF$. Then we have a Cartesian diagram of stacks over $\ov k$
\begin{equation}
\xymatrix{\cZ'[D_{\bu}]\ar[d]\ar[r] & \cH'[D_{\bu}]\ar[d]^{(p_{0}, p_{r})}\\
\cM'[D_{0}]\ar[r]^-{(\Id,\Frob)} & \cM'[D_{0}]\times \cM'[D_{r}]}
\end{equation}
Note since $D_{r}={}^{\t}D_{0}$, the Frobenius morphism sends $\cM'[D_{0}]$ to $\cM'[D_{r}]$. We claim that the map
\begin{equation}
\Pi[D_{\bu}]: \cH'[D_{\bu}]\to \cM'[D_{r}]\times {X'_{\ol{k}}}^{I_{0}}
\end{equation}
is smooth of relative dimension $(2n-3)|I_{0}| + (2n-2)|I_{+}| + (n-2) |I_-| +(n-1) |I_{\pm}|$.  Then by \cite[Lemma 9.3]{FYZ}, the fibers of $\cZ'[D_{\bu}]\to (X'_{\ol{k}})^{I_0} $, which are fibers of $\pi'[N;I_{\bu}]$, have  dimension $\le (2n-3)|I_{0}| + (2n-2)|I_{+}| + (n-2) |I_-| +(n-1) |I_{\pm}|$.

For $0\le j\le r$, let $\cH'_{\ge j}$ be the moduli stack defined similarly to $\cH'[D_{\bu}]$ but classifying only $x'_{\bu}\in {X'_{\ol{k}}}^{I_0}$ and saturated maps $\{t_{i}: \cL(D_{i}) \to \cF_{i}\}_{j\le i\le r}$ (and $\cF_{i}$ and $\cF_{i+1}$ are still related to each other by elementary modifications at $x'_{i+1}$ for $j\le i<r$). We can factorize $\Pi[D_{\bu}]$ as
\begin{equation}
\Pi[D_{\bu}]: \cH'[D_{\bu}]=\cH'_{\ge0} \xr{ \Pi_{1}} \cH'_{\ge1}\xr{\Pi_{2} }\cdots \xr{ \Pi_{r}}\cH'_{\ge r}=\cM'[D_{r}]\times {X'_{\ol{k}}}^{I_0}.
\end{equation}
The smoothness claim follows after we establish the following four statements:
\begin{enumerate}
\item[($H0$)] If $i\in I_{0}$, then $\Pi_{i}$ exhibits $\cH'_{\ge i-1}$ as an open substack in a $\PP^{n-1}$-bundle over a $\PP^{n-2}$-bundle over $\cH'_{\ge i}$.
\item[($H+$)] If $i\in I_{+}$, then $\Pi_{i}$ exhibits $\cH'_{\ge i-1}$ as an open substack in a $\PP^{n-1}$-bundle over a $\PP^{n-1}$-bundle over $\cH'_{\ge i}$.
\item[($H-$)] If $i\in I_{-}$, then $\Pi_{i}$ exhibits $\cH'_{\ge i-1}$ as an open substack in a $\PP^{n-2}$-bundle over $\cH'_{\ge i}$.
\item[($H\pm$)] If $i\in I_{\pm}$, then $\Pi_{i}$ exhibits $\cH'_{\ge i-1}$ as an open substack in a $\PP^{n-1}$-bundle over $\cH'_{\ge i}$.
\end{enumerate}

Proof of ($H0$). When $i\in I_{0}$, $D_{i-1}=D_{i}$. We write the modification $\cF_{i-1}\dashrightarrow\cF_{i}$ as 
\begin{equation}\label{eq: flat modification} 
\begin{tikzcd}
\cF_{i-1} & \cF^{\flat}_{i-1/2}\ar[r, "\s x'_{i} ", hook]  \ar[l, "x'_{i}"', hook'] & \cF_{i} 
\end{tikzcd}
\end{equation}
Here both arrows have cokernel of length one supported at the labelled points. Such modifications $\cF^{\flat}_{i-1/2}$ of $\cF_{i}$ are parametrized by a hyperplane $H$ in the fiber $\cF_{i}|_{\s x'_{i}}$ and a line $L$ in the fiber $\cF_{i}|_{x'_{i}}$. 
The lower modifications of $\cF_{i}$ at $\s x'_{i}$ allowed in this case are those for which the map $t_{i}:\cL(D_{i})\to \cF_{i}$ factors through $\cF^{\flat}_{i-1/2}$, which is parametrized by the closed subset of hyperplanes $H \subset \cF_{i}|_{\s x'_{i}}$ containing the line given by the image of $\cL(D_{i})|_{\s x'_{i}}$. The space of choices for $H$ thus form a copy of $\PP^{n-2}$. The upper modifications of $\cF_{i}$ at $x'_{i}$ allowed in this case are those for which the map $t_{i-1}:\cL(D_{i})\to \cF^{\flat}_{i-1/2}\to \cF_{i-1}$ is saturated, which is parametrized by the open subset of those lines  $L \subset  \Cal{F}_i|_{x_i'}$ not equal to the image of $t_{i}(x'_{i})$. The space of such choices of $L$  thus form a copy of $\PP^{n-1}-\{\pt\}$.

This argument globalizes in the evident way as $(\{\cL(D_{j})\xr{t_{j}}\cF_{j}\}_{i\le j\le r}, \{x'_{i}\}_{i\in I_{0}}) $ moves over $\cH'_{\ge i}$, exhibiting $\Pi_i$ as an open substack in a $\PP^{n-2}\times\PP^{n-1}$-bundle. This applies similarly for the analogous arguments below for the other cases, so we focus on analyzing the fibers.

Proof of ($H+$). When $i\in I_{+}$, we have $D_{i-1}=D_{i}-\s x'_{i} $. We use the same notation $(H,L)\in \PP^{\vee}(\cF_{i}|_{\s x'_{i}})\times \PP(\cF_{i}|_{x'_{i}})$ as in the $(H0)$ case. This time the allowable lower modifications of $\cF_{i}$ at $\s x'_{i}$ are parametrized by the open subset of $H \subset \PP(\cF_{i}|_{\s(x'_{i})})$ that do not contain the image of $\cL(D_{i})|_{\s(x'_{i})}$. This forms a copy of $\PP^{n-1}-\PP^{n-2}$. The allowable upper modifications of $\cF_{i}$ at $x'_{i}$ are again parametrized by those $L$ not equal to the image of $t_{i}(x'_{i})$.  This is a copy of $\PP^{n-1}-\{\pt\}$. So the fibers of $\Pi_{i}$ in this case are isomorphic to $(\PP^{n-1}-\PP^{n-2})\times(\PP^{n-1}-\{\pt\})$.


Proof of ($H-$). When $i\in I_{-}$, we have $D_{i-1}=D_{i}+x'_{i}$. This time the allowable lower modifications of $\cF_{i}$ at $\s x'_{i}$ are parametrized by the closed subset of $H \subset \PP(\cF_{i}|_{\s(x'_{i})})$ that contain the image of $\cL(D_{i})|_{\s(x'_{i})}$. This forms a copy of $\PP^{n-2}$. The allowable upper modifications of $\cF_{i}$ at $x'_{i}$ are parametrized by a single point where  $L$ is equal to the image of $t_{i}(x'_{i})$.  So the fibers of $\Pi_{i}$ in this case are isomorphic to $\PP^{n-2}$.


Proof of ($H\pm$). When $i\in I_{\pm}$, we have $D_{i}+x'_{i}=D_{i-1}+\s x'_{i}$. This time the allowable lower modifications of $\cF_{i}$ at $\s x'_{i}$ are parametrized by the open subset of $H \subset \PP(\cF_{i}|_{\s(x'_{i})})$ that do not contain the image of $\cL(D_{i})|_{\s(x'_{i})}$. This forms a copy of $\PP^{n-1}-\PP^{n-2}$. The allowable upper modifications of $\cF_{i}$ at $x'_{i}$ are parametrized by a single point where  $L$ is equal to the image of $t_{i}(x'_{i})$.  So the fibers of $\Pi_{i}$ in this case are isomorphic to $\PP^{n-1}-\PP^{n-2}$.


 (3) By (2) and Lemma \ref{l:Dqf} we have
\begin{eqnarray*}
\dim \cZ'[N,I_{\bu}]&\le& |I_{0}|+\dim \frD(N;I_{\bu})+ (2n-3)|I_{0}|+(2n-2)|I_{+}|+(n-2)|I_{-}|+(n-1)|I_{\pm}|\\
&\le& |I_{0}|+|I_{+}|+ (2n-3)|I_{0}|+(2n-2)|I_{+}|+(n-2)|I_{-}|+(n-1)|I_{\pm}|\\
&=&(2n-2)|I_{0}|+\frac{3n-3}{2}|I_{+}|+\frac{3n-3}{2}|I_{-}|+(n-1)|I_{\pm}|.
\end{eqnarray*}
Here we use that $|I_{+}|=|I_{-}|$. Since $(3n-3)/2\le 2n-2$ and $n-1\le 2n-2$, we conclude that the last quantity above is $\le (2n-2)(|I_{0}|+|I_{+}|+|I_{-}|+|I_{\pm}|)=(2n-2)r$. Moreover, if equality holds, then we must have $|I_{+}|= |I_{-}|=|I_{\pm}|=0$, i.e., $I_{0}=\{1,2,\cdots, r\}$. 
\end{proof}

\part{Evidence} 

For the whole of Part 3, we assume $X'/X$ is a geometrically nontrivial double cover. 

\section{Nonsingular Fourier coefficients for unitary similitude groups}
\label{s:Eis GU}

In this section we extend the main result of \cite{FYZ} to the case of unitary similitude groups. One advantage of doing this is that we get central derivative formulas for the Siegel-Eisenstein series when the sign of the functional equation is $-1$ (when $n$ is odd). 

\subsection{Siegel--Eisenstein series on unitary groups with similitudes}\label{ss:Eis} 
We extend the result from \cite[\S2]{FYZ} to the case of unitary groups with similitudes.
For any one-dimensional $F$-vector space $L$, let $\Herm_{n}(F,L)$ be the $F$-vector space of $F'/F$-Hermitian forms $h: F'^{n}\times F'^{n}\to L\ot_{F}F'$ (with respect to the involution $1\ot\s$ on $L\ot_{F}F'$). For any $F$-algebra $R$, $\Herm_{n}(R, L):=\Herm_{n}(F,L)\ot_{F}R$ is the set of $L\ot_{F}R'$-valued $R'/R$-Hermitian forms on $R'^{n}$, where $R'=R\ot_{F}F'$.  When $L=F$ we write $\Herm_{n}(F)=\Herm_{n}(F,F)$ and $\Herm_{n}(R)=\Herm_{n}(F)\ot_F R$ for any $F$-algebra $R$.

Let $W$ be the standard split $F'/F$-skew-Hermitian space of dimension $2n$. Let $H_n=U(W)$ be the unitary group, and let $\wt H_n=GU(W)$ be the unitary group with similitudes, both as algebraic groups over $F$. 
Let $c:\wt H_n\to \BG_m$ denote the similitude character. Let $P_n(\BA)=M_n(\BA)N_n(\BA)$ be the standard Siegel parabolic subgroup of $H_n(\BA)$, where
\begin{align*}
  M_n(\BA)&=\left\{m(\alpha)=\begin{pmatrix}\alpha & 0\\0 &{}^t\bar \alpha^{-1}\end{pmatrix}: \alpha\in \GL_n(\BA_{F'})\right\},\\
  N_n(\BA) &= \left\{n(\beta)=\begin{pmatrix} 1_n & \beta \\0 & 1_n\end{pmatrix}: \beta\in\Herm_n(\BA)\right\}.
\end{align*}
Similarly, let $\wt P_n(\BA)=\wt M_n(\BA)N_n(\BA)$ be the standard Siegel parabolic subgroup of $\wt H_n(\BA)$, where
\begin{align*}
 \wt M_n(\BA)&=\left\{m(\alpha,c)=\begin{pmatrix}\alpha & 0\\0 &c\, {}^t\bar \alpha^{-1}\end{pmatrix}:c\in \BA^\times, \alpha\in \GL_n(\BA_{F'})\right\} \cong M_{n}(\BA)\times \BA^{\times}.
\end{align*}

Let $\eta: \BA^\times/F^\times\rightarrow \mathbb{C}^\times$ be the quadratic character associated to $F'/F$ by class field theory. Since $X'/X$ is \'etale, the character descends to $\y: \Pic_{X}(k)/\Pic_{X'}(k)\to \{\pm1\}$. Fix $\chi: \BA_{F'}^\times/F'^\times\rightarrow \mathbb{C}^\times$ a character such that $\chi|_{\BA_F^\times}=\eta^n$. We may view $\chi$ as a character on $M_n(\BA)\simeq \GL_{n}(\BA_{F'})$ by $\chi(\alpha)=\chi(\det(\alpha))$ and extend it to $P_n(\BA)$ trivially on $N_n(\BA)$. Fix a character $\chi_0:\BA^\times/F^\times\rightarrow \mathbb{C}^\times$. Define the \emph{degenerate principal series} of $\wt H_n(\BA)$ to be the unnormalized smooth induction 
$$
 I_n(s,(\chi,\chi_0))=\Ind_{\wt P_n(\BA)}^{ \wt H_n(\BA)}(\chi |\det|_{F'}^{s+n/2},\chi_0 |\cdot|_{F}^{-n(s+n/2)}),\quad s\in \mathbb{C}.
$$
In other words, its sections $\Phi(-,s)$ satisfy 
$$
\Phi(m(\alpha,c)n(\b)g ,s)=\chi(\alpha)\chi_0(c) |c|_{F}^{-n(s+n/2)} |\det\alpha|_{F'}^{s+n/2} \Phi(g,s).
$$

 For a standard section $\Phi(-, s)\in I_n(s,\chi)$, define the associated \emph{Siegel--Eisenstein series} 
$$E(g,s, \Phi)= \sum_{\gamma\in P_n(F)\backslash H_n(F)}\Phi(\gamma g, s),\quad g\in\wt H_n(\BA),
$$ which converges for $\Re(s)\gg 0$ and admits meromorphic continuation to $s\in \mathbb{C}$. Here we have used $P_n(F)\backslash H_n(F)\simeq \wt P_n(F)\backslash \wt H_n(F)$.

Notice that $E(g,s,\Phi)$ depends on the choice of $\chi$. 
In this paper, we will assume both $\chi$ and $\chi_{0}$ are unramified everywhere. Eventually it will be convenient to take $\chi_{0}=\y^{n}$ but we do not make this assumption until \S\ref{ss:comparison n=1}. Then  $I_n(s,(\chi,\chi_{0}))$ is unramified and we fix $\Phi(-, s)\in I_n(s,(\chi,\chi_{0}))$  as the unique $K=\wt H_n(\wh\cO)$-invariant section normalized by
\begin{align}\label{eq:unr sec}
\Phi(1_{2n}, s)=1.
\end{align}
Similarly we normalize $\Phi_v \in I_n(s,(\chi_v,\chi_{0}))$ for every $v\in|X|$ and we then have a factorization $\Phi=\bigotimes_{v\in|X|}\Phi_v$.

\subsection{Fourier expansion} 
Let $\om_{F}$ be the generic fiber of the canonical bundle $\om_X$ of $X$, and $\BA_{\om_{F}}=\BA\ot_{F}\om_{F}$. Let $\Herm_{n}(\BA, \om_{F})$ (resp. $\Herm_{n}(\BA)$) denote the Hermitian forms on $\BA^n$ valued in $ \BA_{\om_{F}}$ (resp.  $\BA$). The residue pairing $\Res: \BA_{\om_{F}}\times \BA\to k$ induces a pairing  
$$\j{\cdot,\cdot}: \Herm_{n}(\BA, \om_{F})\times \Herm_{n}(\BA)\to k$$
given by $\j{T,b}=\Res(-\Tr(Tb))$. Composing this pairing with the fixed nontrivial additive character $\psi_{0}: k\to \BC^\times$  exhibits $\Herm_{n}(\BA,\om_{F})$ as the Pontryagin dual of $\Herm_{n}(\BA)$. Moreover, it exhibits $\Herm_{n}(F,\om_{F})$ as the Pontryagin dual of $\Herm_{n}(F)\bs \Herm_{n}(\BA)=N_{n}(F)\bs N_{n}(\BA)$. The global residue pairing is the sum of local residue pairings $\j{\cdot,\cdot}_{v}: \Herm_{n}(F_{v}, \om_{F_{v}})\times \Herm_{n}(F_{v})\to k$ defined by $\j{T,b}_{v}=\tr_{k_v/k}\Res_{v}(-\Tr(Tb))$.


We have a Fourier expansion 
$$E(g,s,\Phi)=\sum_{T\in\Herm_n(F,\om_{F})}E_T(g,s,\Phi), \quad g\in\wt H_n(\BA),$$ 
where
 $$E_T(g,s,\Phi)=\int_{N_n(F)\backslash N_n(\mathbb{A})} E(n(b)g,s,\Phi)\psi_{0}(\j{T,b})\,\rd n(b),$$ 
and the Haar measure $\rd n(b)$ is normalized such that $N_{n}(F)\bs N_{n}(\BA)$ has volume $1$. 

When $T$ is nonsingular, for a factorizable $\Phi=\bigotimes_{v\in|X|}\Phi_v$ we have a factorization of the Fourier coefficient into a product (cf. \cite[\S 4]{Kudla1997})
\begin{equation}\label{Eis factor Wh}
E_T(g,s,\Phi)=|\om_{X}|_{F}^{-n^{2}/2}\prod_v W_{T,v}(g_v, s, \Phi_v),
\end{equation}
where the \emph{local (generalized) Whittaker function} is defined by 
$$W_{T,v}(g_v, s, \Phi_v)=\int_{N_n(F_{v})}\Phi_v(w_n^{-1}n(b)g_v,s) \psi_{0}(\j{T,b}_{v})
\, \rd_{v} n(b),\quad w_n=
\begin{pmatrix}
0  & 1_n\\
  -1_n & 0\\
\end{pmatrix}$$ 
and has an analytic continuation to $s\in \mathbb{C}$. Here the local Haar measure $\rd_{v} n(b)$ is normalized so that the volume of $N_{n}(\cO_{v})$ is $1$. The factor $|\om_{X}|_{F}^{-n^{2}/2}$ is the ratio between the global measure $\rd n$ and the product of the local measures $\prod_{v}\rd_{v}n$.

Note that for  $m(\a,c)\in \wt M_{n}(F_{v})$,
\begin{equation}\label{eq:ch base}
W_{T,v}(m(\alpha,c), s, \Phi_v)=\chi(\a)(\chi_0\eta^n)(c) |c|_{F_v}^{-n(-s+n/2)}|\det(\alpha)|^{-s+n/2}_{F'_v}W_{c^{-1}\,^t\bar\alpha\,T\alpha,v}(1, s, \Phi_v).
\end{equation}

We define the {\em regular} part of the Eisenstein series to be
\begin{equation}\label{eq:def reg}
E^{\rm reg}(g,s,\Phi)=\sum_{T\in\Herm_n(F,\om_{F})\atop \rank T=n}E_T(g,s,\Phi),\quad g\in\wt H_n(\BA).
\end{equation}
Analogous to \cite[\S2.6]{FYZ} we view $E^{\rm reg}$ as a function on
$$\wt M_n(F)\bs \wt M_n(\BA)/\wt M_n(\wh\cO)\simeq  \Bun_{\GL(n)'}(k)\times \Pic_X(k).
$$

For $(\cE,\LL)\in \Bun_{\GL(n)'}(k)\times \Pic_{X}(k)$ and $a:\cE\to \s^{*}\cE^\vee\ot\nu^{*}\LL$, we can define the $a$-th Fourier coefficient $E_{a}(m(\cE,\LL), s,\Phi)$ (similar to what is done in \cite[\S2.6]{FYZ}).

\begin{thm}\label{th:Eis Den} 
Let $\cE$ be a vector bundle over $X'$ of rank $n$.  
Then 
\begin{equation}\label{eq:reg}
E^{\rm reg}(m(\cE,\LL),s,\Phi)=\sum_{a}E_{a}(m(\cE,\LL), s,\Phi)
\end{equation}
where the sum runs over all injective Hermitian maps  $a: \cE\to \s^{*}\cE^\vee\ot\nu^{*}\LL $, and
\begin{align*}
E_{a}(m(\cE,\LL),s,\Phi)=&(\chi_0\eta^n)(\LL)\chi(\det\cE)q^{-(\deg(\cE)-n\deg(\LL))( s-n/2)-\frac{1}{2}n^2\deg(\omega_X) }\\
&\times \sL_n(s) ^{-1} \Den(q^{-2s}, \coker(a)).
\end{align*}
Here 
\begin{equation}\label{eq: L-factor}
\sL_n(s)=\prod_{i=1}^nL(i+2s, \eta^i).
\end{equation}
The density function $ \Den(q^{-2s}, \coker(a))$ (see \cite[\S2.6, \S5.1]{FYZ}) is  a polynomial in $q^{-s}$ of degree 
$$
\len (\coker(a))=\deg (\s^{*}\cE^\vee\ot\LL)-\deg(\cE)=2 n(\deg\LL+\deg\omega_X)-2\deg(\cE).
$$

\end{thm}
\begin{proof} 
By  \eqref{Eis factor Wh} and \eqref{eq:ch base} we have
\begin{align*}
E_T(m(\alpha,c), s, \Phi_v)=&\chi(\a)(\chi_0\eta^n)(c) |c|_{F}^{-n(-s+n/2)}|\det(\alpha)|^{-s+n/2}_{F'}  |\om_{X}|^{-\frac{1}{2}n^2}\\
&\times \prod_{v\in|X|} W_{c^{-1}\,^t\bar\alpha\,T\alpha,v}(1, s, \Phi_v).
\end{align*}
Note 
$$
|\det(\alpha)|_{F'}=q^{\deg(\cE)}, \quad |c|_{F}=q^{\deg(\LL)}.
$$
The rest is the same as  (the proof of) \cite[Thm.~2.7, Thm.~5.1]{FYZ}.
\end{proof}

 \subsection{Non-singular terms with similitudes}\label{ssec: similitude higher KR} 
Now we can state a generalization of the main result of \cite{FYZ} to Hermitian shtukas with similitudes. 
%
%
%
%
 \begin{thm}
\label{th:ns} 
Let $\cE$ be a vector bundle of rank $n$ on $X'$, and let $d=-\deg(\cE)+n(\deg\LL+\deg\omega_X)$. Let $a: \cE\to \s^{*}\cE^\vee\ot\nu^{*}\LL $ be an injective (i.e. non-singular) Hermitian map. Then
 $$
 \deg [ \cZ_{\cE,\LL}^{r}(a)]=\frac{1}{(\log q)^r}  \left(\frac{d}{ds}\right)^r\Big|_{s=0} \left( q^{ds} \Den(q^{-2s},\coker(a))\right).
 $$
 \end{thm} 
 \begin{remark}
Here we note that  $\Den(q^{-2s}, \coker(a))$ is a polynomial in $q^{-s}$ of degree $2d=-2\deg\cE+2(\deg\omega_X+\deg\LL)$. The right hand side of the above formula is symmetric up to the sign $\eta^n(\LL)$ with respect to the substitution $s\mapsto -s$. Therefore the right side vanishes if $(-1)^{r}\ne \eta^{n}(\LL)$.

On the other hand, by Lemma \ref{lem: when non-empty}, $\Sht^{r}_{U(n),\LL}$ is empty when $(-1)^{r}\ne \eta^{n}(\LL)$, so in that case the identity in the theorem holds trivially. The theorem is nontrivial only when $(-1)^r=\eta^{n}(\LL)$. 

\end{remark}
 \begin{proof} The proof is similar to \cite[Thm.~12.1]{FYZ}. We introduce a generalization of the moduli stack of torsion Hermitian sheaves $\Herm_{2d}(X'/X,\LL)$ that classifies $(\cQ,h)$ where $\cQ$ is a torsion coherent sheaf on $X'$ of length $2d$, and $h$ is an isomorphism $\cQ\isom \s^{*}\cQ^{\vee}\ot \nu^{*}\LL$ such that $\s^{*}h^{\vee}=h$. 
 
 The arguments in {\it loc.\ cit.\ } also show
 \begin{enumerate}
\item there is a graded virtual perverse sheaf on $\Herm_{2d}(X'/X,\LL)$
\begin{equation*}
\cK^{\Eis}_{d}=\bigoplus_{i=0}^{d}\cK^{\Eis}_{d,i}
\end{equation*}such that
\begin{equation}\label{intro KEis}
\Den (q^{-2s}, \cQ)= \sum_{i=0}^{d}\Tr(\Frob_{\cQ}, (\cK^{\Eis}_{d,i})_{\cQ})q^{-2is}, \quad \mbox{ for } \cQ\in \Herm_{2d}(X'/X,\LL)(k).
\end{equation}
\item there is a graded virtual perverse sheaf on $\Herm_{2d}(X'/X,\LL)$
\begin{equation*}
\cK^{\Int}_{d}=\bigoplus_{i=0}^{d}\cK^{\Int}_{d,i}
\end{equation*}such that
\begin{equation}\label{intro KInt}
\deg [\cZ_{\cE,\LL}^{r}(a)]=\sum_{i=0}^{d}\Tr(\Frob_{\cQ}, (\cK^{\Int}_{d,i})_{\cQ})\cdot (d-2i)^{r}.
\end{equation}Here $\cQ=\coker(a)\in \Herm_{2d}(X'/X,\LL)(k)$.
\end{enumerate}
By the same proof of {\it loc. cit.}, $\cK^{\Eis}_{d}$ and $\cK^{\Int}_{d}$ are virtual linear combinations of isotypic summands of the Hermitian-Springer sheaf $\Spr^{\Herm}_{2d}$ on $\Herm_{2d}(X'/X,\LL)$ under the action of $W_{d}=(\Z/2\Z)^{d}\rtimes S_{d}$. The same proof of \cite[Prop.12.3]{FYZ} again shows
\begin{equation}\label{KK}
\cK^{\Eis}_{d}\cong \cK^{\Int}_{d}
\end{equation}
as graded virtual perverse sheaves on $\Herm_{2d}(X'/X,\LL)$, and the proof is complete. \end{proof}

\begin{remark} When $\y^{n}(\frL)=-1$ (so $n$ is necessarily odd), $q^{ds}\Den(q^{-2s}, \coker(a))$ is an odd function in $s$. Theorem \ref{th:ns} then gives a geometric interpretation of odd order central derivative of nonsingular Fourier coefficients of the normalized Eisenstein series in terms of degrees of special cycles. This complements the even derivative case treated in \cite{FYZ}.
\end{remark}

\subsection{A refinement of non-singular coefficients}
In certain cases, the special cycles $\cZ^{r}_{\cE,\LL}(a)$ can be further decomposed into the union of two open-closed parts. We will prove a refinement of Theorem \ref{th:ns} that calculates the degree of the $0$-cycle on each part.


Let $\LL^{(n)}:=\om_X^{\ot n-1}\ot \LL^{\ot n}$. Taking determinant induces a map
\begin{equation}
\det: \Sht^{r}_{U(n),\LL}\to \Sht^{r}_{U(1), \LL^{(n)}}.
\end{equation}

Below we consider the case $\y^{n}(\LL)=1$.  In this case, $\Sht_{U(n), \LL}^{r}=\vn$ unless $r$ is even, by Lemma \ref{lem: when non-empty}. So we also assume $r$ is even.

Let $\frN=\om_{X}^{\ot n}\ot \LL^{\ot n}$.
Since $\frN$ is a norm (as $\y(\LL^{\ot n})=1$ by assumption and  $\om_{X}$ is known to be a norm), Lemma \ref{lem: 2 connected components} implies that the set $\Irr(\Prym_{\frN})$ of irreducible components of $\Prym_{\frN}$ (defined over $k$) is $\Z/2\Z$-torsor. For $\e\in \Irr(\Prym_{\frN})$, let $\Prym^{\e}_{\frN}$ be the corresponding component. Let $p: \Sht^{r}_{U(1), \LL^{(n)}}\to \Prym_{ \frN}$ be the map recording $\cF_{0}$. Let $\Sht^{r,\e}_{U(1), \LL^{(n)}}$ be the preimage of $\Prym^{\e}_{ \frN}$ under $p$, and let $\Sht^{r,\e}_{U(n),\LL}$ be the further preimage under $\det$. Define $\cZ^{r,\e}_{\cE,\LL}(a)$ to be the preimage of $\Sht^{r,\e}_{U(n),\LL}$ under the map $\z: \cZ^{r}_{\cE,\LL}(a)\to \Sht^{r}_{U(n), \LL}$.

\begin{thm}\label{thm: deg equal on components} Assume that $\y^{n}(\LL)=1$ and $r$ even and $r>0$.  Let $\cE$ be a vector bundle of rank $n$ on $X'$, and $d=-\deg(\cE)+n(\deg\LL+\deg\omega_X)$. Let $a: \cE\to \s^{*}\cE^\vee\ot\nu^{*}\LL $ be an injective (i.e. non-singular) Hermitian map. Then for any $\e\in \Irr(\Prym_{\frN})$ we have:
 $$
 \deg [ \cZ_{\cE,\LL}^{r,\e}(a)]=\frac{1}{2}\deg [ \cZ_{\cE,\LL}^{r}(a)]=\frac{1}{2(\log q)^r}  \left(\frac{d}{ds}\right)^r\Big|_{s=0} \left( q^{ds} \Den(q^{-2s},\coker(a))\right).
 $$
\end{thm}
\begin{proof}
By Theorem \ref{th:ns}, it suffices to show that $\deg [ \cZ_{\cE,\LL}^{r,\e}(a)]=\frac{1}{2}\deg [ \cZ_{\cE,\LL}^{r}(a)]$.

Define $X_d$ to be the $d^{\mrm{th}}$ symmetric power of $X$, $X_d'$ similarly, and $\nu_d \co X_d' \rightarrow X_d$ to be the map induced by $\nu$. Consider the moduli stack $\cP_{d}$ that classifies $(D,\cF,\io)$ where $D$ is an effective divisor on $X$ of degree $d$,  $\cF\in \Pic_{X'}$ and $\io$ is an isomorphism $\Nm(\cF)\cong \cO_{X}(D)$. The map $p: \cP_{d}\to X_{d}$ given by forgetting $\cF$ is a torsor for $\Prym$. Let $\cP_{d}\to X^{\sh}_{d}\xr{\mu_{d}} X_{d}$ be the Stein factorization of $p$. Since $\Prym$ has two geometrically connected components, $\mu_{d}: X^{\sh}_{d}\to X_{d}$ is an \'etale double cover. Consider the  map $\a: X'_{d}\to \cP_{d}$ (over $X_{d}$) sending $D'\in X'_{d}$ to $(\nu_{d}(D'), \cO_{X'}(D'),\io)$ where $\io$ is the canonical isomorphism $\Nm(\cO_{X'}(D'))\cong \cO_{X}(D)$. It induces a map $\nu^{\sh}_{d}: X'_{d}\to X^{\sh}_{d}$ such that $\nu_{d}: X'_{d} \to X_{d}$ factorizes as
\begin{equation}
\nu_{d}: X'_{d}\xr{\nu^{\sh}_{d}} X^{\sh}_{d}\xr{\mu_{d}} X_{d}.
\end{equation}
If we base change $\mu_{d}$ along the symmetrization map $X^{d}\to X_{d}$ we get an \'etale double cover $X^{d,\sh}$ of $X^{d}$. We claim that the \'etale double cover $X^{d,\sh}\to X^{d}$ is given by the homomorphism
\begin{equation}
\pi_{1}(X^{d})\to  (\Z/2\Z)^{d}\xr{\textup{sum}}\Z/2\Z
\end{equation}
where the first map classifies the $(\Z/2\Z)^{d}$-torsor  $X'^{d}/X^{d}$. Indeed we may consider the map $\wt\a: X'^{d}\to X'_{d}\xr{\a} \cP_{d}$; it sends $(x'_{1},\cdots, x'_{d})$ to $(\nu(x'_{1})+\cdots+\nu(x'_{d}), \cO(x'_{1}+\cdots+x'_{d}))$. Each time we change $x'_{i}$ to $\s(x'_{i})$ the resulting point under $\wt\a$ moves to a different component of the corresponding fiber of $\cP_{d}\to X_{d}$, hence the resulting map $X'^{d}\to X^{\sh}_{d}$ factors through the quotient of $X'^{d}$ by the subgroup $\ker(\textup{sum}: (\Z/2\Z)^{d}\to \Z/2\Z)$.

To summarize, $\mu_{d}: X^{\sh}_{d}\to X_{d}$ is the double cover attached to the local system $\y_{d}$ on $X_{d}$ (see \cite[\S11.4]{FYZ}).

Let $\cM^{\ns}_{d}$ be the open-closed substack of $\cM^{\ns}_{\GL(n)', U(n), \LL}=\cM^{\c}_{\GL(n)', U(n), \LL}$ where $\deg\cE=n(\deg\LL+\deg\om_{X})-d$. Let $\cA^{\ns}_{d}$ be the corresponding Hitchin base and $f_d: \cM^{\ns}_{d}\to \cA^{\ns}_{d}$ the Hitchin map. Write $\Herm_{2d}:=\Herm_{2d}(X'/X, \LL)$. Let  $\Lagr_{2d}:=\Lagr_{2d}(X'/X, \LL)$  be the stack classifying pairs $(\cL\subset\cQ)$ where $\cQ\in\Herm_{2d}$ and $\cL\subset \cQ$ is a Lagrangian subsheaf with the natural map $\ups_{2d}:\Lagr_{2d}\to \Herm_{2d} $ sending $(\cL\subset\cQ)$ to $\cQ$. Recall from \cite[Lemma 8.8]{FYZ} we have a commutative diagram where the left side square is Cartesian
\begin{equation*}
\xymatrix{\cM^{\ns}_{d}\ar[d]^{f_{d}}\ar[r]^-{\e'_{d}} & \Lagr_{2d}\ar[d]^{\ups_{2d}} \ar[r]^-{} & X'_{d}\ar[d]^{\nu_{d}}\\
\cA^{\ns}_{d} \ar[r]^-{\e_{d}} & \Herm_{2d}\ar[r] & X_{d} }
\end{equation*}Here
the map $\Herm_{2d}\to X_{d}$ is the descent of the divisor of the Hermitian torsion sheaf, and $\Lagr_{2d}\to X'_{d}$ records the divisor of the Lagrangian subsheaf. 

Let $\mu_{d}^{\Herm}: \Herm^{\sh}_{2d}\to \Herm_{2d}$ and $\mu^{\cA}_{d}: \cA^{\ns,\sh}_{d}\to \cA^{\ns}_{d} $ be the base changes of the double cover $\mu_{d}: X^{\sh}_{d}\to X_{d}$. Then we have a diagram where all squares are Cartesian
\begin{equation*}
\xymatrix{\cM^{\ns}_{d}\ar[d]^{f^{\sh}_{d}}\ar[r]^-{\e'_{d}} & \Lagr_{2d}\ar[d]^{\ups^{\sh}_{2d}}\\
\cA^{\ns,\sh}_{d}\ar[d]^{\mu^{\cA}_{d}}\ar[r]^-{\e^{\sh}_{d}} & \Herm^{\sh}_{2d}\ar[d]^{\mu^{\Herm}_{d}}\\
\cA^{\ns}_{d} \ar[r]^-{\e_{d}} & \Herm_{2d}
}
\end{equation*}

We claim that the double cover $\mu^{\cA}_{d}: \cA^{\ns,\sh}_{d}\to \cA^{\ns}_{d}$ is trivial. To show this, we prove a stronger statement: the $\Prym$-torsor $\cP_{d}$ becomes trivial when pulled back along $\cA^{\ns}_{d}\to X_{d}$. Indeed, by assumption, $\LL^{\ot n}$ is a norm, hence $(\om_{X}\ot \LL)^{\ot n}$ is also a norm since $\om_{X}$ is known to be a norm. Say $(\om_{X}\ot \LL)^{\ot -n}=\Nm(\frM)$ for some $\frM\in \Pic_{X'}(k)$.  This means $\frM$ carries a Hermitian form $a_{\frM}: \frM\isom \s^{*}\frM^{\ot -1}\ot \nu^{*}(\om_X\ot\LL)^{\ot -n}$. We define a map $\b_{\frM}: \cA^{\ns}_{d}\to \cP_{d}$ sending $(\cE,a)$ to $(\div(a), (\det\cE)^{\ot -1}\ot \frM^{\ot -1}, \io)$. Here $\io: \Nm((\det\cE)^{\ot -1}\ot \frM^{\ot -1})\cong \cO(\div(a))$ is defined as follows. Note that $\det(a)$ is a Hermitian map $\det\cE\to \s^{*}(\det\cE)^{\ot -1}\ot \nu^{*}(\om_{X}\ot \LL)^{\ot n}$. Then $\det(a)\ot a_{\frM}$ gives a Hermitian map $\det\cE\ot \frM\to \s^{*}(\det\cE\ot \frM)^{\ot -1}$ whose divisor descends to $\div(a)\in X_{d}$. Therefore it induces  a  canonical Hermitian isomorphism $\det\cE\ot \frM\isom \s^{*}(\det\cE\ot \frM)^{\ot -1}\ot \nu^{*}\cO(-\div (a))$, which gives the desired isomorphism $\io$. The map $\b_{\frM}$ then trivializes the $\Prym$-torsor $\cP_{d}|_{\cA^{\ns}_{d}}$ (but not canonically, since the trivialization depends on the choice of $\frM$). 



Let $(\cE, a^{\sh}_{1})$ and $(\cE, a^{\sh}_{2})\in \cA^{\ns,\sh}_{d}(k)$  be the two preimages of $(\cE,a)$.  We have a canonical map $\cZ^{r}_{\cE,\frL}(a)\to \cM^{\ns}_{d}$ by recording $\cE\xr{t_{0}}(\cF_{0},h_{0})$. We have a further map $\cM^{\ns}_{d}\to \cP_{d}$ mapping $(\cE\xr{t}(\cF,h))$ to $(\div(a), \det\cF\ot(\det\cE)^{\ot -1})$ where $a=\s^{*}t^{\vee}\c h\c t$. By construction, for each $\e\in \Irr(\Prym_{\frN})$, $\cZ^{r,\e}_{\cE,\frL}(a)$ maps exactly to the one (out of two) components of the fiber  $p^{-1}(\div(a))\subset \cP_{d}$ over $\div(a)\in X_{d}(k)$. Therefore, the decomposition of $\cZ^{r}_{\cE,\frL}(a)$ according to $\e$ is the same as the decomposition given by the map $\cZ^{r}_{\cE,\frL}(a)\to f_{d}^{-1}(\cE, a)\to \mu^{\cA,-1}_{d}(\cE, a)=\{(\cE, a^{\sh}_{1}), (\cE, a_{2}^{\sh})\}$. Let $\cZ^{r}_{\cE,\frL}(a^{\sh}_{i})$ be the fiber over $(\cE,a^{\sh}_{i})$. We thus reduce to showing
\begin{equation}\label{Zdeg half}
\deg[\cZ^{r}_{\cE,\frL}(a^{\sh}_{1})]=\deg[\cZ^{r}_{\cE,\frL}(a^{\sh}_{2})].
\end{equation}
Here $ [\cZ^{r}_{\cE,\frL}(a^{\sh}_{i})]$ is defined as the restriction of $[\cZ^{r}_{\cE,\frL}(a)]$ to the open-closed substack $\cZ^{r}_{\cE,\frL}(a^{\sh}_{i})$, $i=1,2$.

Let $\cQ_{i}^{\sh}\in \Herm^{\sh}_{2d}(k)$ be the image of $(\cE,a^{\sh}_{i})$. Define $\cZ^{r}_{\cQ^{\sh}_{i}}$ and its fundamental class using the Shtuka construction for $\Lagr_{2d}$ (see \cite[\S11.1]{FYZ}). Then the smoothness of $\e_{d}$ implies $\deg[\cZ^{r}_{\cE,\frL}(a^{\sh}_{i})]=\deg[\cZ^{r}_{\cQ^{\sh}_{i}}]$.  Recall the self-correspondence $\Hk^{1}_{\Lagr_{2d}}$ of  $\Lagr_{2d}$ over $\Herm_{2d}$. Applying the Lefschetz trace formula \cite[Prop. 11.8]{FYZ} to the same situation as \cite[Corollary 11.9]{FYZ} except that we are now working over the base $\Herm^{\sh}_{2d}$ rather than $\Herm_{2d}$, we get 
\begin{equation}\label{degZsh}
\deg[\cZ^{r}_{\cQ^{\sh}_{i}}]=\Tr\left([\Hk_{\Lagr_{2d}}^{1}]^{r}\c\Frob_{\cQ^{\sh}_{i}}, (R(\ups_{2d}^{\sh})_{*}\Qlbar)_{\cQ_{i}^{\sh}}\right), \quad i=1,2.
\end{equation}
The only caveat here is that the self-correspondence $\Hk_{\Lagr_{2d}}^{1}$ of $\Lagr_{2d}$ is not over $\Herm^{\sh}_{2d}$ but only over $\Herm_{2d}$: the two compositions $\ups_{2d}^{\sh}\c\pr_{0}$ and $\ups_{2d}^{\sh}\c\pr_{1}: \Hk_{\Lagr_{2d}}^{1}\to \Herm^{\sh}_{2d}$ differ by the involution $\s_{d}$ of $\Herm^{\sh}_{2d}/\Herm_{2d}$ (base changed from the involution of $X_{d}^{\sh}/X_{d}$). Therefore $[\Hk_{\Lagr_{2d}}^{1}]$ induces a map $R(\ups_{2d}^{\sh})_{*}\Qlbar\to \s_{d}^{*}R(\ups_{2d}^{\sh})_{*}\Qlbar$. Since $r$ is even, $[\Hk_{\Lagr_{2d}}^{1}]^{r}$ induces an endomorphism of $R(\ups_{2d}^{\sh})_{*}\Qlbar$. The proof of the Lefschetz trace formula in {\em loc.cit} adapts to this twisted situation.

Let $W_{d}=(\Z/2\Z)^{d}\rtimes S_{d}$ and $\chi_{d}: W_{d}\to \Z/2\Z$ be trivial on $S_{d}$ and nontrivial on each factor of  $\Z/2\Z$. Let $W'_{d}=\ker(\chi_{d})$. Recall the map $\pi^{\Herm}_{2d}: \wt\Herm_{2d}\to\Herm_{2d}$ (\cite[\S4.2]{FYZ}), and it factors through $\pi^{\sh}_{2d}: \wt\Herm_{2d}\to\Herm^{\sh}_{2d}$. This is a small map that is generically a $W'_{d}$-torsor. Therefore $\Spr^{\sh}_{2d}:=R(\pi^{\sh}_{2d})_{*}\Qlbar$ is a middle extension perverse sheaf that carries an action of $W'_{d}$.  For a representation $\r$ of $W'_{d}$, let $\Spr^{\sh}_{2d}[\r]=\Hom_{W'_{d}}(\r, \Spr^{\sh}_{2d})=(\r^{\vee}\ot R(\pi^{\sh}_{2d})_{*}\Qlbar )^{W'_{d}}$.  Then, analogous to \cite[Lemma 11.4, Corollary 11.5]{FYZ},  $R(\ups_{2d}^{\sh})_{*}\Qlbar\cong \Spr^{\sh}_{2d}[\Ind_{S_{d}}^{W'_{d}}\one]$. Recall $\{\r_{i}\}_{0\le i\le d}$ are irreducible direct summands of $\Ind_{S_{d}}^{W_{d}}\one$ (as $W_{d}$-representation) defined in \cite[\S11.2]{FYZ}. We have a corresponding decomposition of $W'_{d}$-representations
\begin{equation*}
\Ind_{S_{d}}^{W'_{d}}\one=\bigoplus_{0\le i\le [d/2]} \r'_{i}
\end{equation*}
where $\r'_{i}\in\Irr(W'_{d})$ is characterized by:
\begin{equation}\label{ind r'}
\Ind_{W'_{d}}^{W_{d}}\r'_{i}\cong \begin{cases}\r_{i}\op \r_{d-i}, & i<d/2;  \\  \r_{d/2}, & i=d/2.\end{cases}
\end{equation}
Let $\cK^{\sh}_{d,i}=\Spr^{\sh}_{2d}[\r'_{i}]$ for $0\le i\le [d/2]$, then 
\begin{equation*}
R(\ups_{2d}^{\sh})_{*}\Qlbar\cong\bigoplus_{0\le i\le [d/2]}\cK^{\sh}_{d,i}.
\end{equation*}
The action of $[\Hk_{\Lagr_{2d}}^{1}]^{r}$ preserves each $\cK^{\sh}_{d,i}$ and acts on it by the scalar $(d-2i)^{r}$. Therefore \eqref{degZsh} implies
\begin{equation*}
\deg[\cZ^{r}_{\cQ^{\sh}_{1}}]=\sum_{0\le i<d/2} (d-2i)^{r}\Tr(\Frob_{\cQ^{\sh}_{1}},(\cK^{\sh}_{d,i})_{\cQ^{\sh}_{1}} ).
\end{equation*}
Here we have omitted the term $i=d/2$  because its coefficient $(d-2i)^{r}=0$ (using that $r>0$).

The same formula is true when $\cQ^{\sh}_{1}$ is replaced by $\cQ^{\sh}_{2}$. Therefore to show \eqref{Zdeg half} it suffices to show that $\cK^{\sh}_{d,i}$ has the same Frobenius trace at conjugate points under $\s_{d}$. Better, we will show that $\cK^{\sh}_{d,i}$ descends to $\Herm_{2d}$. Indeed, \eqref{ind r'} implies $\r_{i}|_{W'_{d}}\cong  \r'_{i}$ for $0\le i<d/2$, hence $\Spr^{\sh}_{2d}[\r_{i}|_{W'_{d}}]\cong \Spr^{\sh}_{2d}[\r'_{i}]$. Note that $\Spr^{\sh}_{2d}[\r_{i}|_{W'_{d}}]\cong \mu_{d}^{\Herm, *}\cK^{\Int}_{d,i}$, hence $\cK^{\sh}_{d,i}$ is the pullback of $\cK^{\Int}_{d,i}$. The proof is complete.

\end{proof}

\begin{remark}
Interestingly, the statement of Theorem \ref{thm: deg equal on components} does not hold for $r=0$. For example, take $\LL=\cO_{X}$, and take any Hermitian bundle $(\cF,h)$ over $X'$. Let $x'\in |X'|$ be inert over $X$. Let $\cE\subset \cF$ be a subsheaf such that $\cF/\cE\cong k_{x'}$. Such an $\cE$ corresponds to a hyperplane $H$ in the fiber $V_{x'}$ of $\cF$ at $x'$. We choose $H$ to be non-degenerate under the Hermitian form $h_{x'}$ on $V_{x'}$. With this choice of $H$, the cokernel $\cQ$ of $a: \cE\to \cF\xr{h}\s^{*}\cF^{\vee}\to\s^{*}\cE^{\vee}$ is isomorphic to $k_{x'}[\varpi]/\varpi^{2}$, where $\varpi$ is a uniformizer at $x=\nu(x')$. In particular, the Hermitian torsion sheaf $\cQ$ has only one Lagrangian $\varpi k_{x'}[\varpi]/\varpi^{2}$. Hence the special cycle $\cZ^{0}_{\cE}(a)$ (a discrete scheme over $k$) has only one geometric point, and cannot distribute evenly over the two possibilities of $\e$.

Reviewing the proof above, the reason it does not work for $r=0$ is that the sheaf $\cK^{\sh}_{d,d/2}$ occurring in the proof does not descend to $\Herm_{2d}$, so its Frobenius trace cannot be the same at all conjugate rational points for all $\F_q$. 
\end{remark}

\subsection{Normalized Eisenstein series}
As a preparation for the next two sections, we introduce a normalized version of the Siegel--Eisenstein series. 

There is an intertwining operator 
$$M(s): I(s,(\chi,\chi_0))\to I(-s,(\chi,\chi_0\eta^n))=I(-s,(\chi,\chi_0))\otimes (\eta^n\circ c)$$
where $c:\wt H_n(\BA)\to \BA^\times$ is the similitude factor.
The image of the unramified section is
$$
M(s)\Phi(s,g)=q^{-2ns\deg \omega_X}\frac{\sL_n(-s)}{\sL_n(s)} \eta^n(c(g))\Phi(-s,g),
$$
which follows from the computation in \cite[Prop. 2.1 and p.170]{Tan99} and the functional equation (cf. \eqref{FE L eta} below).
Later in \S\ref{ss:E0 n=1} we will recall the well-known computation when $n=1$.

We define a normalized Eisenstein series
\begin{align}\label{eq:nor Eis}
\wt E(g,s,\Phi)=q^{n\deg \omega_X s} \sL_n(s)E(g,s,\Phi).
\end{align}
Then it satisfies a functional equation
\begin{align}\label{eq:Eis FE}
\wt E(g,s,\Phi)=\eta^n(c(g))\wt E(g,-s,\Phi), \quad g\in\wt H_n(\BA).
\end{align}
Note that when $n$ is odd and $\eta(c(g))=-1$, the sign of the functional equation is $-1$.

By Theorem \ref{th:Eis Den}, for injective $a: \cE\to \s^{*}\cE^\vee\ot\nu^{*}\LL $,
the $a$-th Fourier coefficient (expanded at $g=m(\cE,\LL)$)  then has a very simple form
\begin{align}\label{eq:nor Eis a}
\wt E_{a}(m(\cE,\LL),s,\Phi)=&(\chi_0\eta^n)(\LL)\chi(\det\cE)q^{d( s-n/2)} \Den(q^{-2s}, \coker(a)),
\end{align}
where
\begin{align}\label{eq:deg n}
d=n(\deg\LL+\deg\omega_X)-\deg\cE
\end{align}
is the half of the degree of $ \Den(q^{-2s}, \coker(a))$ (as a polynomial in $q^{-s}$). Note that $d$ depends on $\cE$ via its degree. This normalization differs from \cite{FYZ} in that here we do not absorb the trivial terms.  

\section{Modularity: the case of $U(1)$}\label{sec: constant term}

In this section we prove the modularity Conjecture \ref{c:mod} for $n=m=1$. We show in Corollary \ref{th:n=1 mod} that the modularity can be checked after taking the degrees of special cycles (on each connected component, if there are multiple). The degrees of nonzero terms in the generating series in this case are taken care of by Theorem \ref{thm: deg equal on components}. The bulk of this section is devoted to the calculation of the degree of the $0$-th term in the generating series, which we relate to the higher derivatives of an $L$-function, completing the higher Siegel-Weil formula in this case.

\subsection{The constant term of the Eisenstein series}\label{ss:E0 n=1}

  From now on,  $\cE$  denotes a line bundle on $X'$. We compute the constant Fourier coefficient of the Siegel-Eisenstein series for $\wt H_1=GU(2)$. We use notations from \S\ref{s:Eis GU}.

By definition, the constant term is equal to 
$$E_0(g,s,\Phi)=\Phi(g,s)+M(s)\Phi(g,s),
$$
where  $M(s)$ is the intertwining operator 
$$\xymatrix{M(s): I(s,(\chi,\chi_0))\ar[r]& I(-s,(\chi,\chi_0\eta^n))}
$$
defined by
$$
M(s)\Phi(g,s):=\int_{N(\BA)}\Phi(w^{-1} n(b)g,s )\, \rd n(b).
$$Since our section $\Phi$ is unramified (see \eqref{eq:unr sec}),  so is $M(s)\Phi(g,s)$. Therefore it suffices to determine the value of $M(s)\Phi(g,s)$ at $g=1$. 
By \cite[Lem.~4.3.2]{Shahidi}, translating into the current context and noting that ${\rm vol}(\wh \cO)=q^{-\frac{1}{2}\deg \omega_X}$ for the self-dual measure,  we obtain $$
M(s)\Phi(1,s)=q^{-\frac{1}{2}\deg \omega_X} \frac{L(2s,\eta)}{L(2s+1,\eta)}
$$
and hence
 $$
M(s)\Phi(g,s)=\eta(c(g))q^{-\frac{1}{2}\deg \omega_X} \frac{L(2s,\eta)}{L(2s+1,\eta)}\Phi(g,-s).
$$
Therefore 
$$
E_0(g,s,\Phi)=\Phi(g,s)+\eta(c(g))q^{-\frac{1}{2}\deg \omega_X} \frac{L(2s,\eta)}{L(2s+1,\eta)}\Phi(g,-s).
$$

\begin{remark}Note that the formula for the constant term is consistent with the functional equation
$$
E(g,s,\Phi)=E(g,-s,M(s)\Phi),
$$or equivalently $M(-s)M(s)\Phi=\Phi$. 
In fact, by the above formula on $M(s)\Phi$, we have 
$$M(-s)M(s)\Phi=
q^{-\frac{1}{2}\deg \omega_X} \frac{L(-2s,\eta)}{L(-2s+1,\eta)} \cdot q^{-\frac{1}{2}\deg \omega_X} \frac{L(2s,\eta)}{L(2s+1,\eta)} \Phi.
$$
Then $M(-s)M(s)\Phi=\Phi$ follows from the functional equation
\begin{align}\label{FE L eta}
q^{\frac{1}{2}\deg \omega_X s }L(s,\eta)=q^{\frac{1}{2}\deg \omega_X (1-s) } L(1-s,\eta),
\end{align} 
where we note that $L(s,\eta)$ is a polynomial in $q^{-s}$ of degree $\deg\omega_X$.
\end{remark}

Now we evaluate the constant term at $g=m(\cE,\LL)$ for line bundles $\cE\in\Pic_{X'}(k),\LL\in\Pic_X(k)$, 
\begin{align}\label{eq:E0}
E_0(m(\cE,\LL),s,\Phi)=&\chi(\cE) \chi_0(\LL)q^{(\deg\cE-\deg \LL)(s+1/2)}\notag
\\&+\eta(\LL) \chi(\cE) \chi_0(\LL) q^{-\frac{1}{2}\deg \omega_X} \frac{L(2s,\eta)}{L(2s+1,\eta)}q^{(\deg\cE-\deg \LL)(-s+1/2)}.
\end{align}

The normalized Eisenstein series \eqref{eq:nor Eis},  specialized to the case $n=1$, gives
$$\wt E(g,s,\Phi)=q^{\deg \om_X s} L(2s+1,\eta)E(g,s,\Phi)
$$
and \eqref{eq:deg n} becomes 
\begin{align}\label{eq:deg 1}
d=\deg\LL+\deg\omega_X-\deg\cE.
\end{align}

By \eqref{eq:E0}, when $a=0$,
\begin{align}\label{eq:wt E0}
&\wt E_0(m(\cE,\LL),s,\Phi)\notag
\\
=& (\chi_0\eta)(\LL)\chi(\cE) q^{d(s-1/2)}  L(2s,\eta) +\chi_{0}(\LL) \chi(\cE) q^{-d(s+1/2)}q^{\deg \omega_X (1/2+2s) }  L(2s+1,\eta) .
\end{align}
By the functional equation \eqref{FE L eta},  the two summands in \eqref{eq:wt E0} are switched (up to the sign $\eta(\LL)$) with respect to the substitution $s\mapsto -s$.

\subsection{The constant term of the generating series}\label{ssec: result for constant term} Fix a line bundle $\LL\in \Pic_X(k)$.
For $1\le i\le r$, let $\ell_{i}$ be the line bundle on $\Sht^{r}_{U(1),\LL}$ whose fiber at $(\{x_{i}\},\{\cF_{i}\})$ is the fiber $\cF_{i}|_{\s x_{i}}$. According to Definition \ref{def:cap Chern} and Definition \ref{def:special cycle classes}, the proposed constant term for the generating series is a sum of two terms. When $n=1$, one of them vanishes and hence we have
\begin{equation}\label{eq:Z(0)}
[\cZ^{r}_{\cE}(0)]=[\cZ^{r}_{\cE}[\cE](0)]:=\prod_{i=1}^{r}c_{1}(p_{i}^{*}\s^{*}\cE^{-1}\ot\ell_{i})\in \Ch^{r}(\Sht^{r}_{U(1),\LL})=\Ch_{0}(\Sht^{r}_{U(1),\LL}).
\end{equation}
Note that on the left hand side we have suppressed the dependence on $\LL$, for brevity.  

The goal now is to calculate the degree of $[\cZ^{r}_{\cE}(0)]$ in terms of higher derivatives of the $L$-function $L(s,\y)$.
We have $L(s,\y)=\z_{X'}(s)/\z_{X}(s)$, and it is a polynomial in $q^{-s}$ (because of our assumption that $X'/X$ is non-split) of degree $2g-2$.

\begin{thm}\label{th:U1 const term} Let $r\in\BZ_{\geq 0}$ be such that $(-1)^r=\eta(\LL)$. Then we have
\begin{equation}\label{eq:deg}
\deg[\cZ^{r}_{\cE}(0)]=2(\log q)^{-r}\frac{d^{r}}{ds^{r}}\Big|_{s=0}\left(q^{ds}L(2s,\y)\right)
\end{equation}
where $d$ is defined in \eqref{eq:deg 1}.
\end{thm}

For the brevity of notation, we will denote  
\begin{eqnarray*}
\Sht^r_{U(1)} : =\Sht^{r}_{U(1),\LL},\quad  \LN: =\LL\ot \omega_X,
\end{eqnarray*}
and then $d=\deg\LN-\deg\cE$. 
%

\subsection{Calculation of the Chern classes}
Recall the Cartesian diagram from \eqref{eq: prym cartesian square} 
\begin{equation}
\xymatrix{\Sht^{r}_{U(1)}\ar[r]^-{p} \ar[d]_{p_{[1,r]}:=(p_{1},\cdots, p_{r})}& \Prym_{ \LN}\ar[d]^{\Lang}\\
X'^{r}\ar[r]^{\AJ^{r}} & \Prym^{\epsilon(\LN)}}
\end{equation}

 Let $p_{ab}=(p_{a},p_{b}): \Sht^{r}_{U(1)}\to X'\times X'$. Let $\D\subset X'\times X'$ be the diagonal and $\D^{-}\subset X'\times X'$ be the anti-diagonal consisting of points $(x,\s x)$.

Let $\cP$ be the Poincar\'e line bundle over $X'\times \Prym_\LN$. 

In the following, all Chern classes lie in $\ell$-adic cohomology groups. Also, when we write $\cohog{i}{Z, \Qlbar}$ or $\cohog{i}{Z}$ for a stack $Z$ over $k$ we mean $\cohog{i}{Z_{\ov k}, \Qlbar}$.

\begin{lemma}\label{l:chern term a} For $1\le a\le r$ we have an equality in $\cohog{2}{\Sht^{r}_{U(1)},\Qlbar(1)}$:
\begin{equation}
c_{1}(p_{a}^{*}\s^*\cE^{-1}\ot\ell_{a})=(\s p_{a}, p)^{*}c_{1}(\cP)+ \sum_{b<a}p_{ba}^{*}c_{1}(\cO(\D-\D^{-}))-p_{a}^{*}c_{1}(\cE\ot \om_{X'}).
\end{equation}
\end{lemma}
\begin{proof}
We have
\[
\ell_{a}|_{(\{x'_{i}\}, \{\cF_{i}\})} = \cF_a/\cF^{\flat}_{a-1/2} =  (\cF_a)_{\s x'_a} =  \cF_{0}(\s x'_{1}+\cdots+\s x'_{a}-x'_{1}-\cdots- x'_{a})|_{\s x'_{a}}.
\]
Therefore 
\[
\ell_{a}\cong (\s p_{a},p)^{*}\cP\ot \left(\ot_{1\le b< a}p_{ba}^{*}\cO(\D-\D^{-}) \right)\ot p_{a}^{*}(\cO(\D-\D^{-})|_{\D}).
\]
Since $\cO(\D)|_{\D}\cong \om_{X'}^{-1}$, $\cO(\D^{-})|_{\D}\cong \cO_{X'}$ and $c_1(\s^{*}\cE)=c_1(\cE)$, we obtain the desired formula.
\end{proof}

Denote $V=\cohog{1}{X',\Qlbar}^{\s=-1}$ as a $\Frob$-module. Denote the action of $\Frob$ on $V$ by $\phi$. 
Let $\xi\in \cohog{2}{X',\Qlbar(1)}$ be the fundamental class of any closed point on $X'_{\ov k}$ and use it to identify $\cohog{2}{X',\Qlbar(1)}\simeq \Qlbar$. Let \begin{equation}
\j{\cdot,\cdot}: V\times V\to \Qlbar(-1)
\end{equation}
be the symplectic pairing on $V$ induced from the cup product, i.e.,
\begin{equation}
v\cup v'=\j{v,v'}\xi\in \cohog{2}{X',\Qlbar}, \quad v,v'\in V.
\end{equation}
For dual bases $\{v_{i}\}$ and $\{v^{i}\}$ under this symplectic pairing, i.e.,  $\sum_{i}\j{v_{i},\a}v^{i}=\a$ for all $\a\in V$. Let
\begin{equation}
\b=\sum_{i}v_{i}\ot v^{i}\in \wedge^{2}(V)\subset V\ot V.
\end{equation}

\begin{lemma}\label{l:c1 D-D}
We have  $c_{1}(\cO(\D-\D^{-}))=-2\b\in V\ot V\subset \cohog{2}{X'\times X'}$.\end{lemma}
\begin{proof} Note that the group $\Aut(X'/X)\times \Aut(X'/X)$ acts on $\cohog{2}{X'\times X'}$ and $V\ot V$ is exactly the isotypic subspace for the character $\chi: \Aut(X'/X)\times \Aut(X'/X)\to \{\pm1\}$ that is nontrivial on both factors. A straightforward computation shows that $c_{1}(\cO(\D-\D^{-}))$ must be in this isotypic subspace, therefore $c_{1}(\cO(\D-\D^{-}))\in V\ot V$.

For any class $\g\in \cohog{2}{X'\times X'}$, let $\g^{\hs}$ be the projection of $\g$ to the $\chi$-isotypic subspace $V\ot V\subset \cohog{2}{X'\times X'}$. 
We {\em claim} that $c_{1}(\cO(\D))^{\hs}=-\b$. Note that $c_{1}(\cO(\D))$ is the cycle class $cl(\D)$ of the diagonal $\D$  in $X'\times X'$,  so we need to show that $cl(\D)^{\hs}=-\b$. Under the Kunneth decomposition and the Poincar\'e duality,  $cl(\D)$ corresponds to the identity endomorphism of $\cohog{*}{X'}$. In particular, 
\begin{equation}
cl(\D)^{\hs}\cup (\a\ot 1)=\xi\ot \a\in \cohog{2}{X'}\ot V, \quad\forall \a\in V.
\end{equation}
(Here $1\in H^0(X')$ is the fundamental class of $X'$.) This property characterizes $cl(\D)^{\hs}$. To show $cl(\D)^{\hs}=-\b$, it suffices to check that
\begin{equation}
-\b\cup (\a\ot 1)=\xi\ot \a, \quad\forall \a\in V.
\end{equation}
This holds because $-(\sum_{i}v_{i}\ot v^{i})\cup (\a\ot 1)=\sum_{i} (v_{i}\cup \a)\ot v^{i}$ (by Koszul sign convention), which is $\sum_{i} \j{v_{i}, \a}\xi\ot v^{i}=\xi\ot(\sum_{i}\j{v_{i},\a}v^{i})=\xi\ot \a$.

Finally, by $\cO(\D^-)\simeq (\sigma,1)^{*}\cO(\D)$ we know that  $c_{1}(\cO(\D^-))^{\hs}=-c_{1}(\cO(\D))^{\hs}=\beta$. Therefore $c_{1}(\cO(\D-\D^{-}))=c_{1}(\cO(\D-\D^{-}))^{\hs}=-2\b$.

\end{proof}

The Abel-Jacobi map 
\begin{equation}\begin{aligned}
\AJ_{1}: X'& \to \Prym^{1}\\
 x'& \mapsto \cO(\sigma x'-x')
\end{aligned}
\end{equation}
induces an injective map $\AJ_{1}^{*}$ on $\upH^{1}$ and identifies the image:
\begin{equation}
\cohog{1}{\Prym^{1}}\isom \cohog{1}{X'}^{\s=-1}=V.
\end{equation}

We claim that  $\cohog{1}{\Prym^{0}}$ as well as $\cohog{1}{\Prym^{\e}_{\frN}}$ for any $\frN\in \Pic_{X}(k)$ are {\em canonically} identified with $V$. Indeed, if $A$ is a (geometrically) connected group scheme over $k$ and $A_{1}$ is any $A$-torsor over $k$, then any choice of $b\in A_{1}(\ov k)$ identifies $A_{\ov k}$ with $A_{1,\ov k}$ hence gives an isomorphism $\cohog{*}{A}\cong \cohog{*}{A_{1}}$. Different choices of $b$ give the same isomorphism because $b$ varies in $A_{1}$ which is geometrically connected.  Applying this principle to $A=\Prym^{0}$ and $A$-torsors $\Prym^{\e}_{\frN}$ and $\Prym^{1}$, we see that there are canonical isomorphisms 
\begin{equation}\label{isomV}
\cohog{1}{\Prym^{\e}_{\frN}}\cong \cohog{1}{\Prym^{0}}\cong \cohog{1}{\Prym^{1}}\cong V.
\end{equation}

\begin{lemma}\label{l:c1P} 
For $\e\in\Irr(\Prym_{\LN})$, we have
$$c_{1}(\cP)|_{X'\times\Prym_\LN^{\ep}}=2\b+\deg_X\LN \,( \xi\ot 1)$$ where  $\b\in V\ot V\subset \cohog{1}{X'}\ot \cohog{1}{\Prym_\LN^{\ep}}\subset \cohog{2}{X'\times\Prym_\LN^{\ep}}$ and $\xi\ot 1\in \cohog{2}{X'}\ot \cohog{0}{\Prym_\LN^{\ep}}\subset \cohog{2}{X'\times\Prym_\LN^{\ep}}$.
\end{lemma}
\begin{proof} Choose $\LN'\in\Prym_\LN^{\ep}(\bar k)$ (in particular $\Nm(\LN')\cong \LN$). Pulling back by $(\Id_{X'}, \AJ_{1})$ to $X'\times X'$, $\cP$ becomes $\cO(\D^--\D)\otimes \pr_1^\ast \LN'$ where $\pr_1:X'\times X'\to X'$ is the  projection to the first factor. Then note that $c_1( \pr_1^\ast \LN')=\deg_{X'}\LN' \,(\xi\ot 1)=\deg_{X}\LN\,(\xi\ot 1)\in  \cohog{2}{X'}\ot \cohog{0}{X'}$. Now the lemma follows from Lemma \ref{l:c1 D-D}.
\end{proof}

\begin{lemma}\label{l:pullback c1P} For $1\le a\le r$, we have
\begin{align*}
-\frac{1}{2}(\s p_{a},p)^{*}c_{1}(\cP)=&-\frac{1}{2}\deg_X \LN \cdot \ p_{a}^{*}\xi+ \Tr((\phi-1)^{-1}|V)p_{a}^{*}\xi   \notag\\
&+\sum_{b<a}p^{*}_{ba}(((\phi-1)^{-1}\ot 1)\b)+\sum_{b>a}p^{*}_{ab}((1\ot (\phi-1)^{-1})\b).
\end{align*}
\end{lemma}
\begin{proof}
We first have 
\begin{align}\label{eq:pb xi}
(\s p_{a},p)^{*} (\xi\ot 1)=p_a^* \xi.
\end{align}
Next we use the commutative diagram
\begin{equation}
\xymatrix{\Sht^{r}_{U(1)}\ar[r]^-{(p_{a},p)}\ar[d]_{p_{[1,r]}} & X'\times \Prym_\LN\ar[d]^{\Id\times\Lang}\\
X'^{r}\ar[r]^-{(\pr_{a}, \AJ^{r})} & X'\times \Prym^{\ep(\LN)}}
\end{equation}
Here $\pr_{a}: X'^{r}\to X'$ is the $a$-th projection. The pullback along the Lang map  $\cohog{1}{\Prym^{\ep(\LN)}}\to \cohog{1}{\Prym_\LN^{\ep}}$ is the isomorphism $\phi-1$ of $V$ under the isomorphisms \eqref{isomV}, for each component $\epsilon\in \Irr(\Prym_{\LN})$. Therefore, 
\begin{equation}
\b=(\Id_{X'}\times\Lang)^{*}(1\ot (\phi-1)^{-1})\b.
\end{equation}
Here we view $(1\ot (\phi-1)^{-1})\b\in V\ot V$ as an element of $\cohog{2}{X'\times \Prym^{\ep(\LN)}}$. Hence by the above commutative diagram
\begin{align}
(p_{a},p)^{*}\b=&(p_{a},p)^{*}(\Id_{X'}\times\Lang)^{*}(1\ot (\phi-1)^{-1})\b\notag\\
\label{pullback Poin}=&p_{[1,r]}^{*}(\pr_{a}, \AJ^{r})^{*}(1\ot (\phi-1)^{-1})\b.
\end{align}
Since $\AJ^{r}$ can be decomposed into a composition $X'^{r}\xr{\AJ_{1}^{r}}(\Prym^{1})^{r}\xr{m}\Prym^{\ep(\LN)}$ (the map $m$ is multiplication), and 
$$
m^{*}(v)=1\ot \cdots \ot1\ot v+\cdots+v\ot1\ot\cdots\ot 1
$$
for $v\in V=\cohog{1}{\Prym^{\ep(\LN)}}$, we get $\AJ^{r*}v=\sum_{1\le b\le r}\pr_{b}^{*}v\in \cohog{1}{X'^{r}}$. Hence
\begin{equation}\label{sum ab}
p_{[1,r]}^{*}(\pr_{a}, \AJ^{r})^{*}(1\ot (\phi-1)^{-1})\b=\sum_{i}p_{a}^{*}v_{i}\ot\left(\sum_{b=1}^{r}p_{b}^{*}(\phi-1)^{-1}v^{i}\right).
\end{equation}
Note that when $b<a$, the term $p_{a}^{*}v_{i}\ot p_{b}^{*}(\phi-1)^{-1}v^{i}$ is $-p_{ba}^{*}((\phi-1)^{-1}v^{i}\ot v_{i})$; but after summing over $i$, using that $\sum_{i}v^{i}\ot v_{i}=-\b$,  we obtain 
$$-p_{ba}^{*}((\phi-1)^{-1}\ot 1)(-\b)=p_{ba}^{*}((\phi-1)^{-1}\ot 1)(\b).$$ 
When $b>a$, the corresponding term in \eqref{sum ab} is $p_{ab}^{*}(1\ot (\phi-1)^{-1})(\b)$. When $b=a$, the corresponding term in \eqref{sum ab} is $$\sum_{i}p_{a}^{*}(v_{i}\cup (\phi-1)^{-1}v^{i})=\Tr((\phi-1)^{-1}|V)p_{a}^{*}\xi.$$ 
Here we are using that the bases $\{v_{i}\}$ and $\{v^{j}\}$ satisfy $v_{i}\cup v^{j}=\d_{ij}\xi$. Combining these with \eqref{pullback Poin} and \eqref{sum ab} we get
\begin{equation}\label{papb}
(p_{a},p)^{*}\b=\Tr((\phi-1)^{-1}|V)p_{a}^{*}\xi+\sum_{b<a}p^{*}_{ba}(((\phi-1)^{-1}\ot 1)\b)+\sum_{b>a}p^{*}_{ab}((1\ot (\phi-1)^{-1})\b).
\end{equation}

Since the action of $\s$ on $V$ is by $-1$, we have 
\begin{align}\label{eq:s beta}
(\s p_{a},p)^{*} (\beta)=(p_{a},p)^{*}(\sigma, \Id)^{*}\beta=-(p_{a},p)^{*} (\beta).
\end{align}
Combining this with \eqref{papb}, \eqref{eq:pb xi} and Lemma \ref{l:c1P}, we get the desired identity. 
\end{proof}

\subsection{Taylor expansion of $L(s,\eta)$}
\begin{lemma}\label{l:fn} Let $\a\in \CC$.  Write the $n$-th derivative of $\log(1-\a q^{-s})$ as
\begin{equation}
(-1)^{n-1}(\log q)^{n}\frac{f_{n}(\a q^{-s})}{(1-\a q^{-s})^{n}}, \quad n=1,2,\cdots.
\end{equation}
where $f_{n}(x)$ is a polynomial in $x$. Then
\begin{equation}
f_{n}(x)=\sum_{c\in C_{n}} x^{\d(c)}.
\end{equation}
Here $C_{n}$ is the set of cyclic permutations (i.e., with only one cycle) on $\{1,2,\cdots, n\}$; for $c\in C_{n}$, $\d(c)$ is the number of $1\le i\le n$ such that $c(i)\le i$ (when $n=1$, $\d(c)=1$).
\end{lemma}
\begin{proof}
From the definition we get a recursive relation:
\begin{equation}
f_{n+1}(x)=nxf_{n}(x)+x(1-x)f'_{n}(x).
\end{equation}
Also $f_{1}(x)=x$. From this it is easy to see that $\deg f_{n}(x)\le n$, 
and $f_{n}(0)=0$.
Write $f_{n}(x)=a^{(n)}_{1}x+a^{(n)}_{2}x^{2}+\cdots +a^{(n)}_{n}x^{n}$. Then
\begin{equation}\label{ani}
a^{(n+1)}_{i}=ia_{i}^{(n)}
+(n+1-i)a_{i-1}^{(n)}
, \quad i=1,\cdots, n.
\end{equation}
On the other hand, let $C_{n,i}$ be the set of $c\in C_{n}$ such that $\d(c)=i$. We must show that $|C_{n,i}|=a^{(n)}_{i}$. We do this by checking that $|C_{n,i}|$ satisfies the same recursive relation \eqref{ani}. 

For $c\in C_{n+1}$, let $1\le i_{c}, j_{c}\le n$ be defined by $c(i_{c})=n+1$ and $c(n+1)=j_{c}$. 
We have a map $\pi: C_{n+1}\to C_{n}$ sending $c\in C_{n+1}$ to $c'\in C_{n}$ defined by $c'(i)=c(i)$ if $i\ne i_{c}$ and $c'(i_{c})=j_{c}$. We decompose $C_{n+1,i}=C'_{n+1,i}\sqcup C''_{n+1,i}$, where $C'_{n+1,i}$ is the set of $c\in C_{n,i}$ such that $i_{c}>j_{c}$. Then $\pi$ restricts to a $i$ to $1$ map $\pi': C'_{n+1,i}\to C_{n,i}$ (the preimage of $c'$ are in bijection with $i$ such that $c'(i)<i$) and an $(n+1-i)$ to $1$ map $\pi'': C''_{n+1,i}\to C_{n,i-1}$ (the preimage of $c'\in C_{n,i-1}$ are in bijection with $i$ such that $c'(i)>i$). This shows
\begin{equation}
|C_{n+1,i}|=|C'_{n+1,i}|+|C''_{n+1,i}|=i|C_{n,i}|+(n+1-i)|C_{n,i-1}|.
\end{equation}
This shows that $|C_{n,i}|$ satisfies the same recursive relation \eqref{ani} as $(a^{(n)}_{i})$. Since the initial values match $|C_{1,1}|=1=a^{(1)}_{1}$, the lemma follows.
\end{proof}

\begin{cor}\label{c:Taylor log} The Taylor expansion of $\log L(s,\y)$ at $s=0$ is:
\begin{equation}
\log L(s,\y)=\log L(0,\y)-\sum_{\ell\ge1}\sum_{c\in C_{\ell}}\Tr\left(\frac{\phi^{\d(c)}}{(1-\phi)^{\ell}}\Big |V\right)\frac{(\log q)^{\ell} (-s)^{\ell}}{\ell!}.
\end{equation}
\end{cor}
\begin{proof} Let $\{\a_{i}\}$ be the multiset of eigenvalues of $\phi$ acting on $V$. By Lemma \ref{l:fn} evaluated at $s=0$, we get
\begin{equation}
\log(1-\a_{i}q^{-s})=\log(1-\a_{i})-\sum_{\ell\ge1}\sum_{c\in C_{\ell}} \frac{\a_{i}^{\d(c)}}{(1-\a_{i})^{\ell}}\frac{(\log q)^{\ell}(- s)^{\ell}}{\ell!}. 
\end{equation}
Taking sum over $\a_{i}$, noting that $L(s,\y)=\prod_{i}(1-\a_{i}q^{-s})$ and $\sum_{i}\a_{i}^{\d(c)}/(1-\a_{i})^{\ell}=\Tr(\frac{\phi^{\d(c)}}{(1-\phi)^{\ell}}|V)$, we get the desired formula.
\end{proof}

\subsection{Proof of Theorem \ref{th:U1 const term}}

Combining Lemma \ref{l:chern term a}, Lemma \ref{l:c1 D-D} and Lemma \ref{l:pullback c1P}, we get
\begin{eqnarray}\label{c1a}
&&- \frac{1}{2}c_{1}(p_{a}^{*}\s^*\cE^{-1}\ot\ell_{a})\\
\notag&=&
p_{[1,r]}^{*}\left(\sum_{b<a}\pr^{*}_{ba}((\phi(\phi-1)^{-1}\ot 1)\b)+\sum_{b>a}\pr^{*}_{ab}((1\ot (\phi-1)^{-1})\b)\right)\\
\notag&+&p_{[1,r]}^{*}\left(\left(\Tr(\phi(\phi-1)^{-1}|V)-d/2\right)\pr_{a}^{*}\xi\right).
\end{eqnarray}
Here we are using \eqref{eq:deg 1} and
\begin{eqnarray*}
-c_{1}(\cE\ot \om_{X'})-2\Tr((\phi-1)^{-1}|V)\xi = -2\left(\frac{\deg\cE}{2}+\frac{\deg\om_{X'}}{2}+\Tr((\phi-1)^{-1}|V)\right)\xi.
\end{eqnarray*}
Note that $\frac{\deg\om_{X'}}{2}=\deg\om_{X}=\dim V$, hence the last two terms combine to give $\Tr(1+(\phi-1)^{-1}|V)=\Tr(\phi(\phi-1)^{-1}|V)$.

Taking the product of $-\frac{1}{2}c_{1}(p_{a}^{*}\s^{*}\cE^{-1}\ot\ell_{a})$ over all $1\le a\le r$ (the order does not matter because these classes have even degree),  using \eqref{c1a} and extracting the coefficient of $p_{[1,r]}^{*}\xi^{r}$ we get
\begin{equation}\label{prod c1a}
\prod_{a=1}^{r}-\frac{1}{2}c_{1}(p_{a}^{*}\s^{*}\cE^{-1}\ot\ell_{a})=(\sum_{g\in S_{r}}A_{g})p_{[1,r]}^{*}(\xi^{\ot r})
\end{equation}
where $A_{g}\in \Qlbar$ is defined as
\begin{align}\label{eq:def Ag}
A_{g} \,\xi^{\ot r}=&\prod_{g(a)<a}(\phi(\phi-1)^{-1}\ot 1)\b _{g(a)a}\prod_{g(a)>a}(1\ot (\phi-1)^{-1})\b_{ag(a)}\\
\notag &\quad \times \prod_{g(a)=a}\left(\Tr(\phi(\phi-1)^{-1}|V)-d/2\right)\xi_{a}.
\end{align}
Here  we use the  abbreviations  $(-)_{ba}=\pr_{ba}^{*}(-), (-)_{a}=\pr_{a}^{*}(-)$.  When $r=0$ we understand the sum $\sum_{g\in S_{r}}A_{g}$ as $1$.   

We form the generating series of $A_{g}$ for $g\in S_{r}$ for all $r\ge0$. Our aim is to show 
\begin{equation}\label{Ag der}
\sum_{r\ge0, g\in S_{r}}A_{g}\frac{(\log q)^{r}(-s)^{r}}{r!}=q^{ds/2}L(s,\y)L(0,\y)^{-1}.
\end{equation}
Indeed if this holds, then making a change of variables 
$s\mapsto 2s$ and extracting the coefficient of $s^{r}$ we get
\begin{equation}
(-2\log q)^{r}\sum_{g\in S_{r}}A_{g}=L(0,\y)^{-1}\frac{d^{r}}{ds^{r}}\left(q^{ds}L(2s,\y)\right)\Big|_{s=0}.
\end{equation}
Taking the degrees of both sides of \eqref{prod c1a} we get
\begin{equation}\label{Zr Ag}
\deg[\cZ^{r}_{\cE}(0)]=(-2)^r |\Prym(\F_{q})|\sum_{g\in S_{r}}A_{g}.
\end{equation}
Here the factor $|\Prym(\F_{q})|$ is the degree of $p_{[1,r]}:\Sht^{r}_{U(1)}\to X'^{r}$. Using \eqref{Zr Ag} and the fact that $|\Prym(\F_{q})|=2L(0,\y)$ we get
\begin{equation}
\deg[\cZ^{r}_{\cE}(0)]=2L(0,\y)\cdot (-2)^{r}\sum_{g\in S_{r}}A_{g}=2(\log q)^{-r}\frac{d^{r}}{ds^{r}}\left(q^{ds}L(2s,\y)\right)\Big|_{s=0}.
\end{equation}
This is exactly \eqref{eq:deg} and hence Theorem \ref{th:U1 const term} is proved.

Now it remains to prove \eqref{Ag der}. Let $C(g)$ be the set of cycles of $g$. For each $c\in C(g)$, let $\ell(c)$ be its length and recall $\d(c)$ is the number of $1\le i\le n$ such that $c(i)$ is defined and $c(i)\le i$. We claim that we can write $A_{g}=\prod_{c\in C(g)}A_{g,c}$ where
\begin{equation}\label{eq:Agc}
A_{g,c}=\begin{cases}-\Tr(\phi(1-\phi)^{-1}|V)-d/2,& \ell(c)=1,\\
-\Tr(\phi^{\d(c)}(1-\phi)^{-\ell(c)}|V), & \ell(c)>1.\end{cases}
\end{equation}
Indeed, suppose a cycle $c=(a,g(a),\cdots, g^{\ell-1}(a))$ has length $\ell=\ell(c)$. If $\ell=1$ then $a$ is a fixed point of $g$ and the factor corresponding to such a fixed point can be directly read from the definition of $A_{g}$ in \eqref{eq:def Ag}. If $\ell>1$, write $a_{s}=g^{s-1}(a)$ for $s=1,2,\cdots, \ell$. We assume that $a$ is the largest element in the cycle. If $a_{s}>a_{s+1}$, the corresponding factor $(\phi(\phi-1)^{-1}\ot 1)\b _{a_{s+1}a_{s}}=\sum_{i}\phi(\phi-1)^{-1}v_{i, a_{s+1}}\ot v^{i}_{a_{s}}$ (recall that $v^{i}_{a_{s}}$ means $v^{i}$ put in the $a_{s}$-th factor of $\cohog{*}{X'}^{\ot r}$). We rewrite it as
\begin{equation}\label{as dec}
\sum_{i}v^{i}_{a_{s}} \ot\phi(1-\phi)^{-1}v_{i, a_{s+1}}
\end{equation}
where switching the terms produces a minus sign which cancels with the change from $\phi-1$ to $1-\phi$. Similarly, if $a_{s}<a_{s+1}$, writing $\b=-\sum_{i}v^{i}\ot v_{i}$, the corresponding factor $(1\ot (\phi-1)^{-1})\b _{a_{s}a_{s+1}}$ is
\begin{equation}\label{as inc}
\sum_{i}v^{i}_{a_{s}}\ot(1-\phi)^{-1}v_{i, a_{s+1}}.
\end{equation}
For $1\le s\le \ell$, let $T^{(s)}\in \End(V)$ be defined as
\begin{equation*}
T^{(s)}=\begin{cases} \phi(1-\phi)^{-1},  & \mbox{if }a_{s}>a_{s+1}\\
(1-\phi)^{-1}, & \mbox{if }a_{s}<a_{s+1}.
\end{cases}
\end{equation*}
Then the product of the terms \eqref{as dec} or \eqref{as inc} for $s=1,\cdots, \ell$ (in that order) is
\begin{equation}\label{mult trace}
\sum_{i_{1},i_{2},\cdots, i_{\ell}}v^{i_{1}}_{a_1}\ot (T^{(1)}v_{i_{1}}\cup v^{i_{2}})_{a_{2}}\ot (T^{(2)}v_{i_{2}}\cup v^{i_{3}})_{a_{3}}\ot\cdots\ot (T^{(\ell-1)}v_{i_{\ell-1}}\cup v^{i_{\ell}})_{a_{\ell}}\ot (T^{(\ell)}v_{i_{\ell}})_{a_\ell}.
\end{equation}

For any endomorphism $T$ of $V$,  we have $Tv_{i}\cup v^{j}=T_{ji}\xi$ where $T_{ji}$ is the $(i,j)$-entry of the matrix of $T$ under the basis $\{v_{i}\}$. The product in \eqref{mult trace} is then a multiple of $\xi^{\ot r}$, and the multiple is the negative of the trace of the operator  $\prod_{s=1}^{\ell} T^{(s)}\in \End(V)$. The negative sign comes from taking the cup product of the first and the last factor since $v^{i_{1}}_{a}\cup T^{(\ell)}v_{i_{\ell},a}=-T^{(\ell)}_{i_{1},i_{\ell}}\xi$. Since
\begin{equation*}
\prod_{s=1}^{\ell}T^{(s)}=(\phi(1-\phi)^{-1})^{\d(c)}((1-\phi)^{-1})^{\ell(c)-\d(c)}=\phi^{\d(c)}(1-\phi)^{-\ell(c)},
\end{equation*}
the identity \eqref{eq:Agc} follows.


The formula \eqref{eq:Agc} depends only on the cyclic permutation $c$ on an {\em ordered} set. We write $A_{g,c}$ as $B_{c}$ with the understanding that the ordered set on which $c$ operates is a subset of $\N$. Now we re-organize the sum over $g\in S_{r}$ by grouping first according to the partitions of the set $\{1,2,\cdots, r\}$ and then according to the conjugacy classes. We have surjections
\begin{equation}
\pi: S_{r}\xr{\pi_{1}} \Pi_{r} \xr{\pi_{2}} P_{r}
\end{equation}
where $\Pi_{r}$ is the set of partitions of the set $\{1,2,\cdots, r\}$, and $P_{r}$ is the set of partitions of $r$. The map $\pi_{1}$ takes $g\in S_{r}$ to its cycles, and $\pi_{2}$ takes the lengths of the cycles. For $I_{\bu}\in \Pi_{r}$, corresponding to a partition $\{I_{\a}\}$ of $\{1,\cdots, r\}$, the contribution of $\pi_{1}^{-1}(I_{\bu})$ to $\sum_{g\in S_{r}}A_{g}x^{r}/r!$ is
\begin{equation}
\Sigma_{I_{\bu}}:=\frac{\prod_{\a}|I_{\a}|!}{r_{!}}\prod_{\a}\left(\sum_{c\in C(I_{\a})}B_{c}\frac{x^{|I_{\a}|}}{|I_{\a}|!}\right)
\end{equation}
where the sum is over the set $C(I_{\a})$ of cyclic permutations of $I_{\a}$. Clearly the sum $\sum_{c\in C(I_{\a})}B_{c}$ depends only on the cardinality $|I_{\a}|$ and not on the ordering of $I_{\a}$. Denote
\begin{equation}
\Gamma_{\ell}:=\sum_{c\in C_{\ell}}B_{c}.
\end{equation}
(Here recall from Lemma \ref{l:fn} that $C_{\ell}$ is the set of cyclic permutations on $\{1,2,\cdots, \ell\}$.) Write $\l:=\pi_{2}(I_{\bu})\in P_{r}$ as $\l_{1}^{m_{1}}\cdots \l_{t}^{m_{t}}$, where $\l_{1}>\cdots>\l_{t}$, and $m_{i}$ is the multiplicity of $\l_{i}$, then 
\begin{equation}
\Sigma_{I_{\bu}}=\frac{\prod (\l_{i}!)^{m_{i}}}{r_{!}}\prod_{i}\left(\Gamma_{\l_{i}}\frac{x^{\l_{i}}}{\l_{i}!}\right)^{m_{i}}.
\end{equation}
In particular, $\Sigma_{I_{\bu}}$ depends only on the partition $\pi_{2}(I_{\bu})\in P_{r}$.  Therefore the contribution of $\l\in P_{r}$ to $\sum_{g\in S_{r}}A_{g}x^{r}/r!$ is
\begin{equation}
|\pi_{2}^{-1}(\l)|\frac{\prod (\l_{i}!)^{m_{i}}}{r_{!}}\prod_{i}\left(\Gamma_{\l_{i}}\frac{x^{\l_{i}}}{\l_{i}!}\right)^{m_{i}}=\prod_{i}\frac{1}{ m_{i}!}\left(\Gamma_{\l_{i}}\frac{x^{\l_{i}}}{\l_{i}!}\right)^{m_{i}}.
\end{equation}
Here we are using $|\pi_{2}^{-1}(\l)|=|O_{\l}|/|\pi_{1}^{-1}(I_{\bu})|$ (where $O_{\l}\subset S_{r}$ is the conjugacy class corresponding to $\l$), $|O_{\l}|=r!/(\prod \l_{i}^{m_{i}} m_{i}! )$ and $|\pi_{1}^{-1}(I_{\bu})|=\prod_{i}((\l_{i}-1)!)^{m_{i}}$. 

Summing over all partitions of $r$ and then over all $r\ge 0$, we get
\begin{equation}\label{Ag exp}
\sum_{r\ge0,g\in S_{r}}A_{g}\frac{x^{r}}{r!}=\prod_{\ell\ge1}\sum_{m\ge0}\frac{1}{m!}\left(\Gamma_{\ell}\frac{x^{\ell}}{\ell!}\right)^{m}=\exp\left(\sum_{\ell\ge1}\Gamma_{\ell}\frac{x^{\ell}}{\ell!}\right).
\end{equation}

Using the formula \eqref{eq:Agc} for $A_{g,c}$ we have
\begin{equation}
\Gamma_{\ell}=\begin{cases}-\Tr\left(\frac{\phi}{1-\phi}\Big|V\right)-d/2,  & \ell=1\\ \sum_{c\in C_{\ell}}-\Tr\left(\frac{\phi^{\d(c)}}{(1-\phi)^{\ell}}\Big|V\right), & \ell>1.\end{cases}
\end{equation}

Plugging into \eqref{Ag exp} we get
\begin{equation}\label{exp sum}
\sum_{r\ge0,g\in S_{r}}A_{g}\frac{x^{r}}{r!}=\exp\left(\left(-\Tr\left(\frac{\phi}{1-\phi}\Big|V\right)-d/2\right)x-\sum_{\ell\ge2}\sum_{c\in C_{\ell}}\Tr\left(\frac{\phi^{\d(c)}}{(1-\phi)^{\ell}}|V\right)\frac{x^{\ell}}{\ell!}\right).
\end{equation}
By Corollary \ref{c:Taylor log}, and letting $x=-(\log q)s$, we have
\begin{eqnarray*}
&&\left(-\Tr\left(\frac{\phi}{1-\phi}\Big|V\right)-d/2\right)(\log q)(-s)-\sum_{\ell\ge2}\sum_{c\in C_{\ell}}\Tr\left(\frac{\phi^{\d(c)}}{(1-\phi)^{\ell}}\Big|V\right)\frac{(\log q)^{\ell}(-s)^{\ell}}{\ell!}\\
&=&\log L(s,\y)-\log L(0,\y)+d/2\,(\log q)s.
\end{eqnarray*}
Taking the exponential, and plugging into \eqref{exp sum} we get
\begin{equation*}
\sum_{r\ge0,g\in S_{r}}A_{g}\frac{(\log q)^{r}(-s)^{r}}{r!}=q^{ds/2}L(s,\y)\cdot L(0,\y)^{-1},
\end{equation*}
which is exactly \eqref{Ag der}. \qed

\subsection{The complete comparison}\label{ss:comparison n=1}
We now take the definition of Eisenstein series in \S\ref{ss:Eis}. We make the following choices of characters:
\begin{itemize}
\item $\chi_0=\eta$;
\item $\chi$ is any character on $\Pic_{X'}(k)$ such that $\chi|_{\Pic_{X}(k)}=\y$.
\end{itemize}

Recall from Definition \ref{def: gen series before descent} and \eqref{Zm function} the generating series 
$$\wt Z^{r}_{1}: \Bun_{\wt P_{1}}(k)=\wt P_{1}(F)\bs \wt H_{1}(\BA)/\wt H_{1}(\wh \cO)\to \Ch_{0,c}(\Sht^{r}_{GU(1)}).$$
Note $\wt Z^{r}_{1}(g)$ is compactly supported because its support is contained in $\Sht^{r}_{U(1), \frL}$ (where $\frL=c(g)$) which is proper. 

\begin{thm}\label{th:n=1} We have for all $g\in\wt H_1(\BA)$,
\begin{equation}\label{eq: intro higher derivatives}
\frac{1}{(\log q)^r}  \left(\frac{d}{ds}\right)^r\Big|_{s=0} \left( \wt  E(g,s,\Phi) \right)  =  \deg \wt Z^{r}_{1}(g).
\end{equation}
\end{thm}

\begin{proof} Since both sides are $\wt H_1(\wh \cO)$-invariant, it suffices to show the Fourier expansions at $g=m(\cE,\LL)$ match term-wise: 
\begin{equation}\label{eq:der norm Eis}
\frac{1}{(\log q)^r}  \left(\frac{d}{ds}\right)^r\Big|_{s=0} \wt E_{a}(m(\cE,\LL),s,\Phi)= \chi(\cE)q^{-d /2} \deg [\cZ^{r}_{\cE,\LL}(a)]
\end{equation}
for every $\cE\in \Pic_{X'}(k), \frL\in \Pic_{X}(k)$ and $a\in \cA_{\cE,\LL}(k)$. Here $d$ is as in \eqref{eq:deg 1}. We may further assume that $(-1)^r=\eta(\LL)$, since otherwise both sides vanish.

When $a\neq 0$, by \eqref{eq:nor Eis a} specialized to $\chi_0=\eta$, we have
\begin{align*}
\wt E_a(m(\cE,\LL),s,\Phi)=\chi(\cE)q^{-d/2} q^{d s }\Den(q^{-2s}, \coker(a)).
\end{align*}
Then \eqref{eq:der norm Eis} follows from Theorem \ref{th:ns} specialized to $n=1$, which relates the degree of $\cZ^{r}_{\cE,\LL}(a)$ to the local density.

It remains to consider the case $a=0$.  By \eqref{eq:wt E0} specialized to $\chi_0=\eta$ and the symmetry with respect to $s\mapsto -s$, we have
\begin{align}\label{eq:wt E0 eta}
\left(\frac{d}{ds}\right)^r\Big|_{s=0} \wt E_0(m(\cE,\LL),s,\Phi)=2 \chi(\cE)  q^{-d/2}\left(\frac{d}{ds}\right)^r\Big|_{s=0}  q^{ds}L(2s,\eta) .
\end{align}
On the geometric side, since $\cE$ is a line bundle, there are two terms in the decomposition \eqref{decomp Z kernel}:
\begin{equation*}
\cZ^r_{\cE, \LL}(0) =\cZ_{\cE, \LL}^r[\cE](0)  \coprod \cZ^r_{\cE, \LL}(0)^\circ,
\end{equation*}
where the first term is isomorphic to $\Sht_{U(1), \LL}^r$. Correspondingly, in Definition \ref{def:special cycle classes}, there are two terms in $[\cZ^{r}_{\cE, \LL}( 0)]$ in this case. Since the rank of $\cF_i$ is $n=1$, an injective $\cE\to \cF_i$ must give rise to non-zero $a: \cE\to \s^*\cE^\vee\ot\nu^*\LL$. It follows that the stack $\cZ^{r}_{\cE, \LL}( 0)^{\c}$ is empty. Hence there is only one term left, i.e.,  $[\cZ^{r}_{\cE, \LL}( 0)]=[\cZ^r_{\cE, \LL}[\cE](0)]$, which is defined by the Chern classes of the tautological line bundles. This term has the desired degree  
by Theorem \ref{th:U1 const term} and \eqref{eq:wt E0 eta}.
This completes the proof.
 
\end{proof}

\begin{cor}\label{th:n=1 mod} 
The generating series $g\in \wt P_{1}(F)\bs \wt H_1(\BA)/\wt H_{1}(\wh\cO)\mapsto \wt Z^{r}_{1}(g)$ is automorphic, i.e., it is left $\wt H_1(F)$-invariant and hence descends to a map 
$$Z^{r}_{1}: \Bun_{GU(2)}(k)\to \Ch_{0}(\Sht^{r}_{GU(1)}).$$
In other words, Conjecture \ref{c:mod} holds for $n=m=1$.
\end{cor}

\begin{proof} The case $r=0$ is classical and follows from the modularity of theta functions (proved by Poisson summation).

Now consider the case $r>0$. Let $g\in \wt H_{1}(\BA)$ with similitude $c(g)\in \BA^{\times}$ corresponding to $\frL\in\Pic_{X}(k)$. Then $\wt Z^{r}_{1}(g)\in \Ch_{0}(\Sht^{r}_{U(1),\frL})$. By Corollary \ref{cor: chow of Sht} below, the (component-wise) degree map induces an isomorphism $\Ch_0(\Sht^{r}_{U(1),\LL}) \xrightarrow{\sim} \Q^{\pi_0 (\Sht^{r}_{U(1),\LL})}$. Hence it suffices to show that $Z^{r}_{1}$ is automorphic after composing with component-wise degree.

Assume that $\Sht^r_{U(1),\LL}$ is non-empty (otherwise the statement is vacuously true). According to Lemma \ref{lem: 2 connected components}, $\Sht^r_{U(1),\LL}$ has two connected components if $r$ is even (and positive) and one connected component when $r$ is odd. When $r$ is odd, Theorem \ref{th:n=1} implies that $\deg \cZ_{\cE, \LL}^r(a)$ is equal to $\frac{1}{(\log q)^r}  \left(\frac{d}{ds}\right)^r\Big|_{s=0} \left( \wt  E(g,s,\Phi) \right)$, which is automorphic in $g$. 

When $r>0$ is even, we claim that $\deg \cZ_{\cE, \LL}^r(a)  =  \frac{1}{2} \frac{1}{(\log q)^r}  \left(\frac{d}{ds}\right)^r\Big|_{s=0} \left( \wt  E(g,s,\Phi) \right)$ on both components of $\Sht^{r}_{U(1),\frL}$, hence is also automorphic. For $a \neq 0$, it follows from Theorem \ref{thm: deg equal on components}. For $a=0$, it is immediate from the calculation of Theorem \ref{th:U1 const term} that the degrees of $\cZ_{\cE, \LL}^r(0)$ on both components of $\Sht^r_{U(1),\LL}$ are the same.

\end{proof}

\subsubsection{Chow groups of zero-cycles}\label{sssec: chow groups of zero-cycles}

For a stack $\cY$ over $k$, we denote by $\Ch_0(\cY)^{\deg 0}$ the subgroup of $\Ch_0(\cY)$ whose degree on each \emph{proper} connected component of $\cY$ vanishes. We will show that $\Ch_0(\cY)^{\deg 0}$ vanishes for DM stacks satisfying mild conditions.

\begin{lemma}\label{lem: path connected}
Let $Y$ be a quasi-compact connected scheme of finite type over a field. Then any zero-cycle on $Y$ lies on a connected (but possibly reducible) curve contained in $Y$. 
\end{lemma}

\begin{proof}
If $Y$ is quasi-projective, then the result follows from \cite[Corollary 1.9]{CP16}. In general, we may cover $Y$ by a finite number of affine varieties $U_1, \ldots, U_s$. Without loss of generality, we may enlarge our zero-cycle $D$ so that whenever $U_i \cap U_j$ is non-empty, then $D \cap (U_i \cap U_j) $ is also non-empty. By the quasi-projective case, for each $i$ we may find a connected curve $C_i$ containing $D \cap U_i$. Then $\bigcup C_i$ is a connected curve containing $D$. 
\end{proof}

\begin{lemma}\label{lem: quasiprojective chow 0}
Let $Y$ be a quasi-compact separated scheme of finite type over $\F_q$. Then $\Ch_0(Y)^{\deg 0} = 0$. 
\end{lemma}

\begin{proof}
We immediately reduce to the case where $Y$ is connected. Next we will reduce to the proper connected case. A compactification $\ol{Y}$ of $Y$ exists, by Nagata's Theorem. Then the map $\Ch_0(\ol{Y})^{\deg 0} \rightarrow \Ch_0(Y)^{\deg 0}$ is surjective, since to a zero-cycle on $Y$ we may add an appropriate (rational) multiple of any closed point on the boundary point of $\ol{Y}$ so that the sum has degree $0$ on $\ol{Y}$. Hence it suffices to show that $\Ch_0(\ol{Y})^{\deg 0} = 0$. 

So we may and do assume that $Y$ is proper and connected. Let $D \in \Ch_0(Y)^{\deg 0}$. By Lemma \ref{lem: path connected}, we may find a connected curve $C$ in $Y$ containing $D$. Since $Y$ is proper we may furthermore assume that $C$ is proper by replacing it with its closure if necessary. 

We next reduce to the case where $C$ is irreducible. Indeed, suppose that $C = \bigcup C_i$ is the union of irreducible components and that $\Ch_0(C_i)^{\deg 0} = 0$ for each $i$. Then any zero cycle in $\Ch_0(C_i) $ is equivalent to one concentrated at a single point (with $\Q$-coefficients); applying this repeatedly, any zero-divisor on $C$ is equivalent to one supported on a single $C_i$.

So we may assume that $C$ is proper and irreducible, and let $\wt{C} \rightarrow C$ be its normalization. Any $D \in \Ch_0(C)^{\deg 0}$ is the image of $\wt{D} \in \Ch_0(\wt{C})^{\deg 0} = \Pic^0_C(\F_q) \otimes_{\Z} \Q$, which vanishes by the finiteness of $ \Pic^0_C(\F_q) $. 

\end{proof}

\begin{cor}\label{cor: chow of mild DM}
Suppose $\cY$ is a finite type separated Deligne-Mumford stack over a field, admitting a Zariski cover by open substacks that each has a finite flat atlas from a quasi-projective scheme. Then $\Ch_0(\cY)^{\deg 0} = 0$. 
\end{cor}

\begin{proof}
By the Keel-Mori Theorem \cite{KM97} (as explained in \cite[Theorem 1.1]{Con}), $\cY$ has a coarse moduli space $Y$. The hypothesis implies that the conditions in \cite[\S 3]{Con} hold. In particular, $Y$ is a scheme and $\cY \rightarrow Y$ is a proper universal homeomorphism \cite[Theorem 3.1]{Con}, so it induces a bijection of connected components that matches proper components with proper components. Applying \cite[Theorem 6.8]{Gil84} to each of the connected components of $\cY$, we obtain $\Ch_0(\cY) \xrightarrow{\sim} \Ch_0(Y)$. As the proper components are also in bijection, this isomorphism takes $\Ch_0(\cY)^{\deg 0} \xrightarrow{\sim} \Ch_0(Y)^{\deg 0} $, which vanishes by Lemma \ref{lem: quasiprojective chow 0}.
\end{proof}

\begin{cor}\label{cor: chow of Sht}
We have $\Ch_0(\Sht^r_{U(1), \LL})^{\deg 0} = 0$. 
\end{cor}

\begin{proof}
The hypotheses of Corollary \ref{cor: chow of mild DM} are satisfied by (a variant with identical proof of) \cite[Proposition 2.16]{Var04}. 
\end{proof}

\section{The corank one case: testing against CM cycles}
\label{sec: CM}
We provide further evidence for the modularity in the corank one case, by intersecting against a certain class of CM cycles constructed in Example \ref{ex:CM}.
In the number field case, an analogous problem was studied by Howard \cite{How12}.

\subsection{Setup}
Let $Y$ be another  smooth projective curve and $\th: Y\rightarrow X$ be a map of degree $n$,  and let $Y'=X'\times_X Y$ their fiber product:
\[
\begin{tikzcd}
& Y'\ar[dl, "2","\nu'"'] \ar[dd, "n", "\th'"'] \\
 Y \ar[dd, "n", "\th"']  & \\
& X' \ar[dl, "2", "\nu"'] \\
X 
\end{tikzcd}
\]
Abusing notation, we will let $\s$ denote the nontrivial involution on both $Y'/Y$ and $X'/X$.
We allow $Y$ to be disconnected and ramified over $X$; but we will assume that the cover $Y'/Y$ remains geometrically non-split over every component (i.e., for every connected component $Y_\alpha$ of $Y$, $Y_\alpha\times_X X'$ is geometrically connected).  

For a line bundle $\LL$ over $X$, let $\Sht^r_{U(1)/Y, \th^\ast\LL}$ be the moduli stack constructed in Example \ref{ex:CM} (see also Example \ref{ex: CM pullback special cycle}). The non-split hypothesis ensures that $\Sht^r_{U(1)/Y, \th^\ast\LL}$ is proper. Taking direct image $\cF_\bullet \mapsto \th'_\ast\cF_\bullet$ along the map $\th':Y'\to X'$ induces a finite unramified morphism
\begin{align}\label{eq:Th def}
\xymatrix{\Th \co \Sht^r_{U(1)/Y, \th^\ast\LL} \ar[r]& \Sht_{U(n),\LL}^r.}
\end{align}
This map defines a class  $$\Th_\ast [\Sht^r_{U(1)/Y,\th^\ast\LL}]\in
 \Ch_{r,c}(\Sht^r_{U(n),\LL})$$ in the Chow group of proper cycles on $\Sht^r_{U(n),\LL}$.

 \subsection{Pullback formula}Let $\cE$ be a line bundle on $X'$. Recall that  $\cA_{\cE,\LL}(k)$ is the set of Hermitian maps $a:\cE\to \s^*\cE^\vee\ot \LL$, where $\cE^\vee$ denotes the Serre dual. Previously in  \eqref{eq:tr A} we have defined a trace map 
$$
\xymatrix{\tr\colon \cA _{\th'^*\cE,  \th^*\LL}(k)\ar[r]&\cA_{\cE, \LL}(k) }.
$$

\begin{prop}Let $\cE$ be a line bundle on $X'$ and let $a\in \cA_{\cE, \LL}(k)$. Then there is a natural decomposition into open-closed substacks: 
\begin{align}\label{eq:disj un}
\xymatrix{ 
 \Sht^r_{U(1)/Y, \th^\ast\LL} \times_{  \Sht^r_{U(n),\LL}} \cZ^{r}_{\cE,\LL}(a)\ar[r]^-{\sim} &\coprod_{\wt a} \cZ_{ \th'^\ast \cE,\th^\ast\LL  }^{r}(\wt a) }
\end{align}
where $\wt a$ runs over all elements in $\cA_{\th'^*\cE,  \th^*\LL}(k)$ such that $ \tr (\wt a)=a $, and  the virtual fundamental classes satisfy
\begin{align}\label{eq:pb}
\Th^! [\cZ_{\cE,\LL}^{r}(a)]\Big|_{\cZ_{ \th^\ast \cE, \th^\ast\LL}^{r}(\wt a)}=[\cZ_{ \th'^\ast \cE,\th^\ast\LL}^{ r}(\wt a)].
\end{align}


\end{prop}
\begin{proof}
This follows from
Example \ref{ex: CM pullback special cycle}.
\end{proof}


It follows immediately that, under the intersection pairing
\begin{align}
\xymatrix{\langle-,-\rangle \co\Ch^r(\Sht^r_{U(n),\LL})\times \Ch_{r,c}(\Sht^r_{U(n),\LL}) \ar[r]& \BQ ,}
\end{align}
we have the following pullback formula:
\begin{equation}\label{eq:pullback}
\Big\langle \cZ^{r}_{\cE,\LL}(a) , \Th_*  [\Sht^r_{U(1)/Y,\th^*\LL} ] \Big\rangle =
 \sum_{\substack{\wt a\in\cA_{\th'^*\cE,  \th^*\LL}(k), \\
 \tr (\wt a)=a  }} \deg[\cZ_{ \th'^\ast \cE,\th^\ast\LL}^{r}(\wt a)].
\end{equation}


\begin{remark}If we assume $Y$ is connected, the pullback relation \eqref{eq:pullback} can be proved without using the derived methods behind Example \ref{ex: CM pullback special cycle}. We sketch a direct argument. Since $Y$ is connected, a map $\th'^{*}\cE\to \cF_{\bu}$ is injective if and only if the induced map $\cE\to \th'_{*}\cF_{\bu}$ is injective.  Therefore \eqref{eq:disj un} restricts to the following {\em Cartesian} diagram for the circle loci of the special cycles (see Definition \ref{defn: twisted similitude shtuka definitions}):
\[
\xymatrix{\coprod_{ \wt a } \cZ_{ \th'^\ast \cE,\th^\ast\LL}^{r}(\wt a) ^\circ   \ar[r] \ar[d] &  \cZ^{r}_{\cE,\LL}(a) ^\circ \ar[d] \\
 \Sht^r_{U(1)/Y, \th^\ast\LL} \ar[r]^-{\Th}&  \Sht^r_{U(n),\LL}}
\]
Note that  $\cZ_{\th'^\ast \cE,\th^\ast\LL}^{r}(0)^\circ$ is empty. By Corollary \ref{cor: proper intersection}, 
all terms in the disjoint union have the expected dimension (i.e., every $\cZ_{ \th^\ast \cE,\th^\ast\LL}^ {r}(\wt a) ^\circ$ has dimension zero, and  $\cZ^{r}_{\cE,\LL}(a) ^\circ$ has dimension $r(n-1)$). Since the bottom map is a LCI morphism, and  $ \cZ^{r}_{\cE,\LL}(a) ^\circ$ is LCI by Corollary \ref{cor: proper intersection} (which is much easier to prove in this special case $n=1$), the Gysin pullbacks along the map $\Theta$ of the fundamental classes $[\cZ^{r}_{\cE,\LL}(a) ^\circ]$ are represented by the naive fundamental classes. This almost proves the relation \eqref{eq:pb}, except for the most degenerate term $\cZ_{\cE,\LL}^{r}[\cE](0)$ corresponding to the Chern classes given by Definition \ref{def:cap Chern}. It remains to show 
$$
\Th^![ \cZ_{\cE,\LL}^{ r}[\cE](0)]=  [\cZ_{ \th'^\ast \cE,\th^\ast\LL}^{ r}[\th'^\ast \cE](0)].
$$
Since the tautological line bundles $\ell_{i}$ on $\Sht^{r}_{U(n),\frL}$ pullback to the  tautological line bundles on $\Sht^{r}_{U(1)/Y,\th'^{*}\frL}$ via $\Th$, this identity is easy to check directly.
\end{remark}

\subsection{Evidence for modularity in the corank one case}
Suppose that $Y = \coprod_{\alpha\in \Irr(Y)} Y_{\alpha}$ is the decomposition of $Y$ into connected components. Let $\wt H_{1}(\BA_{Y_\alpha})$ denote the adelic similitude unitary group $GU(2)$ over (the function field of) $Y_\alpha$. Let $\wt H_{1}(\BA_Y):=\prod_{\alpha} \wt H_{1}(\BA_{Y_\alpha})$.

Recall  the generating series $\wt Z^r_{m}$ of corank $m$ special cycles from Definition \ref{def: gen series before descent} and \eqref{Zm function}
\begin{equation*}
\wt H_{m}(\BA)\ni g\mapsto \wt Z^{r}_{m}(g)\in  \Ch_{r(n-m)}(\Sht^{r}_{GU(n)}).
\end{equation*}
 We now specialize it to the corank one case, i.e., $m=1$. We will denote $\wt Z^r_{X'/X}:=\wt Z^r_{m=1}$ and apply similar notation to the double covers $Y'/Y$ and $Y_\alpha' /Y_\alpha$ for irreducible components $Y_\alpha \subset Y$. We want to intersect the cycle class $\wt Z^r(g)\in \Ch_{r(n-1)}(\Sht^{r}_{GU(n)})$ for $g$ satisfying $c(g)=\LL $ with the cycle  $\Th_\ast [\Sht^r_{U(1)/Y,\th^\ast\LL}]\in
 \Ch_{r,c}(\Sht^r_{U(n),\LL})$. 
To make the statement more concise, we now introduce $GU(1)_{Y/X}$ to be the subgroup scheme of $GU(1)_Y$ with similitude line bundle in $\Pic_X$ (so that $GU(1)_{Y/X}$-torsors are the same as $(\cF,\frL, h)\in \Pic_{Y'}\times \Pic_{X}$ where $h$ is a Hermitian isomorphism $h:\cF\isom \s^{*}\cF^{\vee}\ot \nu'^{*}\th^{*}\frL$), and define  $\Sht^{r}_{GU(1)_{Y/X}}$ accordingly. Parallel to \eqref{eq: GU shtuka decomposition}, there is a decomposition 
\[
\Sht^{r}_{GU(1)_{Y/X}}=\coprod_{\LL\in\Pic_X(k)} \Sht^r_{U(1)/Y,\th^\ast\LL} \times B(\Aut(\LL)(k)).
\]
Then we have a finite morphism 
 $$
 \xymatrix{\Th\co\Sht^{r}_{GU(1)_{Y/X}}\ar[r]& \Sht^{r}_{GU(n)}},
 $$ 
which is the union of components \eqref{eq:Th def} indexed by $\LL$.

By Example \ref{ex: CM pullback special cycle}, we have an open-closed partition 
\begin{align}\label{eq:Sht ulr part 3}
\Sht^{r}_{U(1)/Y,\th^*\LL}=\coprod_{\ul{r}} \Sht^{\ul{r}}_{U(1)/Y,\th^*\LL},  
\end{align}
where $\ul{r}=(r_\alpha)_{\alpha}\in\Z^{\Irr(Y)}$ satisfies $ |\ul{r}|:=\sum_{\alpha} r_\alpha=r$, and $ \Sht^{\ul{r}}_{U(1)/Y,\th^*\LL}:= \prod_{\alpha}  \Sht^{r_\alpha }_{U(1)/Y,\th^*\LL }$.
 Note that by our definition the generating series $\wt Z^r_{Y'/Y}$ is a function
\begin{equation}
\wt Z^r_{Y'/Y}: \wt H_{1}(\BA_Y)=\prod_{\alpha} \wt H_{1}(\BA_{Y_\alpha})\to \Ch_{0,c}(\Sht^r_{GU(1)_{Y}})=\bigoplus_{\LN\in\Pic_Y(k)} \Ch_{0,c}(\Sht^r_{U(1)/Y,\LN})
\end{equation}
Viewing $\wt H_{1}(\BA)$  as a subgroup of $\wt H_{1}(\BA_Y)$ via the diagonal embedding, the restriction $ \wt{Z}^r_{Y'/Y}|_{\wt H_1(\BA)}$ takes values in $\Ch_{0,c}(\Sht^r_{GU(1)_{Y/X}})=\bigoplus_{\LL\in\Pic_X(k)} \Ch_{0,c}(\Sht^r_{U(1)/Y,\th^*\LL})$.

On the analytic side, we denote  by $\wt E(g_{\a},s,\Phi_{Y_{\alpha}})$ the normalized Eisenstein series \eqref{eq:nor Eis} in \S\ref{s:Eis GU} for $n=1$, the covering $Y'_\alpha/Y_\alpha$, and the spherical section $\Phi_{Y_\alpha}$. It is an automorphic form on $\wt H_{1}(\BA_{Y_\alpha})$. Let $\Phi_{Y}=\ot\Phi_{Y_{\a}}$. We define for $g=(g_{\a})\in \wt H_{1}(\BA_{Y})$, 
$$
\wt E(g,s,\Phi_{Y})=\prod_{\alpha\in \Irr(Y)}\wt E(g_{\a},s,\Phi_{Y_{\alpha}}),
$$
which is an automorphic form on $\wt H_{1}(\BA_Y)$.

 We have the following result, which provides evidence for the Modularity Conjecture  \ref{c:mod} in the corank $m=1$ case. 
\begin{thm}\label{th:n=1 Y} 
\begin{enumerate}\item
 We have an equality
$$
\Th^! \wt{Z}^r_{X'/X} =  \wt{Z}^r_{Y'/Y}|_{\wt H_1(\BA)}
$$
of functions on $ \wt H_1(\BA)$ with values in $\Ch_{0,c}(\Sht^r_{GU(1)_{Y/X}})$.

\item For  every $g\in \wt H_1(\BA)$, we have
\begin{equation}\label{eq: higher derivatives}
  \Big\langle \wt{Z}^{r}_{X'/X}(g) , \Th_*[\Sht^r_{GU(1)_{Y/X}} ] \Big\rangle\ =\frac{1}{(\log q)^r}  \left(\frac{d}{ds}\right)^r\Big|_{s=0} \wt E(g,s,\Phi_{Y}).
\end{equation}
In particular, the function
$\wt H_1(\BA)  \ni g\mapsto  \Big\langle \wt{Z}^{r}_{X'/X}(g) , \Th_*[\Sht^r_{GU(1)_{Y/X}} ] \Big\rangle$ defines  an automorphic form on $ \wt H_1(\BA)$.
\end{enumerate}
\end{thm}

\begin{proof} To show the first statement,
suppose that $g\in  \wt H_1(\BA)$  has similitude factor $c(g)=\LL$.  Then both sides take values in  $\Ch_{0,c}(\Sht^r_{U(1)/Y,\th^*\LL})$ and the equality follows from the pull back relation \eqref{eq:pb}. 

To show the second statement, for  $\ul{r}=(r_\alpha)_{\alpha}$ satisfying $\sum_{\alpha} r_\alpha=r$, we have 
$$
\deg( \wt Z^r_{Y'/Y}(g) |_{\Sht^{\ul{r}}_{GU(1)/Y}})=\prod_{\alpha\in\Irr(Y)} \deg \wt Z^{r_\alpha}_{Y_\alpha'/Y_\alpha} (g_\alpha), 
$$
for $g=(g_\alpha)\in \wt H_1(\BA_Y)$. Similarly there is a decomposition of the analytic side, by Lebniz's rule,
$$
 \left(\frac{d}{ds}\right)^r\Big|_{s=0} \wt E(g,s,\Phi_{Y})=\sum_{\un r\in \Z^{\Irr(Y)}, |\ul{r}|=r}     \prod_{\alpha} \left(\frac{d}{ds}\right)^{r_\alpha}\Big|_{s=0} \wt E(g_\alpha,s,\Phi_{Y_{\alpha}}).
$$

By the case of modularity when $n=1$, i.e., Theorem \ref{th:n=1}, we have 
$$
\deg \wt Z^{r_\alpha}_{Y_\alpha'/Y_\alpha} (g_\alpha)=\frac{1}{(\log q)^{r_\alpha}}  \left(\frac{d}{ds}\right)^{r_\alpha}\Big|_{s=0} \wt E(g_\alpha,s,\Phi_{Y_{\alpha}})
$$for $g_\alpha\in \wt H_1(\BA_{Y_\alpha})$. The assertion follows by combining these equalities.
\end{proof}

\begin{remark} In view of \eqref{eq:Sht ulr part 3},
 the proof above shows a refinement of \eqref{eq: higher derivatives}, i.e., for any $\un r\in\Z_{\ge0}^{\Irr(Y)}$ such that $|\ul{r}|=r$,  we have 
$$
  \Big\langle \wt Z^{r}_{X'/X}(g) , \Th_*[\Sht^{\ul{r}}_{GU(1)_{Y/X}} ] \Big\rangle\ =\frac{1}{(\log q)^r} \prod_{\alpha} \left(\frac{d}{ds}\right)^{r_\alpha}\Big|_{s=0} \wt E(g,s,\Phi_{Y_{\a}})
$$
as a function of $g\in \wt H_1(\BA)$.
\end{remark}

\begin{remark}
In the number field case, the theorem of Howard \cite{How12} is analogous to our case where $Y$ is connected and $r=1$. It seems that the analog of the case of disconnected $Y$ in the number field case has not been treated.
\end{remark}

\begin{remark}Since $Y$ is allowed to be ramified over $X$, there are infinitely many such covers. We may form the subspace of $\upH^{2(n-1)r}_{c}(\Sht^r_{U(n),\LL})$ spanned by the cycle classes $ \Th_*[\Sht^r_{U(1)/Y,\th^*\LL} ] $ for varying coverings $Y/X$ of degree $n$. It is an interesting question how large this subspace is.
\end{remark}

\appendix

\bibliographystyle{amsalpha}
\bibliography{Bibliography}

\end{document}